\documentclass[a4paper,?,psamsfonts,12pt,reqno,final]{amsbook}   
\pdfoutput=1         

\usepackage{amsmath,amsthm,amssymb,latexsym}

\usepackage{exscale,mathrsfs,dsfont,bold-extra}   

\makeindex         

\newcommand{\st}{\ensuremath{\ast}}

\newcommand{\s}{\mathrm{sp}}

\newcommand{\iu}{\mathrm{i}}

\newcommand{\id}{\mathrm{id}}

\newcommand{\tld}{\widetilde}

\newcommand{\wht}{\widehat}

\newcommand{\sa}{_{\textstyle{sa}}}

\newcommand{\bsa}{_{\textstyle{bsa}}}

\newcommand{\rlambda}{\mathrm{r}_{\textstyle{\lambda}}}

\newcommand{\rsigma}{\mathrm{r}_{\textstyle{\sigma}}}

\newcommand{\univ}{_{\textit{\small{u}}}}

\newcommand{\ra}{_{\textit{\small{ra}}}}

\newcommand{\smu}{_{\text{\small{$\mu$}}}}

\newcommand{\0}{\ensuremath{_{\mathrm{o}}}}

\newcommand{\twiddle}{\ensuremath{^{\sim}}}

\newcommand{\tla}{_{\textstyle{a}}}

\newlength{\skipamount}
\newtheoremstyle{myplain}
{\skipamount}   
{\skipamount}   
{\itshape}   
{}   
{\scshape\bfseries}   
{.}   
{.5em}   
{}   

\newtheoremstyle{mydefinition}
{\skipamount}   
{\skipamount}   
{}   
{}   
{\scshape\bfseries}   
{.}   
{.5em}   
{}   

\swapnumbers
\theoremstyle{myplain}
\newtheorem{theorem}{Theorem}[section]
\newtheorem{proposition}[theorem]{Proposition}
\newtheorem{lemma}[theorem]{Lemma}
\newtheorem{corollary}[theorem]{Corollary}
\newtheorem{reminder}[theorem]{Reminder}
\newtheorem{addendum}[theorem]{Addendum}
\theoremstyle{mydefinition}
\newtheorem{definition}[theorem]{Definition}
\newtheorem{example}[theorem]{Example}
\newtheorem{counterexample}[theorem]{Counterexample}
\newtheorem{introduction}[theorem]{Introduction}
\newtheorem{remark}[theorem]{Remark}

\begin{document}


\frontmatter


\title{Introduction to Normed \st-Algebras \protect\\
         and their Representations, 6th ed.}

\author{Marco Thill}

\maketitle

\clearpage


\phantomsection

\chapter*{Preface}

\vspace{-3mm}
A synopsis of the work is given in the next three paragraphs.

In part \ref{part1}, the by now classical spectral theory of Banach
\st-algebras is developed, including the Shirali-Ford Theorem
\ref{ShiraliFord}, which is less common for an introduction.

In part \ref{part2}, the theory of states on a normed \st-algebra is
examined, stripped of the usual assumptions of completeness,
presence of an \linebreak approximate unit, and isometry of the
involution. One idea is to consider continuity of mappings on the
real subspace $A\sa$ of Hermitian elements only. We also exploit
the idea of Berg, Christensen, Ressel \cite[Theorem 4.1.14]{BCR},
which becomes our items \ref{weaksa} - \ref{sigmaAsa}. This
deduction of contractivity on $A\sa$ from weak continuity on $A\sa$
was published first (in a different form) in the first edition of this work.
Rather than merely proving the Gelfand-Na\u{\i}mark Theorem, we
attach to each normed \st-algebra a C*-algebra, which may be used
to pull down results. The construction of this so-called enveloping
C*-algebra is one of the main goals of this text. We also give an
original proof of the Theorem of Varopoulos, cf.\ \ref{Varopsect}.

In part \ref{part3}, we consider Spectral Theorems using integrals
converging in norm in the bounded case, and pointwise in the
unbounded case. Our Cauchy-Bochner type of approach leads
to easier constructions, simpler proofs, and stronger results than
the more conventional approach via weakly convergent integrals.
We give the Spectral Theorem for non-degenerate representations
in Hilbert spaces $\neq \{0\}$ of perfectly general commutative
\st-algebras, rather than only of commutative Banach \st-algebras,
or even only of commutative C*-algebras. This is possible at little
extra effort, and it does not seem to have permeated the
mathematical folklore that this is possible at all. We have tried to be
didactic in our derivation of the main disintegrations (the abstract
Bochner Theorem as well as the Spectral Theorem), in that we get
these results in an easy step from the corresponding result for
C*-algebras by taking an image measure.

The material is not restricted to unital algebras. Proofs are complete,
and there are plenty of cross-references. \pagebreak

This book is suitable for accompanying a lecture, or for self-study.
One of the primary design goals of this text has been that it should
be able to serve as a reference. It also could function as an
underpinning for an intermediate course in representation theory
leading to the more advanced treatments. Therefore the prerequisite
for reading this work is a first course in functional analysis, which
nowadays can be expected to cover the Commutative
Gelfand-Na\u{\i}mark Theorem as well as its ingredients: the notion
of a C*-algebra, the weak* topology and Alaoglu's Theorem, as well
as the Stone-Weierstrass Theorem. A course in integration theory
covering the Riesz Representation Theorem is needed for the last
part of the book on the spectral theory of representations. The present
text nicely complements the beautiful introductory text by Mosak \cite{Mos}.

What is not included in the text:
H*-algebras and group algebras are not treated for lack of free energy.
We refer the reader to Loomis \cite{LoAHA} for an introduction to these
topics. We have taken a minimalist approach, in that we develop only
so much of the material as is necessary for a good grasp of the basic
theory. Also, we have avoided material on general Banach algebras,
which is not needed for Banach algebras with involution. In this spirit,
the following omissions have been made. The holomorphic functional
calculus is not developed, as we can do with only power series and the
Rational Spectral Mapping Theorem. Similarly the radical is not considered
as we can do with only the \st-radical. Neither touched upon are the
developments of Theodore W.\ Palmer \cite{Palm} concerning spectral
algebras, which surely is a major omission, especially in the present context.

Now for the ``version history'' of the work.

The first four editions were published by Kassel University Press in Germany.

The second edition differs from the first one, in that the first edition
considered normed \st-algebras with isometric involution exclusively. In
this work, the involution may be discontinuous. I gratefully acknow\-ledge
that Theodore W.\ Palmer and Torben Maack Bisgaard brought to my
attention a few crucial mistakes in a draft of the second edition.

Since the third edition, the work is typeset in \LaTeXe, and contains an
index (with currently more than 450 entries). The content of the third
edition differs mainly from that of the second edition, in that the spectral
theory of representations now deals with general commutative \st-algebras,
rather than merely with normed commutative \st-algebras.\pagebreak

Since the fourth edition, the text is typeset in \AmS-\LaTeX. Chapter
\ref{PositiveElements} has been made completely independent of Chapter
\ref{TheGelfandTransformation}, and thereby freed of the use of the full
version of the axiom of choice, cf.\ \ref{Zorn}.

The fifth edition features many minor improvements. Six new paragraphs
have been included: \ref{boundary}, \ref{StoneCech}, \ref{polarfactoris},
\ref{leftspectrum}, \ref{separating}, \ref{scopug}, as well as \hyperlink{remtop}{\S\ 45}
in the appendix. The material depending on measure theory has been
made to work with most, if not all, of the prevalent definitions. Many thanks
to Torben Maack Bisgaard who reported many typos as well as a couple
of errors (one bad), and who made useful suggestions.

The sixth edition fixes some minor matters, many of which were reported
by Torben Maack Bisgaard. Thank you, Torben!

\clearpage


\tableofcontents

\clearpage


\mainmatter


\part{Spectral Theory of Banach \texorpdfstring{$*$-}{\052\055}Algebras}\label{part1}

\chapter{Basic Properties of the Spectrum}


\section{\texorpdfstring{$*$-}{\052\055}Algebras and their Unitisation}

\begin{definition}[algebra]\index{algebra}%
An \underline{algebra} is a complex vector space $A$ equip\-ped
with an associative multiplication $(a,b) \mapsto ab$ satisfying
the following laws:
\begin{gather*}
\lambda (ab) = (\lambda a ) b = a (\lambda b) \\
c(a+b) = ca + cb\\
(a+b)c = ac + bc
\end{gather*}
for all $\lambda$ in $\mathds{C}$ and all $a,b,c$ in $A$.
\end{definition}

Please note that we consider complex algebras only.

\begin{example}[$\mathrm{End}(V)$]
If $V$ is a complex vector space, then the vector space
$\mathrm{End}(V)$ of linear operators on $V$ is an algebra
under \linebreak composition. (Here ``$\mathrm{End}$''
stands for endomorphism.)
\end{example}

\begin{definition}[\protect\st-algebra]\index{adjoint}%
\index{algebra!s-algebra@\protect\st-algebra}%
\index{involution}\index{s01@\protect\st-algebra|see{algebra}}%
A \underline{\st-algebra} is an algebra $A$  together with an
\underline{involution} $\st : A \to A$, $a \mapsto a^*$, such that
for all $\lambda$ in $\mathds{C}$ and all $a, b$ in $A$, one has
\begin{align*}
{(a^*)}^* & = a                      \qquad \qquad \text{(involution)} \\
{(a+b)}^* & = a^* + b^* \\
{(ab)}^* & = b^*a^*\\
{(\lambda a)}^*&  = \overline{\lambda}a^*.
\end{align*}
For $a \in A$, the element $a^*$ is called the \underline{adjoint} of $a$.
\end{definition}

Please note the precise behaviour of the involution as it acts on a product.
\pagebreak

\begin{example}[$B(H)$]\index{B(H)@$B(H)$}%
If $H$ is a complex Hilbert space, then the algebra $B(H)$ of bounded
linear operators on $H$ is a \st-algebra under the operation of taking
the adjoint of an operator.
\end{example}

All pre-Hilbert spaces in this book shall be \underline{complex} pre-Hilbert
spaces. 

\begin{example}[the group ring \hbox{$\mathds{C}[G]$}]%
\index{C1@$\protect\mathds{C}[G]$}\label{CG}%
Let $G$ be a group, and let $G$ be multiplicatively written. One denotes by
$\mathds{C}[G]$ the complex vector space of complex-valued functions $a$
on $G$ of finite carrier: $\mathrm{carr} (a) = \{ \,g \in G : a(g) \neq 0 \,\}$.
One may consider $G$ as a basis for $\mathds{C}[G]$. More precisely,
for $g$ in $G$, one may consider the function
$\delta _g$ in $\mathds{C}[G]$
defined as the function taking the value $1$ at $g$ and vanishing
everywhere else. For $a$ in $\mathds{C}[G]$, we then have
\[ a = \sum _{\text{\small{$g \in G$}}} a(g)\,\delta _{\text{\small{$g$}}}. \]
It follows that $\mathds{C}[G]$ carries a unique structure of
algebra, such that its multiplication extends the one of $G$, namely
\[ ab = \sum _{\text{\small{$g,h \in G$}}} a(g)\,b(h)\,\delta _{\text{\small{$gh$}}}
\qquad (a, b \in \mathds{C}[G]). \]
An involution then is introduced by defining
\[ a^* := \sum _{\text{\small{$g \in G$}}}
\overline{a(g)}\,\delta _{\text{\small{$g$}}^{\text{\footnotesize{\,-1}}}}
\qquad (a \in \mathds{C}[G]). \]
This makes $\mathds{C}[G]$ into a \st-algebra, called the \underline{group ring}
of $G$.
\end{example}

\begin{definition}[Hermitian and normal elements]%
\index{Hermitian!element}\index{normal!element|(}%
\index{self-adjoint!element}%
Let $A$ be a \linebreak \st-algebra. An element $a$ of $A$ is called
\underline{Hermitian}, or \underline{self-adjoint}, if $a = a^*$.
An element $b$ of $A$ is called \underline{normal}, if
$b^*b = b\,b^*$, that is, if $b$ and $b^*$ commute.
Every Hermitian element is normal of course.
\end{definition}

\begin{definition}[$A\protect\sa$]\index{A7@$A\protect\sa$}%
If $A$ is a \st-algebra, one denotes by \underline{$A\sa$} the set
of Hermitian elements of $A$. Here ``sa'' stands for ``self-adjoint''.
\end{definition}

\begin{remark}
One notes that an arbitrary element $c$ of a  \st-algebra
$A$ can be written uniquely as $c = a+\iu b$ with $a,b$ in $A\sa$,
namely $a = (c+c^*)/2$ and $b = (c-c^*)/(2 \iu)$. Furthermore the
element $c$ is normal if and only if $a$ and $b$ commute. Hence:
\pagebreak
\end{remark}

\begin{proposition}\index{normal!element|)}%
A \st-algebra is commutative if and only if each of its elements is normal.
\end{proposition}

\begin{proposition}\label{Hermprod}
If $a, b$ are Hermitian elements of a \st-algebra, then the product
$ab$ is Hermitian if and only if $a$ and $b$ commute.
\end{proposition}

\begin{definition}[unit, unital algebra]%
\index{unit}\index{algebra!unital}\index{unital!algebra}%
Let $A$ be an algebra. A \underline{unit} in $A$ is a non-zero element
$e$ of $A$ such that
\[ ea = ae = a \]
for all $a$ in $A$. One says that $A$ is a \underline{unital} algebra if $A$
has a unit. 
\end{definition}

Please note that a unit is required to be different from the zero-element.
An algebra can have at most one unit, as is easily seen. Hence any unit
in a \st-algebra is Hermitian. We shall reserve the notation ``$e$'' for the unit.

\begin{example}
If $G$ is a group with unit $e$, then $\mathds{C}[G]$ is a unital
\st-algebra with unit $\delta_{\text{\small{$e$}}}$.
\end{example}

\begin{definition}[unitisation]%
\index{unitisation}\index{algebra!unitisation}\index{A8@$\protect\tld{A}$}%
Let $A$ be an algebra. If $A$ is unital, one defines $\tld{A} := A$. Assume
now that $A$ has no unit. One then defines \linebreak
$\tld{A} := \mathds{C} \oplus A$ (direct sum of vector spaces). One imbeds
$A$ into $\tld{A}$ via $a \mapsto (0,a)$. One defines $e := (1,0)$, so that
$(\lambda,a) = \lambda e + a$ for $\lambda \in \mathds{C}$, $a \in A$.
In order for $\tld{A}$ to become a unital algebra with $e$ as unit and with
multiplication extending the one in $A$, the multiplication in $\tld{A}$ must
be \linebreak defined by
\[ ( \lambda e + a ) ( \mu e + b ) =
( \lambda \mu ) e + ( \lambda b + \mu a + ab )
\quad (\lambda, \mu \in \mathds{C}, a, b \in A), \]
and this definition indeed satisfies the requirements. One says that
$\tld{A}$ is the \underline{unitisation} of $A$. If $A$ is a \st-algebra,
one makes $\tld{A}$ into a unital \st-algebra by putting
\[ ( \lambda e + a )^* := \overline{\lambda} e + a^*
\qquad (\lambda \in \mathds{C}, a \in A). \]
\end{definition}

\begin{definition}[algebra homomorphisms]\label{hom}%
\index{homomorphism}\index{algebra!homomorphism}%
\index{algebra!s-algebra@\protect\st-algebra!homomorphism}%
Let $A,B$ be algebras. An \underline{algebra homomorphism}
from $A$ to $B$ is a linear mapping $\pi : A \to B$ such that
$\pi (ab) = \pi (a) \pi(b)$ for all $a,b \in A$.
If $A,B$ are \st-algebras, and if $\pi$ furthermore satisfies
$\pi (a^*) = \pi (a)^*$ for all $a$ in $A$, then $\pi$ is called a
\underline{\st-algebra homomorphism}. \pagebreak
\end{definition}

\begin{definition}[unital algebra homomorphisms]\label{unitalhom}%
\index{unital!homomorphism}\index{homomorphism!unital}%
Let $A,B$ be unital algebras  with units $e_A,e_B$ respectively.
An algebra homomorphism $\pi$ from $A$ to $B$ is called
\underline{unital}, if it satisfies $\pi (e_A) = e_B$.
\end{definition}

Our last topic in this paragraph will be subalgebras.

\begin{definition}[subalgebra, \protect\st-subalgebra]%
\index{subalgebra}\index{algebra!subalgebra}%
\index{algebra!subalgebra!s-subalgebra@\protect\st-subalgebra}%
\index{s06@\protect\st-subalgebra}%
\index{subalgebra!s-subalgebra@\protect\st-subalgebra}%
\index{algebra!s-algebra@\protect\st-algebra!s-subalgebra@\protect\st-subalgebra}%
A \underline{subalgebra} of an algebra $A$ is a complex vector subspace
$B$ of $A$ such that $ab \in B$ whenever $a, b \in B$.
A \underline{\st-subalgebra} of a \st-algebra $A$ is a subalgebra
$B$ of $A$ containing with each element $b$ its adjoint $b^*$ in $A$.
\end{definition}

\begin{definition}[unital subalgebra]\index{unital!subalgebra}%
\index{algebra!subalgebra!*-unital@unital}\index{subalgebra!*-unital@unital}%
\index{s06@\protect\st-subalgebra!*-unital@unital}%
If $A$ is a unital algebra, then a subalgebra of $A$ containing the unit of $A$
is called a \underline{unital} subalgebra of $A$.
\end{definition}

The following two definitions are our own. 

\begin{definition}[\twiddle-unital subalgebra]%
\label{twiddleunital}\index{t-unital@\protect\twiddle-unital}%
\index{subalgebra!t-unital@\protect\twiddle-unital}%
\index{algebra!subalgebra!t-unital@\protect\twiddle-unital}%
\index{unital!t-unital@\protect\twiddle-unital}%
\index{s06@\protect\st-subalgebra!t-unital@\protect\twiddle-unital}%
Let $B$ be a subalgebra of an algebra $A$. We shall say
that $B$ is \underline{\twiddle-unital} in $A$, if either $B$
has no unit, or else the unit in $B$ also is a unit in $A$.
\end{definition}

We next define a sensible imbedding. (The only sensible one.)

\begin{definition}[the canonical imbedding]%
\label{canimb}\index{canonical!imbedding}\index{imbedding!canonical}%
Let $B$ be a subalgebra of an algebra $A$. The
\underline{canonical imbedding} of $\tld{B}$ into
$\tld{A}$ is described as follows. If $B$ has a unit, it is the map
$\tld{B} \to \tld{A}$ which is the identity map on $\tld{B} = B$.
If $B$ has no unit, it is the only linear map $\tld{B} \to \tld{A}$
which is the identity map on $B$ and which maps unit to unit.
\end{definition}

The next statement gives a rationale for the preceding two definitions.

\begin{proposition}\label{twidunitp}%
Let $B$ be a subalgebra of an algebra $A$. Then $B$ is \twiddle-unital
in $A$ if and only if $\tld{B}$ is a unital subalgebra of $\tld{A}$ under
the canonical imbedding.
\end{proposition}

\clearpage


\section{Normed \texorpdfstring{$*$-}{\052\055}Algebras and their Unitisation}

\begin{definition}[algebra norm, normed algebra, Banach algebra]%
\index{norm|(}\index{algebra!norm}\index{algebra!normed algebra}%
\index{norm!algebra norm}\index{normed!algebra|see{algebra}}%
\index{algebra!Banach algebra}\index{Banach!algebra|see{algebra}}%
An \underline{algebra norm} on an algebra
$A$ is a norm $|\cdot|$ on $A$ such that
\[ |\,ab\,| \leq |\,a\,|\cdot|\,b\,| \]
for all $a,b$ in $A$. The pair $(A,|\cdot|)$ then is called a
\underline{normed algebra}. A normed algebra is called
a \underline{Banach algebra} if the underlying normed space is
complete, i.e.\ if it is a Banach space.
\end{definition}

\begin{example}
If $V$ is a complex normed space, we denote by $B(V)$ the
set of bounded linear operators on $V$. It is a normed algebra
under the operator norm. If $V$ furthermore is complete, i.e.\ if
it is a Banach space, then $B(V)$ is a Banach algebra.
\end{example}

\begin{definition}%
[normed \protect\st-algebra, Banach \protect\st-algebra]%
\index{algebra!normed s-algebra@normed \protect\st-algebra}%
\index{normed!s-algebra@\protect\st-algebra|see{algebra}}%
\index{Banach!s-algebra@\protect\st-algebra|see{algebra}}%
\index{algebra!Banach s-algebra@Banach \protect\st-algebra}%
A \st-algebra \linebreak equipped with an algebra norm shall be
called a \underline{normed \st-algebra}. A
\underline{Banach \st-algebra} shall be a normed \st-algebra such
that the underlying normed space is complete.
\end{definition}

\begin{example}[$B(H)$]\index{B(H)@$B(H)$}%
If $H$ is a complex Hilbert space, then the  Banach
algebra $B(H)$ of bounded linear operators on $H$ is a Banach 
\st-algebra under the operation of taking the adjoint of an operator.
\end{example}

\begin{definition}[auxiliary \protect\st-norm]\index{a05@${"|}{"|}\,a\,{"|}{"|}$}%
\index{involution!isometric}\index{s03@\protect\st-norm}%
\index{norm!auxiliary *-norm@auxiliary \protect\st-norm}
\index{auxiliary *-norm@auxiliary \protect\st-norm}%
\index{s03@\protect\st-norm!auxiliary}\label{auxnorm}
Let $(A, | \cdot |)$ be a normed \linebreak \st-algebra. One says
that the involution in $A$ is \underline{isometric} if
\[ |\,a^*\,| = |\,a\,|\quad \text{for all } a \in A. \]
One then also says that the norm $| \cdot |$ is a \underline{\st-norm},
and that $A$ is \underline{\st-normed}.
If the involution is not isometric, one can introduce an auxiliary algebra
norm by putting
\[ \|\,a\,\| := \sup\,\{ \,|\,a\,|,|\,a^*\,| \,\} \qquad (a \in A). \]
In this norm the involution is isometric. We shall call this norm the
\underline{auxiliary \st-norm}.

We shall reserve the notation $\| \cdot \|$ for the auxiliary \st-norm,
for \linebreak \st-norms in general, and for the norms on pre-C*-algebras
in particular, see \ref{preC*alg} below. \pagebreak
\end{definition}

\begin{proposition}\index{involution!continuous}\label{involcont}%
Let $A$ be a normed \st-algebra. If the  involution in $A$ is
continuous, there exists $c > 0$ such that
\[ |\,a^*\,| \leq c\,|\,a\,| \]
for all $a$ in $A$. The auxiliary \st-norm then is equivalent to
the original norm.
\end{proposition}

\begin{example}\index{C1@$\protect\mathds{C}[G]$}\label{CGn}%
Let $G$ be a group. One introduces an algebra norm $| \cdot |$
on $\mathds{C}[G]$ by putting
\[ |\,a\,| := \sum _{g \in G} |\,a(g)\,| \qquad (a \in \mathds{C}[G]), \]
thus making $\mathds{C}[G]$ into a normed \st-algebra with
isometric involution. (Recall \ref{CG}.)
\end{example}

\begin{definition}[pre-C*-algebra, C*-algebra, C*-property]%
\label{preC*alg}\index{algebra!C*-algebra}%
\index{algebra!pre-C*-algebra}\index{pre-C*-algebra}%
\index{C6@C*-algebra}\index{C7@C*-property}%
We \linebreak shall say that a \underline{pre-C*-algebra} is a normed
\st-algebra $( A, \| \cdot \| )$ such that for all $a \in A$ one has
\[ {\|\,a\,\| \,}^{2} = \|\,a^*a\,\|\quad \text{as well as}\quad \|\,a\,\| = \|\,a^*\,\|. \]
The first equality is called the \underline{C*-property}.
A complete pre-C*-algebra is called a \underline{C*-algebra}.
\end{definition}

Please note that pre-C*-algebras have isometric involution.

There is some redundancy in the above definition,
as the next proposition shows.

\begin{proposition}\label{condC*}%
Let $(A,\|\cdot\|)$ be a normed \st-algebra such that
\[ {\|\,a\,\| \,}^{2} \leq \|\,a^*a\,\| \]
holds for all a in A. Then $(A,\|\cdot\|)$ is a pre-C*-algebra.
\end{proposition}

\begin{proof} For $a$ in $A$ we have
${\|\,a\,\| \,}^2 \leq \|\,a^*a\,\| \leq \|\,a^*\,\| \cdot \|\,a\,\|$,
whence $\|\,a\,\| \leq \|\,a^*\,\|$. Hence also
$\|\,a^*\,\| \leq \| \,{(a^*)}^* \,\| = \|\,a\,\|$,
which implies $\|\,a^*\,\| = \|\,a\,\|$.
It follows that $\|\,a^*a\,\| \leq {\|\,a\,\| \,}^2$.
The converse inequality holds by assumption.
\end{proof}

\begin{corollary}\index{B(H)@$B(H)$}%
If $(H, \langle \cdot, \cdot \rangle)$ is a Hilbert space, then $B(H)$
is a \linebreak C*-algebra, and with it every closed \st-subalgebra.
\pagebreak
\end{corollary}

\begin{proof}
For $a$ in $B(H)$ and $x$ in $H$, we have
\[ {\|\,ax\,\| \,}^2 = \langle ax, ax \rangle = \langle a^*ax, x \rangle
\leq \|\,a^*ax\,\| \cdot \|\,x\,\| \leq \|\,a^*a\,\| \cdot {\|\,x\,\| \,}^2, \]
so that upon taking square roots $\|\,a\,\| \leq {\| \,a^*a \,\| \,}^{1/2}$.
\end{proof}

\medskip
The closed \st-sub\-algebras of $B(H)$, with $H$ a Hilbert space,
are the prototypes of C*-algebras, see the Gelfand-Na\u{\i}mark
Theorem \ref{GelfandNaimark}. 

\begin{definition}[${\ell\,}^{\infty}(X), C_b(X), C(K), C\0(\Omega), C_c(\Omega)$]%
\index{C2@$C(K)$}\index{norm!supremum norm}\index{C3@$C\0(\Omega)$}%
\index{C4@$C_c(\Omega)$}\index{C15@$C_b(X)$}\index{l09@${\ell\,}^{\infty}(X)$}%
\label{Ccdense}%
For a set $X \neq \varnothing$, we denote by ${\ell\,}^{\infty}(X)$ the
algebra of bounded complex-valued functions on $X$, with pointwise
operations. It is a C*-algebra when equipped with complex conjugation
as involution and with the supremum norm $|\cdot|_{\infty}$ as norm.
If $X \neq \varnothing$ is a  Hausdorff space, we denote by $C_b(X)$
the algebra of bounded continuous complex-valued functions on $X$.
It is a C*-subalgebra of ${\ell\,}^{\infty}(X)$. If $K \neq \varnothing$ is a
compact Hausdorff space, we have $C(K) = C_b (K)$. Let now
$\Omega \neq \varnothing$ be a locally compact Hausdorff space. One
denotes by $C\0(\Omega)$ the algebra of continuous complex-valued
functions on $\Omega$ vanishing at infinity, i.e.\ the algebra consisting
of all continuous complex-valued functions $f$ on $\Omega$ such that
for given $\varepsilon > 0$, there is a compact subset $K$ of $\Omega$
with the property that $| f | < \varepsilon$ off $K$. Then $C\0(\Omega)$ is
a C*-subalgebra of $C_b(\Omega)$. One denotes by $C_c(\Omega)$ the
algebra of all continuous complex-valued functions $f$ of compact support
$\mathrm{supp}(f) = \overline{\{ \,x \in \Omega: f(x) \neq 0 \,\}}$. Then
$C_c(\Omega)$ is a pre-C*-algebra dense in $C\0(\Omega)$. (If
$f \in C\0(\Omega)$ and $f \geq 0$, consider $(f-1/n)_+ \in C_c(\Omega)$
for all integers $n \geq 1$.) 
\end{definition}

The C*-algebras $C\0(\Omega)$, with $\Omega$ a locally compact Hausdorff
space, are the prototypes of commutative C*-algebras, see the Commutative
Gelfand-Na\u{\i}mark Theorem \ref{commGN}.

We can now turn to the unitisation of normed \st-algebras.
We shall have to treat pre-C*-algebras apart from general
normed \st-algebras.

\begin{proposition}\label{Cstarleftreg}%
For a pre-C*-algebra  $(A,\| \cdot \|)$ and for $a \in A$, we have
\[ \| \,a \,\| = \max\,\{ \,\| \,ax \,\| : x \in A, \,\| \,x \,\| \leq 1\,\}. \]
If $a \neq 0$, the maximum is achieved at $x = a^*/ \| \,a^* \,\|$.
\end{proposition}

\begin{corollary}\label{unitCstar}%
The unit in a unital pre-C*-algebra has norm $1$.\pagebreak
\end{corollary}

\begin{proposition}\label{preCstarunitis}%
\index{unitisation}\index{pre-C*-algebra!unitisation}%
\index{C6@C*-algebra!unitisation}\index{A8@$\protect\tld{A}$}%
\index{algebra!C*-algebra!unitisation}%
Let $(A,\| \cdot \|)$ be a pre-C*-algebra. By defining
\[ \| \,a \,\| = \sup\,\{ \,\| \,ax \,\| : x \in A, \,\| \,x \,\| \leq 1\, \}
\quad (a \in \tld{A}) \]
one makes $\tld{A}$ into a unital pre-C*-algebra. The above
norm extends of course the norm in $A$ by \ref{Cstarleftreg}.
\end{proposition}

\begin{proof}
Assume that $A$ has no unit. Let $a \in \tld{A}$. One defines a linear
operator $L\tla : A \to A$ by $L\tla x := ax $ for $x \in A$. Let
$\lambda \in \mathds{C}$, $b \in A$ with $a = \lambda e + b$. For
$x \in A$ with $\| \,x \,\| \leq 1$ we have
$\| \,L\tla x \,\| \leq | \,\lambda \,| + \| \,b \,\|$,
so that $L\tla$ is a bounded operator. Furthermore, $\| \,a \,\|$ is the
operator norm of $L\tla$. Thus in order to prove that $\| \cdot \|$ is an
algebra norm, it suffices to show that if $L\tla = 0$ then $a = 0$. So
assume that $L\tla = 0$. For all $x$ in $A$ we then have
$0 = L\tla x = \lambda x + b x$. Thus, if $\lambda = 0$,
we obtain $b = 0$ by proposition \ref{Cstarleftreg}. It now suffices to
show that $\lambda = 0$. So assume that $\lambda$ is different from
zero. With $g := -b/ \lambda \in A$, we get $gx = x$ for all $x$ in $A$.
Hence also $xg^* = x$ for all $x$ in $A$. In particular we have
$gg^* = g$, whence, after applying the involution, $g^* = g$. Therefore
$g$ would be a unit in $A$, in contradiction with the assumption made
initially. Now, in order to show that $(\tld{A},\| \cdot \|)$ is a pre-C*-algebra,
it is enough to prove that $\| L\tla \,\| \leq {\| L_{\textstyle{a^*a}} \,\| \,}^{1/2}$,
cf.\ \ref{condC*}. For $x$ in $A$, we have
\[ {\| \,L\tla x \,\| \,}^{2} = {\| \,ax \,\| \,}^{2} = \| \,(ax)^*(ax) \,\|
\leq \| \,x^* \,\| \cdot \| \,a^*ax \,\|
\leq \| \,L_{\textstyle{a^*a}} \,\| \cdot {\| \,x \,\| \,}^{2}. \qedhere \]
\end{proof}

\begin{definition}[unitisation]\label{normunitis}\index{unitisation}%
\index{norm!unitisation}\index{A8@$\protect\tld{A}$}%
Let $(A , | \cdot |)$ be a normed algebra without unit. If $A$ is not
a pre-C*-algebra, one makes $\tld{A}$ into a unital normed algebra
by putting
\[ | \,\lambda e + a \,| := | \,\lambda \,| + | \,a \,| \]
for all $\lambda \in \mathds{C}$, $a \in A$. If $A$ is a pre-C*-algebra,
one makes $\tld{A}$ into a unital pre-C*-algebra as in the preceding
proposition. (See also \ref{C*unitis} below.)
\end{definition}\index{norm|)}

\begin{proposition}\index{unitisation}%
\index{A8@$\protect\tld{A}$}\label{Banachunitis}%
If $A$ is a Banach algebra, then $\tld{A}$ is a Banach algebra as well.
\end{proposition}

\begin{proof}
This is so because $A$ has co-dimension $1$ in $\tld{A}$
if $A$ has no unit, cf.\ the appendix \ref{quotspaces}.
\end{proof}

\clearpage


\section{The Completion of a Normed Algebra}

\begin{proposition}%
If $(A, | \cdot |)$ is a normed algebra, and $a\0,b\0,a,b \in A$, then
\[ |\,a\0 b\0 - a b\,| \leq |\,a\0\,|\cdot|\,b\0-b\,| +
|\,a\0-a\,|\cdot|\,b\0-b\,| + |\,b\0\,|\cdot|\,a\0-a\,|. \]
It follows that multiplication is jointly continuous
and uniformly so on \linebreak bounded subsets.
\end{proposition}

\begin{theorem}[unique continuous extension]\label{continuation}%
Let $f$ be a function defined on a dense subset of a metric space $X$
and taking values in a complete metric space $Y$. Assume that $f$
is uniformly continuous on bounded subsets of its domain of definition.
Then $f$ has a unique continuous extension $X \to Y$.
\end{theorem}

\begin{proof}[\hspace{-3.55ex}\mdseries{\scshape{Sketch of a proof}}]
Uniqueness is clear. Let $x \in X$ and let $(x_{n})$ be a sequence in the
domain of definition of $f$ converging to $x$. Then $(x_{n})$ is a Cauchy
sequence, and thus also bounded. Since uniformly continuous maps take
Cauchy sequences to Cauchy sequences, it follows that $\bigl(f(x_{n})\bigr)$
is a Cauchy sequence in $Y$, hence convergent to an element of $Y$
denoted by $g(x)$, say. (One verifies that $g(x)$ is independent of the
sequence $(x_n)$.) The function $x \mapsto g(x)$ satisfies the requirements.
\end{proof}

\begin{corollary}[completion of a normed algebra]%
\index{algebra!normed algebra!completion}\index{completion}%
Let $A$ be a \linebreak normed algebra. Then the \underline{completion}
of $A$ as a normed space carries \linebreak a unique structure of Banach
algebra such that its multiplication \linebreak extends the one of $A$.
\end{corollary}

\begin{proposition}
If $A$ is a normed \st-algebra with \underline{continuous} \linebreak
\underline{involution}, then the completion of $A$ carries the structure of a
Banach \st-algebra, by continuation of the involution.
\end{proposition}

The same is not necessarily true for normed \st-algebras with \linebreak
discontinuous involution. Indeed, we shall later give an example of a
commutative normed \st-algebra which cannot be imbedded in a Banach
\st-algebra at all. See the remark \ref{counterexbis} below. 

\begin{example}[${\ell\,}^1(G)$]\index{l1@${\ell\,}^1(G)$}\label{l1G}%
If $G$ is a group, then the completion of $\mathds{C}[G]$ is ${\ell\,}^1(G)$.
It is a unital Banach \st-algebra with isometric involution.
(Recall \ref{CG} and \ref{CGn}.) \pagebreak
\end{example}

\begin{proposition}\label{twidunitcompl}%
A normed algebra is a \twiddle-unital subalgebra of its completion.
(Recall \ref{twiddleunital}.) 
\end{proposition}

\begin{proof}
Let $A$ be a normed algebra and let $B$ be the completion
of $A$. If $A$ has a unit $e$, then $e$ also is a unit in $B$ as
is easily seen as follows. Let $b \in B$ and let $(a_n)$ be a
sequence in $A$ converging to $b$. We then have for example
\[ be = \lim _{n \to \infty} a_n\,e = \lim _{n \to \infty} a_n = b. \qedhere \]
\end{proof}

\medskip
We shall need the following technical result.

\begin{proposition}\label{complunit}\index{norm}%
\index{norm!unitisation}\index{unitisation}%
Let $(A, \| \cdot \|)$ be a pre-C*-algebra without unit.
Let $(\tld{A}, \| \cdot \|_{\tld{A}})$ be the unitisation of $A$ as in
\ref{preCstarunitis}. Let $(B, \| \cdot \|_B)$ be the completion of $A$.
Assume that $B$ contains a unit $e$.
The canonical embedding identifies $\tld{A}$ as a unital algebra with
the unital subalgebra $\mathds{C}e + A$ of $B$, cf.\ \ref{canimb}. (See also
\ref{twidunitp} together with either \ref{twiddleunital} or \ref{twidunitcompl}.)
The two norms $\| \cdot \|_{\tld{A}}$ and $\| \cdot \|_B$ then coincide on $\tld{A}$.
In particular, $A$ is dense in $(\tld{A}, \| \cdot \|_{\tld{A}})$.
\end{proposition}

\begin{proof}
For $a \in \tld{A}$, we have
\begin{align*}
\| \,a \,\|_{\tld{A}} & = \sup \,\{ \,\| \,ax \,\| : x \in A, \,\| \,x \,\| \leq1 \,\} \\
 & = \sup \,\{ \,\| \,ax \,\|_B : x \in A, \,\| \,x \,\|_B \leq1 \,\} \\
 & = \sup \,\{ \,\| \,ay \,\|_B : y \in B, \,\| \,y \,\|_B \leq1 \,\} = \| \,a \,\|_B
\end{align*}
because the unit ball of $A$ is dense in the unit ball of $B$, as is easily seen.
It follows that $A$ is dense in $(\tld{A}, \| \cdot \|_B) = (\tld{A}, \| \cdot \|_{\tld{A}})$.
\end{proof}

\clearpage


\section{The Spectrum}

In this paragraph, let $A$ be an algebra.

We start with a few remarks concerning invertible elements.\index{invertible}

\begin{proposition}\label{leftrightinv}%
If $a \in \tld{A}$ has a left inverse $b$ and a right
inverse $c$, then $b = c$ is an inverse of $a$.
\end{proposition}

\begin{proof}
$b = be = b(ac) = (ba)c = ec = c.$
\end{proof}

\medskip
In particular, an element of $\tld{A}$ can have at most one inverse.
If $a,b$ are invertible elements of $\tld{A}$, then $ab$ is
invertible with inverse $b^{-1}a^{-1}$. It follows that the
invertible elements of $\tld{A}$ form a group under multiplication.

\begin{proposition}
Let $a,b \in \tld{A}$. If both of their products $ab$ and $ba$
are invertible, then both $a$ and $b$ are invertible.
\end{proposition}

\begin{proof}
With $c := {(ab)}^{-1}$, $d := {(ba)}^{-1}$, we get
\begin{align*}
 & e = (ab)c = a(bc)   &  & e = c(ab) = (ca)b \\
 & e = d(ba) = (db)a  &  & e = (ba)d = b(ad).  \\
\intertext{The preceding proposition implies that now}
 & bc = db = a^{-1}    &  \text{ and \qquad }  & ca = ad = b^{-1}. \qedhere
\end{align*}
\end{proof}

\medskip
This result is usually used in the following form:

\begin{corollary}\label{comm1}%
If $a,b$ are \underline{commuting} elements of $\tld{A}$,
and if $ab$ is invertible, then both $a$ and $b$ are invertible.
\end{corollary}

\begin{proposition}\label{comm2}%
Let $a$ be an invertible element of $\tld{A}$. An element of $\tld{A}$
commutes with $a$ if and only if it commutes with $a^{-1}$.
\end{proposition}

\begin{proof}
$ab = ba \ \Leftrightarrow \ b = a^{-1}ba
\ \Leftrightarrow \ ba^{-1} = a^{-1}b$.
\end{proof}

\begin{definition}[the spectrum]%
\index{s2@$\protect\s(a)$}\index{spectrum!of an element}%
For $a \in A$ one defines
\[ \s_A(a) := \{ \lambda \in \mathds{C} :
\lambda e - a \in \tld{A} \text{ is not invertible in } \tld{A} \,\}. \]
One says that $\s_A(a)$ is the \underline{spectrum} of the element
$a$ in the algebra $A$. We shall often abbreviate $\s(a) := \s_A(a)$.
\pagebreak
\end{definition}

The next result will be used tacitly in the sequel.

\begin{proposition}\label{specunit}%
For $a \in A$ we have
\[ \s_{\tld{A}} (a) = \s_A (a). \]
\end{proposition}

\begin{proof}
The following statements for $\lambda \in \mathds{C}$ are equivalent.
\begin{align*}
 & \ \lambda \in \s_{\tld{A}} (a), \\
 & \ \lambda e - a  \in \tilde{\tilde{A}}
\text{ is not invertible in } \tilde{\tilde{A}}, \\
 & \ \lambda e - a \in \tld{A}
\text{ is not invertible in } \tld{A}
\qquad \text{(because $\tilde{\tilde{A}} = \tld{A}$)}, \\
 & \ \vphantom{\tilde{\tilde{A}}}\lambda \in \s_A (a). \qedhere
\end{align*}
\end{proof}

\begin{proposition}\label{speczero}%
If $A$ has no unit, then $0 \in \s_A(a)$ for all $a$ in $A$.
\end{proposition}

\begin{proof}
Assume that $0 \notin \s_{A}(a)$ for some $a$ in
$A$. It shall be shown that $A$ is unital. First, $a$
is invertible in $\tld{A}$. So let for example
$a^{-1} = \mu e + b$ with $\mu \in \mathds{C}$, $b \in A$.
We obtain $ e = ( \mu e + b ) a = \mu a + b a \in A$,
so $A$ is unital.
\end{proof}

\medskip
The following Rational Spectral Mapping Theorem
is one of the most widely used results in spectral theory.

\begin{theorem}[the Rational Spectral Mapping Theorem]%
\index{Theorem!Rational Spectral Mapping}%
\index{Rational Spec.\ Mapp.\ Thm.}%
\label{ratspecmapthm}%
Let $a \in \tld{A}$, and let $r(x)$ be a non-constant rational
function without pole on $\s(a)$. Then
\[ \s\bigl(r(a)\bigr) = r\bigl(\s(a)\bigr). \]
\end{theorem}

\begin{proof}
We may write uniquely
\[ r(x) = \gamma \;\prod  _{i\in I} \;( \alpha _i -x )
\;\prod _{j\in J} \;{( \beta _j -x )}^{-1} \]
where $\gamma \in \mathds{C}$ and
$\{ \alpha _i : i \in I \}$, $\{\beta_j : j \in J \}$
are disjoint sets of complex numbers.
The $\beta_j$ $(j \in J) $ are the poles of $r(x)$.
The element $r(a)$ in $\tld{A}$ then is defined by
\[  r(a) := \gamma \;\prod _{i\in I} \;( \alpha_i e-a )
\;\prod _{j\in J} \;{( \beta _j e-a)}^{-1}. \pagebreak \]
Let now $\lambda$ in $\mathds{C}$ be fixed. The function
$\lambda -r(x)$ has the same poles as $r(x)$, occurring with
the same multiplicities. We may thus write uniquely
\[ \lambda - r(x) = \gamma ( \lambda )
\prod _{k \in K( \lambda )} \;\bigl( \delta _k (\lambda ) -x \bigr)
\;\prod _{j \in J} \;{( \beta _j -x )}^{-1}, \]
where $ \gamma ( \lambda ) \neq 0 $ by the assumption that
$r(x)$ is not constant. Hence also
\[ \lambda e - r(a) = \gamma ( \lambda )
\prod _{k \in K( \lambda )} \;\bigl( \delta_k ( \lambda ) e - a\bigr)
\;\prod _{j \in J} \;{( \beta_j e - a)}^{-1}. \]
By \ref{comm1} \& \ref{comm2} it follows that the following
statements are equivalent.
\begin{align*}
 & \lambda \in \s\bigl(r(a)\bigr), \\
 & \lambda e - r(a)
\text{ is not invertible in } \tld{A}, \\
 & \text{there exists } k
\text{ such that } \delta_k( \lambda ) e - a
\text{ is not invertible in } \tld{A}, \\
 & \text{there exists } k
\text{ such that } \delta_k ( \lambda ) \in \s(a), \\
 & \lambda - r(x)
\text{ vanishes at some point } x\0 \in \s(a), \\
 & \text{there exists } x\0 \in \s(a)
\text{ such that } \lambda = r(x\0), \\
 & \lambda \in r\bigl(\s(a)\bigr). \qedhere
\end{align*}
\end{proof}

\begin{theorem}\label{spechom}%
Let $B$ be another algebra, and let $\pi : A \to B$ be an algebra
homomorphism. For $a \in A$ we have
\[ \s_B \bigl(\pi (a)\bigr) \setminus \{0\} \subset \s_A (a) \setminus \{0\}. \]
\end{theorem}

\begin{proof}
Let $a \in A$ and $\lambda \in \mathds{C} \setminus \{0\}$. Assume
that $\lambda \notin \s_{A} (a)$. We will show that
$\lambda \notin \s_{B} \bigl( \pi (a)\bigr)$. Let $p$ be the unit in
$\tld{A}$ and $e$ the unit in $\tld{B}$. Extend $\pi$ to an algebra
homomorphism $\tld{\pi} : \tld{A} \to \tld{B}$ by requiring
that $\tld{\pi}(p) = e$ if $\vphantom{\tld{A}}A$ has no unit.
The element $\lambda p-a$ has an inverse $b$ in $\tld{A}$.
Then $\tld{\pi} (b) + \lambda ^{-1} \bigl(e- \tld{\pi} (p)\bigr)$ is the
inverse of $\lambda e- \pi (a)$ in $\tld{B}$. Indeed we have for example
\begin{align*}
   & \ \bigl[ \,\tld{\pi} (b) + \lambda^{-1} \bigl( e- \tld{\pi} (p) \bigr) \,\bigr]
       \,\bigl[ \,\lambda e- \pi (a) \,\bigr] \\
 = & \ \bigl[ \,\tld{\pi} (b) + \lambda ^{-1} \bigl( e - \tld{\pi} (p) \bigr) \,\bigr]
     \,\bigl[ \,\tld{\pi} ( \lambda p - a ) + \lambda \bigl( e - \tld{\pi} (p) \bigr) \,\bigr] \\
 = & \ \tld{\pi} (p) + 0 + 0 + {\bigl( e - \tld{\pi} (p) \bigr)}^{2}
 = \tld{\pi} (p) + \bigl(e - \tld{\pi} (p)\bigr) = e. \qedhere
\end{align*}
\end{proof}

For the next result, recall \ref{twiddleunital} -- \ref{twidunitp}. \pagebreak

\begin{theorem}\label{specsubalg}%
If $B$ is a subalgebra of $A$, then for an element $b$ of $B$ we have
\[ \s_{A}(b) \setminus \{0\} \subset \s_{B}(b) \setminus\{0\}. \]
If $B$ furthermore is \twiddle-unital in $A$, then
\[ \s_{A}(b) \subset \s_{B}(b). \]
\end{theorem}

\begin{proof}
The first statement follows from the preceding theorem.
The second statement follows from \ref{twidunitp}.
\end{proof}

\medskip
See also \ref{spbdry} below.

\begin{proposition}\label{speccomm}%
For $a, b \in A$, we have
\[ \s(ab) \setminus \{0\} = \s(ba) \setminus \{0\}. \]
\end{proposition}

\begin{proof}
Let $\lambda \in \mathds{C} \setminus \{0\}$. Assume
that $\lambda \notin \s(ab)$. There then exists $c \in \tld{A}$ with
\[ c(\lambda e - ab) = (\lambda e - ab)c = e.  \]
We need to show that $\lambda \notin \s(ba)$. We claim that
$\lambda ^{-1}(e+bca)$ is the inverse of $\lambda e-ba$, i.e.
\[ (e+bca)(\lambda e - ba) = (\lambda e - ba)(e+bca) = \lambda e. \]
Indeed, one calculates
\begin{align*}
 &\ (e+bca)(\lambda e - ba) &\quad &\ (\lambda e - ba)(e+bca) \\
= &\ \lambda e-ba +\lambda bca-bcaba & = &\ \lambda e-ba+\lambda bca-babca \\
= &\ \lambda e-ba+bc(\lambda e-ab)a & = &\ \lambda e-ba+b(\lambda e-ab)ca \\
= &\ \lambda e-ba+ba = \lambda e & = &\ \lambda e-ba+ba = \lambda e. \qedhere
\end{align*}
\end{proof}

\medskip
In this paragraph, only the last result relates to an involution.

\begin{proposition}\label{specstar}%
If $A$ is a \st-algebra and $a \in A$, then
\[ \s(a^*) = \overline{\s(a)}. \]
\end{proposition}

\begin{proof}
This follows from $(\lambda e - a)^* = \overline{\lambda}e - a^*$.
\end{proof}

\medskip
It follows that the spectrum of a Hermitian element is symmetric
with respect to the real axis. It does not follow that the spectrum
of a Hermitian element is real. In fact, a \st-algebra in which each
Hermitian element has real spectrum is called a Hermitian \st-algebra,
cf.\ \ref{secHerm}. \pagebreak

\clearpage


\section[The Spectral Radius Formula: $\protect{\mathrm{r}_\lambda(\,\cdot\,)}$]%
{The Spectral Radius Formula: \texorpdfstring{$\protect{\rlambda(\,\cdot\,)}$}{r lambda}}

\begin{reminder}\label{exponential}%
For $0 < c < \infty$, we have
\[ {c \,}^{1/n} = \exp \,( \,1/n \cdot \ln c \,) \to 1 \qquad (n \to \infty). \]
\end{reminder}

\begin{definition}[$\rlambda(a)$]%
\index{r1@$\protect\rlambda(a)$}\label{rldef}%
If $A$ is a normed algebra, then for $a \in A$ we define
\[ \rlambda(a) := \inf _{n \geq 1} {\bigl|\,{a \,}^n\,\bigr| \,}^{1/n} \leq |\,a\,|. \]
The subscript $\lambda$ is intended to remind of the notation for eigenvalues.
The reason for this will become apparent in \ref{specradform} below.
\end{definition}

\begin{theorem}%
For a normed algebra $A$ and $a \in A$, we have
\[ \rlambda(a) = \lim _{n \to \infty} {\bigl|\,{a \,}^n\,\bigr| \,}^{1/n}. \]
\end{theorem}

\begin{proof}
Let $a \in A$. Let $\varepsilon > 0$ and choose $k$ such that
\[ {\bigl|\,{a \,}^k\,\bigr| \,}^{1/k} \leq \rlambda (a) + \varepsilon / 2. \]
Every positive integer $n$ can be written uniquely in
the form $n = p(n)k + q(n)$ with $p(n), q(n)$
non-negative integers and $q(n) \leq k-1$. Since for
$n \to \infty$, we have
\[ q(n)/n \to 0, \]
it follows that
\[ p(n)k/n \to 1, \]
and so (if $a \neq 0$)
\[ {\bigl|\,{a \,}^n\,\bigr| \,}^{1/n}
\leq {\bigl|\,{a \,}^k\,\bigr| \,}^{p(n)/n} \cdot {|\,a \,| \,}^{q(n)/n}
\to {\bigl|\,{a \,}^{k}\,\bigr| \,}^{1/k} \leq \rlambda(a) + \varepsilon / 2. \]
Thus
\[ {\bigl|\,{a \,}^{n}\,\bigr| \,}^{1/n} < \rlambda(a) + \varepsilon \]
for all sufficiently large values of $n$.
\end{proof}

We shall now need some complex analysis. We refer the reader to
P.\ Henrici \cite{Hen} or W.\ Rudin \cite{RCAC}. We warn the reader
that for us, a power series is an expression of the form
$\sum _{n=0}^{\infty } z^{n} a_{n}$, where the $a_n$ are elements
of a Banach space, and $z$ is a placeholder for a complex variable.%
\pagebreak

\begin{lemma}\label{leftinv}%
Let $A$ be a unital Banach algebra. Let $c$ be a left \linebreak invertible
element of $A$, and let ${c \,}^{-1}$ be a left inverse of $c$.

\noindent For an arbitrary element $a$ of $A$, the geometric series
\[ G(z) = \sum _{n=0}^{\infty} \,{z \,}^n \,{\bigl( a {c \,}^{-1} \bigr) \,}^n \]
has non-zero radius of convergence
\[ 1\,/ \,\rlambda \bigl( a{c \,}^{-1} \bigr). \]
We have that the element ${c \,}^{-1} \,G(z)$ is a left inverse of $c-za$ whenever
$| \,z \,| < 1 \,/ \,\rlambda \bigl( a{c \,}^{-1} \bigr)$.
In particular, for $| \,a \,| < 1 \,/ \,\bigl| \,{c \,}^{-1} \,\bigr|$, the element $c-a$
has a left inverse given by
\[ {c \,}^{-1} \,\sum _{n=0}^{\infty} {\bigl( a{c \,}^{-1} \bigr) \,}^n. \]
Thus the set of elements with a left inverse is open in $A$.

A similar result holds of course for right invertible elements.
\end{lemma}

\begin{proof}
By the Cauchy-Hadamard formula, the radius of convergence of $G(z)$ is
\[ 1 \,/ \,\limsup _{n \to \infty}
\,{\Bigl|{ \,\bigl(a{c \,}^{-1}\bigr) \,}^n \,\Bigr| \,}^{1/n}
= 1 \,/ \,\rlambda \bigl( a{c \,}^{-1} \bigr) > 0. \]
For $z \in \mathds{C}$ with
$| \,z \,| < 1 \,/ \,\rlambda \bigl( a{c \,}^{-1} \bigr)$,
one computes
\[ G(z) \,\bigl( e-za{c \,}^{-1} \bigr) =
\sum _{n=0}^{\infty} \,{z \,}^n \,{\bigl( a{c \,}^{-1} \bigr) \,}^n
- \sum _{n=1}^{\infty} \,{z \,}^n \,{\bigl( a{c \,}^{-1} \bigr) \,}^n = e, \]
whence $G(z) \,(c-za) = c$, so ${c \,}^{-1} \,G(z) \,(c-za) = e$.
Thus ${c \,}^{-1} \,G(z)$ indeed is a left inverse of $c-za$.
\end{proof}

\begin{theorem}\label{GL(A)}%
The group of invertible elements in a unital  Banach algebra
$A$ is open and inversion is continuous on this group.
\end{theorem}

\begin{proof}
Let $c \in A$ be invertible. By \ref{leftrightinv}, for $a \in A$ with
$| \,a \,| < 1 \,/ \,\bigl| \,{c \,}^{-1} \,\bigr|$, the element $c-a$ is
invertible with inverse
\[ {( c-a ) \,}^{-1} =
{c \,}^{-1} \,\sum _{n=0}^{\infty } \,{\bigl( a{c \,}^{-1} \bigr) \,}^{n}. \]
This shows that the set of invertible elements in $A$ is open.
\pagebreak

\noindent Moreover we have
\begin{align*}
\Bigl| \,{( c-a ) \,}^{-1} - {c \,}^{-1} \,\Bigr|
& = \biggl| \ {c \,}^{-1} \,\biggl( \,\sum _{n=0}^{\infty }
\,{\bigl( a{c \,}^{-1} \bigr) \,}^{n} -e \,\biggr) \,\biggr| \\
 & \leq \bigl| \ {c \,}^{-1} \,\bigr| \cdot
\biggl| \ \sum _{n=1}^{\infty } \,{\bigl( a{c \,}^{-1} \bigr) \,}^{n} \,\biggr| \\
 & \leq \bigl| \ {c \,}^{-1} \,\bigr| \cdot \bigl| \,a{c \,}^{-1} \,\bigr|
\cdot \sum _{n=0}^{\infty } \,{\bigl| \,a{c \,}^{-1} \,\bigr| \,} ^{n} \\
 & \leq | \,a \,| \cdot {\bigl| \ {c \,}^{-1}\,\bigr| \,}^2
\cdot \Bigl(\,1 - |\,a\,| \cdot \bigl| \,{c \,}^{-1} \,\bigr| \,\Bigr)^{-1},
\end{align*}
which shows that inversion is continuous.
\end{proof}

\bigskip
The next result is basic for all of spectral theory.

\medskip\begin{theorem}[the Spectral Radius Formula]\label{specradform}%
\index{Theorem!Spectral Radius Formula}\index{spectral!radius formula}%
For a Banach algebra $A$ and $a \in A$, we have:
\begin{itemize}
  \item[$(i)$] $\s(a)$ is a \underline{non-empty compact} set in the complex plane,
 \item[$(ii)$] $\max\,\{\,|\,\lambda \,| : \lambda \in \s(a) \,\} = \rlambda(a)$.
\end{itemize}
Please note that the right side of the equation depends on the norm in $A$,
whereas the left side does not.
\end{theorem}

\begin{proof}
We consider the function $f$ given by $f(\mu):=(e-\mu a)^{-1}$ with its maximal
domain. It is an analytic function: write $e-(\mu+z)a = (e-\mu a)-za$ and apply the
lemma. It is defined at the origin and the radius of convergence of its power
series expansion there is $1/\rlambda(a)$. Hence $e-\mu a$ is invertible for
$|\,\mu\,|<1/\rlambda(a)$. That is, $\lambda \notin \s(a)$ for
$|\,\lambda\,|>\rlambda(a)$. It follows that $\s(a)$ is compact. If
$\rlambda(a) \neq 0$, then $f$ has a singular point $z\0$ with
$|\,z\0\,| = 1/\rlambda(a)$ (the radius of convergence),
cf.\ \cite[Thm.\ 16.2 p.\ 320]{RCAC} or the proof of \cite[Thm.\ 3.3a p.\ 155]{Hen}.
This implies that $e-z\0a$ cannot be invertible, and consequently $1/z\0 \in \s(a)$.
This proves the result for the case $\rlambda(a) \neq 0$. Assume that
$\rlambda(a)=0$. We must show that $0 \in \s(a)$. If $0 \notin \s(a)$, then $a$ is
invertible in $\tld{A}$, and $e={a \,}^n \,{({a \,}^{-1}) \,}^n$, whence
$0 < |\,e\,| \leq \bigl|\,{a \,}^n\,\bigr| \cdot {\bigl|\,{a \,}^{-1}\,\bigr| \,}^n$,
and so $\rlambda(a) \geq {\bigl|\,{a \,}^{-1}\,\bigr| \,}^{-1} > 0$, cf.\ \ref{exponential}.
\end{proof}

\bigskip
We shall use the above result tacitly. \pagebreak

\clearpage


\section[Properties of $\protect{\mathrm{r}_\lambda(\,\cdot\,)}$]%
{Properties of \texorpdfstring{$\protect{\rlambda(\,\cdot\,)}$}{r lambda}}

\begin{theorem}\label{specmodul}%
An element of a \underline{normed} algebra has non-empty
\linebreak spectrum. Indeed, if $A$ is a normed algebra and
$a \in A$, then $\s(a)$ contains a number of modulus $\rlambda(a)$.
\end{theorem}

\begin{proof} Let $B$ denote the completion of $A$. One then
uses the fact that $\s_B (a) \subset \s_A (a)$ by \ref{specsubalg}
as $A$ is \twiddle-unital in $B$ by \ref{twidunitcompl}.
\end{proof}

\begin{theorem}[Gelfand, Mazur]\label{GelfandMazur}%
\index{Theorem!Gelfand-Mazur}\index{Gelfand-Mazur Theorem}%
Let $A$ be a unital normed algebra which is a division ring. The map
\[ \mathds{C} \to A,\ \lambda \mapsto \lambda e \]
then is an algebra isomorphism onto $A$ which also is a homeomorphism.
\end{theorem}

\begin{proof} This map is surjective, because for $a \in A$ there
exists $\lambda \in \s(a)$. Indeed $\lambda e-a$ is not invertible in $A$,
so that $\lambda e-a = 0$, or $a =\lambda e$. The map is injective and
homeomorphic because $|\,\lambda e\,| = |\,\lambda \,|\cdot|\,e\,|$ and
$|\,e\,| \neq 0$.
\end{proof}

\medskip
We shall next give some basic properties of $\rlambda (\,\cdot\,)$.
The first three of them hold in general normed algebras.

\begin{proposition}\label{rlpowers}%
For a normed algebra $A$ and an element $a \in A$, we have
\[ \rlambda \bigl( {a \,}^k \bigr) = {\rlambda (a ) \,}^k \quad \text{for all integers } k \geq 1. \]
\end{proposition}

\begin{proof} One calculates
\[ \rlambda \bigl( {a \,}^k \bigr)
= \lim\limits_{n \to \infty} \,{\Bigl| \,{{ \bigl( {a \,}^k \bigr) \,}^n} \,\Bigr| \,}^{1/n}
= \lim\limits_{n \to \infty} \,{\left( \,{\left|\,{a \,}^{kn}\,\right| \,}^{1/kn} \,\right) \,}^k
= {\rlambda ( a ) \,}^k. \qedhere \]
\end{proof}

\begin{proposition}\label{rlcomm}%
For a normed algebra $A$ and elements $a, b \in A$, we have
\[ \rlambda ( a b ) = \rlambda ( b a ). \]
\end{proposition}

\begin{proof}
We express ${( a b ) \,}^{n+1}$ through ${( b a ) \,}^n$:
\begin{align*}
& \ \rlambda ( a b )
= \lim\limits_{n \to \infty} {\bigl| \,{{( a b ) \,}^{n+1}} \,\bigr| \,}^{1/(n+1)}
= \lim\limits_{n \to \infty} {\bigl| \,a \,{( b a ) \,}^n \,b \,\bigr| \,}^{1/(n+1)} \\
\leq & \ \lim\limits_{n \to \infty} {\bigl( \,{| \,a \,| \cdot | \,b \,|} \,\bigr) \,}^{1/(n+1)} \cdot
\lim\limits_{n \to \infty} {\Bigr( \,{\bigl| \,{{( b a ) \,}^n} \,\bigr| \,}^{1/n} \,\Bigr) \,}^{n/(n+1)}
\leq \rlambda ( b a ).
\end{align*}
This result also follows from \ref{speccomm}. \pagebreak
\end{proof}

\begin{proposition}\label{commrlsub}%
If $a,b$ are \underline{commuting} elements of a normed
\linebreak algebra $A$, then
\begin{gather*}
\rlambda(ab) \leq \rlambda(a) \,\rlambda(b), \\
\rlambda(a+b) \leq \rlambda (a) + \rlambda(b),
\end{gather*}
whence also
\[ |\,\rlambda(a) - \rlambda(b)\,|
\leq \rlambda(a-b) \leq |\,a-b\,|. \]
Thus, if $A$ is commutative, then $\rlambda$ is uniformly continuous.
\end{proposition}

\begin{proof}
By commutativity, we have
\[ { \bigl| \,{( \,a \,b \,) \,}^n \,\bigr| \,}^{1/n} = { \bigl|\,{a \,}^n \,{b \,}^n\,\bigr| \,}^{1/n}
\leq { \bigl| \,{a \,}^n \,\bigr| \,}^{1/n} \cdot { \bigl|\,{b \,}^n \,\bigr| \,}^{1/n}, \]
and it remains to take the limit to obtain the first inequality. In order
to prove the second inequality, let $s, t$ with $\rlambda(a) < s$,
$\rlambda(b) < t$ be arbitrary and define $c := {s \,}^{-1} \,a$, $d := {t \,}^{-1} \,b$.
As then $\rlambda(c)$, $\rlambda(d) <1$, there exists
$0 < \alpha < \infty$ with $|\,{c \,}^n\,|$, $|\,{d \,}^n\,| \leq \alpha$ for all
$n \geq 1$. We then have
\[ \bigl| \,{( \,a+b \,) \,}^n \,\bigr| \leq \sum _{k=0} ^n
\,\bigl({\,^{n\,}_{k\,}}\bigr)
\,{s \,}^k \,{t \,}^{n-k} \,\bigl| \,{c \,}^k \,\bigr| \cdot \bigl| \,{d \,}^{n-k} \,\bigr|
\leq {\alpha \,}^2 { \,( \,s + t \,) \,}^n. \]
It follows that
\[ { \bigl| \,{( \,a+b \,) \,}^n \,\bigr| \,}^{1/n} \leq {\alpha \,}^{2/n} \,( \,s + t \,). \]
By passing to the limit, we get $\rlambda(a+b) \leq s+t$,
cf.\ \ref{exponential}. Since $s$, $t$ with $\rlambda (a) < s$
and $\rlambda (b) < t$ are arbitrary, it follows
$\rlambda (a + b) \leq \rlambda (a) + \rlambda (b)$.
\end{proof}

\medskip
See also \ref{commspeccont} and \ref{specsub} below.

Now four properties of $\rlambda (\,\cdot\,)$ in Banach algebras.

\begin{proposition}\label{rlunit}%
Let $A$ be a Banach algebra without unit. For all
$\mu \in \mathds{C}$, $a \in A$, we then have
\[ \rlambda(\mu e + a) \leq | \,\mu \,| + \rlambda(a)
\leq 3 \,\rlambda(\mu e + a). \]
\end{proposition}

\begin{proof}
The first inequality follows from \ref{commrlsub}. For the same reason
$\rlambda (a) \leq \rlambda (\mu e + a) + | \,\mu \,|$, and so
$| \,\mu \,| + \rlambda (a) \leq \rlambda (\mu e + a) + 2 \,| \,\mu \,|$.
It now suffices to prove that $| \,\mu \,| \leq \rlambda (\mu e + a)$.
We have $0 \in \s (a)$ by \ref{speczero}, and so
$\mu \in \s (\mu e + a)$, which implies
$| \,\mu \,| \leq \rlambda (\mu e + a)$, as required.
\end{proof}

\medskip
In this paragraph, merely the next result relates to an involution.
\pagebreak

\begin{proposition}\label{Brlsymm}%
If $A$ is a Banach \st-algebra, then
\[ \rlambda(a^*) = \rlambda(a) \qquad \text{for all } a \in A. \]
\end{proposition}

\begin{proof}
One uses the Spectral Radius Formula and the fact that
\[ \s(a^*) = \overline{\s(a)}, \]
cf.\ \ref{specstar}.
\end{proof}

\medskip
The above result is not true for normed \st-algebras in general,
cf.\ \ref{counterexbis} below.

\begin{theorem}\label{distance}%
Let $A$ be a Banach algebra, and let $a \in A$.
If $\mu \in \mathds{C} \setminus \s(a)$, then
\[ \mathrm{dist} \,\bigl( \,\mu, \,\s(a) \,\bigr) =
{\rlambda \bigl( \,{(\mu e - a) \,}^{-1} \,\bigr) \,}^{-1}. \]
\end{theorem}

\begin{proof}
By the Rational Spectral Mapping Theorem, we have
\[ \s \,\bigl( \,{(\mu e - a) \,}^{-1} \,\bigr) =
\bigl\{ \,{(\mu - \lambda) \,}^{-1} : \lambda \in \s(a) \,\bigr\}. \]
From the Spectral Radius Formula, it follows that
\begin{align*}
\rlambda \bigl( \,{(\mu e - a) \,}^{-1} \,\bigr)
& = \max \,\bigl\{ \,{| \,\mu - \lambda \,| \,}^{-1} : \lambda \in \s(a) \,\bigr\} \\
& = {\min \,\bigl\{ \,| \,\mu - \lambda \,| : \lambda \in \s(a) \,\bigr\} \,}^{-1} \\
& = {\mathrm{dist} \,\bigl( \,\mu, \,\s(a) \,\bigr) \,}^{-1}. \qedhere
\end{align*}
\end{proof}

\begin{theorem}\label{commspeccont}%
Let $A$ be a Banach algebra. Let $\mathds{D}$ denote the
closed unit disc. For two \underline{commuting} elements
$a, b \in A$, we have
\[ \s(a) \subset \s(b) + \rlambda (b-a) \,\mathds{D}
\subset \s(b) + | \,b-a \,| \,\mathds{D}. \]
If $A$ is commutative, one says that the spectrum function is
uniformly continuous on $A$.
\end{theorem}

\begin{proof}
Suppose the first inclusion is not true. There then exists\vphantom{$\bigr)^{-1}$}
$\mu \in \s(a)$ with $\mathrm{dist}\bigl(\mu, \s(b)\bigr) > \rlambda (b-a)$.
Please note that $\mu e - b$ is invertible. By \ref{distance}, we have
$\rlambda \bigl( \,{(\mu e - b)\,}^{-1} \,\bigr) \,\rlambda (b-a) < 1$. From
\ref{commrlsub}, we get $\rlambda \bigl( \,{(\mu e - b) \,}^{-1} \,(b-a) \,\bigr) < 1$.
We have $\mu e-a = (\mu e - b) + (b-a)$. This implies that also
$\mu e-a = (\mu e - b) \,\bigl[ \,e + {(\mu e - b) \,}^{-1} \,(b-a) \,\bigr]$ is invertible,
a contradiction.\vphantom{$)^{-1}$}
\end{proof}

\medskip
See also \ref{commrlsub} above and \ref{specsub} below. \pagebreak

\clearpage


\section[The Pt\'ak Function: $\protect{\mathrm{r}_\sigma(\,\cdot\,)}$]%
{The Pt\'ak Function: \texorpdfstring{$\protect{\rsigma(\,\cdot\,)}$}{r sigma}}

\begin{definition}[$\rsigma(a)$]%
\index{r2@$\protect\rsigma(a)$}\index{Pt\'ak function}%
If $A$ is a normed \st-algebra, then for $a$ in $A$ we define
\[ \rsigma (a) := {\rlambda (a^*a) \,}^{1/2} \leq {| \,a^*a \,| \,}^{1/2}. \]
The subscript $\sigma$ is intended to remind of the notation for
singular values. The function $a \mapsto \rsigma(a)$ $(a \in A)$
is called the \underline{Pt\'ak function}.
\end{definition}

\begin{proposition}\label{Hermrseqrl}%
For a normed \st-algebra $A$ and a Hermitian element
$a \in A$, we have
\[ \rsigma(a) = \rlambda(a). \]
\end{proposition}

\begin{proof}
By \ref{rlpowers} we have
\[ {\rsigma(a) \,}^2 = \rlambda(a^*a)
= \rlambda \bigl({a \,}^2 \bigr) = {\rlambda(a) \,}^2. \qedhere \]
\end{proof}

\begin{proposition}\label{Bnormalrsleqrl}%
For a \underline{Banach} \st-algebra $A$ and a normal
element $b \in A$, we have
\[ \rsigma(b) \leq \rlambda(b). \]
\end{proposition}

\begin{proof} By \ref{commrlsub} and \ref{Brlsymm} we have
\[ {\rsigma(b) \,}^2 = \rlambda(b^*b) \leq
\rlambda(b^*) \,\rlambda(b) = {\rlambda(b) \,}^2. \qedhere \]
\end{proof}

\begin{proposition}\label{rsC*}
For  a normed \st-algebra $A$ and an element $a \in A$, we have
\begin{gather*}
\rsigma(a^*a) = {\rsigma(a) \,}^2, \\
\rsigma(a^*) = \rsigma(a).
\end{gather*}
\end{proposition}

\begin{proof}
With \ref{rlpowers} we have
\[ \rsigma(a^*a) = {\rlambda(a^*aa^*a) \,}^{1/2}
= {\bigl( \,{\rlambda (a^*a) \,}^2 \,\bigr) \,}^{1/2}
= \rlambda(a^*a)={\rsigma(a) \,}^2. \]
With \ref{rlcomm} we have
\[ \rsigma(a^*) = {\rlambda(aa^*) \,}^{1/2}
= {\rlambda(a^*a) \,}^{1/2} = \rsigma(a). \qedhere \]
\end{proof}

\begin{theorem}\label{C*rs}%
Let $(A,\| \cdot \|)$ be a normed \st-algebra. If $(A,\| \cdot \|)$
is a pre-C*-algebra, then
\[ \| \,a \,\| = \rsigma(a) \]
for all $a$ in $A$, and conversely. \pagebreak
\end{theorem}

\begin{proof}
If $\| \,a \,\| = \rsigma(a)$ for all $a\in A$, then $(A,\| \cdot \| )$ is a
pre-C*-algebra by the preceding proposition.
Conversely, let $(A,\| \cdot \|)$ be a pre-C*-algebra. We then have
\[ \bigl\| \,{( \,a^*a \,) \,}^{( \,{2 \,}^n \,)} \,\bigr\| = { \,\| \,a^*a \,\|\,} ^{( \,{2 \,}^n \,)}
\tag*{$(n \geq 1, a \in A )$}, \]
which is seen by induction as follows. We have
\[ \bigl\| \,{( \,a^*a \,) \,}^2 \,\bigr\| = \| \,{( \,a^*a \,) \,}^* \,( \,a^*a \,) \,\| = {\| \,a^*a \,\|\,}^2, \]
and so
\begin{align*}
\bigl\| \,{( \,a^*a \,) \,}^{( \,{2 \,}^n \,)} \,\bigr\|
 & = \bigl\| \,{\bigl[ \,{( \,a^*a \,) \,}^* \,( \,a^*a \,) \,\bigl] \,}^{( \,{2 \,}^{n-1} \,)} \,\bigr\| \\
 & = {\| \,{( \,a^*a \,) \,}^* \,( \,a^*a \,) \,\| \,} ^{( \,{2 \,}^{n-1} \,)}
 \tag*{by induction hypothesis} \\
 & = {\| \,a^*a \,\|\,}^{( \,{2 \,}^n \,)}.
\end{align*}
It follows that
\begin{align*}
{\rsigma(a) \,}^2
 & = \lim _{n \to \infty} { \bigl\| \,{( \,a^*a \,) \,}^{( \,{2 \,}^n \,)} \,\bigr\|\,}^{1 \,/ \,( \,{2 \,}^n \,)} \\
 & = \lim _{n \to \infty} {\Bigl( \,{\| \,a^*a \,\|\,}^{( \,{2 \,}^n) \,} \,\Bigr) \,}^{1\,/ \,( \,{2 \,}^n \,)}
     = \| \,a^*a \,\| = {\| \,a \,\|\,}^2. \qedhere
\end{align*}
\end{proof}

\begin{theorem}\label{C*rl}%
For a pre-C*-algebra $(A,\| \cdot \|)$ and a normal
element $b$ of $A$, we have
\[ \| \,b \,\| = \rlambda(b). \]
\end{theorem}

\begin{proof}
We can assume that $A$ is a C*-algebra. By \ref{Bnormalrsleqrl},
we get $\| \,b \,\| = \rsigma(b) \leq \rlambda(b) \leq \| \,b \,\|$.
\end{proof}

\begin{corollary}
Let $A$ be a pre-C*-algebra. For a \underline{normal} element $b$
of $A$ and every integer $n \geq 1$ we have
\[ \| \,{b\,}^n \,\| = {\| \,b \,\|\,}^n. \]
\end{corollary}

\begin{proof} This follows now from \ref{rlpowers}. \end{proof}

\begin{theorem}\index{norm}\index{norm!unitisation}\label{C*unitis}%
\index{unitisation}\index{A8@$\protect\tld{A}$}%
Let $(A,\| \cdot \| )$ be a C*-algebra without unit. The \linebreak unitisation
$\tld{A}$ then carries besides its C*-algebra norm $\| \cdot \| $,
the norm $| \,\lambda e+a \,| := | \,\lambda \,| + \|\,a\,\|$
$(\lambda \in \mathds{C}$, $a \in A)$, cf.\ \ref{normunitis}. 
Both norms are equivalent on $\tld{A}$. More precisely, we have
for $b \in \tld{A}$ normal:
\[ \|\,b\,\| \leq |\,b\,| \leq 3\:\|\,b\,\|, \]
respectively $a \in \tld{A}$ arbitrary:
\[ \|\,a\,\| \leq |\,a\,| \leq 6\:\|\,a\,\|. \pagebreak \]
\end{theorem}

\begin{proof}
The statement for normal elements follows directly from \ref{rlunit} and
\ref{C*rl}. The statement for arbitrary elements uses the isometry of the
involution in the C*-algebra norm.
\end{proof}

\medskip
For a later purpose, we record:

\begin{proposition}\label{commrs}%
For elements $a,b$ of a normed \st-algebra such that
$a^*a$ commutes with $b\,b^*$, we have
\[ \rsigma(ab) \leq \rsigma(a)\,\rsigma(b). \]
In particular, the function $\rsigma$ is submultiplicative on a
commutative normed \st-algebra.
\end{proposition}

\begin{proof}
From \ref{rlcomm} and \ref{commrlsub} it follows that
\begin{align*}
 & {\rsigma(ab)}^2 = \rlambda(b^*a^*ab) = \rlambda(a^*abb^*) \\
 \leq\ & \rlambda(a^*a)\,\rlambda(b\,b^*)
 = \rlambda(a^*a)\,\rlambda(b^*b) = {\rsigma(a)}^2\:{\rsigma(b)}^2. \qedhere
\end{align*}
\end{proof}

\clearpage

\vfill


\section{Automatic Continuity}

\begin{definition}[contractive linear map]\index{contractive}%
A linear mapping between (possibly real) normed spaces
is said to be \underline{contractive} if it is bounded with
norm not exceeding $1$.
\end{definition}

\begin{theorem}\label{contraux}%
Let $\pi : A \to B$ be a \st-algebra homomorphism from a Banach
\st-algebra $A$ to a pre-C*-algebra $B$. Then $\pi$ is contractive
in the auxiliary \st-norm $\| \cdot \|$ on $A$. (Cf.\ \ref{auxnorm}.)
\end{theorem}

\begin{proof} We may assume that $B$ is complete. For $a \in A$, we have
\[ \| \,\pi (a) \,\| = \rsigma\bigl(\pi(a)\bigr) \leq \rsigma(a)
\leq {| \,a^*a \,| \,}^{1/2} = {\| \,a^*a \,\| \,}^{1/2} \leq \| \,a \,\|, \]
where the first inequality stems from the Spectral Radius Formula and the
fact that $\s\bigl(\pi (a^*a)\bigr) \setminus \{ 0 \} \subset \s(a^*a) \setminus \{ 0 \}$,
cf.\ \ref{spechom}.
\end{proof}

\begin{corollary}%
A \st-algebra homomorphism from a Banach \st-alge\-bra with isometric
involution to a pre-C*-algebra is contractive.
\end{corollary}

In particular we have:

\begin{corollary}\label{Cstarcontr}%
A \st-algebra homomorphism from a C*-algebra to a \linebreak
pre-C*-algebra is contractive.
\end{corollary}

\begin{corollary}\label{Cstarisom}%
\index{isomorphism of C*-algebras}\index{C6@C*-algebra!isomorphism}%
\index{algebra!C*-algebra!isomorphism}%
A \st-algebra isomorphism between two C*-algebras is isometric,
and shall therefore be called a \underline{C*-algebra isomorphism}
or an \underline{isomorphism of C*-algebras}.
\end{corollary}

See also \ref{C*injisometric} \& \ref{reminjiso} below.

Our next aim is the automatic continuity result \ref{automcont}.
The following two lemmata will find a natural reformulation at a later place.
(See \ref{reform2} and \ref{reform1} below.)

\begin{lemma}\label{kerclosed}%
Let $\pi : A \to B$ be a \st-algebra homomorphism from a Banach \st-algebra
$A$ to a pre-C*-algebra $B$. Then $\ker \pi$ is closed in $A$.\pagebreak
\end{lemma}

\begin{proof} We may assume that $B$ is complete. We shall
use the fact that then
\[ \rlambda\bigl(\pi(a)\bigr) \leq \rlambda(a) \]
holds for all $a \in A$, cf.\ \ref{spechom}. Let $a$ belong to the closure of
$\ker \pi$. Let $x \in A$ be arbitrary. Then $xa$ also is in the closure of
$\ker \pi$. For a sequence $(k_n)$ in $\ker \pi$ converging to $xa$, we have
\[ \rlambda\bigl(\pi (xa)\bigr) = \rlambda\bigl(\pi (xa-k_n)\bigr)
\leq \rlambda(xa-k_n) \leq | \,xa-k_n \,| \to 0, \]
which implies
\[ \rlambda\bigl(\pi(xa)\bigr) = 0. \]
Choosing $x := a^*$, we obtain
\[ 0 = \rsigma\bigl(\pi(a)\bigr) = \| \,\pi (a) \,\|, \]
i.e.\ $a$ is in $\ker \pi$.
\end{proof}

\begin{lemma}\label{injclosgraph}%
Let $\pi : A \to B$ be a \st-algebra homomorphism from a normed
\st-algebra $A$ to a pre-C*-algebra $B$, which is continuous on
$A\sa$. If $\pi$ is injective, then the involution in $A$ has closed
graph.
\end{lemma}

\begin{proof}
Let $a = \lim a_n$, $b = \lim a_n^*$ in $A$. Since $\pi$ is bounded
on $A\sa$, there exists $c \geq 0$, such that for all $n$, one has
\begin{align*}
 {\| \,\pi (b-a_n^*) \,\| \,}^2 & = \| \,\pi \bigl((b-a_n^*)^*(b-a_n^*)\bigr) \,\| \\
 & \leq c \cdot | \,(b-a_n^*)^*(b-a_n^*) \,| \\
 & \leq c \cdot | \,b^*-a_n \,| \cdot | \,b-a_n^* \,|.
\end{align*}
Since $(a_n)$ is bounded in $A$, it follows that
\[ \pi (b) = \lim \pi (a_n^*). \]
Similarly, we get
\[ \pi (a^*) = \lim \pi (a_n^*). \]
One concludes that $\pi (b) = \pi (a^*)$, which implies
$b = a^*$ if $\pi$ is injective.
\end{proof}

We now have the following automatic continuity result.

\begin{theorem}\label{automcont}%
A \st-algebra homomorphism from a Banach \st-algebra
to a pre-C*-algebra is continuous. \pagebreak
\end{theorem}

\begin{proof}
Let $\pi : A \to B$ be a \st-algebra homomorphism from a Banach \linebreak
\st-algebra $A$ to a pre-C*-algebra $B$. The kernel of $\pi$ is closed,
by \ref{kerclosed}. This implies that $C := A \,/ \ker \pi$ is a Banach \st-algebra
in the quotient norm, cf.\ the appendix \ref{quotspaces}. Furthermore, the map
$\pi$ factors to an injective \st-algebra homomorphism $\pi_1$ from $C$ to $B$.
The mapping $\pi_1$ is contractive in the auxiliary \st-norm on $C$,
cf.\ \ref{contraux}. This implies that the involution in $C$ has closed graph,
cf.\ \ref{injclosgraph}. But then the involution is continuous, as $C$ is complete.
Hence the auxiliary \st-norm on $C$ and the quotient norm on $C$ are equivalent,
cf.\ \ref{involcont}. It follows that $\pi_1$ is continuous in the quotient norm on $C$,
which is enough to prove the theorem. 
\end{proof}

\begin{proposition}\label{bdedlambcontract}%
Let $A,B$ be normed algebras and let $\pi : A \to B$ be a continuous algebra
homomorphism. We then have
\[ \rlambda\bigl(\pi(a)\bigr) \leq \rlambda(a) \qquad \text{for all } a \in A. \]
\end{proposition}

\begin{proof}
Let $c > 0$ be a bound of $\pi$. Let $a \in A$. For all
integers $n \geq 1$, we have
\[ \rlambda\bigl(\pi(a)\bigr) \leq { \bigl| \,{\pi(a) \,}^n \,\bigr| \,}^{1/n}
= {\bigl| \,\pi \bigl( {a \,}^n \bigr) \,\bigr| \,}^{1/n}
\leq {c \,}^{1/n} { \,\bigl| \,{a \,}^n \,\bigr| \,}^{1/n} \to \rlambda(a). \qedhere \]
\end{proof}

\begin{corollary}\label{commbdedcontr}%
Let $\pi$ be a continuous \st-algebra homomorphism from a normed
\st-algebra $A$ to a pre-C*-algebra $B$. For a \underline{normal}
\linebreak element $b$ of $A$, we then have
\[ \| \,\pi(b) \,\| \leq \rlambda(b) \leq | \,b \,|. \]
It follows that if $A$ is commutative then $\pi$ is contractive.
\end{corollary}

\begin{proof}
This follows now from \ref{C*rl}.
\end{proof}

\medskip
Hence the following enhancement of \ref{automcont}.

\begin{corollary}\label{commBcontr}%
A \st-algebra homomorphism from a commutative
Banach \st-algebra to a pre-C*-algebra is contractive.
\end{corollary}

The following two paragraphs prepare the ground for the ensuing
\linebreak \ref{secHerm} on Hermitian Banach \st-algebras.

\clearpage


\section{Square Roots of Invertible Elements}

\begin{theorem}[Mertens]\label{Mertens}%
\index{Theorem!Mertens}\index{Mertens' Theorem}%
Let $A$ be a Banach algebra. Let
\[ \sum _{n \geq 0} a_n\quad \text{and}\quad \sum _{n \geq 0} b_n \]
be two convergent series in $A$ and assume that the first of them
converges \underline{absolutely}. Put
\[ c_{n} :=\sum _{k=0}^n a_k b_{n-k} \]
for every integer $n \geq 0$. The series
\[ \sum _{n \geq 0} c_n \]
is called the \underline{Cauchy product} of the above two series. We have
\[ \sum _{n \geq 0} c_n = \biggl( \,\sum _{n \geq 0} a_n \biggr)
\cdot \biggl( \,\sum _{n \geq 0} b_n \biggr). \]
\end{theorem}

\begin{proof} Let
\begin{alignat*}{3}
A_n & := \sum _{k=0}^n a_k,
 & \qquad B_n & := \sum _{k=0}^n b_k,
 & \qquad C_n & := \sum _{k=0}^n c_k, \\
s & := \sum _{n \geq 0} a_n,
 & \qquad t & := \sum _{n \geq 0} b_n,
 & \qquad d_n & := B_n - t.
\end{alignat*}
Rearrangement of terms yields
\begin{align*}
C_n & = a\0 b\0 + (a\0 b_1 + a_1 b\0) + \cdots + ( a\0 b_n + a_1 b_{n-1} + \cdots + a_n b\0) \\
 & = a\0 B_n + a_1 B_{n-1} + \cdots + a_n B\0 \\
 & = a\0 ( t + d_n ) + a_1 ( t + d_{n-1} ) + \cdots + a_n ( t + d\0 ) \\
 & = A_n t + a\0 d_n +a_1 d_{n-1} + \cdots + a_n d\0.
\end{align*}
Put $h_n := a\0 d_n + a_1 d_{n-1} + \cdots + a_n d\0$. We have to show that
$C_n \to st.$ Since $A_n t \to st$, it suffices to show that $h_n \to 0$. So
let $\varepsilon > 0$ be given. Let $\alpha := \sum _{n \geq 0} | \,a_n \,|$,
which is finite by assumption. We have $d_n \to 0$ because $B_n \to t$.
Hence we can choose $N \geq 0$ such that $| \,d_n \,| \leq \varepsilon / (1 + \alpha)$
for all integers $n > N$, and for such $n$ we have
\begin{align*}
| \,h_n \,| & \leq | \,a_n d\0 + \cdots + a_{n-N} d_N \,|
+ | \,a_{n-(N+1)} d_{N+1} + \cdots + a\0 d_n \,| \\
 & \leq | \,a_n d\0 + \cdots + a_{n-N} d_N \,| + \varepsilon.
\end{align*}
Letting $n \to \infty$, we get $\limsup _{n \to \infty} | \,h_n \,| \leq \varepsilon$
by $a_n \to 0$. The statement follows because $\varepsilon > 0$ was arbitrary.
\pagebreak
\end{proof}

\begin{proposition}%
Let $x,y$ be two commuting elements of a ring with the same square.
If $x+y$ is not a left divisor of zero then $x=y$.
\end{proposition}

\begin{proof} One computes $(x+y)(x-y) = {x \,}^2-xy+yx-{y \,}^2 = 0$, so that
the assumption implies $x-y = 0$. \end{proof}

\begin{theorem}\label{sqroot}\index{square root|(}%
Let $A$ be a Banach algebra and let $a \in A$ with $\rlambda(a) < 1$. The series
\[ \sum\limits_{n=1}^{\infty} \,\bigl( _{\ \text{\small{$n$}}}^{1/2} \bigr) \,{a \,}^{n} \]
then converges absolutely to an element $b$ of $A$ with ${(e+b) \,}^2 = e+a$
and $\rlambda(b) < 1 $. Furthermore, there is no other element $c$ of $\tld{A}$
with ${(e+c) \,}^2 = e+a$ and $\rlambda(c) \leq 1$.
\end{theorem}

\begin{proof} If the above series converges absolutely, it follows by the \linebreak
Theorem of Mertens \ref{Mertens} concerning the Cauchy product
of series that ${(e+b) \,}^2 = e+a$. Assume that $\rlambda(a) < 1$. The series for
$b$ converges absolutely by comparison with the geometric series. It shall
next be shown that $\rlambda(b) < 1$. For $k \geq 1$, we have by \ref{rlpowers}
and \ref{commrlsub}:
\[ \rlambda \biggl( \,\sum\limits_{n=1}^k
\,\bigl( _{\ \text{\small{$n$}}}^{1/2} \bigr) \,{a \,}^n \biggr) \leq \sum\limits_{n=1}^k
\ \Bigl| \bigl( _{\ \text{\small{$n$}}}^{1/2} \bigr) \Bigr| \ {\rlambda ( a ) \,}^n \]
whence, by applying the last part of \ref{commrlsub} to the closed
subalgebra of $A$ generated by $a$, we get
\begin{align*}
\rlambda (b) \leq \sum\limits_{n=1}^{\infty}
\ \Bigl| \bigl( _{\ \text{\small{$n$}}}^{1/2} \bigr) \Bigr| \ {\rlambda(a) \,}^n
 & = -\sum\limits_{n=1}^{\infty} \ \bigl( _{\ \text{\small{$n$}}}^{1/2} \bigr)
        \,{\bigl(-\rlambda (a)\bigr) \,}^n \\
 & = 1-\sqrt{1-\rlambda (a)} < 1.
\end{align*}
Let $c$ now be an element of $\tld{A}$ with $\rlambda(c) \leq 1$ and
${(e+c) \,}^2 = e+a$. It shall be shown that $c = b$. The elements
$x := e+b, y := e+c$ have the same square. Furthermore, we have
$y(e+a) = {(e+c) \,}^3 = (e+a)y$, i.e.\ the element $y$ commutes with $a$. By
continuity, it follows that $y$ also commutes with $b$, hence with $x$.
Finally, $x+y = 2e+(b+c)$ is invertible as
$\rlambda(b+c) \leq  \rlambda(b) + \rlambda(c) < 2$, cf.\ \ref{commrlsub}.
It follows from the preceding proposition that $x = y$, or in other words,
that $b = c$. \pagebreak
\end{proof}

\begin{lemma}[Ford's Square Root Lemma]%
\index{Ford's Square Root Lemma}\index{Lemma!Ford's Square Root}%
\index{Theorem!Ford's Square Root Lemma}\label{Ford}%
If $a$ is a Hermitian \linebreak element of a Banach \st-algebra,
satisfying $\rlambda(a) < 1$, then the square root $e+b$ of $e+a$
(with $b$ as in the preceding theorem) is Hermitian. 
\end{lemma}

\begin{proof} This is obvious if the involution is continuous.
Now for the general case.
Consider the element $b$ of the preceding proposition.
We have to show that $b^* = b$. The element $b^*$ satisfies
${(e+b^*) \,}^2 = (e+a)^* = e+a$ as well as $\rlambda(b^*) = \rlambda(b) < 1$
by \ref{Brlsymm}. So $b = b^*$ by the uniqueness statement of the
preceding theorem.
\end{proof}

\begin{theorem}\label{possqroot}%
Let $A$ be a Banach algebra, and let $x$ be an element of $\tld{A}$ with
$\s(x) \subset\ ]0,\infty [$. (Please note: \ref{speczero}.) Then $x$ has a
unique square root $y$ in $\tld{A}$ with $\s(y) \subset\ ]0, \infty [$. Moreover
$x$ has no other square root with spectrum in $[0, \infty [$. The element $y$
belongs to the closed subalgebra of $\tld{A}$ generated by $x$. If $A$ is a
Banach \st-algebra, and if $x$ is Hermitian, so is $y$.
\end{theorem}

\begin{proof} We may assume that $\rlambda(x) \leq 1$. The
element $a := x-e$ then has $\s(a) \subset \ ]-1,0]$, whence
$\rlambda(a) < 1$. Consider now $b$ as in \ref{sqroot}. Then
$y := e + b$ is a square root of $e + a = x$, by \ref{sqroot}.
Furthermore, $\s(y) \subset \mathds{R} \setminus \{ 0 \}$ by the Rational
Spectral Mapping Theorem. Also $\rlambda(b) < 1$, which implies that
$\s(y) \subset\ ]0,\infty [$. Let now $z$ be another square root of $x$ with
$\s(z) \subset\ [0,\infty [$. Then $\s(z) \subset \,[0,1]$ by the Rational Spectral
Mapping Theorem. Thus, the element $c:= z-e$ \linebreak has $\s(c) \subset \,[-1,0]$,
so $\rlambda(c) \leq 1$. Since ${(e+c) \,}^2 = x = e+a$, theorem \ref{sqroot}
implies that $c = b$, whence $z = y$. Let now $B$ denote the closed
subalgebra of $\tld{A}$ generated by $x$. Clearly $y$ belongs to the closed
subalgebra of $\tld{A}$ generated by $x$ and $e$, which is
$\mathds{C}e + B$, cf.\ the proof of \ref{Banachunitis}. So let $y = \mu e + f$
with $f \in B$ and $\mu \in \mathds{C}$. Then $x = {\mu \,}^2 e + g$ with
$g \in B$. So if $\mu \neq 0$, then $e \in B$, whence $y \in B$. If $\mu = 0$,
then $y = f \in B$ as well. The last statement follows from the preceding lemma.
\end{proof}
\index{square root|)}

\medskip
For C*-algebras, a stronger result holds, cf.\ \ref{C*sqroot}.

\clearpage


\section{The Boundary of the Spectrum}%
\label{boundary}

We proved in \ref{specsubalg} that if $B$ is a subalgebra of an algebra
$A$, then
\[ \s _A (b) \setminus \{0\} \subset \s _B (b) \setminus \{0\} \]
holds for all $b \in B$. Also, if $B$ furthermore is \twiddle-unital in $A$, then
\[ \s _A (b) \subset \s _B (b) \]
holds for all $b \in B$. We set out to show a partial converse inclusion,
see \ref{spbdry} below.

\begin{proposition}\label{invbded}%
Let $A$ be a unital Banach algebra. Assume that $a \in A$ is the limit
of a sequence $\{ \,a_n \,\}$ of invertible elements of $A$. If the sequence
$\{ \,{a_n \,}^{-1} \,\}$ is bounded, then it converges in $A$ to an inverse
of $a$.
\end{proposition}

\begin{proof}
Let $c > 0$ with $| \,{a_n \,}^{-1} \,| \leq c$ for all $n$.
For all $m, n$ we then have
\[ | \,{a_m \,}^{-1} - {a_n \,}^{-1} \,| =
| \,{a_m \,}^{-1} \,( \,a_n - a_m \,) \,{a_n \,}^{-1} \,|
\leq c^{\,2} \,| \,a_n - a_m \,|. \]
This shows that $\{ \,{a_n \,}^{-1} \,\}$ is a Cauchy sequence,
and thus convergent to some limit $b \in A$. The element
$b$ is an inverse of $a$ by continuity of multiplication.
\end{proof}

\begin{definition}[topological divisors of zero]%
\index{topological divisors of zero}%
Let $A$ be a normed algebra. An element $a \in A$ is called a
\underline{left topological divisor of zero}, if there exists a sequence
$\{ \,s_n \,\}$ in the unit sphere of $A$, such that
\[ \lim _{n \to \infty} a \,s_n = 0. \]
Similarly, one defines a \underline{right topological divisor of zero}.
An element $a \in A$ is called a \underline{joint topological divisor of zero},
if there exists a sequence $\{ \,s_n \,\}$ in the unit sphere of $A$, such that
\[ \lim _{n \to \infty} a \,s_n = \lim _{n \to \infty} s_n \,a = 0. \]
\end{definition}

\begin{theorem}
Let $A$ be a unital Banach algebra. The boundary of the group of
invertible elements of $A$ then consists of joint topological divisors of zero.
\pagebreak
\end{theorem}

\begin{proof}
Let $a \in A$ belong to the boundary of the group of invertible elements
of $A$, and let $\{ \,a_n \,\}$ be a sequence of invertible elements of $A$
converging to $a$. Since the group of invertible elements of $A$ is open
by \ref{GL(A)}, the element $a$ cannot be invertible. The preceding
proposition \ref{invbded} implies that the sequence $\{ \,{a_n\,}^{-1} \,\}$
must be unbounded. By going over to a subsequence, we can assume
that $| \,{a_n\,}^{-1} \,| \geq n$ for all $n \geq 1$. The elements
$s_n := {| \, {a_n \,}^{-1} \,|\,}^{-1} \,{a_n\,}^{-1}$
belong to the unit sphere of $A$, and
\[ a \,s_n = ( \,a-a_n \,) \,s_n + a_n \,s_n
= ( \,a - a_n \,) \,s_n + {| \,{a_n \,}^{-1} \,| \,}^{-1} \,e \ \to \ 0. \]
Thus $\lim _{n \to \infty} a \,s_n = 0$, and similarly $\lim _{n \to \infty} s_n \,a = 0$.
\end{proof}

\begin{theorem}
Let $B$ be a subalgebra of a unital normed algebra $A$. Then a left
topological divisor of zero in $B$ is not left invertible in $A$. A similar
result holds with ``right'' instead of ``left''.
\end{theorem}

\begin{proof}
Assume for example that $b$ is a left topological divisor of zero in $B$, and
that $b$ is left invertible in $A$. There then exists a sequence $\{ \,s_n \,\}$
in the unit sphere of $B$, such that $\lim _{n \to \infty} b \,s_n = 0$, and there
exists an element $a$ of $A$ with $a \,b = e$. We arrive at the contradiction
\[ 1 = | \,s_n \,| = | \, e \,s_n \,| = | \,a \,b \,s_n \,| \leq | \,a \,| \cdot | \,b \,s_n \,|
\ \to \ 0. \qedhere \]
\end{proof}

\begin{corollary}\label{bdrynotinv}%
Let $B$ be a complete unital subalgebra of a unital normed algebra $A$.
Then the elements of the boundary of the group of invertible elements of
$B$ are neither left nor right invertible in $A$.
\end{corollary}

\begin{theorem}[\v{S}ilov]\label{spbdry}%
\index{Silov@\v{S}ilov}\index{Theorem!Silov@\v{S}ilov}%
Let $B$ be a complete subalgebra of a normed algebra $A$,
and let $b \in B$. Then
\[ \partial \,\s _B (b) \subset \partial \,\s _A (b), \]
where $\partial$ denotes ``the boundary of'' (relative to the complex plane).
\end{theorem}

\begin{proof}
Assume first that $B$ is \twiddle-unital in $A$. Then $\tld{B}$ is a unital
\linebreak
subalgebra of $\tld{A}$, cf.\ \ref{twidunitp}. For our purpose it is allowed
to equip $\tld{B}$ with the norm inherited from $\tld{A}$. Then $\tld{B}$ is
complete, cf.\ the proof of \ref{Banachunitis}. Let $\lambda \in \partial \,\s _B (b)$
and let $(\lambda _n)$ be a sequence in $\mathds{C} \setminus \s _B (b)$
converging to $\lambda$. Then each $c_n := \lambda _n e - b$ is invertible
in $\tld{B}$, while $c := \lambda e - b$ is not. That is, $c$ is in the boundary
of the group of invertible elements of $\tld{B}$. So the preceding corollary
\ref{bdrynotinv} yields that $c$ is not invertible in $\tld{A}$. This says that
$\lambda \in \s _A (b)$, and hence $\lambda \in \partial \,\s _{A} (b)$, by
\ref{specsubalg}. This settles the \twiddle-unital case. \pagebreak

Assume next that $B$ is not \twiddle-unital in $A$. Then $B$ must have
a unit, $p$ say, which is not a unit in $A$, cf.\ \ref{twiddleunital}. Let $e$
denote the unit of $\tld{A}$, and consider $C := \mathds{C} e+B$, which
is a complete subalgebra of $\tld{A}$, cf.\ the proof of \ref{Banachunitis}.
Please note that $C$ is unital in $\tld{A}$, and so \twiddle-unital in $\tld{A}$.
Let now $\lambda \in \partial \,\s _B (b)$, and consider $d := \lambda e - b \in C$.
We have $\lambda \in \s _C (b)$. Indeed, suppose $d$ were invertible
in $C$. Let $g \in C$ be an inverse of $d$ in $C$. Then $gd = dg = e$.
Multiplication from the left and from the right with $p$ gives
$pgdp = pdgp = {p \,}^2$, which is the same as
$pg(\lambda p - b) = (\lambda p - b)gp = p$ (as $p$ is the unit in $B$),
so that $pg = gp \in B$ would be an inverse of $\lambda p - b$ in $B$,
cf.\ \ref{leftrightinv}. This contradiction shows that $\lambda \in \s _C (b)$,
but then also $\lambda \in \partial \,\s _C (b)$. Indeed, for $\lambda \neq 0$,
this follows from \ref{specsubalg}. If $\lambda = 0$, it follows from
\ref{specsubalg} and the fact that then $0 \in \s _B (b)$ and $0 \in \s _C (b)$.
By the preceding paragraph, we now get $\lambda \in \partial \,\s _A (b)$,
and the proof is complete.
\end{proof}

\medskip
This result has two corollaries guaranteeing equality of spectra \linebreak
under certain conditions. The first condition is on the spectrum in the
subalgebra, the second one is on the spectrum in the whole algebra.

\begin{corollary}\label{nointerior}%
Let $B$ be a complete subalgebra of a normed \linebreak algebra $A$,
and let $b \in B$ such that $\s _B (b)$ has no interior (e.g.\ is real). Then
\[ \s _B (b) \setminus \{0\} = \s _A (b) \setminus \{0\}. \]
If $B$ furthermore is \twiddle-unital in $A$, then
\[ \s _B (b) = \s _A (b). \]
\end{corollary}

\begin{proof}
This follows from the preceding theorem and \ref{specsubalg}.
\end{proof}

\begin{corollary}%
Let $B$ be a \twiddle-unital closed subalgebra of a Banach
algebra $A$, and let $b \in B$ such that $\s _A (b)$ does not
separate the complex plane (e.g.\ is real). Then
\[ \s _B (b) = \s _A (b). \]
\end{corollary}

\begin{proof}
It suffices to prove that $\s _B (b) \subset \s _A (b)$, cf.\ \ref{specsubalg}.
So suppose that there exists $\lambda \in \s _B (b) \setminus \s _A (b)$.
Since $\s _A (b)$ does not separate the complex plane, there exists a
continuous path disjoint from $\s _A (b)$, that connects $\lambda$ with
the point at infinity. It is easily seen that this path meets the boundary of
$\s _B (b)$ in at least one point, which then also must belong to $\s _A (b)$,
by \ref{spbdry}, a contradiction. \pagebreak
\end{proof}

\clearpage


\section{Hermitian Banach \texorpdfstring{$*$-}{\052\055}Algebras}\label{secHerm}

\begin{definition}[Hermitian \protect\st-algebras]%
\index{Hermitian!s-algebra@\protect\st-algebra}%
\index{algebra!Banach s-algebra@Banach \protect\st-algebra!Hermitian}
\index{algebra!s-algebra@\protect\st-algebra!Hermitian}
A \st-algebra is said to be \linebreak \underline{Hermitian} if each of its
Hermitian elements has real spectrum.
\end{definition}

\begin{proposition}%
A \st-algebra $A$ is Hermitian if and only if its \linebreak unitisation $\tld{A}$ is.
\end{proposition}

\begin{proof}
For $a \in A\sa$ and $\lambda \in \mathds{R}$, we have
\[ \s_{\tld{A}} \,(\lambda e+a) = \lambda + \s_{\tld{A}} \,(a) = \lambda + \s_A \,(a), \]
cf.\ \ref{specunit}.
\end{proof}

\begin{theorem}\label{fundHerm}%
For a Banach \st-algebra $A$ the following statements are equivalent.
\begin{itemize}
   \item[$(i)$] \ $A$ is Hermitian,
  \item[$(ii)$] \ $\iu \notin \s(a)$ whenever $a \in A$ is Hermitian,
 \item[$(iii)$] \ $\rlambda(a) \leq \rsigma(a)$ for all $a \in A$,
 \item[$(iv)$] \ $\rlambda(b) = \rsigma(b)$ for all normal $b \in A$,
  \item[$(v)$] \ $\rlambda(b) \leq {|\,b^*b\,| \,}^{1/2}$ for all normal $b \in A$.
\end{itemize}
\end{theorem}

\begin{proof} (i) $\Rightarrow$ (ii) is trivial.

(ii) $\Rightarrow$ (iii): Assume that (ii) holds and that
$\rlambda(a) > \rsigma(a)$ for some $a \in A.$ Upon replacing $a$
with a suitable multiple, we can assume that $1 \in \s(a)$ and
$\rsigma(a) < 1$. There then exists a Hermitian element $b = b^*$ of
$\tld{A}$ such that ${b\,}^2 = e-a^*a$ (Ford's Square Root Lemma \ref{Ford}).
Since $b$ is invertible (Rational Spectral Mapping Theorem), it follows that
\[(e+a^*) \,(e-a) = {b\,}^2 + a^* -a = b \,\bigl( e + b^{-1} \,(a^* - a) \,b^{-1} \bigr) \,b. \]
Now the element $-\iu \,b^{-1} \,(a^*-a) \,b^{-1}$ is a Hermitian element of $A$,
so that $\iu \,\bigl( e+b^{-1} \,(a^*-a) \,b^{-1} \bigr)$ is invertible by assumption.
This implies that $e-a$ has a left inverse. One can apply a similar argument to
\[ (e-a) \,(e+a^*) = c \,\bigl( e+c^{-1} \,(a^*-a) \,c^{-1} \bigr) \,c, \]
where $c$ is a Hermitian element of $\tld{A}$ with ${c\,}^2 = e-aa^*$
(by $\rsigma(a^*) = \rsigma(a) < 1$, cf.\ \ref{rsC*}). Thus $e-a$ also has a
right inverse, so that $1$ would not be in the spectrum of $a$ by
\ref{leftrightinv}, a contradiction.

(iii) $\Rightarrow$ (iv) follows from the fact that
$\rsigma(b) \leq \rlambda(b)$ for all normal
$b \in A$, cf.\ \ref{Bnormalrsleqrl}.

(iv) $\Rightarrow$ (v) is trivial. 

(v) $\Rightarrow$ (i): Assume that (v) holds. Let $a$ be a Hermitian element
of $A$ and let $\lambda \in \s(a)$. We have to show that $\lambda$ is real.
For this purpose, we can assume that $\lambda$ is non-zero. Let $\mu$ be an
arbitrary real number, and let $n$ be an integer $\geq 1$. Define
$u := \lambda^{-1}a$ and $b := {(a + \iu \mu e) \,}^n \,u$. Then
${(\lambda + \iu \mu) \,}^n$ is in the spectrum of $b$
(Rational Spectral Mapping Theorem). This implies
\[ {|\,\lambda + \iu \mu \,|\,}^{2n} \leq {\rlambda( \,b \,) \,}^2 \leq |\,b^*b\,|
\leq \bigl| \,{ \bigl( \,{a \,}^2 + {\mu \,}^2 \,e \,\bigr) \,}^n \,\bigr| \cdot |\,u^*u\,|. \]
It follows that
\[ {| \,\lambda + \iu \mu \,|\,}^2
\leq {\bigl| \,{ \bigl( \,{a \,}^2 + {\mu \,}^2 \,e \,\bigr) \,}^n \,\bigr| \,}^{1/n}
\cdot {| \,u^*u \,| \,}^{1/n} \]
for all integers $n \geq 1$, so that by \ref{commrlsub}, we have
\[ {| \,\lambda + \iu \mu \,|\,}^2 \leq \rlambda \bigl( \,{a \,}^{2} + {\mu \,}^{2} \,e \,\bigr)
\leq \rlambda \bigl( \,{a \,}^2 \,\bigr)+ {\mu \,}^2. \]
Decomposing $\lambda = \alpha + \iu \beta $ with $\alpha $ and $\beta $ real,
one finds
\[ {\alpha \,}^2 + {\beta \,}^2 + 2 \beta \mu \leq \rlambda \bigl( \,{a \,}^2 \,\bigr) \]
for all real numbers $\mu$, so that $\beta$ must be zero, i.e.\ $\lambda$
must be real.
\end{proof}

\medskip
We now have the fundamental result:

\begin{theorem}\label{C*Herm}%
A C*-algebra is Hermitian.
\end{theorem}

\begin{proof}
For an element $a$ of a C*-algebra we have by \ref{rldef} and \ref{C*rs}:
\[ \rlambda(a) \leq \| \,a \,\| = \rsigma(a). \qedhere \]
\end{proof}

\begin{theorem}\label{Herminher}%
A closed \st-subalgebra of a Hermitian Banach
\linebreak \st-algebra is Hermitian as well.
\end{theorem}

\begin{proof}
One uses the fact that a Banach \st-algebra $C$ is
Hermitian if and only if for every $a \in C$ one has
$\rlambda(a) \leq \rsigma(a)$.
\end{proof}

\begin{lemma}\label{lspecHermincl}%
Let $B$ be a complete Hermitian unital \st-subalgebra of a unital
normed \st-algebra $A$. If an element of $B$ is invertible in $A$,
then it is invertible in $B$.
\end{lemma}

\begin{proof}
Let $b \in B$ be invertible in $A$ $\bigl($with inverse ${b}^{-1}$ in
$A\,\bigr)$. The element ${b}^{*}b$ then also is invertible in $A$
$\bigl($with inverse ${b}^{-1}{({b}^{-1})}^{*}$ in \pagebreak
$A\bigr)$.\linebreak Since $B$ is Hermitian, the Hermitian element
${b}^{*}b$ is the limit of the invertible elements
$\frac1n \iu e + {b}^{*}b$ $(n \geq 1)$. If ${b}^{*}b$ were not invertible
in $B$, it would follow that ${b}^{*}b$ were in the boundary of the
group of invertible elements of $B$, and thus not invertible in $A$,
cf. \ref{bdrynotinv}. Thus ${b}^{*}b$ has an inverse ${({b}^{*}b)}^{-1}$
in $B$. Then ${({b}^{*}b)}^{-1}{b}^{*}$ is a left inverse of $b$ in $B$.
Similarly, $b$ has a right inverse in $B$, and so $b$ is invertible in $B$,
by \ref{leftrightinv}.
\end{proof}

\begin{theorem}\label{specHermincl}%
Let $B$ be a complete Hermitian \st-subalgebra
of a \linebreak normed \st-algebra $A$.
For an element $b$ of $B$, we then have
\[ \s_B (b) \subset \s_A (b). \]
\end{theorem}

\begin{proof}
The proof of the \twiddle-unital case follows from the above lemma,
cf.\ \ref{twidunitp}. (When $\tld{B}$ is equipped with the norm inherited
from $\tld{A}$, then it is complete, cf.\ the proof of \ref{Banachunitis}.)
Assume next that $B$ is not \twiddle-unital in $A$. Then $B$ contains
a unit $p$, say, which is not a unit in $\tld{A}$, cf.\ \ref{twiddleunital}.
Let $e$ denote the unit in $\tld{A}$, and consider $C := \mathds{C} e + B$,
which is a complete \st-subalgebra of $\tld{A}$, cf.\ the proof of
\ref{Banachunitis}. It can be seen that $C$ is Hermitian by using
\ref{specsubalg}. Please note that $C$ is unital in $\tld{A}$. Assume that
$\lambda \notin \s_A (b)$ for some $b \in B$ and $\lambda \in \mathds{C}$.
The element $d := \lambda e - b \in C$ then is invertible in $\tld{A}$.
It follows that $d$ is invertible in $C$, by the preceding lemma.
So let $g \in C$ be an inverse of $d$ in $C$. We then have $gd = dg = e$.
Multiplication from the left and from the right with $p$ yields
$pgdp = pdgp = p^2$, which is the same as
$pg(\lambda p - b) = (\lambda p - b)gp = p$ (as $p$ is the unit in $B$), so that
$pg = gp \in B$ is an inverse of $\lambda p - b$ in $B$, cf.\ \ref{leftrightinv}.
This says that $\lambda \notin \s_B (b)$, as was to be shown. Please note that
the present proof closely parallels the second paragraph of the proof of \ref{spbdry}.
\end{proof}

\medskip
From \ref{specsubalg} it follows now

\begin{theorem}\label{specHermsubalg}%
Let $B$ be a complete Hermitian \st-subalgebra
of a \linebreak normed \st-algebra $A$.
For an element $b$ of $B$ we then have
\[ \s_A (b) \setminus \{ 0 \} = \s_B (b) \setminus \{ 0 \}. \]
If furthermore $B$ is \twiddle-unital in $A$, then
\[ \s_A (b) = \s_B (b). \]
\end{theorem}

It is clear that \ref{lspecHermincl}, \ref{specHermincl}, and \ref{specHermsubalg}
are most of the time used in conjunction with \ref{Herminher}. \pagebreak

Now for unitary elements.

\begin{definition}[unitary elements]\index{unitary element}%
An element $u$ of a unital \linebreak \st-algebra is called
\underline{unitary} if $u^*u = uu^* = e$. The unitary elements
of a unital \st-algebra form a group under multiplication. 
\end{definition}

\begin{proposition}[Cayley transform]\index{Cayley transform}%
Let $A$ be a Banach \st-al\-gebra and let $a$ be a Hermitian
element of $\tld{A}$. Let $\mu > \rlambda(a)$. Then
\[ u := (a - \iu\mu e){(a + \iu\mu e)}^{-1} \]
is a unitary element of $\tld{A}$. One says that $u$ is a
\underline{Cayley transform} of $a$.
\end{proposition}

\begin{proof}
Please note first that $a + \iu\mu e$ is invertible as $\mu > \rlambda(a)$
by assumption, so $u$ is well defined. By commutativity, it follows that
\begin{align*}
u^* & = \bigl[{(a + \iu\mu e)}^{-1}\bigr]^*(a - \iu\mu e)^* \\
       & = {(a - \iu\mu e)}^{-1}(a + \iu\mu e) \\
       & = (a + \iu\mu e){(a - \iu\mu e)}^{-1},
\end{align*}
so that $u^*u = uu^* = e$.
\end{proof}

\begin{theorem}\label{specunitary}%
A Banach \st-algebra $A$ is Hermitian if and only if the spectrum of every
unitary element of $\tld{A}$ is contained in the unit circle.
\end{theorem}

\begin{proof}
Assume that $A$ is Hermitian and let $u$ be a unitary element of $\tld{A}$.
We then have $\rlambda(u) = \rsigma(u) = 1$ by \ref{fundHerm}
(i) $\Rightarrow$ (iv). It follows that the spectrum of $u$ is contained in the
unit disk. Since also $u^{-1} = u^*$ is unitary, it follows from
${\s(u)}^{-1} = \s(u^{-1})$ (Rational Spectral Mapping Theorem) that also
${\s(u)}^{-1}$ is contained in the unit disk. This makes that $\s(u)$ is
contained in the unit circle.

Conversely, let $a = a^*$ be a Hermitian element of $A$ and let
$\mu > \rlambda(a)$. Then the Cayley transform
$u := (a - \iu\mu e){(a + \iu\mu e)}^{-1}$ is a unitary element of $\tld{A}$.
By the Rational Spectral Mapping Theorem, we have $\s(u) = r\bigl(\s(a)\bigr)$
with $r$ denoting the Moebius transformation
$r(z) = (z - \iu\mu){(z + \iu\mu)}^{-1}$. Thus, if the spectrum of $u$ is
contained in the unit circle, it follows that $\s(a) \subset \mathds{R}$.
\end{proof}

\begin{corollary}\label{C*unitary}
The spectrum of a unitary element of a unital \linebreak
C*-algebra is contained in the unit circle. \pagebreak
\end{corollary}

\clearpage


\section{An Operational Calculus}

We next give a version of the operational calculus with an elementary
proof avoiding Zorn's Lemma. (See also \ref{opcalc}.)

\begin{theorem}[the operational calculus, weak form]%
\index{operational calculus}\label{weakopcalc}%
Let $A$ be a \linebreak C*-algebra, and let $a = a^*$ be a Hermitian
element of $A$. Denote by \linebreak $C\bigl(\s(a)\bigr)\0$ the C*-subalgebra
of $C\bigl(\s(a)\bigr)$ consisting of the continuous complex-valued functions
on $\s(a)$ vanishing at $0$. There is a unique \st-algebra homomorphism
\[ C\bigl(\s(a)\bigr)\0 \to A \]
which maps the identity function on $\s(a)$ to $a$. This mapping is an
isomorphism of C*-algebras from $C\bigl(\s(a)\bigr)\0$ onto the closed
subalgebra of $A$ generated by $a$. (Which also is the C*-subalgebra
of $A$ generated by $a$.) This mapping is called the
\underline{operational calculus} for $a$, and it is denoted by $f \mapsto f(a)$.
\end{theorem}

\begin{proof}
Let $\mathds{C}[x]\0$ denote the set of complex polynomials without
\linebreak constant term. For $p \in \mathds{C}[x]\0$, we shall denote by
$p|_{\s(a)} \in C\bigl(\s(a)\bigr)\0$ the restriction to $\s(a)$ of the corresponding
polynomial function. We can define an isometric mapping
$p|_{\text{\small{$\s(a)$}}} \mapsto p(a)$
$(p \in \mathds{C}[x]\0)$. Indeed, for $p \in \mathds{C}[x]\0$,
the element $p(a) \in A$ is normal, and one calculates
\begin{align*}
\| \,p(a) \,\| & = \rlambda\,\bigl( \,p(a) \,\bigr) \tag*{by \ref{C*rl}} \\
 & = \max \ \bigl\{ \,| \,\lambda \,| : \lambda \in \s \,\bigl( \,p(a) \,\bigr) \,\bigr\}
\tag*{by \ref{specradform}} \\
 & = \max \ \bigl\{ \,| \,p(\mu) \,| : \mu \in \s(a) \,\bigr\} \tag*{by \ref{ratspecmapthm}} \\
 & = \bigl| \,p|_{\text{\small{$\s(a)$}}} \,\bigr| _{\infty}.
\end{align*}
Thus, if $p|_{\text{\small{$\s(a)$}}} = q|_{\text{\small{$\s(a)$}}}$ for
$p,q \in \mathds{C}[x]\0$, then
$\bigl| \,(p-q)|_{\text{\small{$\s(a)$}}} \,\bigr| _{\infty} = 0$,
whence $\| \,(p-q) (a) \,\| = 0$, so $p(a) = q(a)$, and the mapping in
question is well defined. The mapping in question is isometric as well.

The Theorem of Weierstrass implies that the set of complex \linebreak
polynomial functions vanishing at $0$ is dense in $C\bigl(\s(a)\bigr)\0$.
This \linebreak together with the continuity result \ref{Cstarcontr} already
implies the uniqueness statement. The above mapping has a unique
extension to a continuous mapping from $C\bigl(\s(a)\bigr)\0$ to $A$. This
continuation then is a C*-algebra isomorphism from $C\bigl(\s(a)\bigr)\0$
onto the closed subalgebra of $A$ generated by $a$. It is used that the
continuation is isometric, which implies that its image is complete, hence
closed in $A$. \pagebreak
\end{proof}

\begin{corollary}\label{calcplus}%
Let $a = a^*$ be a Hermitian element of a \linebreak C*-algebra $A$.
If $f \in C\bigl(\s(a)\bigr)\0$ satisfies $f \geq 0$ on $\s(a)$ then $f(a)$
is Hermitian and $\s\bigl(f(a)\bigr) \subset [0,\infty[$.
\end{corollary}

\begin{proof}
Let $g \in C\bigl(\s(a)\bigr)\0$ with $f = {g \,}^2$ and $g \geq 0$ on $\s(a)$.
Then $g$ \linebreak is Hermitian in $C\bigl(\s(a)\bigr)\0$ as $g$ is real-valued.
Hence $g(a)$ is Hermitian, \linebreak and thus
$\s \bigl(g(a)\bigr) \subset \mathds{R}$. So $f(a) = {g(a) \,}^2$ is Hermitian
and $\s \bigl(f(a)\bigr) = \s {\bigl(g(a)\bigr) \,}^2 \subset [0,\infty[$ by the Rational
Spectral Mapping Theorem.
\end{proof}

\begin{theorem}\index{square root}\label{C*sqroot}%
Let $A$ be a C*-algebra. Let $a$ be a Hermitian element of $A$ with
$\s(a) \subset [0,\infty[$. Then $a$ has a unique Hermitian square root
${a \,}^{1/2}$ in $A$ with $\s \,\bigl( \,{a \,}^{1/2} \,\bigr) \subset [0,\infty[$. The
element ${a \,}^{1/2}$ belongs to the closed subalgebra of $A$ generated
by $a$. (Which also is the C*-subalgebra of $A$ generated by $a$.)
See also \ref{C*sqrootrestate}.
\end{theorem}

\begin{proof}
For existence, consider the function $f$ given by $f(x) := \sqrt{x}$ for
$x \in \s(a) \subset [0, \infty[$, and put ${a \,}^{1/2} := f(a)$. Now for
uniqueness. Let $b$ be a Hermitian square root of $a$ with
$\s(b) \subset [0,\infty[$. We have to show that $b = {a \,}^{1/2}$. Let $c$ be a
Hermitian square root of ${a \,}^{1/2}$ in the closed subalgebra of $A$
generated by ${a \,}^{1/2}$. Let $d$ be a Hermitian square root of $b$. Then
$d$ commutes with $b = {d \,}^2$, hence also with $a = {b \,}^2$. But then $d$
also commutes with ${a \,}^{1/2}$ and $c$. This is so because ${a \,}^{1/2}$ lies
in the closed subalgebra of $A$ generated by $a$ (by construction of
${a \,}^{1/2}$), and because $c$ is in the closed subalgebra of $A$ generated
by ${a \,}^{1/2}$ (by assumption on $c$). We conclude that the Hermitian
elements $a, {a \,}^{1/2}, b, c, d$ all commute. One now calculates
\begin{align*}
0 = & \,\bigl( a - {b \,}^2 \bigr) \,\bigl( {a \,}^{1/2} - b \bigr) \\
  = & \,\bigl( {a \,}^{1/2} - b \bigr) \,{a \,}^{1/2}
         \,\bigl( {a \,}^{1/2} - b \bigr) + \bigl( {a \,}^{1/2}  -b \bigr) \,b \,\bigl( {a \,}^{1/2} - b \bigr) \\
  = & \,{\bigl[ \,\bigl( {a \,}^{1/2} - b \bigr) \,c \,\bigr] \,}^2+{\bigl[ \,\bigl( {a \,}^{1/2} - b \bigr) \,d \,\bigr] \,}^2.
\end{align*}
Both terms in the last line have non-negative spectrum. Indeed, both terms
are squares of Hermitian elements. Please note \ref{Hermprod} here.
Since the two terms are opposites one of another, both must have spectrum
$\{0\}$, and thence must both be zero, cf.\ \ref{C*rl}. The difference of these
terms is
\[ 0 = \bigl( {a \,}^{1/2} - b \bigr) \,{a \,}^{1/2} \,\bigl({a \,}^{1/2} - b \bigr)
- \bigl( {a \,}^{1/2} - b \bigr) \,b \,\bigl( {a \,}^{1/2} - b \bigr)
= { \bigl( {a \,}^{1/2} - b \bigr) \,}^3. \]
With the C*-property \ref{preC*alg} it follows that
\[ 0 = \| \,{\bigl( {a \,}^{1/2} - b \bigr) \,}^4 \,\| = {\| \,{ \bigl( {a \,}^{1/2} - b \bigr) \,}^2 \,\|\,}^2
= {\| \,{a \,}^{1/2} - b \,\|\,}^4. \pagebreak \qedhere \]
\end{proof}

\clearpage


\section{Odds and Ends: Questions of Imbedding}%
\label{questimbed}

This paragraph is concerned with imbedding for example
a normal element of a normed \st-algebra in a closed
commutative \st-subalgebra. The problem is posed by
possible discontinuity of the involution.

\begin{definition}[\st-stable (or self-adjoint) and normal subsets]%
\label{selfadjointsubset}\index{normal!subset}\index{self-adjoint!subset}%
\index{s054@\protect\st-stable}\index{subset!self-adjoint}%
\index{subset!s-stable@\protect\st-stable}\index{subset!normal}%
Let $S$ be a subset of a \st-algebra. One says that $S$ is
\underline{\st-stable} or \underline{self-adjoint} if with each
element $a$, it also contains the adjoint $a^*$. One says
that $S$ is \underline{normal}, if it is \st-stable and if its elements
commute pairwise.
\end{definition}

\begin{definition}[the commutant]\index{commutant}%
Let $S$ be a subset of an algebra $A$.
One denotes by $\underline{S'}$ the
\underline{commutant of $S$ in $A$}, which
is defined as the set of those elements of $A$
which commute with every element in $S$.
\end{definition}

The reader will easily prove the next two propositions.

\begin{proposition}\label{commutantp1}%
Let $S$ be a subset of an algebra $A$.
The commutant $S'$ enjoys the following properties:
\begin{itemize}
   \item[$(i)$] $S'$ is a subalgebra of $A$,
  \item[$(ii)$] if $e$ is a unit in $A$, then $e \in S'$,
 \item[$(iii)$] if $A$ is a normed algebra then $S'$ is closed,
 \item[$(iv)$] if $A$ is a \st-algebra, and if $S$ is \st-stable,
                       then $S'$ is a \st-subalgebra of $A$.
\end{itemize}
\end{proposition}

\begin{proposition}\label{commutantp2}%
Let $S,T$ be subsets of an algebra. We then have
\begin{itemize}
   \item[$(i)$] $S \subset T \Rightarrow T' \subset S'$,
  \item[$(ii)$] $S \subset T \Rightarrow S'' \subset T''$,
 \item[$(iii)$] $S \subset S''$, and $S'= S'''$,
 \item[$(iv)$] the elements of $S$ commute $\Leftrightarrow S \subset S'$.
\end{itemize}
\end{proposition}

\begin{proof}
The second statement of (iii): from $S \subset S''$ it follows
by (i) that $(S'')' \subset S'$. Also $S' \subset (S')''$.
\end{proof}

\begin{definition}[the second commutant]%
\index{commutant!second}\label{scndcomm}%
Let $S$ be a subset of an algebra $A$. The subset $S''$ is called the
\underline{second commutant} of $S$. It is a subalgebra of $A$ containing
$S$. It follows from the preceding proposition that $S''$ is commutative if
the elements in $S$ commute. \pagebreak
\end{definition}

\begin{proposition}\label{inv2comm}%
If $a$ is an invertible element of a unital algebra $A$, then its inverse
${a\,}^{-1}$ belongs to the second commutant $\{ \,a \,\}''$ in $A$.
\end{proposition}

\begin{proof}
See \ref{comm2}.
\end{proof}

\medskip
The following result is sometimes useful.

\begin{proposition}\label{scndcommspec}%
Let $S$ be a subset of an algebra $A$. Let $B$ denote the second
commutant of $S$ in $\tld{A}$. For $b \in B \cap A$ we then have
\[  \s_{B}(b) = \s_{A}(b). \]
\end{proposition}

\begin{proof}
Let $b \in B \cap A$. We know that $\s_A(b) = \s_{\tld{A}}(b)$,
cf.\ \ref{specunit}. We also have $\s_{\tld{A}}(b) \subset \s_B(b)$,
cf.\ \ref{specsubalg}. So assume that (for some $\lambda \in \mathds{C}$)
the element $\lambda e - b$ is invertible in $\tld{A}$ with inverse
$d \in \tld{A}$. It is enough to prove that $d \in B$. Since $b \in B$,
we have $d \in \{ \,b \,\}'' \subset B'' = B$, cf.\ \ref{inv2comm} and
\ref{commutantp2} (ii) \& (iii), as $B$ is a (second) commutant.
\end{proof}

\medskip
The preceding item is often used in the following form:

\begin{theorem}[Civin \& Yood]\index{Civin}\label{CivinYood}%
\index{Theorem!Civin \& Yood}\index{Yood}%
Let $A$ be a normed algebra, and let $S$ be a subset consisting of
commuting elements of $A$. There then exists a commutative unital
closed subalgebra $B$ of $\tld{A}$ containing $S$ such that for all
$b \in B \cap A$ one has $\s_{B}(b) = \s_{A}(b)$. If $A$ is a normed
\st-algebra, and if $S$ is normal, then $B$ can be chosen to be a
\st-subalgebra.
\end{theorem}

\begin{proof}
\ref{scndcomm}, \ref{commutantp1} (ii) \& (iii) \& (iv), and \ref{scndcommspec}.
\end{proof}

\medskip
This often permits one to reduce proofs to the case of commutative unital
Banach algebras. See for example \ref{specsub} and \ref{instability} below.

\begin{proposition}\label{normalscndcomm}%
Let $A$ be a normed \st-algebra, and let $S$ be a normal subset
of $A$. The second commutant of $S$ in $A$ then is a closed
commutative \st-subalgebra of $A$ containing $S$.
\end{proposition}

\begin{proof}
\ref{scndcomm} and \ref{commutantp1} (iii) \& (iv).
\end{proof}

\begin{definition}%
Let $S$ be a subset of a normed \st-algebra $A$.
The \underline{closed \st-subalgebra of $A$ generated by $S$}
is defined as the intersection of all closed \st-subalgebras of $A$
containing $S$. By the preceding proposition, it is commutative if
$S$ is normal. \pagebreak
\end{definition}

Please note that when the involution in $A$ is continuous, then the
closed \st-subalgebra of $A$ generated by a subset $S$ is the closure
of the \st-subalgebra of $A$ generated by $S$. This holds in particular
for \linebreak C*-algebras, and we shall make tacit use of this fact in the
sequel.

\begin{remark}\label{Zorn}
So far we did not use the full version of the axiom of choice. (We
did use the countable version of the axiom of choice, however.
Namely in using the fact that if $S$ is a subset of a metric space,
then a point in the closure of $S$ is the limit of a sequence in $S$.)
The next chapter \ref{TheGelfandTransformation} on the Gelfand
transformation will make crucial use of Zorn's Lemma, which is
equivalent to the full version of the axiom of choice. The ensuing
chapter \ref{PositiveElements} on positive elements, however, is
completely independent of chapter \ref{TheGelfandTransformation},
and does not require Zorn's Lemma. From a logical point of view,
one might well skip the following chapter \ref{TheGelfandTransformation}
if the priority is to use Zorn's Lemma as late as possible.
\end{remark}

\clearpage


\chapter{The Gelfand Transformation}%
\label{TheGelfandTransformation}


\setcounter{section}{13}

\section{Multiplicative Linear Functionals}

\begin{definition}[multiplicative linear functionals, $\Delta(A)$]%
\index{multiplicative!linear functional}%
\index{D(A)@$\Delta(A)$}\index{functional!multiplicative}%
If $A$ is an algebra, then a \underline{multiplicative linear functional}
on $A$ is a \underline{non-zero} linear functional $\tau$ on $A$
such that $\tau (ab) = \tau (a) \,\tau (b)$ for all $a,b \in A$. The set of
multiplicative linear functionals on $A$ is denoted by $\underline{\Delta (A)}$.
Multiplicative linear functionals are also called
\underline{complex homomorphisms}.

Please note that if $A$ is unital with unit $e$, then $\tau (e) = 1$ for all
$\tau \in \Delta (A)$ by $\tau (e) \,\tau (a) = \tau (a)$ for some $a \in A$ with
$\tau (a) \neq 0$.
\end{definition}

\begin{proposition}\label{mlfbounded}%
A multiplicative linear functional $\tau$ on a Banach
algebra $A$ is continuous with $| \,\tau \,| \leq 1$.
\end{proposition}

\begin{proof} Otherwise there would exist $a \in A$ with
$| \,a \,| < 1$ and $\tau (a) = 1$. With $b := \sum _{n=1}^{\infty } a^{n}$,
we would then get $a+ab = b$, and so $\tau (b) = \tau (a+ab)
= \tau (a) + \tau (a) \tau (b) = 1+ \tau (b)$, a contradiction.
\end{proof}

\begin{definition}[the Gelfand transform of an element]%
\index{Gelfand!transform}\index{a07@$\protect\wht{a}$}%
Let $A$ be an algebra. For $a \in A$, one defines
\[ \wht{a}(\tau) := \tau(a) \qquad \bigl(\tau \in \Delta (A)\bigr). \]
The \underline{Gelfand transform of $a$} is the function
\begin{align*}
\wht{a} : \Delta (A) & \to      \mathds{C}  \\
                 \tau    & \mapsto  \wht{a}(\tau).
\end{align*}
\end{definition}

\begin{proposition}[the Gelfand transformation]%
\index{Gelfand!transformation}\label{GelfTransf}%
Let $A$ be a Banach algebra with $\Delta(A) \neq \varnothing$.
The map
\begin{align*}
A & \to     {\ell\,}^{\infty} \bigl(\Delta (A)\bigr) \\
a & \mapsto \wht{a}
\end{align*}
then is a contractive algebra homomorphism.
It is called the \underline{Gelfand} \linebreak \underline{transformation}.
\pagebreak
\end{proposition}

\begin{proof}
Let $a$, $b \in A$, $\tau \in \Delta (A)$. We then have
\[ \wht{(ab)} ( \tau ) = \tau (ab) = \tau (a) \cdot \tau (b)
= \wht{\vphantom{b}a}\,(\tau) \cdot \wht{b}\,(\tau)
= \bigl( \,\wht{\vphantom{b}a} \cdot \wht{b} \,\bigr) (\tau). \]
Also $\bigl| \,\wht{a} ( \tau ) \,\bigr| = | \,\tau (a) \,| \leq | \,a \,|$,
cf.\ \ref{mlfbounded}, whence
$\bigl| \,\wht{a} \,\bigr| _{\infty} \leq | \,a \,|$.
\end{proof}

\begin{definition}[ideals]\index{ideal!left}%
\index{ideal}\index{ideal!proper}\index{ideal!maximal}%
Let $A$ be an algebra. A \underline{left ideal} in $A$ is a vector subspace
$I$ of $A$ such that $aI \subset I$ for all $a \in A$. A left ideal $I$ of $A$
is called \underline{proper} if $I \neq A$. A left ideal of $A$ is called a
\underline{maximal left ideal} or simply \underline{left maximal} if it is
proper and not properly contained in any other proper left ideal of $A$.
Likewise with ``right'' or ``two-sided'' instead of ``left''. If $A$ is commutative,
we drop these adjectives.
\end{definition}

\begin{lemma}\label{notinvideal}%
In a unital algebra $A$ with unit $e$, the following statements hold:
\begin{itemize}
   \item[$(i)$] a left ideal $I$ in $A$ is proper if and only if $e \notin I$,
  \item[$(ii)$] an element of $A$ is not left invertible in $A$ if and only if it
lies in some proper left ideal of $A$.
\end{itemize}
If $A$ is a unital Banach algebra, then furthermore:
\begin{itemize}
 \item[$(iii)$] if a left ideal $I$ of $A$ is proper, so is its closure $\overline{I}$,
 \item[$(iv)$] a maximal left ideal of $A$ is closed.
\end{itemize}
\end{lemma}

\begin{proof}
(i) is obvious. (ii) Assume that $a$ lies in a proper left ideal $I$ of $A$.
If $a$ had a left inverse $b$, one would have $e = ba \in I$, so that $I$
would not be proper. Conversely, assume that $a$ is not left invertible.
Consider the left ideal $I$ given by $I := \{ ba : b \in A \}$. Then
$e \notin I$, so that $a$ lies in the proper left ideal $I$. (iii) follows from
the fact that the set of left invertible elements of a unital Banach algebra
is an open neighbourhood of the unit, cf.\ \ref{leftinv}. (iv) follows from (iii).
\end{proof}

\begin{lemma}\label{inmaxideal}%
A proper left ideal in a unital algebra is contained in a maximal left ideal.
\end{lemma}

\begin{proof}
Let $I$ be a proper left ideal in a unital algebra $A$. Let $Z$ be the set of
proper left ideals in $A$ containing $I$, and order $Z$ by inclusion. It shall
be shown that $Z$ is inductively ordered. So let $C \neq \varnothing$ be a
chain in $Z$. Let $J := \bigcup \,C$. It shall be shown that $J \in Z$. The left
ideal $J$ is proper. Indeed one notes that $J$ does not contain the unit of
$A$, because otherwise some element of $Z$ would contain the unit. Now
Zorn's Lemma yields a maximal element of $Z$. \pagebreak
\end{proof}

\begin{lemma}\label{field}%
Let $I$ be a proper two-sided ideal in a unital algebra $A$.
The unital algebra $A/I$ then is a division ring if and only if $I$
is both a maximal left ideal and a maximal right ideal. 
\end{lemma}

\begin{proof}
Please note first that the algebra $A/I$ is unital because the two-sided ideal $I$
is proper, cf.\ \ref{notinvideal} (i). For $b \in A \setminus I$, consider the left ideal
\[ Ab+I := \{ \,ab+i : a \in A,\ i \in I \,\}. \]
It contains $I$ properly. (Because $b \in Ab+I$ as $A$ is unital,
and because $b \notin I$ by the assumption on the element $b$.)
It is the smallest left ideal $J$ in $A$ with $I \subset J$ and $b \in J$.
The following statements are equivalent.
(We let $\underline{a} := a + I \in A/I$ for $a \in A$.)
\par each non-zero element $\underline{b}$ of $A/I$ is left invertible
(please note \ref{leftrightinv}),
\par for all $\underline{b} \neq \underline{0}$ in $A/I$
there exists $\underline{c}$ in $A/I$ with
$\underline{c} \,\underline{b} = \underline{e}$,
\par for all $b \in A \setminus I$ there exist $c \in A$,
$i \in I$ with $cb + i = e$,
\par for all $b \in A \setminus I$ one has $e \in Ab+I$,
\par for all $b \in A \setminus I$ one has $Ab+I = A$,
\par for each left ideal $J$ containing $I$ properly, one has $J = A$,
\par the proper two-sided ideal $I$ is a maximal left ideal.
\end{proof}

\begin{theorem}\label{taumaxid}%
Let $A$ be a unital Banach algebra. The map $\tau \mapsto \ker \tau$
establishes a bijective correspondence between the set $\Delta (A)$
of multiplicative linear functionals $\tau$ on $A$ and the set of proper
two-sided ideals in $A$ which are both left maximal and right maximal.
\end{theorem}

\begin{proof}
If $\tau \in \Delta (A)$, then $\ker \tau$ has co-dimension $1$ in
$A$, so the two-sided ideal $\ker \tau$ is both left and right maximal.
Conversely, for a two-sided ideal $I$ in $A$ which is both left and
right maximal, $A/I$ is a Banach algebra by \ref{notinvideal} (iv),
and a division ring by \ref{field}, so isomorphic to $\mathds{C}$ by
the Gelfand-Mazur Theorem \ref{GelfandMazur}.
Let $\pi$ be such an isomorphism, and put $\tau (a) := \pi (a+I)$
$(a \in A)$. Then $\tau \in \Delta (A)$ and $\ker \tau =I$. To see
that the map $\tau \mapsto \ker \tau$ is injective, let
$\tau _1,\tau _2 \in \Delta (A)$, $\tau _1 \neq \tau _2$.
Let $a \in A$ with $\tau _1 (a) = 1$, $\tau _2 (a) \neq 1$.
With $b := a^2 -a$ we have
$\tau _i (b) = \tau _i (a) \cdot \bigl( \tau _i (a) -1 \bigr)$.
Thus $\tau _1 (b) = 0$, and either $\tau _2 (a) = 0$ or
else $\tau _2 (b) \neq 0$. So either
$b \in \ker \tau _1 \setminus \ker \tau _2$ or
$a \in \ker \tau _2 \setminus \ker \tau _1$.
\end{proof}

\medskip
In the remainder of this paragraph, we shall mainly consider
\linebreak \underline{commutative} Banach algebras. For
emphasis, we restate: \pagebreak

\begin{corollary}\label{taumaxidcomm}%
Let $A$ be a \underline{unital commutative} Banach algebra.
The map $\tau \mapsto \ker \tau$ establishes a bijective
correspondence between the set $\Delta (A)$ of multiplicative
linear functionals $\tau$ on $A$ and the set of maximal ideals
in $A$.
\end{corollary}

\begin{theorem}[abstract form of Wiener's Theorem]%
\label{abstrWien}%
\index{Theorem!Wiener's!abstract form}%
\index{Wiener's Theorem!abstract form}%
\index{abstract!form of Wiener's Thm.}%
Let $A$ be a \underline{unital commutative} Banach
algebra. An element $a$ of $A$ is not invertible in $A$ if and
only if $\wht{a}$ vanishes at some $\tau$ in $\Delta (A)$.
\end{theorem}

\begin{proof}
\ref{notinvideal}(ii), \ref{inmaxideal}, \ref{taumaxidcomm}.
\end{proof}

\begin{theorem}\label{rangetld}%
Let $A$ be a commutative Banach algebra. For $a \in \tld{A}$ we
then have $\s(a) = \wht{a}\,\bigl(\Delta ( \tld{A} )\bigr)$.
\end{theorem}

\begin{proof}
For $\lambda \in \mathds{C}$, the following statements are equivalent.
\begin{align*}
& \lambda \in \s(a), \\
& \lambda e-a \text{ is not invertible in } \tld{A}, \\
& \wht{ ( \lambda e - a ) } \,\bigl( \tld{ \tau } \bigr) = 0
\text{ for some } \tld{ \tau } \in \Delta (\tld{A}),
\quad \text{by \ref{abstrWien}} \\
& \lambda = \wht{a} \,\bigl( \tld{ \tau } \bigr)
\text{ for some } \tld{ \tau } \in \Delta (\tld{A}). \qedhere
\end{align*}
\end{proof}

\begin{theorem}%
If $A$ is a commutative Banach algebra without unit, we have
$\s(a) = \wht{a}\,\bigl(\Delta(A)\bigr) \cup \{ 0 \}$ for all $a \in A$.
Please note \ref{speczero}.
\end{theorem}

\begin{proof}
For $\lambda \in \mathds{C}$, the following statements are equivalent.
\begin{align*}
& \lambda \in \s(a), \\
& \lambda = \wht{a} \,\bigl( \tld{ \tau } \bigr)
\text{ for some } \tld{ \tau } \in \Delta (\tld{A}),
\quad \text{by \ref{rangetld}} \\
& \text{either } \lambda = 0 \text{ or } \lambda = \wht{a} \,( \tau )
\text{ for some } \tau \in \Delta (A), \\
& \text{(according as } \tld{ \tau } \text{ vanishes on all of } A
\text{ or not).} \qedhere
\end{align*}
\end{proof}

\begin{corollary}\label{rangeGT}%
For a commutative Banach algebra $A$ and $a \in A$, we have
\[ \underline{\s(a) \setminus \{ 0 \}
= \wht{a} \,\bigl( \Delta (A) \bigr) \setminus \{ 0 \}}. \]
\end{corollary}

\begin{theorem}\label{normGTrl}%
For a commutative Banach algebra $A$ and $a \in A$, we have
\[ \bigl| \,\wht{a} \,\bigr|_{\infty} = \rlambda(a). \pagebreak \]
\end{theorem}

\begin{proof}
This follows from the preceding corollary by applying the
Spectral Radius Formula.
\end{proof}

\begin{corollary}\label{C*normGT}%
For a commutative C*-algebra $A$ and $a \in A$, we have
\[ \bigl| \,\wht{a} \,\bigr| _{\infty} = \|\,a\,\|. \]
\end{corollary}

\begin{proof}
By \ref{C*rl} we have $\rlambda (a) = \|\,a\,\|$.
\end{proof}

\begin{corollary}\label{C*GTisom}%
If $A$ is a commutative C*-algebra $\neq \{ 0 \}$, then
$\Delta (A) \neq \varnothing$ and the Gelfand transformation
$A \to {\ell\,}^{\infty}\bigl(\Delta(A)\bigr)$ is iso\-metric. 
\end{corollary}

\begin{definition}[Hermitian linear functionals]%
\label{Hermfunct}\index{functional!Hermitian}\index{Hermitian!functional}%
A linear functional $\varphi$ on a \st-algebra $A$ is called
\underline{Hermitian} if $\varphi (a^*) = \overline{\varphi (a)}$
for all $a \in A$.
\end{definition}

\begin{proposition}\label{mlfHerm}%
A commutative Banach \st-algebra $A$ is Hermitian
precisely when each $\tau \in \Delta (A)$ is Hermitian.
\end{proposition}

\begin{proof}
Let $a \in A$. We then have
$\{ \tau (a) : \tau \in \Delta (A) \} \setminus \{ 0 \} = \s(a) \setminus \{ 0 \}$.
Thus, if every $\tau$ is Hermitian, then $\s(a)$ is real for $a$ Hermitian,
so that $A$ is Hermitian. Conversely, if $A$ is Hermitian, then each
$\tau$ assumes real values on Hermitian elements. Given any $c \in A$,
however, we can write $c = a + \iu b$ with $a,b$ Hermitian, so that for
$\tau \in \Delta (A)$ we have
$\tau (c^*) = \tau (a) - \iu \tau (b) = \overline{\tau (a) + \iu \tau (b)}
= \overline{\tau (c)}$.
\end{proof}

\begin{corollary}\label{Hermstarhom}%
Let $A$ be a \underline{Hermitian} commutative Banach \linebreak
\st-algebra with $\Delta(A) \neq \varnothing$. The Gelfand transformation
then is a \linebreak \st-algebra homomorphism. In particular, the range of
the Gelfand transformation then is a \st-subalgebra of
${\ell\,}^{\infty}\bigl(\Delta(A)\bigr)$ (not only a subalgebra).
\end{corollary}

\begin{proposition}%
Let $A$ be a commutative Banach algebra. For two elements $a,b$
of $A$ we have
\begin{align*}
\s(a+b) \subset & \ \s(a)+\s(b) \\
\s(ab) \subset & \ \s(a) \,\s(b).
\end{align*}
\end{proposition}

\begin{proof}
This is an application of \ref{rangetld}. Indeed, if $\lambda \in \s(ab)$,
there exists $\tld{\tau} \in \Delta (\tld{A})$ such that
$\tld{\tau} (ab) = \lambda$. But then $\lambda = \tld{\tau} (a) \,\tld{\tau} (b)$
where $\tld{\tau} (a) \in \s(a)$ and $\tld{\tau} (b) \in \s(b)$. \pagebreak
\end{proof}

\begin{theorem}\label{specsub}%
Let $A$ be a Banach algebra. If $a$, $b$ are \underline{commuting}
elements of $A$, then
\begin{align*}
\s(a+b) \subset & \ \s(a)+\s(b) \\
\s(ab) \subset & \ \s(a) \,\s(b).
\end{align*}
\end{theorem}

\begin{proof}
This follows from the preceding proposition by an application
of the Theorem of Civin \& Yood, cf.\ \ref{CivinYood}.
Indeed, this result says that there exists a commutative
closed subalgebra $B$ of $\tld{A}$ containing $a$ and $b$,
such that $\s _A (c) = \s _B (c)$ for $c \in \{ \,a, b, a+b, ab \,\}$.
\end{proof}

\medskip
For the preceding theorem, see also \ref{commrlsub} and
\ref{commspeccont}. 

Finally an application to harmonic analysis:

\begin{example}%
Consider the Banach \st-algebra ${\ell\,}^1 (\mathds{Z})$,
cf.\ example \ref{l1G}. We are interested in knowing the set of
multiplicative linear functionals of ${\ell\,}^1(\mathds{Z})$. To this
end we note that a multiplicative linear functional on
${\ell\,}^1(\mathds{Z})$ is determined by its value at $\delta _1$.
So let $\tau$ be a multiplicative linear functional on
${\ell\,}^1(\mathds{Z})$ and put $\alpha := \tau (\delta _1)$.
We then have $|\,\alpha \,| \leq 1$ by \ref{mlfbounded}. But as
$\delta _{-1}$ is the inverse of $\delta _1$, we also have
$\tau (\delta _{-1}) = \alpha ^{-1}$ and $|\,\alpha ^{-1}\,| \leq 1$.
This makes that $\alpha$ belongs to the unit circle:
$\alpha = {\mathrm{e} \,}^{\iu t}$ for some
$t \in \mathds{R} / 2 \pi \mathds{Z}$. For an integer $k$, we obtain
$\tau ( \delta _k ) = \bigl( \tau ( \delta _1 ) \bigr)^k = {\mathrm{e} \,}^{\iu kt}$.
This implies that for $a \in {\ell\,}^1(\mathds{Z})$, we have
\[ \wht{a}(\tau) = \sum _{k \in \mathds{Z}} a(k) \,{\mathrm{e} \,}^{\iu kt}. \]
Conversely, by putting for some $t \in \mathds{R} / 2 \pi \mathds{Z}$
\[ \tau (a) := \sum _{k \in \mathds{Z}} a(k) \,{\mathrm{e} \,}^{\iu kt} \qquad (a \in A), \]
a multiplicative linear functional $\tau$ on ${\ell\,}^1(\mathds{Z})$ is defined.
Indeed, we compute:
\begin{align*}
\tau (ab) & = \sum _{k \in \mathds{Z}} (ab)(k) \,{\mathrm{e} \,}^{\iu kt}
            = \sum _{k \in \mathds{Z}} \,\sum _{l+m=k} a(l) \,b(m) \,{\mathrm{e} \,}^{\iu kt} \\
          & = \biggl( \,\sum _{l \in \mathds{Z}} a(l) \,{\mathrm{e} \,}^{\iu lt} \biggr) \cdot
              \biggl( \,\sum _{m \in \mathds{Z}} b(m) \,{\mathrm{e} \,}^{\iu mt} \biggr)
            = \tau (a) \cdot \tau (b).
\end{align*}
As a multiplicative linear functional is determined by its
value at $\delta _1$, it follows that
$\Delta \bigl({ \ell\,}^1 (\mathds{Z})\bigr)$ can be identified with
$\mathds{R} / 2 \pi \mathds{Z}$. Also, a function on
$\mathds{R} / 2 \pi \mathds{Z}$ is the Gelfand transform
of an element of ${\ell\,}^1 (\mathds{Z})$
if and only if it has an absolutely convergent Fourier expansion.
\pagebreak
\end{example}

From the abstract form of Wiener's Theorem \ref{abstrWien}, it follows now:

\begin{theorem}[Wiener]%
\index{Theorem!Wiener's}\index{Wiener's Theorem}%
Let $f$ be a function on $\mathds{R} / 2 \pi \mathds{Z}$
which has an absolutely convergent Fourier expansion. If
$f(t) \neq 0$ for all \linebreak $t \in \mathds{R} / 2 \pi \mathds{Z}$,
then also $1/f$ has an absolutely convergent Fourier
\linebreak expansion.
\end{theorem}

This easy proof of Wiener's Theorem was an early success
of the theory of Gelfand.

Using \ref{mlfHerm}, the reader can now check that
${\ell\,}^1 (\mathds{Z})$ is a Hermitian Banach \st-algebra.

We shall return to the non-commutative theory in \ref{leftspectrum}.

\clearpage


\section{The Gelfand Topology}%
\label{Gtopo}

In this paragraph, let $A$ be a Banach algebra.

\begin{definition}[the Gelfand topology, the spectrum]%
\index{Gelfand!topology}\index{topology!Gelfand}\label{Gelftopo}%
\index{spectrum!of an algebra}\index{D(A)@$\Delta(A)$}%
One imbeds the set $\Delta (A)$ in the unit ball, cf.\ \ref{mlfbounded},
of the dual space of $A$ equipped with the weak* topology, cf.\ the
appendix \ref{weak*top}. The relative topology on $\Delta (A)$ is
called the \underline{Gelfand topology}. When equipped with the
Gelfand topology, $\Delta (A)$ is called the \underline{spectrum} of $A$.
\end{definition}

\begin{proposition}%
We have:
\begin{itemize}
  \item[$(i)$] the closure $\overline{\Delta(A)}$ is a compact Hausdorff space.
 \item[$(ii)$] each point in $\overline{\Delta (A)}$ is either zero or else a
multiplicative linear functional.
\end{itemize}
\end{proposition}

\begin{proof}
(i) follows from Alaoglu's Theorem, cf.\ the appendix \ref{Alaoglu}.
(ii): every adherence point of $\Delta (A)$ is a linear and multiplicative functional.
\end{proof}

\begin{theorem}\label{mlf0compact}%
Either $\Delta(A)$ is compact or
$\overline{\Delta(A)} = \Delta(A) \cup \{ \tau _{\infty} \}$
where $\tau_{\infty} := 0$, in which case $\Delta (A)$ is locally compact and
$\overline{\Delta (A)}$ is the one-point compactification of $\Delta (A)$.
\end{theorem}

\begin{corollary}\label{unitalcomp}%
If $A$ is unital, then $\Delta (A)$ is compact.
\end{corollary}

\begin{proof}
We have $\tau (e) = 1$ for all $\tau \in \Delta (A)$,
whence $\sigma (e) = 1$ for all $\sigma \in \overline{\Delta (A)}$ because
the property in question is preserved by pointwise convergence. On the other
hand we have $\tau_{\infty} (e) = 0$, so that
$\tau_{\infty} \notin \overline{\Delta (A)}$, whence it follows that
$\Delta (A)$ is compact.
\end{proof}

\begin{lemma}[the abstract Riemann-Lebesgue Lemma]%
\index{Theorem!abstract!Riemann-Lebesgue Lemma}%
\index{Theorem!Riem.-Leb.\ Lemma, abstract}%
\index{abstract!Riemann-Lebesgue Lemma}%
\index{Lemma!abstract Riemann-Lebesgue}%
\index{Riem.-Leb.\ Lemma, abstract}%
If $\Delta (A) \neq \varnothing$, then for all $a \in A$, we have
\[ \wht{a} \in C\0\bigl(\Delta (A)\bigr). \]
\end{lemma}

\begin{proof}
If $\Delta (A)$ is not compact, we have
\[ \lim _{\tau \to \tau _{\infty}} \wht{a} ( \tau )
= \lim _{\tau \to \tau _{\infty}} \tau (a)
= \tau_{\infty} (a) = 0. \qedhere \]
\end{proof}

\medskip
From \ref{GelfTransf} it follows: \pagebreak

\begin{corollary}%
If $\Delta (A) \neq \varnothing$, then the
Gelfand transformation $a \mapsto \wht{a}$
is a contractive algebra homomorphism
\[ A \to C\0\bigl(\Delta (A)\bigr). \]
\end{corollary}

\begin{proposition}\label{Hermdense}%
Let $A$ be a \underline{Hermitian} commutative Banach \linebreak
\st-algebra with $\Delta (A) \neq \varnothing$. The range of the
Gelfand  transformation then is dense in $C\0\bigl(\Delta (A)\bigr)$.
\end{proposition}

\begin{proof}
One applies the Stone-Weierstrass Theorem, see the appendix
\linebreak \ref{StW}. This requires \ref{Hermstarhom}.
\end{proof}

The remaining statements in this paragraph concern commutative
\linebreak C*-algebras, as well as their dual objects: locally compact
Hausdorff spaces.

\begin{theorem}[the Commutative Gelfand-Na\u{\i}mark Theorem]%
\label{commGN}\index{Theorem!Gelfand-Na\u{\i}mark!Commutative}%
\index{Gelfand-Naimark@Gelfand-Na\u{\i}mark!Theorem!Commutative}
If $A$ is a commutative C*-algebra $\neq \{ 0 \}$, then the
Gelfand transformation establishes a C*-algebra
isomorphism from $A$ onto $C\0\bigl(\Delta (A)\bigr)$.
\end{theorem}

\begin{proof} It follows from \ref{C*GTisom} that
$\Delta (A) \neq \varnothing$ and that the Gelfand transformation is
isometric, and therefore injective. The Gelfand transformation is a
\st-algebra homomorphism by \ref{Hermstarhom}. It suffices now
to show that it is surjective. Its range is dense \ref{Hermdense} and
on the other hand complete (because $a \mapsto \wht{a}$ is
isometric) and thus closed. \end{proof}

Also of great importance is the following result which identifies the
spectrum of a normal element $b$ of a C*-algebra $A$
with the spectrum \ref{Gelftopo} of the C*-subalgebra of $\tld{A}$
generated by $b$, $b^*$, and $e$. Please note that the
C*-subalgebra in question is commutative as $b$ is normal.

\begin{theorem}\label{spechomeom}%
Let $A$ be a C*-algebra and let $b$ be a normal element of $A$.
Denote by $B$ the C*-subalgebra of $\tld{A}$ generated by
$b$, $b^*$, and $e$. Then $\wht{b}$ is a homeomorphism from
$\Delta (B)$ onto $\s_A (b) = \s_B (b)$.
\end{theorem}

\begin{proof}
We have $\s_A (b) = \s_B (b)$ by \ref{specHermsubalg}. Furthermore,
the map $\wht{b}$ is continuous and maps $\Delta (B)$ surjectively
onto $\s_B (b)$, cf.\ \ref{rangetld}. It is also injective: if
$\wht{b} \,( \tau _1 ) = \wht{b} \,( \tau _2 )$ then
$\tau _1 \bigl(p(b,b^*)\bigr) = \tau _2 \bigl(p(b,b^*)\bigr)$ for all
$p \in \mathds{C} [ z , \overline{z} ]$. Since the polynomial functions
in $b$ and $b^*$ are dense in $B$, it follows that $\tau _1 = \tau _2$.
The continuous bijective function $\wht{b}$ is a homeomorphism
because its domain $\Delta (B)$ is compact and its range $\s_B (b)$
is a Hausdorff space, cf.\ the appendix \ref{homeomorph}. \pagebreak
\end{proof}

\medskip
The preceding two theorems together yield the so-called
``operational calculus'':

\begin{corollary}[the operational calculus]%
\index{operational calculus}\label{opcalc}%
Let $A$ be a  C*-algebra and let $b$ be a normal element of $A$.
Let $B$ be the C*-subalgebra of $\tld{A}$ generated by
$b$, $b^*$, and $e$. For $f \in C\bigl(\s(b)\bigr)$, one denotes by $f(b)$
the element of $B$ satisfying
\[ \wht{f(b)} = f \circ \wht{b} \in C\bigl(\Delta(B)\bigr). \]
The mapping
\begin{alignat*}{2}
C\bigl(\s(b)\bigr) \ & & \to & \ B \\
f \ & & \mapsto & \ f(b)
\end{alignat*}
is a C*-algebra isomorphism, called the
\underline{operational calculus} for $b$. It is
the only \st-algebra homomorphism from
$C\bigl(\s(b)\bigr)$ to $B$ mapping the identity
function to $b$ and the constant function $1$ to $e$.
The following \underline{Spectral Mapping Theorem} holds:
\[ \s\bigl(f(b)\bigr) = f\bigl(\s(b)\bigr) \qquad \Bigl( f \in C\bigl(\s(b)\bigr) \Bigr). \]
\end{corollary}

\begin{proof}
The uniqueness statement follows from \ref{Cstarcontr}
and the Stone-\linebreak Weierstrass Theorem, cf.\ the appendix \ref{StW}.
\end{proof}

\medskip
See also \ref{weakopcalc}.

Now for locally compact Hausdorff spaces.

If $\Omega$ is a locally compact Hausdorff space, then what is the spectrum
of the commutative C*-algebra $C\0(\Omega)$? Precisely $\Omega$, as the
next theorem shows.

\begin{theorem}\label{nonews}
Let $\Omega$ be a locally compact Hausdorff space $\neq \varnothing$.
Then $\Delta \bigl( C\0(\Omega) \bigr)$ is homeomorphic to $\Omega$.

Indeed, for each $\omega \in \Omega$, consider the evaluation
$\varepsilon_{\omega}$ at $\omega$ given by
\begin{alignat*}{2}
\varepsilon_{\omega} : \ & & C\0(\Omega) \to \ & \mathds{C} \\
                                          \ & & f \mapsto  \ & \varepsilon_{\omega} (f) := f(\omega).
\end{alignat*}
Then the map $\omega \mapsto \varepsilon_{\omega}$ is a homeomorphism
$\Omega \to \Delta \bigl(C\0(\Omega)\bigr)$, and
$\wht{f}(\varepsilon_{\omega}) = f(\omega)$ for all $f \in C(\Omega)$ and all
$\omega \in \Omega$. Thus, upon identifying
$\varepsilon_{\omega} \in \Delta \bigl(C\0(\Omega)\bigr)$ with $\omega \in \Omega$,
the Gelfand transformation on $C\0(\Omega)$ becomes the identity map.
\pagebreak
\end{theorem}

Before going to the proof, we note the following corollary.

\begin{corollary}
Two locally compact Hausdorff spaces $\Omega$, $\Omega' \neq \varnothing$
are homeomorphic whenever the C*-algebras $C\0(\Omega)$, $C\0(\Omega')$
are isomorphic as \st-algebras. (Cf.\ \ref{Cstarisom}.) In this sense a locally
compact Hausdorff space $\Omega$ is determined up to a homeomorphism
by the \linebreak \st-algebraic structure of the C*-algebra $C\0(\Omega)$ alone. 
\end{corollary}

\begin{proof}[\hspace{-3.55ex}\mdseries{\scshape{Proof of \protect\ref{nonews}}}]%
Assume first that $\Omega$ is compact, so $C\0(\Omega) = C(\Omega)$.
Clearly each $\varepsilon_{\omega}$ is a multiplicative linear functional
on $C(\Omega)$. The map $\omega \mapsto \varepsilon_{\omega}$ is
injective since $C(\Omega)$ separates the points of $\Omega$ (by
Urysohn's Lemma). This map is continuous in the weak* topology because
if $\omega_{\alpha} \to \omega$, then $f(\omega_{\alpha}) \to f(\omega)$
for all $f \in C(\Omega)$. It now suffices to show that this map is surjective
onto $\Delta \bigl(C(\Omega)\bigr)$ because a continuous bijection from
a compact space to a Hausdorff space is a homeomorphism, cf.\ the appendix
\ref{homeomorph}. By corollary \ref{taumaxidcomm}, the sets
$M_{\omega} := \{ \,f \in C(\Omega) : f(\omega) = 0 \,\}$ $(\omega \in \Omega)$
are maximal ideals in $C(\Omega)$, and it is enough to show that each
maximal ideal in $C(\Omega)$ is of this form. It is sufficient to prove that
each maximal ideal $I$ in $C(\Omega)$ is contained in some $M_{\omega}$.
So let $I$ be a maximal ideal in $C(\Omega)$, and suppose that for each
$\omega \in \Omega$ there exists $f_{\omega} \in I$ with
$f_{\omega} (\omega) \neq 0$. The open sets $\{ \,f_{\omega} \neq 0 \,\}$
$(\omega \in \Omega)$ cover $\Omega$, so by passing to a finite subcover,
we obtain $f_1, \ldots , f_n \in I$ such that $g:= \sum_{i=1}^{n} \overline{f_i} f_i$
is $> 0$ on $\Omega$ . Then $g \in I$ would be invertible in $C(\Omega)$,
contradicting proposition \ref{notinvideal} (ii).

Next assume that $\Omega$ is not compact, and consider the one-point
compactification $K$ of $\Omega$. Let $\infty \in K \setminus \Omega$
be the corresponding point at infinity. We then may identify $C\0(\Omega)$
with the set of continuous complex-valued functions on $K$ vanishing at
$\infty$, and hence we may identify $C(K)$ with
$\mathds{C} 1_{\Omega} \oplus C\0(\Omega)$. That is, $C(K)$ ``is'' the
unitisation of $C\0(\Omega)$. This implies that $\Delta \bigl( C(K) \bigr)$
may be identified with $\Delta \bigl( C\0(\Omega) \bigr) \cup \{ \tau_\infty \}$,
where $\tau_\infty := 0$. Hence $\Delta \bigl( C(K) \bigr)$ ``is'' the one-point
compactification of $\Delta \bigl( C\0(\Omega) \bigr)$, cf.\ \ref{mlf0compact}.
(It is used that $\Delta \bigl( C\0(\Omega) \bigr)$ is not compact, for else
$C\0(\Omega)$ would be unital by the commutative Gelfand-Na\u{\i}mark
Theorem \ref{commGN}.) From the preceding paragraph, we have that the
map $\omega \mapsto \varepsilon_\omega$ is a homeomorphism from $K$
onto $\Delta \bigl( C(K) \bigr)$, and thus restricts to a homeomorphism from
$\Omega$ onto $\Delta \bigl( C\0(\Omega) \bigr)$.
\end{proof}

\medskip
It also follows that there is a bijective correspondence between the
equivalence classes of commutative C*-algebras modulo \st-algebra
isomorphism and the equivalence classes of locally compact Hausdorff
spaces modulo homeomorphism. \pagebreak

\clearpage


\section{The Stone-\texorpdfstring{\v{C}}{C}ech Compactification}\label{StoneCech}

(This paragraph can be skipped at first reading. Although this is the right
time in a mathematician's education to acquire the present material.)

We first have to fix some terminology.

\begin{definition}[completely regular spaces]\index{completely regular!space}%
A topological space $\Omega$ is called \underline{completely regular}
if it is Hausdorff, and if for every non-empty proper closed subset $C$ of
$\Omega$ and every point $\omega \in \Omega \setminus C$ there exists
a continuous function $f$ on $\Omega$ taking values in $[0, 1]$ such that
$f$ vanishes on $C$ and $f(\omega) = 1$.
\end{definition}

Please note that every subspace of a completely regular Hausdorff space
is itself a completely regular Hausdorff space.

For example every normal Hausdorff space is completely regular
by Urysohn's Lemma. In particular, every compact Hausdorff space is
completely regular. It follows that every locally compact Hausdorff space
is completely regular (by considering the one-point compactification).

\begin{definition}[compactification]\index{compactification}%
A \underline{compactification} of a Hausdorff space $\Omega$ is a pair
$(K, \phi)$, where $K$ is a compact Hausdorff space, and $\phi$ is a
homeomorphism of $\Omega$ onto a dense subset of $K$. We shall
identify $\Omega$ with its image $\phi(\Omega) \subset K$, and simply
speak of ``the compactification $K$ of $\Omega$''.
\end{definition}

Please note that at most completely regular Hausdorff spaces can have
a compactification.

For example $([-1, 1], \tanh)$ is a compactification of $\mathds{R}$.

\begin{definition}[equivalence of compactifications]%
\index{equivalent!compactifications}%
Any  two compactifications $(K_1,\phi_1)$, $(K_2,\phi_2)$ of a
Hausdorff space $\Omega$ are said to be \underline{equivalent}
in case there exists a homeomorphism $\theta : K_1 \to K_2$ with
$\theta \circ \phi_1 = \phi_2$.
\end{definition}

\begin{definition}[completely regular algebras]\index{completely regular!algebra}%
Let $\Omega \neq \varnothing$ be a Hausdorff space, and let $A$ be a unital
C*-subalgebra of $C_b(\Omega)$. One says that $A$ is \underline{completely regular}
if for every non-empty proper closed subset $C$ of $\Omega$ and for every point
$\omega \in \Omega \setminus C$, there is a function $f \in A$ taking values in
$[0, 1]$ such that $f$ vanishes on $C$ and $f(\omega) = 1$. In particular $\Omega$
is completely regular if and only if the unital C*-algebra $C_b(\Omega)$ is
completely regular.
\end{definition}

We can now characterise all compactifications of a completely
\linebreak regular Hausdorff space: 

\begin{theorem}[all compactifications]\label{allcomp}%
Let $\Omega \neq \varnothing$ be a completely regular Hausdorff space.
The following statements hold.
\begin{itemize}
   \item[$(i)$]   If $A$ is a completely regular unital C*-subalgebra of
                         $C_b(\Omega)$ then $(\Delta(A), \varepsilon)$ is a
                         compactification of $\Omega$, where $\varepsilon$ is the map
                         \linebreak
                         associating with each $\omega \in \Omega$ the evaluation
                         $\varepsilon_{\omega}$ at $\omega$ of the functions in $A$.
                         For $f \in A$, there exists a (necessarily unique) function
                         $g \in C\bigl( \Delta (A) \bigr)$ with $f = g \circ \varepsilon$.
                         It is interpreted as the ``continuation of $f$ to the compactification''.
                         It is given as the Gelfand transform of $f$.
  \item[$(ii)$]   If $(K,\phi)$ is a compactification of $\Omega$, then
                         \[ A_K := \{ \,g \circ \phi : g \in C(K) \,\} \subset C_b(\Omega) \]
                         is a completely regular unital C*-subalgebra of $C_b(\Omega)$,
                         and the \linebreak ``restriction map'' $C(K) \to A_K$,
                         $g \mapsto g \circ \phi$ is a C*-algebra \linebreak isomorphism.
                         The compactification $(K,\phi)$ then is equivalent to the
                         compactification associated with $A := A_K$ as in $(i)$.
 \item[$(iii)$]   Two compactifications $K$, $K'$ of $\Omega$ are equivalent
                         if and only if $A_K = A_{K'}$.
\end{itemize}
In other words, the equivalence classes of compactifications $K$
of $\Omega$ are in bijective correspondence with the spectra of
the completely regular unital C*-subalgebras $A_K$ of
$C_b(\Omega)$. The $A_K$ consisting of those functions in
$C_b(\Omega)$ which have a continuation to $K$.
\end{theorem}

\begin{proof}
(i): Let $A$ be a completely regular unital C*-subalgebra of
$C_b(\Omega)$. Then $\Delta(A)$ is a compact Hausdorff space by
\ref{unitalcomp}.
For $\omega \in \Omega$, the evaluation $\varepsilon_{\omega}$ at
$\omega$ of the functions in $A$ is a multiplicative linear functional on
$A$ because it is non-zero as $A$ is unital in $C_b(\Omega)$. The map
$\varepsilon : \omega \mapsto \varepsilon_{\omega}$ thus takes
$\Omega$ to $\Delta (A)$.
Moreover the Gelfand transform $g$ of a function $f$ in $A$ satisfies
$f = g \circ \varepsilon$ as
$f (\omega) = \varepsilon_{\omega} (f) = g (\varepsilon_{\omega})$
for all $\omega \in \Omega$.
The set $\varepsilon (\Omega)$ is dense in $\Delta (A)$. For otherwise
there would exist a non-zero continuous function $g$ on $\Delta (A)$
vanishing on $\varepsilon (\Omega)$. (By Urysohn's Lemma.) The
Commutative Gelfand-Na\u{\i}mark Theorem \ref{commGN} would imply
that there exists a unique function $f \in A$ whose Gelfand transform is $g$,
and since $g \neq 0$, one would have $f \neq 0$, contradicting
$f = g \circ \varepsilon = 0$.
The map $\varepsilon : \Omega \to \Delta (A)$ is injective because $A$
separates the points of $\Omega$. (As $A$ is completely regular and
$\Omega$ is Hausdorff.) This map is continuous in the weak* topology
because if $\omega_{\alpha} \to \omega$ then 
$f(\omega_{\alpha}) \to f(\omega)$ for all $f \in A$.
To see that $\varepsilon$ is a homeomorphism
from $\Omega$ onto $\varepsilon (\Omega)$, it remains to show that
$\varepsilon : \Omega \to \varepsilon (\Omega)$ is open. So let $U$
be a non-empty proper open subset of $\Omega$. For $\omega \in U$,
let $f \in A$ such that $f$ vanishes on $\Omega \setminus U$
and $f(\omega) = 1$. Denote by $g$ the Gelfand transform of $f$.
The set $V := g^{-1}(\mathds{C} \setminus \{0\})$ is an open
neighbourhood of $\varepsilon_{\omega}$ in $\Delta (A)$. Now the set
$\varepsilon (U)$ contains $V \cap \varepsilon (\Omega)$ and so
is a neighbourhood of $\varepsilon_{\omega}$ in $\varepsilon (\Omega)$.
Therefore $\varepsilon (U)$ is open in $\varepsilon (\Omega)$. 

(ii): Let $(K,\phi)$ be a compactification of $\Omega$. The map from $C(K)$
to $A_K$ given by $g \mapsto g \circ \phi$, is interpreted as ``restriction to
$\Omega$''. It clearly is a \st-algebra homomorphism. It is isometric (and
thus injective) by density of $\phi(\Omega)$ in $K$. Hence the range
$A_K$ of this map is complete, and thus a unital and completely regular
C*-subalgebra of $C_b(\Omega)$. Since the above map is surjective by
definition, it is a C*-algebra isomorphism. It shall next be shown that
$(K,\phi)$ is equivalent to the compactification associated with $A := A_K$
as in (i). For $x \in K$, consider $\theta_x : f \to g(x)$, $(f \in A_K)$, where
$g$ is the unique function in $C(K)$ such that $f = g \,\circ \,\phi$, see above.
Then $\theta_x$ is a multiplicative linear functional on $A_K$ because
it is non-zero as $A_K$ is unital in $C_b(\Omega)$. Now the map
$\theta : K \to \Delta(A_K)$, $x \mapsto \theta_x$ is the desired
homeomorphism $K \to \Delta(A_K)$. Indeed, this map is a
homeomorphism by theorem \ref{nonews} because $A_K$ is essentially
the same C*-algebra as $C(K)$, see above. It is also clear that
$\theta \circ \phi = \varepsilon$.

(iii) follows essentially from the last statement of (ii).
\end{proof}

\medskip
For example if $\Omega \neq \varnothing$ is a locally compact Hausdorff
space which is not compact, then the one-point compactification of $\Omega$
is associated with $\mathds{C} 1_{\Omega} \oplus C\0(\Omega)$, as is easily seen.

\begin{definition}[Stone-\v{C}ech compactification]%
\index{Stone-Cech@Stone-\v{C}ech compactification}%
\index{compactification!Stone-Cech@Stone-\v{C}ech}%
Let $(K,\phi)$ be a compactification of a Hausdorff space $\Omega$.
Then $(K,\phi)$ is called a \underline{Stone-\v{C}ech compactification}
of $\Omega$ if every $f \in C_b(\Omega)$ is of the form $f = g \circ \phi$
for some $g \in C(K)$. The function $g$ then is uniquely determined by
density of $\phi (\Omega)$ in $K$.

In other words, a Stone-\v{C}ech compactification of a Hausdorff space
$\Omega$ is a compact Hausdorff space $K$ containing $\Omega$ as
a dense subset, such that every function $f \in C_b(\Omega)$ has a
(necessarily unique) continuation to $K$. \pagebreak
\end{definition}

\begin{corollary}%
Let $\Omega$ be a completely regular Hausdorff space. There then
exists a Stone-\v{C}ech compactification of $\Omega$,
and it is unique within equivalence. It is denoted by
$(\beta\/\Omega, \varepsilon)$, or simply by $\beta\/\Omega$.
\end{corollary}

\begin{proof}
Apply theorem \ref{allcomp} above with $A := C_b(\Omega)$
if $\Omega \neq \varnothing$.
\end{proof}

\medskip
Hence the following memorable result.

\begin{corollary}
A Hausdorff space is completely regular if and only if it is homeomorphic
to a subspace of a compact Hausdorff space. 
\end{corollary}

\begin{theorem}[the universal property]%
\index{universal!property!Stone-Cech@of Stone-\v{C}ech compactif.}
Let $\Omega \neq \varnothing$ be a completely regular Hausdorff space.
The Stone-\v{C}ech compactification $(\beta\/\Omega,\varepsilon)$ of
$\Omega$ enjoys the following \underline{universal property}. Let $K$
be a compact Hausdorff space, and $\phi$ be a continuous function
$\Omega \to K$. Then $\phi$ has a (necessarily unique) continuation
to $\beta\/\Omega$. That is, there exists a (necessarily unique) continuous
function $\tld{\phi} : \beta\/\Omega \to K$ with $\tld{\phi} \circ \varepsilon = \phi$.
If $(K,\phi)$ is a compactification of $\Omega$, then $\tld{\phi}$ is surjective.
This shows among other things that the Stone-\v{C}ech compactification
$\beta\/\Omega$ of $\Omega$ is the ``largest'' compactification of $\Omega$
in the sense that any other compactification of $\Omega$ is a continuous
image of $\beta\/\Omega$.
\end{theorem}

\begin{proof}
Put $A := C_b(\Omega)$ and $B := C(K)$. Let $(\beta K, i)$ be the Stone-\v{C}ech
compactification of $K$. Then $\beta K$ is homeomorphic to $K$ as $K$ is
compact. The map $\phi$ induces an ``adjoint'' map $\pi : B \to A$ by \linebreak
putting $\pi (g) := g \circ \phi$, $(g \in B)$. The map $\pi$ is a unital \st-algebra
homomorphism, and hence induces in turn a ``second adjoint'' map
$\pi ^* : \Delta (A) \to \Delta (B)$, by letting $\pi ^* (\tau) := \tau \circ \pi$
$\bigl( \tau \in \Delta (A) \bigr)$. The map $\pi ^*$ is continuous by the universal
property of the weak* topology, cf.\ the appendix \ref{weak*top}. We have
$\pi ^* \circ \varepsilon = i \circ \phi$ as for $\omega \in \Omega$, $g \in B$,
one computes
\begin{align*}
& \ \bigl( (\pi ^* \circ \varepsilon) (\omega) \bigr) (g)
= \bigl( \pi ^* (\varepsilon_{\omega}) \bigr) (g)
= ( \varepsilon_{\omega} \circ \pi ) (g)
= \varepsilon_{\omega} \bigl( \pi (g) \bigr)
= \varepsilon_{\omega} (g \circ \phi) \\
= & \ (g \circ \phi) (\omega)
= g \bigl( \phi (\omega) \bigr)
= \Bigl( i \bigl( \phi (\omega) \bigr) \Bigr) (g)
= \bigl( ( i \circ \phi ) (\omega) \bigr) (g).
\end{align*}
Now let $\tld{\phi} := i^{-1} \circ \pi ^*$. This map is continuous because
$i$ is open. Also
$\tld{\phi} \circ \varepsilon = i^{-1} \circ \pi ^* \circ \varepsilon = \phi$
as required. Finally, if $(K,\phi)$ is a compacti\-fication, then $\phi(\Omega)$
is dense in $K$. Then $\tld{\phi}(\beta\/\Omega)$ is dense in $K$ the more,
and also compact, so that $\tld{\phi}(\beta\/\Omega) = K$. \pagebreak
\end{proof}

\clearpage


\chapter{Positive Elements}%
\label{PositiveElements}


\setcounter{section}{16}

\section{The Positive Cone in a Hermitian Banach \texorpdfstring{$*$-}{\052\055}Algebra}

\begin{reminder}[cone]\index{cone}%
A \underline{cone} in a real vector space $V$ is a non-empty subset $C$
of $V$ such that for $c \in C$ and $\lambda > 0$ also $\lambda c \in C$.
\end{reminder}

\begin{definition}[positive elements]%
\index{A9@$A_+$}\index{positive!element}\index{a1@$a \geq 0$}%
Let $A$ be a \st-algebra. We shall denote by
\[ A_+ := \{ \,a \in A\sa : \s _A (a) \subset [0, \infty[ \,\} \]
the set of \underline{positive} elements in $A$. To indicate that an
element $a \in A$ is positive, we shall also write $\underline{a \geq 0}$.
We stress that positive elements are required to be \underline{Hermitian}.
The set of positive elements in $A$ is a cone in $A\sa$.
\end{definition}

\begin{proposition}\label{hompos}%
Let $\pi : A \to B$ be a \st-algebra homomorphism from a \st-algebra
$A$ to a \st-algebra $B$. If $a \in A_+$, then $\pi (a) \in B_+$.
\end{proposition}

\begin{proof}
$\s\bigl(\pi (a)\bigr) \setminus \{ 0 \} \subset \s(a) \setminus \{ 0 \} $,
cf.\ \ref{spechom}.
\end{proof}

\begin{proposition}\label{plussubal}%
Let $A$ be a normed \st-algebra. Let $B$ be a Hermitian complete
\st-subalgebra of $A$. (Cf.\ \ref{Herminher}.) For $b \in B$, we have
\[ b \in B_+ \Leftrightarrow b \in A_+.\]
\end{proposition}

\begin{proof}
$\s_B (b) \setminus \{ 0 \} = \s_A (b) \setminus \{ 0 \}$,
cf.\ \ref{specHermsubalg}.
\end{proof}

\begin{theorem}[(in-)stability of $A_+$ under multiplication]\label{instability}%
Let $A$ be a Hermitian Banach \st-algebra. If $a,b \in A_+$, then
the product $ab$ is in $A_+$ if and only if $a$ and $b$ commute.
\end{theorem}

\begin{proof}
For the ``only if'' part, see \ref{Hermprod}. Assume now that
$a, b \in A_+$ commute. We can then assume that $A$ is
commutative and unital by an application of the Theorem of
Civin \& Yood, cf.\ \ref{CivinYood}. Indeed, this result says that
there exists an automatically Hermitian \ref{Herminher}
commutative unital closed \st-subalgebra $B$ of $\tld{A}$
containing $a$ and $b$, such that $\s_A (c) = \s_B (c)$ for
$c \in \{ \,a, b, ab \,\}$. Now the elements $\frac{1}{n}e+a$
$(n \geq 1)$ have Hermitian square roots $a_n = {a_n}^*$,
cf.\ \ref{possqroot}. The elements $\frac{1}{n}e+b$ $(n \geq 1)$
have Hermitian square roots $b_n = {b_n}^*$. The elements
$a_n b_n$ are Hermitian by commutativity, and so have real
spectrum. Hence $\s \bigl({(a_n b_n)}^2\bigr) \subset \mathds{R}_+$
for all $n$ (Rational Spectral Mapping Theorem). \linebreak
The continuity of the spectrum function in a commutative
Banach \linebreak algebra \ref{commspeccont} shows that
\[ ab = \lim _{n \to \infty} {a_n}^2 \,{b_n}^2 =
\lim _{n \to \infty} {(a_n b_n) \,}^2 \geq 0. \qedhere \]
\end{proof}

\smallskip
Our next aim is the Shirali-Ford Theorem \ref{ShiraliFord}.
On the way, we shall see that the set of positive elements
of a Hermitian Banach \st-algebra $A$ is a closed convex
cone in $A\sa$.

\begin{proposition}%
Let $A$ be a Hermitian Banach \st-algebra. For two Hermitian
elements $a = a^*, b = b^*$ of $A$, we have
\[ \rlambda(ab) \leq \rlambda(a) \,\rlambda(b).\]
\end{proposition}

\begin{proof} 
By \ref{fundHerm} (i) $\Rightarrow$ (iii), we have
\[ \rlambda (ab) \leq \rsigma (ab) =
{\rlambda (baab) \,}^{1/2} = {\rlambda \bigl( \,{a \,}^2 \,{b\,}^2 \,\bigr) \,}^{1/2} \]
where we have made use of the fact that
$\rlambda (cd) = \rlambda (dc)$ for $c,d \in A$, cf.\ \ref{rlcomm}.
By induction it follows that
\[ \rlambda (ab) \leq
{\rlambda \Bigl( \,{a \,}^{( \,{2 \,}^{n})} \,{b \,}^{(\,{2\,}^{n})} \,\Bigr) \,}^{1 \,/ \,( \,{{2\,}^{n}} \,)}
\leq {\bigl|\,{a \,}^{( \,{2 \,}^{n})}\,\bigr| \,}^{1 \,/ \,( \,{{2 \,}^{n}} \,)}
\cdot {\bigl|\,{b \,}^{( \,{2 \,}^{n} \,)}\,\bigr| \,}^{1 \,/ \,( \,{{2 \,}^{n}} \,)} \]
for all $n \geq 1$. It remains to take the limit for $n \to \infty$. \end{proof}

\begin{theorem}\label{plusconvexcone}%
If $A$ is a Hermitian Banach \st-algebra
then $A_+$ is a convex cone in $A\sa$.
\end{theorem}

\begin{proof}
Let $a,b$ be positive elements of $A$. It shall be shown that
$a+b$ is positive. Upon replacing $a,b$ by suitable multiples, it
suffices to show that $e+a+b$ is invertible. Since $e+a$ and
$e+b$ are invertible, we may define
$c := a \,{(e+a)}^{-1}, d := b \,{(e+b)}^{-1}$. The Rational Spectral
Mapping Theorem yields $\rlambda(c) < 1$ and $\rlambda(d) < 1$.
\pagebreak
The preceding proposition then gives $\rlambda(cd)<1$, so $e-cd$
is invertible. We now have $(e+a)\,(e-cd)\,(e+b) = e+a+b$, so that
$e+a+b$ is invertible.
\end{proof}

\begin{proposition}\label{Hermrlcontsa}%
Let $A$ be a Hermitian Banach \st-algebra. Let $a = a^*, b = b^*$
be Hermitian elements of $A$. We then have
\[ \rlambda(a+b) \leq \rlambda(a)+\rlambda(b), \]
whence also
\[ | \,\rlambda(a)-\rlambda(b) \,| \leq \rlambda(a-b) \leq | \,a-b \,|.\]
It follows that $\rlambda$ is uniformly continuous on $A\sa$.
\end{proposition}

\begin{proof}
We have
\[ \rlambda(a) \,e \pm a \geq 0, \quad \rlambda(b) \,e \pm b \geq 0, \]
and so, by \ref{plusconvexcone}
\[ \bigl(\rlambda(a)+\rlambda(b)\bigr)\,e \pm (a+b) \geq 0, \]
whence
\[ \rlambda(a+b) \leq \rlambda(a)+\rlambda(b). \qedhere \]
\end{proof}

\begin{lemma}\label{critpos}%
Let $A$ be a Banach algebra, and let $a$ be an element of
$A$ with real spectrum. For $\tau \geq \rlambda(a)$, we have
\[ \min \,\s(a) = \tau - \rlambda ( \tau e - a ). \]
(Please note that $\s(\tau e-a) \subset [0, \infty[$.) Hence the
following statements are equivalent.
\begin{itemize}
  \item[$(i)$] $\s(a) \subset [0, \infty[$,
 \item[$(ii)$] $\rlambda(\tau e-a) \leq \tau$.
\end{itemize}
\end{lemma}

\begin{proof}
For $\tau \geq \rlambda(a)$, we obtain
\begin{align*}
\rlambda(\tau e-a) & = \max\,\{\,| \,\alpha \,| : \alpha \in \s(\tau e-a)\,\} \\
 & = \max \,\{ \,| \,\alpha \,| : \alpha\in\tau-\s(a) \,\} \\
 & = \max \,\{ \,| \,\tau-\lambda \,| : \lambda\in \s(a) \,\} \\
 & = \max \,\{ \,\tau-\lambda:\lambda\in \s(a) \,\} \\
 & = \tau - \min \,\s(a). \qedhere
\end{align*}
\end{proof}

\begin{theorem}\label{plusclosed}%
If $A$ is a Hermitian Banach \st-algebra then $A_+$ is closed in $A\sa$.
\pagebreak
\end{theorem}

\begin{proof}
Let $(a_n)$ be a sequence in $A_+$ converging to an
element $a$ of $A\sa$. There then exists $\tau \geq 0$ such
that $\tau \geq | \,a_n \,| \geq \rlambda(a_n)$ for all $n$. It
follows that also $\tau \geq | \,a \,| \geq \rlambda(a)$.
By \ref{Hermrlcontsa} and \ref{critpos}, we have
\[ \rlambda(\tau e-a) = \lim_{n \to \infty} \rlambda(\tau e-a_n)
\leq \tau. \qedhere \]
\end{proof}

\begin{theorem}%
If $A$ is a Hermitian Banach \st-algebra then
$A_+$ is a closed convex cone in $A\sa$.
\end{theorem}

\begin{proof}
This is the theorems \ref{plusconvexcone} and \ref{plusclosed}.
\end{proof}

\begin{theorem}\label{preShiraliFord}%
Let $A$ be a Banach \st-algebra. Let $\alpha \geq 2$ and put
$\beta := \alpha+1$. The following statements are equivalent.
\begin{itemize}
   \item[$(i)$] $A$ is Hermitian,
  \item[$(ii)$] for each $a \in A$ one has $\s(a^*a) \subset [-\lambda ,\lambda]$, \\
where $\lambda := \max \,\bigl( \,( \,\s(a^*a) \cap \mathds{R}_+ \,) \cup \{\,0\,\} \,\bigr)$,
 \item[$(iii)$] for all $a \in A$ with $\rsigma(a) \leq 1$, one has $\rsigma(c) \leq 1$, \\
where $c := {\alpha \,}^{-1} ( \beta a-aa^*a )$,
 \item[$(iv)$] $\s(a^*a)$ does not contain any $\lambda < 0$ whenever $a \in A$.
  \item[$(v)$] $e+a^*a$ is invertible in $\tld{A}$ for all $a \in A$.
\end{itemize}
\end{theorem}

\begin{remark}
In their definition of a C*-algebra, Gelfand and Na\u{\i}mark
had to assume that $e+a^*a$ is invertible in $\tld{A}$ for all
$a \in A$. The above theorem reduces this to showing that a
C*-algebra is Hermitian \ref{C*Herm}. Before going to the
proof, we note the following consequence.
\end{remark}

\begin{theorem}[the Shirali-Ford Theorem]\label{ShiraliFord}%
\index{Theorem!Shirali-Ford}\index{Shirali-Ford Theorem}%
A  Banach \st-algebra $A$ is Hermitian if and only if $a^*a \geq 0$ for all $a \in A$.
\end{theorem}

\begin{proof}[\hspace{-3.55ex}\mdseries{\scshape{Proof of \protect\ref{preShiraliFord}}}]%
The proof is a variant of the one in the book of Bonsall and
Duncan \cite{BD}, and uses a polynomial rather than a rational function.

(i) $\Rightarrow$ (ii). Assume that $A$ is Hermitian, let $a \in A$, and
let $\lambda$ be as in (ii). We obtain
$\s(a^*a) \subset \ ] - \infty , \lambda \,]$ as $a^*a$ is Hermitian. Let
now $a=b + \iu c$ where $b,c \in A\sa$. We then have
$aa^*+a^*a = 2({b\,}^2+{c\,}^2)$, whence
$\lambda e+aa^*= 2({b\,}^2+{c\,}^2)+(\lambda e-a^*a) \geq 0$ by convexity of
$\tld{A}_+$ \ref{plusconvexcone}, so that $\s(aa^*) \subset [ - \lambda , \infty [$.
One concludes that $\s(a^*a) \subset [ - \lambda , \lambda ]$ as
$\s(a^*a) \setminus \{ 0 \} = \s(aa^*) \setminus \{ 0 \}$, cf.\ \ref{speccomm}.
\pagebreak

Intermediate step. Consider the polynomial $p(x) = {\alpha \,}^{-2} \,x \,{( \beta -x) \,}^2$.
If $c$ is as in (iii), then $c^*c = p(a^*a)$, whence $\s(c^*c) = p\bigl(\s(a^*a)\bigr)$.
Please note that $p(1) = 1$ and hence $p(x) \leq 1$ for all $x \leq 1$ by monotony.

(ii) $\Rightarrow$ (iii). Assume that (ii) holds. Let $a \in A$ with
$\rsigma(a) \leq 1$. Then $\s(a^*a) \subset [-1,1]$ by (ii). Let $c$
be as in (iii). Then $\s(c^*c) \subset \,]-\infty,1]$ by the intermediate
step. Again by (ii), it follows that $\s(c^*c) \subset [-1,1]$.

(iii) $\Rightarrow$ (iv). Assume that (iii) holds and that (iv) does not hold. Let
$\delta := - \inf \,\bigcup _{\text{\smaller{$\,a \in A, \,\rsigma (a) \leq 1$}}}
\,\s(a^*a) \cap \mathds{R}$. We then have $\delta > 0$. It follows that there
exists $a \in A$, $\gamma > 0$, such that $\rsigma (a) \leq 1$,
$- \gamma \in \s(a^*a)$, and $\gamma \geq {( \alpha / \beta ) \,}^2 \,\delta$.
With $c$ as in (iii), we then have $\rsigma (c) \leq 1$ by assumption. The
definition of $\delta$ and the intermediate step then imply
$p( - \gamma ) \geq - \delta$, that is
${\alpha \,}^{-2} \,\gamma \,{( \beta + \gamma ) \,}^2 \leq \delta$. Hence, using
$\gamma > 0$, we obtain the contradiction $\gamma < {( \alpha / \beta ) \,}^2 \delta$.

(iv) $\Rightarrow$ (v) is trivial.

(v) $\Rightarrow$ (i). Assume that (v) holds. Let $a=a^*$ be a
Hermitian element of $A$. Then $-1 \notin \s({a \,}^2)$, so that
$i \notin \s(a)$ by the Rational Spectral Mapping Theorem. It
follows via \ref{fundHerm} (ii) $\Rightarrow$ (i) that $A$ is Hermitian.
\end{proof}

\begin{lemma}\label{rsHerm}%
Let $\pi : A \to B$ be a \st-algebra homomorphism from a Banach
\st-algebra $A$ to a Banach \st-algebra $B$. If $B$ is Hermitian, and if
\[ \rsigma \bigl(\pi (a)\bigr) = \rsigma(a)\quad \text{for all } a \in A, \]
then $A$ is Hermitian.
\end{lemma}

\begin{proof}
This follows from \ref{preShiraliFord} (i) $\Leftrightarrow$ (iii).
\end{proof}

\begin{theorem}[stability\,of $A_+$\,under \st-congruence]%
\label{stcongr}%
If $A$ is a Hermitian Banach \st-algebra, and $a \in A_+$,
then $c^*ac \in A_+$ for every $c \in \tld{A}$.
\end{theorem}

\begin{proof}
The elements $\frac{1}{n}e+a$ $(n \geq 1)$ have Hermitian square
roots $b_n = {b_n}^*$ in $\tld{A}\sa$, cf.\ \ref{possqroot}. It follows that
\[ c^*ac = \lim _{n \to \infty} c^* {b_n}^2 c =
\lim _{n \to \infty} {(b_n c)}^*{(b_n c)} \geq 0, \]
cf.\ \ref{ShiraliFord} and \ref{plusclosed}.
\end{proof}

\medskip
We also give the next result, in view of its strength, exhibited by
the ensuing corollary. \pagebreak

\begin{proposition}\label{sapreShiraliFord}%
Let $A$ be a Banach \st-algebra. Let $\alpha \geq 2$ and put
$\beta := \alpha+1$. The following statements are equivalent.
\begin{itemize}
   \item[$(i)$] $A$ is Hermitian,
  \item[$(ii)$] for each $a \in A\sa$ one has
  $\s \bigl({a \,}^2 \bigr) \subset [- \lambda , \lambda ]$, \\
  where $\lambda := \max \,\bigl( \,( \,\s \bigl({a \,}^2 \bigr) \cap \mathds{R}_+ \,) \cup \{\,0\,\} \,\bigr)$,
 \item[$(iii)$] for all $a \in A\sa$ with $\rlambda (a) \leq 1$, one has $\rlambda(c) \leq 1$, \\
  where $c := {\alpha \,}^{-1} \bigl( \beta a - {a \,}^3 \bigr)$,
 \item[$(iv)$] for all $a \in A\sa$ with $\rlambda(a) \leq 1$, and all $\lambda \in \s(a)$, \\
  one has $\bigl| \,{\alpha \,}^{-1} \lambda \,\bigl( \beta -{\lambda \,}^2 \bigr) \,\bigr| \leq 1$,
  \item[$(v)$] $\s \bigl({a \,}^2 \bigr)$ does not contain any $\lambda < 0$ for all $a \in A\sa$.
 \item[$(vi)$] $e + {a \,}^2$ is invertible in $\tld{A}$ for all $a \in A\sa$.
\end{itemize}
\end{proposition}

\begin{proof}
The proof follows the lines of the proof of theorem \ref{preShiraliFord},
cf.\ \ref{Hermrseqrl}.
\end{proof}

\begin{corollary}%
A Banach \st-algebra $A$ is Hermitian if (and only if) there exists a
compact subset $S$ of the open unit disc, such that $a \in A\sa$ and
$\rlambda (a) \leq 1$ together imply $\s(a) \subset \,[ \,-1,1 \,] \ \cup \ S$.
\end{corollary}

\begin{proof}
Let $A$ be a Banach \st-algebra. Assume that $S$ is a compact subset
of the open unit disc such that $a \in A\sa$ and $\rlambda (a) \leq 1$
together imply $\s(a) \subset \,[ \,-1,1 \,] \ \cup \ S$. To show that $A$ is
Hermitian, one applies (iv) of the preceding proposition with $\alpha \geq 2$
so large that
\[ \alpha/(\alpha+2)\geq\max\,\{\,|\,\lambda\,| : \lambda \in S \,\}. \]
Indeed, with $\beta := \alpha +1$ we have the following.
Let $a \in A\sa$ with $\rlambda (a) \leq 1$, and let $\lambda \in \s(a)$.
If $\lambda$ does not belong to $[ \,-1,1 \,]$, then it belongs to $S$, whence
\[|\,\lambda\,|\leq\alpha/(\alpha+2)=\alpha/(\beta+1)
\leq \alpha / \bigl| \,\beta - {\lambda \,}^2 \,\bigr|,\]
Also, if $\lambda$ belongs to $[ \,-1,1 \,]$, then
\[ \bigl| \,{\alpha \,}^{-1} \lambda \,\bigl( \beta - {\lambda \,}^2 \bigr) \,\bigr| \leq 1 \]
by differential calculus.
\end{proof}

\begin{remark}\label{remark2}
Some of the preceding results have easy proofs in pre\-sence of
commutativity, using multiplicative linear functionals.
Also, in \ref{instability}, the assumption that $A$ be Hermitian
can be dropped, cf.\ \ref{specsub}. Our proofs leading to these
facts however use the Lemma of Zorn, cf.\ the remark \ref{Zorn}.
\pagebreak
\end{remark}

\clearpage


\section{The Polar Factorisation of Invertible Elements}%
\label{polarfactoris}

\begin{definition}%
[polar factorisation]\label{polfactdef}\index{polar factorisation}%
If $A$ is a Hermitian Banach \linebreak \st-algebra, and if
$a \in \tld{A}$ is \underline{invertible}, then a pair
$(u, p) \in \tld{A} \times \tld{A}$ shall be called a
\underline{polar factorisation} of $a$, when
\[ a = u p, \quad \text{with } u \text{ unitary in } \tld{A},
\text{ and } p \text{ positive in } \tld{A}. \]
We stress that this is a good definition only for invertible elements.
\end{definition}

\begin{theorem}[existence and uniqueness of a polar factorisation]%
\label{polfactexuni}\index{absolute value}\index{a2@${"|}\,a\,{"|}$}%
\index{factorisation!polar}%
Let $A$ be a Hermitian Banach \st-algebra. An \underline{invertible}
element $a$ of $\tld{A}$ has a unique polar factorisation $(u, | \,a \,|)$.
The element $| \,a \,|$ is the unique positive square root of $a^*a$ in
$\tld{A}$. It is called the \underline{absolute value} of $a$. (Even though
the notation clashes with the notation for the norm.)
\end{theorem}

\begin{proof}
If $(u,p)$ is a polar factorisation of $a$, then $a^*a = {p\,}^2$.
Indeed $a^*a = (up)^*(up) = pu^*up = {p\,}^2$ since $p = p^*$
(as $p$ is positive), and because $u^*u = e$ (as $u$ is unitary).
In other words $p$ must be a positive square root of $a^*a$.
We shall show next existence and uniqueness of a positive
square root of $a^*a$ in $\tld{A}$.

The Shirali-Ford Theorem \ref{ShiraliFord} says that since $\tld{A}$
is Hermitian, we have $a^*a \geq 0$. In particular
$\s (a^*a) \subset [0, \infty[$. Also, with $a$, the elements $a^*$ and
$a^*a$ are invertible (with respective inverses ${\bigl({a \,}^{-1}\bigr)}^*$
and ${a \,}^{-1}{\bigl({a \,}^{-1}\bigr)}^*$). Thus $\s(a^*a) \subset \ ]0, \infty[$.
It follows from theorem \ref{possqroot} that $a^*a$ has a unique
square root with spectrum contained in $]0, \infty[$. Furthermore
this square root is Hermitian, cf.\ the same theorem. We shall denote
this square root with $| \,a \,|$. Then $| \,a \,|$ is an invertible positive
element of $\tld{A}$ with ${| \,a \,| \,}^2 = a^*a$. Please note also that
this square root also is unique within all positive elements of $\tld{A}$,
cf.\ the same theorem.

Since $| \,a \,|$ is invertible, the equation $a = u \,| \,a \,|$ implies
$u = a \,{| \,a \,| \,}^{-1}$. In particular, $u$ is uniquely determined,
and so there is at most one polar factorisation of $a$. We can
also use this formula to define $u$. We note that
$u := a \,{| \,a \,| \,}^{-1}$ satisfies $u^*u = e$ because
\[ u^*u = {\bigl({ \,| \,a \,| \,}^{-1} \,\bigr)}^*\,a^*\,a\,{| \,a \,| \,}^{-1}
=  {\bigl({ \,| \,a \,| \,}^{-1} \,\bigr)}^*\,{| \,a \,| \,}^2\,{| \, a \,| \,}^{-1} = e. \]
We also note that $u$ is invertible, as $a$ is invertible by assumption.
So $u$ has $u^*$ as a left inverse, as well as some right inverse.
It follows from \ref{leftrightinv} that $u^*$ is an inverse of $u$,
i.e.\ the element $u$ is unitary. In other words, a polar factorisation
of $a$ exists, and the proof is complete.
\end{proof}

The reader will recognise the terms below from Hilbert space theory.\pagebreak

\begin{definition}[idempotents and involutory elements]%
\index{involutory}\index{idempotent}%
Let $A$ be an algebra. One says that $p \in A$ is an
\underline{idempotent}, if ${p\,}^2 = p$. One says that
$a \in \tld{A}$ is \underline{involutory}, if ${a\,}^2 = e$,
that is, if $a = a^{\,-1}$.
\end{definition}

\begin{definition}[projections and reflections]%
\index{reflection}\index{projection}%
\index{orthogonal!projections}\index{projections!orthogonal}%
\index{complementary projections}\index{projections!complementary}%
Let $A$ be a \st-algebra.

One says that $p \in A$ is a \underline{projection},
if $p$ is a Hermitian idempotent: $p^* = p = {p\,}^2$.

Two projections $p$ and $q$ in $A$ are called
\underline{orthogonal}, if $p \,q = 0$. Then also
$q \,p = 0$ by taking adjoints.

Two projections $p$ and $q$ in $\tld{A}$ are
called \underline{complementary}, if they are
orthogonal, and if $p + q = e$.

One says that $u \in \tld{A}$ is a \underline{reflection},
if $u$ is both Hermitian and unitary: $u = u^* = u^{\,-1}$.
\end{definition}

\begin{example}%
If $p$ and $q$ are two complementary projections
in a unital \st-algebra, then $p - q$ is a reflection.
Indeed, one computes:
\[ (p - q)^*(p - q) = (p - q)(p - q)^* = {(p - q)\,}^2
= {p\,}^2 + {q\,}^2 - p \,q - q \,p = p + q = e. \]
\end{example}

\begin{theorem}[the structure of reflections]\label{structrefl}%
Let $u$ be a reflection in a unital \st-algebra. Then the
following statements hold.
\begin{itemize}
   \item[$(i)$] $u$ is involutory: ${u \,}^2 = e$,
  \item[$(ii)$] $\s(u) \subset \{ -1, 1\}$,
 \item[$(iii)$] we have $u = p - q$, where the elements
                       $p := \frac12 \,( \,e + u \,)$ and $q:= \frac12 \,( \,e - u \,)$
                       are complementary projections.
\end{itemize}
\end{theorem}

\begin{proof}
The element $u$ is involutory because $u = u^* = u^{\,-1}$ implies that
${u\,}^2 = e$. It follows from the Rational Spectral Mapping Theorem that
$\s(u) \subset \{ -1, 1\}$. Next, the elements $p$ and $q$ as in (iii)
obviously are Hermitian, and we have $p - q= u$ as well as $p + q = e$.
To show that $p$ and $q$ are idempotents, we note that
\[ {\biggl(\frac{e \pm u}{2}\biggr)}^2 = \frac{{e\,}^2 + {u\,}^2 \pm 2u}{4}
= \frac{2e \pm 2u}{4} = \frac{e \pm u}{2}. \]
To show that the projections $p$ and $q$ are orthogonal, we calculate:
\[ p \,q = \frac12 \,( \,e + u \,) \cdot \frac12 \,( \,e - u \,)
= \frac14 \,( \,{e\,}^2 - {u\,}^2 \,) = 0. \qedhere \]
\end{proof}

\begin{corollary}%
In a unital \st-algebra, the reflections are precisely the differences
of complementary projections. \pagebreak
\end{corollary}

\begin{theorem}[orthogonal decomposition]\label{anorthdeco}%
\index{orthogonal!decomposition}\index{decomposition!orthogonal}%
\index{a3@$a_+$}\index{a4@$a_-$}
If $A$ is a Hermitian Banach \st-algebra, and if $a$ is an
\underline{invertible} Hermitian element of $\tld{A}$,
there exists a unique decomposition $a = a _{+} - a _{-}$,
with $a _{+}$ and $a _{-}$ positive elements of
$\tld{A}$ satisfying $a _{+} \,a _{-} = a _{-} \,a _{+} = 0$,
namely $a _{+} = \frac12 \,\bigl( \,| \,a \,| + a \,\bigr)$ and
$a _{-} = \frac12 \,\bigl( \,| \,a \,| -a \,\bigr)$.
\end{theorem}

\begin{proof}
To show uniqueness, let $0 \leq a _{+}$, $a _{-} \in \tld{A}$
with $a = a _{+} - a _{-}$, and $a _{+} \,a _{-} = a _{-} \,a _{+} = 0$.
Then ${a\,}^2 = {(a _{+} - a _{-}) \,}^ 2 = {(a _{+} + a _{-}) \,}^2$, whence
$a _{+} + a _{-} \geq 0$ \ref{plusconvexcone} is a positive square
root of ${a \,}^2 = a^*a$. This makes that $a _{+} + a _{-} = | \,a \,|$,
the absolute value of $a$, cf.\ theorem \ref{polfactexuni}. The
elements $a _{+}$ and $a _{-}$ now are uniquely determined
by the linear system $a _{+} - a _{-} = a$, $a _{+} + a _{-} = | \,a \,|$,
namely $a _{+} = \frac12 \,\bigl( \,| \,a \,| + a \,\bigr)$,
$a _{-} = \frac12 \,\bigl( \,| \,a \,| -a \,\bigr)$.

Now for existence.
By theorem \ref{polfactexuni}, there exists a unique unitary
$u \in \tld{A}$ such that $a = u \,| \,a \,|$. It has also been
noted in the statement of this theorem that $| \,a \,|$ is the
unique positive square root of $a^*a = {a\,}^2$ in $\tld{A}$.
Please note that with $a$ and $u$, also the elements ${a \,}^2$
and $| \,a \,| = {u \,}^{-1} \,a$ are invertible. (The set of invertible
elements of $\tld{A}$ forming a group.)

We need to show that the elements $a$, $| \,a \,|$, ${| \,a \,| \,}^{-1}$, $u$
commute pairwise. Theorem \ref{possqroot} says
that $| \,a \,|$ is in the closed subalgebra of $\tld{A}$ generated
by $a^*a = {a\,}^2$. It follows that if $c \in \tld{A}$ commutes
with ${a\,}^2$, then $c$ also commutes with $| \,a \,|$, and thus
also with ${| \,a \,| \,}^{-1}$, cf.\ \ref{comm2}. Choosing $c := a$,
we see that $a$ commutes with $| \,a \,|$ and ${| \,a \,| \,}^{-1}$.
Hence $a$ also commutes with $u = a \,{| \,a \,| \,}^{-1}$. Choosing
$c := u$, we see that $u$ commutes with $| \,a \,|$ and
${| \,a \,| \,}^{-1}$. The elements $a$, $| \,a \,|$, ${| \,a \,| \,}^{-1}$, $u$
thus commute pairwise.

Our intuition says that $u$ should be Hermitian because $a$ is so.
Indeed $u = a { \,| \,a \,| \,}^{-1}$ is Hermitian by the commutativity
of $a$ and ${| \,a \,| \,}^{-1}$, cf.\ \ref{Hermprod}.
This makes that $u$ is a Hermitian unitary element, that is:
a reflection, and thus of the form $u = p - q$ with
$p$ and $q$ two complementary projections in
$\tld{A}$, cf.\ \ref{structrefl}. We can now put
$a _{+} := p \,| \,a \,|$, $a _{-} := q \,| \,a \,|$. Please note that a
projection in $\tld{A}$ is positive by the Rational Spectral
Mapping Theorem. So, to see that the elements $a _{+}$ and
$a _{-}$ are positive, we need to show that our two
projections commute with $| \,a \,|$, cf.\ theorem \ref{instability}.
This however follows from the commutativity of $u$ and $| \,a \,|$
by the special form of our projections, see theorem \ref{structrefl} (iii).
The commutativity of our two projections with $| \,a \,|$ is used
again in the following two calculations:
$a _{+} \,a _{-} = p \,| \,a \,| \,q \,| \,a \,| =  p \,q \,{| \,a \,| \,}^2 = 0$, and
$a _{-} \,a _{+} = q \,| \,a \,| \,p \,| \,a \,| = q \,p \,{| \,a \,| \,}^2 = 0$ by
the orthogonality of our two projections.
Also $a = u \,| \,a \,| = (p - q) \,| \,a \,| = a _{+} - a _{-}$. \pagebreak
\end{proof}

A stronger result holds in C*-algebras, cf.\ theorem \ref{orthdec} below.

We have again avoided the Gelfand transformation,
see our remarks \ref{Zorn} and \ref{remark2}.

\clearpage


\section{The Order Structure in a C*-Algebra}

In this paragraph, let $(A,\| \cdot \|)$ be a C*-algebra.

\begin{theorem}%
We have $A_+ \cap (-A_+) = \{ 0 \}$.
\end{theorem}

\begin{proof}
For $a \in A_+ \cap (-A_+)$, one has $\s(a) = \{ 0 \}$, whence
$\| \,a \,\| = \rlambda(a) = 0$, cf.\ \ref{C*rl}, so that $a = 0$.
\end{proof}

\begin{definition}[proper convex cone]%
A convex cone $C$ in a real vector space is called
\underline{proper} if $C \cap (-C) = \{ 0 \}$.
\end{definition}

The positive cone $A_+$ thus is a proper closed convex cone in $A\sa$.

\begin{definition}[ordered vector space]%
Let $B$ be a real vector space and let $C$
be a proper convex cone in $B$. By putting
\[ a \leq b :\Leftrightarrow b-a \in C \tag*{$(a,b \in B)$,} \]
one defines an order relation in $B$. In this way $B$ becomes
an \underline{ordered} \underline{vector space} in the following
sense:
\begin{align*}
    (i) &\quad a \leq b \Rightarrow a+c \leq b+c
    \tag*{$(a,b,c \in B)$} \\
   (ii) &\quad a \leq b \Rightarrow \lambda a \leq \lambda b
   \tag*{$(a,b \in B, \lambda \in \mathds{R}_+)$.}
\end{align*}
\end{definition}

The Hermitian part $A\sa$ thus is an ordered vector space.

\begin{proposition}\label{ordernorm}%
For $a \in A\sa$ we have $- \| \,a \,\| \,e \leq a \leq \| \,a \,\| \,e$.
\end{proposition}

\begin{theorem}[${a \,}^{1/2}$]%
\index{square root}\label{C*sqrootrestate}%
A positive element $a \in A$ has a unique positive square root,
which is denoted by $\underline{{a \,}^{1/2}}$. It belongs to the
closed subalgebra of $A$ generated by $a$. (Which also is the
 C*-subalgebra of $A$ generated by $a$.)
\end{theorem}

\begin{proof}
This is a restatement of theorem \ref{C*sqroot}.
\end{proof}

\begin{theorem}\label{sqrtShiraliFord}
\index{positive!element!in C*-algebra}%
For $a \in A$, the following statements are equivalent.
\begin{itemize}
   \item[$(i)$] $a \geq 0$,
  \item[$(ii)$] $a = {b\,}^2$ for some $b$ in $A_+$,
 \item[$(iii)$] $a = c^*c$ for some $c$ in $A$.
\end{itemize}
\end{theorem}

\begin{proof}
(iii) $\Rightarrow$ (i): the Shirali-Ford Theorem \ref{ShiraliFord}. \pagebreak
\end{proof}

\begin{theorem}%
Let $\pi : A \to B$ be a \st-algebra homomorphism from
$A$ to a \st-algebra $B$. Then $\pi$ is injective if and
only if $a \in A_+$ and $\pi (a) = 0$ imply $a = 0$. 
\end{theorem}

\begin{proof}
If $0 \neq b \in \ker \pi$, then $0 \neq b^*b \in \ker \pi \cap A_+$
by the C*-property \ref{preC*alg} and the preceding theorem.
\end{proof}

\begin{definition}\label{absval}%
\index{absolute value}\index{a2@${"|}\,a\,{"|}$}%
The \underline{absolute value} of $a \in A$ is defined as
\[ |\,a\,| := {(a^*a) \,}^{1/2}. \]
This is compatible with theorem \ref{polfactexuni} above.
\end{definition}

\begin{theorem}\label{homabs}%
If $\pi$ is a \st-algebra homomorphism
from $A$ to another C*-algebra, then
\[ \pi \,\bigl(\,|\,a\,|\,\bigr) = \bigl| \,\pi (a) \,\bigr| \qquad \text{for all } a \in A. \]
\end{theorem}

\begin{proof}
It suffices to note that $\pi \,\bigl(\,|\,a\,|\,\bigr)$ is the positive square root of
\linebreak
${\pi (a)}^*\pi (a)$. Indeed $\pi \,\bigl(\,|\,a\,|\,\bigr) \geq 0$ as $|\,a\,| \geq 0$,
cf.\ \ref{hompos}, and ${\pi \,\bigl(\,|\,a\,|\,\bigr) \,}^2 = \pi \,\bigl(\,{|\,a\,| \,}^2\,\bigr)
= \pi (a^*a) = {\pi (a)}^*{\pi (a)}$.
\end{proof}

\begin{proposition}\label{isonorm}%
For $a \in A$, we have
\[ \|\,|\,a\,|\,\| = \| \,a \,\| \quad \text{ as well as } \quad
\rsigma(a) = \rlambda \bigl(\,|\,a\,|\,\bigr). \]
\end{proposition}

\begin{proof}
One calculates with \ref{rlpowers}:
\[ {\|\,|\,a\,|\,\|\,}^2 = \,{\rlambda \bigl(\,|\,a\,|\,\bigr) \,}^2
            = \,\rlambda \bigl(\,{|\,a\,|\,}^2\,\bigr)
            = \,\rlambda (a^*a)
            = \,{\rsigma (a) \,}^2 = {\| \,a \,\|\,}^2. \qedhere \]
\end{proof}

\begin{proposition}\label{monotnorm}%
For $a$, $b$ $\in A\sa$, if $0 \leq a \leq b$, then $\|\,a\,\| \leq \|\,b\,\|$.
\end{proposition}

\begin{proof}
If $0 \leq a \leq b$ then $b \leq \|\,b\,\|\,e$ and
$b \pm a \geq 0$. It follows that $\|\,b\,\| \,e \pm a \geq 0$,
whence $\|\,b\,\| \geq \rlambda(a) = \|\,a\,\|$.
\end{proof}

\begin{corollary}%
For $a$, $b$ $\in A$, if $|\,a\,| \leq |\,b\,|$, then also $\|\,a\,\| \leq \|\,b\,\|$.
\end{corollary}

\begin{proof}
If $|\,a\,| \leq |\,b\,|$, then $\|\,|\,a\,|\,\| \leq \|\,|\,b\,|\,\|$, cf.\ \ref{monotnorm},
which is the same as $\|\,a\,\| \leq \|\,b\,\|$, cf.\ \ref{isonorm}. \pagebreak
\end{proof}

\begin{theorem}[the orthogonal decomposition]\label{orthdec}%
\index{orthogonal!decomposition}\index{decomposition!orthogonal}%
\index{a3@$a_+$}\index{a4@$a_-$}%
For a  Hermitian element $a = a^*$ of $A$, one defines
\[ a_+ := \frac{1}{2} \,\bigl(\,|\,a\,|+a\,\bigr)\quad and 
\quad a_- := \frac{1}{2} \,\bigl(\,|\,a\,|-a\,\bigr). \]
We then have
\begin{itemize}
   \item[$(i)$] $a_+, a_- \in A_+$,
  \item[$(ii)$] $a = a_+ - a_-$,
 \item[$(iii)$] $a_+ \,a_- = a_ - \,a_+ = 0$,
\end{itemize}
and these properties characterise $a_+$, $a_-$ uniquely.
\end{theorem}

This is compatible with theorem \ref{anorthdeco} above. 

\medskip
\begin{proof} We shall first show uniqueness. So let $a_1, a_2 \in A_+$
with $a = a_1-a_2$ and $a_1\,a_2 = a_2 \,a_1= 0$. We obtain
${|\,a\,| \,}^2 = {a \,}^2 = {(a_1-a_2) \,}^2 = {(a_1+a_2) \,}^2$, whence
$|\,a\,| = a_1+a_2$ by the uniqueness of the positive square root in $A$,
and so $a_1 = a_+$, $a_2 = a_-$.

We shall next prove that $a_+$ and $a_-$ do enjoy these properties.
(ii) is clear. (iii) holds because for example
$a_+ \,a_- = \frac{1}{4} \,\bigl(\,|\,a\,|+a\,\bigr) \,\bigl(\,|\,a\,|-a\,\bigr)
= \frac{1}{4} \,\bigl(\,{|\,a\,| \,}^2-{a \,}^2\,\bigr)= 0$.
For the proof of (i) we make use of the weak form of the operational
calculus \ref{weakopcalc}.
Let ``$\mathrm{abs}$'' denote the function $t \to | \,t \,|$ on $\s(a)$,
and let ``$\id$'' denote the identity function on $\s(a)$.
We have $a = \id (a)$ by a defining property of the operational
calculus. We also have
$| \,a \,| = {(a^*a) \,}^{1/2} = {\bigl({ \,a \,}^2 \,\bigr) \,}^{1/2} = \mathrm{abs}(a)$
by \ref{calcplus} together with the uniqueness of the positive
square root. Now consider the functions on $\s(a)$ given by
$f_+ := \frac{1}{2} \,( \mathrm{abs} + \id )$
and $f_- := \frac{1}{2} \,( \mathrm{abs} - \id )$.
We then have $f_+ (a), f_- (a) \in A_+$ by \ref{calcplus}.
From the homomorphic nature of the operational calculus, we have
$a = f_+ (a) - f_- (a)$ as well as $f_+ (a) f_- (a) = f_- (a) f_+(a) = 0$.
It follows from the uniqueness property shown above that
$f_+ (a) = a_+$ and $f_- (a) = a_-$, and the statement follows.
\end{proof}

\begin{corollary}\label{generic}%
The positive part $A_+$ of $A$ is a generic
convex cone in $A\sa$, i.e.\ $A\sa = A_+ +(-A_+)$.
\end{corollary}

\begin{corollary}\label{sandwich}%
For a Hermitian element $a = a^*$ of $A$, we have
\[ -| \,a \,| \leq a \leq | \,a \,|. \pagebreak \]
\end{corollary}

\begin{corollary}\label{product}%
Every Hermitian element of the C*-algebra $A$ is the
product of two Hermitian elements of $A$. Please note
that the factors necessarily commute, cf.\ \ref{Hermprod}.
\end{corollary}

\begin{proof} Let $a = a^*$ be a Hermitian element of $A$.
We can assume that $A$ is commutative by going over to the
C*-subalgebra of $A$ generated by $a$, which is commutative
as $a$ is Hermitian. We then compute
\begin{align*}
a = a_+ - a_-
   & = \Bigl[ \,{{( \,a_+ \,)}^{\,1/2}} \,\Bigr]^{\,2} - \Bigl[ \,{{( \,a_- \,)}^{\,1/2}} \,\Bigr]^{\,2} \\
   & = \Bigl[ \,{{( \,a_+ \,)}^{\,1/2} + {( \,a_- \,)}^{\,1/2}} \,\Bigr]
\cdot \Bigl[ \,{{( \,a_+ \,)}^{\,1/2} - {( \,a_- \,)}^{\,1/2}} \,\Bigr].
\qedhere
\end{align*}
\end{proof}

\medskip
We now turn to special C*-algebras.

\begin{proposition}\label{Coabs}%
For a locally compact Hausdorff space $\Omega \neq \varnothing$,
\linebreak
consider the C*-algebra $A := C\0(\Omega)$. Let $f \in A$. Then the
two meanings of $| f |$ coincide. (One being the absolute value in the
C*-algebra $A$, the other being the absolute value of the function.)
\end{proposition}

\begin{proof}
Let $| f |$ be the absolute value of the function $f$. Then $| f | \in A_+$
because $\left| f \right| = \sqrt{ \left| f \right| } \sqrt{ \left| f \right| }$.
(Please note that $\sqrt{ \left| f \right| }$ is a  Hermitian element of
$C\0(\Omega)$.) From ${| \,f \,| \,}^2 = \overline{f}f$, one sees that $| f |$
is a positive square root of $\overline{f}f$. The statement now follows
from the uniqueness of the positive square root.
\end{proof}

\begin{corollary}\label{pointwiseorder}
For a locally compact Hausdorff space $\Omega \neq \varnothing$,
the order in $C\0(\Omega)\sa$ is the pointwise order.
\end{corollary}

For the C*-algebra $B(H)$ with $H$ a Hilbert space, we have:

\begin{theorem}\label{posop}%
\index{positive!bounded operator}\index{operator!positive bounded}%
Let $H$ be a Hilbert space and let $a \in B(H)$.
The following properties are equivalent.
\begin{itemize}
  \item[$(i)$] $a \geq 0$,
 \item[$(ii)$] $\langle ax, x \rangle \geq 0$ for all $x \in H$.
\end{itemize}
\end{theorem}

\begin{proof} Each of the above two properties implies that $a$ is Hermitian. \\
(i) $\Rightarrow$ (ii). If $a \geq 0$, there exists
$b \in B(H)\sa$ such that $a = {b\,}^2$, and so
\[ \langle ax, x \rangle = \langle {b\,}^2x, x \rangle =
\langle bx, bx \rangle \geq 0\quad \text{for all } x \in H. \pagebreak \]
(ii) $\Rightarrow$ (i). Assume that $\langle ax, x \rangle \geq 0$ for
all $x \in H$ and let $\lambda <  0$. Then
\[ {\| \,(a-\lambda \mathds{1})x \,\| \,}^2 =
\ {\| \,ax \,\| \,}^2 - 2 \,\lambda \,\langle ax, x \rangle + {\lambda \,}^2 \,{\| \,x \,\| \,}^2
\geq \,{\lambda \,}^2 \,{\| \,x \,\| \,}^2. \]
This implies that $a-\lambda \mathds{1}$ has a bounded left inverse. 
In particular, the null space of $a-\lambda \mathds{1}$ is $\{ 0 \}$.
It follows that the range of $a-\lambda \mathds{1}$ is dense in $H$,
as $a-\lambda \mathds{1}$ is Hermitian. So the left inverse of
$a-\lambda \mathds{1}$ has a continuation to a bounded linear
operator defined on all of $H$. By taking adjoints, it follows that
$a-\lambda \mathds{1}$ also is right invertible, and thus invertible
by \ref{leftrightinv}.
\end{proof}

\begin{theorem}[the order completeness of $B(H)\sa$]\label{ordercompl}%
Let $H$ be a Hilbert space. Let $(a_i)_{i \in I}$ be an increasing
net of Hermitian operators in $B(H)\sa$ which is bounded above
in $B(H)\sa$. There then exists in $B(H)\sa$ the supremum $a$
of the $a_i$ $(i \in I)$. We have
\[ ax = \lim _{i \in I} a_ix\quad \text{for all } x \in H. \]
\end{theorem}

\begin{proof} For $x \in H$, the net
$\bigl( \langle a_ix,x \rangle \bigr)_{i \in I}$
is increasing and bounded above, so that
$\lim _{i \in I} \,\langle a_ix,x \rangle =
\sup _{i \in I} \,\langle a_ix,x \rangle$
exists and is finite. By polarisation, we may define
a Hermitian sesquilinear form $\varphi$ on $H$ via
\[ \varphi (x,y) = \lim _{i \in I} \,\langle a_ix,y \rangle
\quad (x,y \in H). \]
This Hermitian form is seen to be bounded in the sense that
\[ \sup _{\|\,x\,\|,\,\|\,y\,\| \,\leq \,1} \,|\,\varphi \,(x,y)\,|
< \infty. \]
Hence there is a unique $a \in B(H)\sa$ with
\[ \langle ax, y \rangle = \varphi (x,y) = \lim _{i \in I}
\,\langle a_ix,y \rangle\quad \text{for all } x,y \in H. \]
For $x \in H$ we get
$\langle ax,x \rangle = \sup _{i \in I} \,\langle a_ix,x \rangle$,
and it follows by the preceding theorem that $a$ is the
supremum of the set $\{a_i:i\in I\}$. It remains to be shown that
the net $(a_i)_{i \in I}$ converges pointwise to $a$. For this
purpose, we use that $a-a_i \geq 0$ for all $i \in I$.
By \ref{stcongr} we have
\begin{align*}
{(a-a_i) \,}^2 = & \ {(a-a_i) \,}^{1/2} \,(a-a_i) \,{(a-a_i) \,}^{1/2} \\
\leq & \ {(a-a_i) \,}^{1/2} \,\| \,a-a_i\,\| \,{(a-a_i) \,}^{1/2} = \|\,a-a_i\,\| \,(a-a_i).
\end{align*}
Whence, for $x \in H$,
\begin{align*}
{\|\,(a-a_i) \,x\,\| \,}^2 = & \ \langle {(a-a_i) \,}^2 \,x, x \rangle \\
 \leq & \ \|\,a-a_i\,\| \,\langle (a-a_i) \,x, x \rangle.
\end{align*}
The statement follows because any section
$( \,\|\,a-a_i\,\| \,)_{\text{\small{$i \geq i\0$}}}$
is decreasing by \ref{monotnorm}, and thus bounded. \pagebreak
\end{proof}

\medskip
We return to abstract C*-algebras with the following theorem.
It will play a decisive role in the next paragraph.

\begin{theorem}[monotony of the inverse operation]%
\label{monotinv}%
Let $a,b \in \tld{A}\sa$ with $0 \leq a \leq b$ and with $a$ invertible.
Then $b$ is invertible with $0 \leq {b \,}^{-1} \leq {a \,}^{-1}$.
\end{theorem}

\begin{proof}
Since $a \geq 0$ is invertible, we have $a \geq \alpha e$ for some
$\alpha > 0$, so that also $b \geq \alpha e$, from which it follows that
$b \geq 0$ is invertible. Also ${b \,}^{-1} \geq 0$ by the Rational Spectral
Mapping Theorem. On one hand, we have
\[ 0 \leq a \leq b = {b \,}^{1/2} \,e \,{b \,}^{1/2} \]
whence, by \ref{stcongr},
\[ 0 \leq {b \,}^{-1/2} \,a \,{b \,}^{-1/2} \leq e \]
so that
\[ \bigl\| \,{b \,}^{-1/2} \,a \,{b \,}^{-1/2} \,\bigr\|
= \rlambda \bigl( \,{b \,}^{-1/2} \,a \,{b \,}^{-1/2} \,\bigr) \leq 1. \]
On the other hand, we have
\[ \bigl\| \,{a \,}^{1/2} \,{b \,}^{-1} \,{a \,}^{1/2} \,\bigr\|
\leq \bigl\| \,{a \,}^{1/2} \,{b \,}^{-1/2} \,\bigr\|
\cdot \bigl\| \,{b \,}^{-1/2} \,{a \,}^{1/2} \,\bigr\|. \]
By the isometry of the involution and the C*-property \ref{preC*alg},
it follows that
\[ \bigl\| \,{a \,}^{1/2} \,{b \,}^{-1} \,{a \,}^{1/2} \,\bigr\|
\leq {\bigl\| \,{a \,}^{1/2} \,{b \,}^{-1/2} \,\bigr\| \,}^{2} =
\bigl\| \,{b \,}^{-1/2} \,a \,{b \,}^{-1/2} \,\bigr\| \leq 1. \]
This implies
\[ 0 \leq {a \,}^{1/2} \,{b \,}^{-1} \,{a \,}^{1/2} \leq e \]
which is the same as $b^{-1} \leq a^{-1}$, by \ref{stcongr} again.
\end{proof}

\medskip
And now for something completely different:

\begin{definition}[positive linear functionals]%
A linear  functional $\varphi$ on a \st-algebra is called
\underline{positive} if $\varphi (a^*a) \geq 0$ for all elements $a$.
\end{definition}

A linear functional $\varphi$ on the C*-algebra $A$ is positive
if and only if $\varphi (a) \geq 0$ for all $a \in A_+$,
cf.\ \ref{sqrtShiraliFord}.

\begin{theorem}\label{plfC*bded}%
\index{positive!functional!on C*-algebra}%
Every positive linear functional on the  C*-algebra $A$ is
Hermitian and continuous. \pagebreak
\end{theorem}

\begin{proof} Let $\varphi$ be a positive linear functional on $A$.
The fact that $\varphi$ is Hermitian \ref{Hermfunct}
follows from the orthogonal decomposition \ref{orthdec}. It shall
next be shown that $\varphi$ is bounded on the positive part
of the unit ball of $A$. If it were not, one could choose a sequence
$(a_n)$ in the positive part of the unit ball of $A$ such that
$\varphi (a_n) \geq n^2$ for all $n$. The elements
\[ a := \sum _n \frac{1}{n^2} \,a_n,\quad
   b_m := \sum _{n \leq m} \frac{1}{n^2} \,a_n,\quad
   a - b_m = \sum _{n \geq m+1} \frac{1}{n^2} \,a_n \]
would all lie in the closed convex cone $A_+$. It would follow that
\[ \varphi (a) \geq \varphi (b_m) =
    \sum _{n \leq m} \frac{\varphi (a_n)}{n^2} \geq m \]
for all $m$, a contradiction. We obtain
\[ \mu := \sup \,\bigl\{ \,\varphi (a) : a \in A_+, \| \,a \,\| \leq 1 \,\bigr\}
< \infty. \]
For $b$ in $A\sa$, we have $-|\,b\,| \leq b \leq |\,b\,|$,
cf.\ \ref{sandwich}, whence
\[ | \,\varphi (b) \,| \leq \varphi \bigl(\,|\,b\,|\,\bigr)
\leq \mu \,\| \,|\,b\,| \,\| = \mu \,\|\,b\,\|. \]
For $c$ in $A$ and $c = a + \iu b$ with $a,b \in A\sa$, we have
$\|\,a\,\|, \|\,b\,\| \leq \|\,c\,\|$, and so $|\,\varphi (c) \,| \leq 2 \,\mu \,\| \,c \,\|$.
\end{proof}

\clearpage


\section{The Canonical Approximate Unit in a C*-Algebra}

\begin{definition}%
A \underline{two-sided approximate unit} in a normed
algebra $A$ is a net $(e_j)_{j \in J}$ in $A$, such that
\[ \lim _{\text{\footnotesize{$j \in J$}}} \,e_j \,a
= \lim _{\text{\footnotesize{$j \in J$}}} \,a \,e_j
= a \quad \text{for all } a \in A. \]
\end{definition}

\begin{theorem}\label{C*canapprunit}%
\index{unit!approximate!canonical}\index{J(A)@$J(A)$}%
\index{canonical!approximate unit}\index{approximate unit!canonical}%
Let $A$ be a C*-algebra. We put
\[ J(A) := \{\,j \in A_+ : \|\,j\,\| < 1\,\}. \]
Then $(j)_{j \in J(A)}$ is a two-sided approximate unit in $A$.
It is called the \linebreak \underline{canonical approximate unit}.
\end{theorem}

\begin{proof}
It shall first be shown that $J(A)$ is upward directed. Let $a,b \in J(A)$.
The elements
\[ a_1 := a \,{( \,e - a \,) \,}^{-1}, \quad b_1 := b \,{( \,e - b \,) \,}^{-1} \]
then belong to $A_+$ by the Rational Spectral Mapping Theorem.
Now consider
\[ c := ( \,a_1 + b_1 \,) \,{\bigl[ \,e + ( \,a_1 + b_1 \,) \,\bigr] \,}^{-1}. \]
Then $c$ is in $J(A)$. Moreover we have
\begin{align*}
 & \ e-{\bigl[ \,e + ( \,a_1 + b_1 \,) \,\bigr] \,}^{-1} \\
 = & \ \bigl[ \,e + ( \,a_1 + b_1 \,) \,\bigr]
 \,{\bigl[ e + ( \,a_1 + b_1 \,) \,\bigr] \,}^{-1} - {\bigl[ \,e + ( \,a_1 + b_1 \,) \,\bigr] \,}^{-1} \\
 = & \ \bigl[ \,e + ( \,a_1 + b_1 \,) -e \,\bigr] \,{\bigl[ \,e + ( \,a_1 + b_1 \,) \,\bigr] \,}^{-1} = c.
\end{align*}
From the monotony of the inverse operation \ref{monotinv} it follows that
\begin{align*}
c = e-{\bigl[ \,e + ( \,a_1 + b_1 \,) \,\bigr] \,}^{-1} & \geq e - {( \,e + a_1 \,) \,}^{-1} \\
& \geq e - {\bigl[ \,e + a \,{( \,e - a \,) \,}^{-1} \,\bigr] \,}^{-1} \\
& \geq e - {\bigl[ \,( \,e - a + a \,) \,{( \,e - a \,) \,}^{-1} \,\bigr] \,}^{-1} \\
& \geq e - ( \,e - a \,) = a,
\end{align*}
and similarly $c \geq b$, which proves that $J(A)$ is upward directed.

It shall be shown next that the net $(j)_{j \in J(A)}$ is a two-sided approximate
unit in $A$. We claim that it suffices to show that for each $a \in A_+$ one has
\[ \lim _{\text{\footnotesize{$j \in J(A)$}}} \| \,a \,( \,e - j \,) \,a \,\| = 0. \pagebreak \]
To see this, we note the following. For $a \in A\sa$, and $j \in J(A)$, we have
\begin{align*}
    & \ {\| \,a - j \,a \,\| \,}^2 = {\| \,( \,e - j \,) \,a \,\| \,}^2 \\
 = & \ \| \,a \,{( \,e - j \,) \,}^2 \,a\,\| = \| \,a \,{( \,e - j \,) \,}^{1/2}
 \,( \,e - j \,) \,{( \,e - j \,) \,}^{1/2} \,a \,\| \\
\leq & \ \| \,a \,{( \,e - j \,) \,}^{1/2} \,{( \,e - j \,) \,}^{1/2} \,a \,\| = \| \,a \,( \,e - j \,) \,a \,\|.
\end{align*}
One also uses the fact that $A$ is linearly spanned by $A_+$,
cf.\ \ref{generic}. The isometry of the involution is used to pass
from multiplication on the left to multiplication on the right. 

So let $a \in A_+$, $\varepsilon > 0$. Put
\[ \alpha := \bigl( \,1+ \|\,a\,\| \,\bigr) / \varepsilon > 0, \]
and consider
\[ j\0 := \alpha a \,{( \,e + \alpha a \,) \,}^{-1}, \]
which belongs to $J(A)$ by the Rational Spectral Mapping Theorem.
For $j \in J(A)$ with $j \geq j\0$, we get
\begin{align*}
0 \leq a \,( \,e - j \,) \,a \leq a \,( \,e - j\0 \,) \,a
= & \ {a \,}^2 \,\bigl[ \,e - \alpha a \,{( \,e + \alpha a \,) \,}^{-1} \,\bigr] \\
 = & \ {a \,}^2 \,\bigl[ \,( \,e + \alpha a - \alpha a \,) \,{( \,e + \alpha a \,) \,}^{-1} \,\bigr] \\
 = & \ {a \,}^2 \,{( \,e + \alpha a \,) \,}^{-1} \\
 = & \ {\alpha \,}^{-1}\,a \cdot {\alpha \,a} \,{( \,e + \alpha a \,) \,}^{-1} \\
 \leq & \ {\alpha \,}^{-1}\,a,
\end{align*}
where the last inequality follows from the
Rational Spectral Mapping Theorem. Hence
\[ \|\,a \,( \,e - j \,) \,a \,\| \leq {\alpha \,}^{-1} \,\| \,a \,\|
= \varepsilon \cdot \frac{\| \,a \,\|}{1 + \|\,a\,\|} < \varepsilon. \qedhere \]
\end{proof}

\clearpage


\section{Quotients and Images of C*-Algebras}

\begin{reminder}[unital homomorphism]%
Let $A, B$ be unital algebras. Recall from \ref{unitalhom} that an algebra
homomorphism $A \to B$ is called \underline{unital} if it maps unit to unit.
\end{reminder}

\begin{proposition}[unitisation of a homomorphism]%
\label{homunit}\index{homomorphism!unitisation}%
\index{unitisation!of homomorphism}%
Let $A, B$ be normed algebras, and let $\pi : A \to B$ be a
non-zero algebra homomorphism with dense range. The mapping
$\pi$ then has a (necessarily unique) extension to a unital algebra
homomorphism $\tld{\pi} : \tld{A} \to \tld{B}$. We shall say that $\tld{\pi}$
is the \underline{unitisation} of $\pi$. If $\pi$ is injective, so is $\tld{\pi}$.
\end{proposition}

\begin{proof}
If $A$ has a unit $e_A$, then $\pi (e_A) \neq 0$ is a unit in $B$, by
density of $\pi (A)$ in $B$. If $A$ has no unit, we extend $\pi$ to a
unital algebra homomorphism $\tld{A} \to \tld{B}$ by requiring that
unit be mapped to unit. Injectivity then is preserved because the unit
in $\tld{B}$ is linearly independent of $\pi (A)$, as $\pi (A)$ has no
unit by injectivity.
\end{proof}

\begin{lemma}\label{injpos}%
Let $A$, $B$ be C*-algebras.
Let $\pi : A \to B$ be an injective \st-algebra homomorphism.
For $a \in A$, we then have
\[ \pi (a) \geq 0 \Leftrightarrow a \geq 0. \]
\end{lemma}

\begin{proof}
We know that $a \geq 0 \Rightarrow \pi (a) \geq 0$, cf.\ \ref{hompos}.
Let now $a \in A$ with $\pi (a) \geq 0$. Then
$\pi(a) = \bigl| \,\pi(a) \,\bigr| = \pi \,\bigl( \,| \,a \,| \,\bigr)$,
cf.\ \ref{homabs}. Hence $a = | \,a \,| \geq 0$ by injectivity of $\pi$,
as was to be shown.
\end{proof}

\begin{lemma}\label{posposrs}%
Let $A$ be a unital Banach \st-algebra. Let $B$ be a
Hermitian unital Banach \st-algebra. Let $\pi : A \to B$
be a unital \st-algebra homomorphism. Assume that
\[ \pi (b) \geq 0 \Leftrightarrow b \geq 0
\qquad \text{for every } b \in A\sa. \]
Then
\[ \rsigma \bigl( \pi (a) \bigr) = \rsigma (a)
\qquad \text{for every } a \in A.  \]
\end{lemma}

\begin{proof}
We have $\rlambda \bigl(\pi (a)\bigr) \leq \rlambda (a)$ for all $a \in A$ because
\[ \s_B \bigl(\pi(a)\bigr) \setminus \{ 0 \} \subset \s_A (a) \setminus \{ 0 \}, \]
cf.\ \ref{spechom}.
\pagebreak

Let now $b$ be a Hermitian element of $A$. Then
\[ \rlambda \bigl(\pi(b)\bigr) \,e_B \pm \pi(b) \geq 0 \]
as $B$ is Hermitian. Our assumption then implies
\[ \rlambda \bigl(\pi(b)\bigr) \,e_A \pm b \geq 0, \]
whence $\rlambda (b) \leq \rlambda \bigl(\pi(b)\bigr)$. Collecting inequalities,
we find  $\rlambda \bigl(\pi(b)\bigr) = \rlambda(b)$ for every Hermitian
element $b$ of $A$, which is enough to prove the statement.
\end{proof}

\begin{theorem}\label{C*injisometric}%
An injective \st-algebra homomorphism $\pi$  from a \linebreak
C*-algebra $A$ to a pre-C*-algebra $B$ is isometric.
As a consequence, $\pi (A)$ is complete.
\end{theorem}

\begin{proof}
We can assume that $B$ is a C*-algebra, and that $\pi$ has dense
range and is non-zero. Then the unitisation $\tld{\pi}$ is well-defined
and injective by proposition \ref{homunit}. The preceding two lemmata
imply that $\rsigma \bigl(\pi(a)\bigr) = \rsigma(a)$ for all $a \in A$,
which is the same as $\| \,\pi(a) \,\| = \| \,a \,\|$ for all $a \in A$, cf.\ \ref{C*rs}.
\end{proof}

\begin{remark}\label{reminjiso}
This is a strong improvement of \ref{Cstarisom} in so far
as it allows us to conclude that $\pi(A)$ is a C*-algebra.
\end{remark}

Next an expression for the quotient norm.
(Cf.\ the appendix \ref{quotspace}.)
Here we use the canonical approximate unit,
see the preceding paragraph.

\begin{lemma}[quotient norm]\label{quotexpr}%
Let $I$ be a \st-stable \ref{selfadjointsubset} closed \linebreak
two-sided ideal in a C*-algebra $A$. For $a \in A$ we then have
\[ \inf _{\text{\footnotesize{$c \in I$}}} \,\| \,a + c \,\|
= \lim _{\text{\footnotesize{$j \in J(I)$}}} \,\| \,a - a \,j \,\|. \]
\end{lemma}

\begin{proof}
Let $a \in A$. Denote $\alpha := \inf _{\,\text{\footnotesize{$c \in I$}}} \,\| \,a + c \,\|$.
For $\varepsilon > 0$, take $c \in I$ with
$\| \,a + c \,\| < \alpha + \varepsilon / 2$. Then choose
$j\0 \in J(I)$ such that $\| \,c - c \,j \,\| < \varepsilon / 2$
whenever $j \in J(I)$, $j \geq j\0$. Under
these circumstances we have
\begin{align*}
\| \,a - a \,j \,\| \leq & \ \| \,( \,a + c \,) \,( \,e - j \,)\,\| + \|\,c - c \,j \,\| \\
      \leq & \ \| \,a + c \,\| + \|\,c - c \,j \,\| < \alpha + \varepsilon. \pagebreak \qedhere
\end{align*}
\end{proof}

\begin{theorem}[factorisation]%
\index{quotient!of C*-algebra}\label{C*quotient}%
Let $I$ be a \st-stable \ref{selfadjointsubset} closed \linebreak
two-sided ideal in a C*-algebra $A$. Then $A/I$ is a C*-algebra.
\end{theorem}

\begin{proof}
It is clear that $A/I$ is a Banach \st-algebra, cf.\ the appendix
\ref{quotspaces}. The preceding lemma \ref{quotexpr} gives 
the following expression for the square of the quotient norm:
\begin{align*}
{| \,a + I \,| \,}^2 &
= \inf _{\text{\footnotesize{$c \in I$}}} {\| \,a + c \,\| \,}^2 \\
 & = \lim _{\text{\footnotesize{$j \in J(I)$}}} \,{\| \,a \,( \,e - j \,)\,\| \,}^2 \\
 & = \lim _{\text{\footnotesize{$j \in J(I)$}}} \,\| \,( \,e - j \,) \,a^*a \,( \,e - j \,)\,\| \\
 & \leq \lim _{\text{\footnotesize{$j \in J(I)$}}} \,\| \,a^*a \,( \,e - j \,)\,\|
 = \inf _{\text{\footnotesize{$c \in I$}}} \| \, a^*a+c \,\|
 = |\,a^*a+I\,|.
\end{align*}
The statement follows now from \ref{condC*}.
\end{proof}

\begin{theorem}%
Let $\pi$ be a \st-algebra homomorphism from a C*-algebra
$A$ to a pre-C*-algebra. Then $A \,/ \ker \pi$ is a C*-algebra.
The \st-algebra isomorphism
\begin{alignat*}{2}
 & A \,/ \ker \pi & \ \to & \ \pi (A)  \\
 & a+\ker \pi & \ \mapsto & \ \pi (a)
\end{alignat*}
is isometric, and hence the range $\pi (A)$ is a C*-algebra.
\end{theorem}

\begin{proof} It follows from \ref{Cstarcontr} that $\ker \pi$
is closed, so that $A \,/ \ker \pi$ is a \linebreak C*-algebra by
the preceding result. The statement then follows from
\ref{C*injisometric}.
\end{proof}

\medskip
In other words:

\begin{theorem}[factorisation of homomorphisms]\label{homfact}%
Let $\pi$ be a \linebreak \st-algebra homomorphism from a C*-algebra
$A$ to a pre-C*-algebra. Then $A \,/ \ker \pi$ and $\pi(A)$ are
C*-algebras as well. The \st-algebra homomorphism $\pi$ factors
to an isomorphism of C*-algebras from $A \,/ \ker \pi$ onto $\pi(A)$.
\end{theorem}

In particular, the C*-algebras form a subcategory of the category of
\st-algebras, whose morphisms are the \st-algebra homomorphisms.
\pagebreak

\clearpage


\part{Representations and Positive Linear Functionals}\label{part2}

\chapter{General Properties of Representations}


\setcounter{section}{21}

\section{Continuity Properties of Representations}

For a complex vector space $V$, we shall denote by
$\mathrm{End}(V)$ the algebra of linear operators on $V$.
Here ``$\mathrm{End}$'' stands for endomorphism.

\begin{definition}[representation in a pre-Hilbert space]%
\index{representation}\index{R3@$R(\pi)$}\index{representation!space}%
If $A$ is a \linebreak \st-algebra, and if $(H, \langle \cdot , \cdot \rangle)$
is a \underline{pre}-Hilbert space, then a \underline{representation} of $A$
in $H$ is an algebra homomorphism
\[ \pi : A \to \mathrm{End}(H) \]
such that
\[ \langle \pi (a)x,y \rangle = \langle x, \pi (a^*)y \rangle \]
for all $x,y \in H$ and for all $a \in A$. One then defines the \underline{range}
of $\pi$ as
\[ R(\pi) := \pi (A) = \{ \pi (a) \in \mathrm{End}(H) : a \in A \}. \]
One says that $(H, \langle \cdot , \cdot \rangle)$ is the
\underline{representation space} of $\pi$. (Recall that all the pre-Hilbert
spaces in this book shall be \underline{complex} pre-Hilbert spaces.)
\end{definition}

\begin{theorem}\label{HellToepl}
Let $\pi$ be a representation of a \st-algebra $A$ in a \underline{Hilbert}
space $H$. The operators $\pi(a)$ $(a \in A)$ then all are bounded. \linebreak
We may thus consider $\pi$ as a \st-algebra homomorphism $A \to B(H)$.
\end{theorem}

\begin{proof}
Decompose $\pi (a) = \pi (b) + \iu \pi (c)$, where $b,c \in A\sa$ and use the
Hellinger-Toeplitz Theorem, which says that a symmetric operator, defined
on all of a Hilbert space, is bounded, cf.\ e.g.\ Kreyszig \cite[10.1-1]{Krey}.
\end{proof}

\medskip
We shall now be interested in a criterion that guarantees boundedness
of the operators $\pi (a)$ for a representation $\pi$ acting in a pre-Hilbert
space merely: \ref{weaksa} - \ref{sigmaAsa}. \pagebreak

\begin{definition}[positive linear functionals]%
\index{functional!positive}\index{positive!functional}%
A linear functional $\varphi$ on a \st-algebra $A$ is called
\underline{positive} if $\varphi (a^*a) \geq 0$ for all $a \in A$.
\end{definition}

If $\pi$ is a representation of a \st-algebra $A$ in a pre-Hilbert space $H$,
then for all $x \in H$, the linear functional given by
\[ a \mapsto \langle \pi(a) x, x \rangle \qquad (a \in A) \]
is a positive linear functional. Indeed, for $x \in H$ and $a \in A$, we have
\[ \langle \pi (a^*a) x, x \rangle = \langle \pi (a) x, \pi (a) x \rangle \geq 0. \]

\begin{definition}[weak continuity]%
\index{representation!weakly continuous!on Asa@on $A\protect\sa$}%
\index{weakly continuous!As@on $A\protect\sa$}%
\index{representation!weakly continuous}\index{weakly continuous}%
Let $\pi$ be a representation of a \linebreak normed \st-algebra $A$ in a
pre-Hilbert space $H$. We shall say that $\pi$ is
\underline{weakly continuous on a subset} $B$ of $A$, if for all $x \in H$,
the positive linear functional
\[ a \mapsto \langle \pi (a) x, x \rangle \]
is continuous on $B$. We shall say that the representation $\pi$ is
\underline{weakly} \underline{continuous}, if $\pi$ is weakly continuous
on all of $A$.
\end{definition}

\begin{theorem}\label{weaksa}
Let $\pi$ be a representation of a normed \st-algebra $A$ in a
pre-Hilbert space $H$. Assume that $\pi$ is weakly continuous on $A\sa$.
The operators $\pi (a)$ $(a \in A)$ then all are bounded with
\[ \|\,\pi (a) \,\| \leq \rsigma(a)\quad \text{for all } a \in A. \]
\end{theorem}

\begin{proof}
Let $a \in A$ and $x \in H$ with $\|\,x\,\| \leq 1$ be fixed. For $n \geq 0$, put
\[ r_n := {\Bigl({ \,\| \,\pi (a) \,x \,\| \,}^2 \,\Bigr) \,}^{({ \,{2 \,}^n \,})}
\qquad \text{and} \qquad
s_n := \langle \,\pi \,\Bigl( { \,( \,a^*a \,) \,}^{( \,{2 \,}^n \,)} \,\Bigr) \,x , \,x \,\rangle. \]
We then have $r_n \leq s_n$ for all $n \geq 0$. Indeed, this obviously holds
for $n = 0$. Let now $n$ be any integer $\geq 0$. The Cauchy-Schwarz
inequality yields
\begin{align*}
{| \,s_n \,| \,}^2 \,= &
\ { \Bigl|\,\langle \,\pi \,\Bigl( { \,( \,a^*a \,) \,}^{( \,{2 \,}^n \,)} \,\Bigr) \,x, \,x \,\rangle \,\Bigr| \,}^2 \\
    \leq & \ \langle \,\pi \,\Bigl( { \,( \,a^*a \,) \,}^{( \,{2 \,}^n \,)} \,\Bigr) \,x,
    \,\pi \,\Bigl( { \,( \,a^*a \,) \,}^{( \,{2 \,}^n \,)} \,\Bigr) \,x \,\rangle {\,\| \,x \,\| \,}^2 \\
    \leq & \ \langle \,\pi \,\Bigl( { \,( \,a^*a \,) \,}^{( \,{2 \,}^{n+1} \,)} \,\Bigr) \,x, \,x \,\rangle = s_{n+1}.
\end{align*}
Thus, if $r_n \leq s_n$, then $s_n = |\,s_n\,|$, and hence
\[ r_{n+1} = {r_n \,}^2 \leq {s_n \,}^2 = {|\,s_n\,| \,}^2 \leq s_{n+1}. \pagebreak \]
By induction it follows that $r_n \leq s_n$ for all $n \geq 0$.
Let now $c > 0$ be such that
\[ |\, \langle \pi (b)x, x \rangle \,| \leq c \,|\,b\,| \]
for all $b \in A\sa$. For $n \geq 0$, we get
\begin{align*}
\|\,\pi (a) \,x \,\| \ & \leq
\ {\Biggl[ { \,\langle \,\pi \,\Bigl({ \,( \,a^*a \,) \,}^{( \,{2 \,}^n \,) \,}\Bigr) \,x,
\,x \,\rangle \,}^{1\,/ \,( \,{{2 \,}^{n}} \,)} \,\Biggr] \,}^{1 \,/ \,2} \\
 & \leq \ {c \,} ^{1 \,/  \,(\,{{2 \,}^{n+1}}\,)}
 \cdot{\Biggl[\,{ \,\Bigl| \,{( \,a^*a \,) \,}^{( \,{2 \,}^n \,)} \,\Bigr| \,}^{1 \,/ \,( \,{{2\,}^{n}} \,)} \,\Biggr] \,}^{1 \,/ \,2}
\to \,\rsigma(a).
\end{align*}
It follows that
\[ \|\,\pi (a) \,x \,\| \leq \rsigma(a). \qedhere \]
\end{proof}

\begin{definition}[$\sigma$-contractive linear maps]%
\label{sigmacontr}\index{s1@$\sigma$-contractive}%
\index{representation!s-contractive@$\sigma$-contractive}%
\index{contractive!s-contractive@$\sigma$-contractive}%
A linear map $\pi$ from a normed \st-algebra $A$ to a normed
space shall be called \underline{$\sigma$-contractive} if
\[ |\,\pi (a)\,| \leq \rsigma(a) \]
holds for all $a \in A$.
\end{definition}

\begin{theorem}\label{sigmaAsa}
Let $\pi$ be a $\sigma$-contractive linear map from a normed \linebreak
\st-algebra $A$ to a normed space. We then have:
\begin{itemize}
   \item[$(i)$] $\pi$ is contractive in the auxiliary \st-norm on $A$, see \ref{auxnorm}.
  \item[$(ii)$] $|\,\pi (a)\,| \leq \rlambda(a) \leq |\,a\,|$ for all Hermitian $a \in A$.
 \item[$(iii)$] $\pi$ is contractive on $A\sa$.
\end{itemize}
\end{theorem}

\begin{proof}
(i) follows from the fact that $\rsigma(a) \leq \|\,a\,\|$ holds for all $a \in A$,
where $\|\,\cdot\,\|$ denotes the auxiliary \st-norm. (ii) follows from the fact
that for a Hermitian element $a$ of $A$, one has $\rsigma(a) = \rlambda(a)$,
cf.\ \ref{Hermrseqrl}.
\end{proof}

\begin{theorem}\label{reprisomcontr}
If $A$ is a normed \st-algebra with isometric involution, then every
$\sigma$-contractive linear map from $A$ to a normed space is contractive.
\end{theorem}

\begin{theorem}\label{commsubBsigmacontr}
If $A$ is a commutative \st-subalgebra of a  Banach \linebreak \st-algebra,
then every $\sigma$-contractive linear map from $A$ to a normed space
is contractive.
\end{theorem}

\begin{proof}
This follows from the fact that in a Banach \st-algebra, one has
$\rsigma(b) \leq \rlambda(b)$ for all normal elements $b$,
cf.\ \ref{Bnormalrsleqrl}. \pagebreak
\end{proof}

\begin{proposition}\label{continvol}
Let $A$ be a normed \st-algebra. For the following statements, we have
\[ (i) \Rightarrow (ii) \Rightarrow (iii) \Rightarrow (iv) \Rightarrow (v) \Rightarrow (vi). \]
\begin{itemize}
   \item[$(i)$] the involution in $A$ is continuous,
  \item[$(ii)$] the involution maps sequences converging to $0$
                       to bounded \linebreak sequences,
 \item[$(iii)$] the map $a \mapsto a^*a$ is continuous at $0$,
 \item[$(iv)$] the function $\rsigma$ is continuous at $0$,
  \item[$(v)$] the function $\rsigma$ is bounded in a
neighbourhood of $0$,
 \item[$(vi)$] every $\sigma$-contractive linear map from $A$
to a normed space is continuous.
\end{itemize}
\end{proposition}

Now for slightly stronger assumptions: weakly
continuous representations in Hilbert spaces.

\begin{theorem}
Let $H, H\0$ be normed spaces, and let $C$ be a subset of the
normed space $B(H,H\0)$ of bounded linear operators $H \to H\0$.
If $H$ is complete, then $C$ is bounded in $B(H,H\0)$ if and only if
\[ \sup _{T \in C} |\,\ell(Tx)\,| < \infty \]
for every $x$ in $H$ and every continuous linear functional $\ell$ on $H\0$.
\end{theorem}

\begin{proof} The uniform boundedness principle is applied twice.
At the first time, one imbeds $H\0$ isometrically into its second dual.
\end{proof}

\begin{theorem}\label{weakcontHilbert}
Let $\pi$ be a \underline{weakly continuous} representation of a normed
\st-algebra $A$ in a \underline{Hilbert} space $H$. Then $\pi$ is continuous.
Furthermore, for a normal element $b$ of $A$, we have
\[ \| \,\pi(b) \,\| \leq \rlambda(b) \leq | \,b \,|. \]
It follows that if $A$ is commutative, then $\pi$ is contractive.
\end{theorem}

\begin{proof}
The operators $\pi (a)$ $(a \in A)$ are bounded by \ref{HellToepl} or \ref{weaksa}.
The linear functionals $A \ni a \mapsto \langle \pi (a) x, y \rangle$ $(x, y \in H)$
are continuous by polarisation. One then applies the preceding theorem
with $C := \{ \,\pi (a) \in B(H) : a \in A, \,|\,a\,| \leq 1 \,\}$ to the effect that $\pi$
is continuous. The remaining statements follow from \ref{commbdedcontr}.
\pagebreak
\end{proof}

Next the case of Banach \st-algebras.

\begin{lemma}
Let $\varphi$ be a positive linear functional on a \underline{unital}
\linebreak Banach \st-algebra $A$. For $a \in A\sa$, we then have
\[ |\,\varphi (a)\,| \leq \varphi (e) \,\rlambda (a) \leq \varphi (e) \,| \,a \,|. \]
In particular, $\varphi$ is continuous on $A\sa$. See also \ref{eBbded}
below.
\end{lemma}

\begin{proof}
Let $b = b^*$ be a Hermitian element of $A$ with $\rlambda(b) < 1$.
There exists by Ford's Square Root Lemma \ref{Ford} a Hermitian square
root $c$ of $e-b$. One thus has $c^*c = c^2 = e-b$, and it follows that
\[ 0 \leq \varphi (c^*c) = \varphi (e-b), \]
whence $\varphi (b)$ is real (because $\varphi (e) = \varphi (e^*e)$ is), and
\[ \varphi (b) \leq \varphi (e). \]
Upon replacing $b$ with $-b$, we obtain
\[ |\,\varphi (b) \,| \leq \varphi (e). \]
Let now $a = a^*$ be an arbitrary Hermitian element of $A$. Let
$\gamma > \rlambda(a)$. From the above it follows that
\[ |\,\varphi (\gamma^{-1}a)\,| \leq \varphi (e), \]
whence
\[ |\,\varphi (a)\,| \leq \varphi (e) \,\gamma, \]
which implies
\[ |\,\varphi (a)\,| \leq \varphi (e) \,\rlambda (a). \qedhere \]
\end{proof}

\begin{theorem}\label{Banachbounded}%
Let $\pi$ be a representation of a Banach \st-algebra $A$ in a pre-Hilbert
space $H$. The operators $\pi(a)$ $(a \in A)$ then all are bounded, and
the representation $\pi$ is continuous.
\end{theorem}

\begin{proof}
Extend $\pi$ to a representation $\tld{\pi}$ of $\tld{A}$ in $H$
by requiring that $\tld{\pi} (e) = \mathds{1}$ if $A$ has no unit (and by
extending linearly). By the preceding lemma, the positive linear functional
\[ a \to \langle \tld{\pi} (a) x, x \rangle \]
is continuous on $\tld{A}\sa$ for each $x \in H$. This implies that $\pi$ is
weakly continuous on $A\sa$. It follows by \ref{weaksa} that the operators
$\pi (a)$ are bounded on $H$. We obtain that $R(\pi)$ is a pre-C*-algebra,
so that $\pi$ is continuous by theorem \ref{automcont}. \pagebreak
\end{proof}

\begin{corollary}
A representation of a commutative Banach \linebreak
\st-algebra in a pre-Hilbert space is contractive.
\end{corollary}

\begin{proof}
This follows now from \ref{weakcontHilbert}.
\end{proof}

\begin{corollary}
A representation of a Banach \st-algebra with \linebreak isometric
involution in a pre-Hilbert space is contractive.
\end{corollary}

\begin{proof}
This follows now from \ref{weaksa} and \ref{reprisomcontr}.
\end{proof}

\medskip
In particular, we have:

\begin{corollary}
A representation of a C*-algebra in a pre-Hilbert space is contractive.
\end{corollary}

\begin{definition}[faithful representation]%
\index{representation!faithful}\index{faithful representation}%
A representation of a \st-al\-gebra in a pre-Hilbert space is
called \underline{faithful} if it is injective.
\end{definition}

\begin{theorem}\label{faithisometric}
A faithful representation of a C*-algebra in a \linebreak Hilbert space is isometric.
\end{theorem}

\begin{proof} \ref{C*injisometric}. \end{proof}

\begin{theorem}\label{rangeC*}
If $\pi$ is a representation of a C*-algebra $A$ in a Hilbert space, then
$A \,/ \ker \pi$ and $R(\pi)$ are C*-algebras as well. The representation
$\pi$ factors to a C*-algebra isomorphism from $A \,/ \ker \pi$ onto $R(\pi)$.
\end{theorem}

\begin{proof} \ref{homfact}. \end{proof}

\medskip
For further continuity results, see \ref{bdedstatesdef} - \ref{kerpsi}.

\begin{remark}[concerning \protect\ref{continvol}]\label{Palmer}%
Theodore W.\ Palmer shows in the second volume of his work
\cite[11.1.4]{Palm} that in a Banach \st-algebra, the function $\rsigma$
is continuous at $0$. This however requires methods above the scope
of this book, and we prove the result only for Hermitian Banach \st-algebras,
cf.\ \ref{weakPalmer}.
\end{remark}

\begin{remark}%
For most purposes of representation theory, the condition of isometry of the
involution can be weakened to requiring that $| \,a^*a \,| \leq {| \,a \,|\,}^2$
holds for all elements. \pagebreak
\end{remark}

\clearpage


\section{Cyclic and Non-Degenerate Representations}%
\label{CyclicNonDeg}

In this paragraph, let $\pi$ be a representation of a \st-algebra $A$
in a Hilbert space $(H, \langle \cdot , \cdot \rangle)$.

\begin{definition}[invariant subspaces]%
\index{invariant}\index{subspace!invariant}%
A subspace $M$ of $H$ is called \underline{invariant} under
$\pi$, if $\pi (a)x \in M$ for every $x \in M$, and every $a \in A$.
\end{definition}

\begin{proposition}%
If $M$ is a subspace of $H$ invariant under $\pi$, so is $M^{\perp}$.
\end{proposition}

\begin{proof}
Let $M$ be a subspace of $H$ invariant under $\pi$. Let $y \in M^{\perp}$,
i.e.\ $\langle y, x \rangle = 0$ for all $x \in M$. For $a \in A$ we then have
\[ \langle \pi (a)y, x \rangle = \langle y,\pi (a^*)x \rangle = 0 \]
for all $x \in M$, so that also $\pi (a)y \in M^{\perp}$.
\end{proof}

\medskip
This is why invariant subspaces are also called \underline{reducing}
subspaces.\index{reducing}\index{subspace!reducing}

\begin{definition}[subrepresentations]\index{subrepresentation}%
\index{representation!subrepresentation}\index{p0@$\pi _M$}%
If $M$ is a closed subspace of $H$ invariant under $\pi$, let $\pi _M$
denote the representation of $A$ in $M$ defined by
\begin{alignat*}{2}
\pi _M : \ & A & \ \to & \  B(M) \\
           & a & \ \mapsto & \ \pi (a)|_M.
\end{alignat*}
The representation $\pi _M$ is called a
\underline{subrepresentation} of $\pi$.
\end{definition}

\begin{definition}[commutants]%
\index{commutant}\index{p1@$\pi '$}\label{repcommutant}%
One denotes by
\[ \pi' := \{ \,c \in B(H) : c \,\pi (a) = \pi (a)c \text{ for all } a \in A \,\} \]
the \underline{commutant} of $\pi$. Similarly, let $S \subset B(H)$.
The set of all operators in $B(H)$ which commute with every operator in
$S$ is called the commutant of $S$ and is denoted by $S'$.
(Cf.\ \ref{questimbed}.)
\end{definition}

\begin{proposition}\label{invarcommutant}%
\index{commutant}\index{subspace!invariant!commutant}%
Let $M$ be a closed subspace of $H$ and let $p$ be the projection on $M$.
Then $M$ is invariant under $\pi$ if and only if $p \in \pi'$. \pagebreak
\end{proposition}

\begin{proof} Assume that $p \in \pi'$. For $x \in M$ and
$a \in A$, we then have $\pi (a)x = \pi (a)px = p \pi (a)x \in M$, so
that $M$ is invariant under $\pi$. Conversely, assume that $M$
is invariant under $\pi$. For $x \in H$ and $a \in A$, we get
$\pi (a)px \in M$, whence $p \pi (a)px = \pi (a)px$, so that
$p \pi (a)p = \pi (a)p$ for all $a \in A$. By applying the involution in
$B(H)$, and replacing $a$ by $a^*$, we obtain that also
$p \pi (a)p = p \pi (a)$, whence $p \pi (a) = \pi (a)p$, i.e.\ $p \in \pi'$.
\end{proof}

\begin{definition}[non-degenerate representations]%
\label{nondegenerate}\index{non-degenerate!representation}%
\index{representation!non-degenerate}%
The representation $\pi$ is called \underline{non-degenerate}
if the closed invariant subspace
\[ N:=\{ \,x \in H : \pi (a)x = 0 \text{ for all } a \in A \,\} \]
is $\{ 0 \}$. Consider the closed invariant subspace
\[ R := \overline{\mathrm{span}} \,\{ \,\pi (a) x \in H: a \in A,\ x \in H \,\}. \]
We have $R^{\perp} = N$, so that $\pi$ is non-degenerate
if and only if $R$ is all of $H$. Also $H = R \oplus N$, or
$\pi = \pi _R \oplus \pi _N$. This says that $\pi$ is the direct sum
of a non-degenerate representation and a null representation. In
particular, by passing to the subrepresentation $\pi _R$, we can
always assume that $\pi$ is non-degenerate.
\end{definition}

\begin{proof}
It suffices to prove that $R^{\perp} = N$. Let $a \in A$ and $x \in H$.
For $y \in N$, we find
\[ \langle \pi (a) x, y \rangle = \langle x, \pi (a^*) y \rangle = 0, \]
so $N \subset R^{\perp}$. On the other hand, for $y \in R^{\perp}$,
we find
\[ \langle \pi (a) y, x \rangle = \langle y, \pi (a^*) x \rangle = 0, \]
which shows that $\pi (a) y \perp H$, or $y \in N$, so
$R^{\perp} \subset N$.
\end{proof}

\begin{definition}[cyclic representations]%
\index{representation!cyclic}\index{cyclic!vector}%
\index{cyclic!representation}\index{vector!cyclic}%
One says that $\pi$ is \underline{cyclic} if $H$
contains a non-zero vector $x$ such that the set
\[ \{ \,\pi (a)x \in H : a \in A \,\} \]
is dense in $H$. The vector $x$ then is called a
\underline{cyclic vector}. In particular, a cyclic representation
is non-degenerate.
\end{definition}

Of great importance is the following observation. \pagebreak

\begin{lemma}[cyclic subspaces]\label{cyclicsubspace}%
\index{cyclic!subspace}\index{subspace!cyclic}%
Assume that $\pi$ is non-degenerate, and that $c$ is a
non-zero vector in $H$. Consider the closed invariant subspace
\[ M := \overline{\pi (A)c}. \]
Then $c$ belongs to $M$. In particular, $\pi _M$ is cyclic with cyclic
vector $c$. One says that $M$ is a \underline{cyclic subspace} of $H$.
\end{lemma}

\begin{proof}
Let $p$ denote the projection upon $M$. We then have $p \in \pi'$.
Thus, for all $a \in A$, we get
\[ \pi(a)(\mathds{1}-p)c = (\mathds{1}-p) \pi (a)c = 0, \]
which implies that $(\mathds{1}-p)c = 0$ as $\pi$ is non-degenerate.
\end{proof}

\begin{definition}[direct sum of Hilbert spaces]%
\index{direct sum!of Hilbert spaces}%
Let $\{\,H_i\,\}_{\,i \in I}$ be a family of Hilbert spaces. The direct sum
$\oplus _{i \in I} \,H_i$ consists of all families $\{\,x_i\,\}_{\,i \in I}$ such
that $x_i \in H_i$ $(i \in I)$ and $\{ \,\| \,x_i \,\| \,\}_{\,i \in I} \in {\ell\,}^2(I)$.
The element $\{\,x_i\,\} _{\,i \in I}$ then is denoted by $\oplus _{i \in I} \,x_i$.
The direct sum $\oplus _{i \in I} \,H_i$ is a Hilbert space with respect to
the inner product
\[ \langle \oplus _{i \in I} \,x_i,\oplus _{i \in I} \,y_i \rangle
:= \sum _{i \in I} \langle x_i, y_i \rangle
\qquad (\oplus _{i \in I} \,x_i, \oplus _{i \in I} \,y_i \in \oplus _{i \in I} \,H_i). \]
\end{definition}

\begin{proof} Let $\oplus _{i \in I} \,x_i$, $\oplus _{i \in I} \,y_i
\in \oplus _{i \in I} \,H_i$. We then have
\begin{align*}
{\biggl( \,\sum _{i \in I} {\| \,x_i+y_i \,\| \,}^2 \biggr)}^{1/2}
\leq & \ {\biggl( \,\sum _{i \in I} {\Bigl( \,\| \,x_i \,\| + \|\, y_i \,\| \,\Bigr) \,}^2 \,\biggr)}^{1/2} \\
\leq & \ {\biggl( \,\sum _{i \in I}{\| \,x_i \,\| \,}^2 \biggr)}^{1/2}
+ {\biggl( \,\sum _{i \in I} {\| \,y_i \,\| \,}^2 \biggr)}^{1/2}
\end{align*}
so that $\oplus _{i \in I} \,H_i$ is a vector space. We also have
\begin{align*}
\sum _{i \in I} | \,\langle x_i ,y_i \rangle \,|
\ & \leq \ \sum _{i \in I}\|\, x_i \,\| \cdot \| \,y_i \,\| \\
  & \leq \ {\biggl( \,\sum _{i \in I} {\| \,x_i \,\| \,}^2 \biggr)}^{1/2}
     {\biggl( \,\sum _{i \in I} {\| \,y_i \,\| \,}^2 \biggr)}^{1/2}, 
\end{align*}
so that $\oplus _{i \in I} \,H_i$ is a pre-Hilbert space.
Let $( \,\oplus _{i \in I} \,x_{i,n} \,)_{\,n}$ be a Cauchy sequence in
$\oplus _{i \in I} \,H_i$. Let $\varepsilon > 0$ be given.
There then is an $n\0 \geq 0$ such that
\[ {\biggl( \,\sum _{i \in I} {\| \,x_{i,m} - x_{i,n} \,\| \,}^2 \biggr)}^{1/2}
\leq \varepsilon\quad \text{for all } m,n \geq n\0. \pagebreak \]
From this we conclude that for each $i \in I$, the sequence
$( \,x_{i,n} \,)_{\,n}$ is a Cauchy sequence in
$H_i$, and therefore converges to an element $x_i$ of $H_i$.
With $F$ a finite subset of $I$, we have
\[ {\biggl( \,\sum _{i \in F} {\| \,x_{i,m} - x_{i,n} \,\| \,}^2 \biggr)}^{1/2}
\leq \varepsilon\quad \text{for all }m,n \geq n\0, \]
so that by letting $m \to \infty$, we get
\[ {\biggl( \,\sum _{i \in F} {\| \,x_i - x_{i,n} \,\| \,}^2 \biggr)}^{1/2}
\leq \varepsilon\quad \text{for all } n \geq n\0, \]
and since $F$ is an arbitrary finite subset of $I$, it follows that
\[ {\biggl( \,\sum _{i \in I} {\| \,x_i - x_{i,n} \,\| \,}^2 \biggr)}^{1/2}
\leq \varepsilon\quad \text{for all } n \geq n\0. \]
This not only shows that $\oplus _{i \in I} \,(x_i -x_{i,n})$, and hence
$\oplus _{i \in I} \,x_i$, belong to $\oplus _{i \in I} \,H_i$, but also that
$( \,\oplus _{i \in I} \,x_{i,n} \,)_{\,n}$ converges in
$\oplus _{i \in I} \,H_i$ to $\oplus _{i \in I} \,x_i$, so that
$\oplus _{i \in I} \,H_i$ is complete. 
\end{proof}

\begin{definition}[direct sum of operators]%
\index{direct sum!of operators}%
Consider a  family of bounded linear operators $\{\,b_i\,\}_{\,i \in I}$
on Hilbert spaces $H_i$ $(i \in I)$ respectively. Assume that
\[ \sup _{i \in I} \,\|\,b_i\,\| < \infty. \]
One then defines
\[ ( \oplus _{i \in I} \,b_i ) ( \oplus _{i \in I} \,x_i ) := \oplus _{i \in I} \,(b_i x_i)
\qquad (\oplus _{i \in I} \,x_i \in \oplus _{i \in I} \,H_i). \]
Then $\oplus _{i \in I} \,b_i$ is a bounded linear operator on
$\oplus _{i \in I} H_i$ with norm
\[ \| \oplus _{i \in I} b_i \,\| = \sup _{i \in I} \,\| \,b_i \,\|.  \]
\end{definition}

\begin{proof} With $c := \sup _{i \in I} \,\| \,b_i \,\|$, we have
\begin{align*}
{\| \,( \oplus _{i \in I} \,b_i ) ( \oplus _{i \in I} \,x_i ) \,\| \,}^2
\,& = \  \sum _{i \in I} \,{\| \,b_i x_i \,\| \,}^2 \\
  & \leq \ \sum _{i \in I} \,{\| \,b_i \,\| \,}^2 \cdot {\| \,x_i \,\| \,}^2 \\
  & \leq \ {c \,}^2 \cdot \,\sum _{i \in I} \,{\| \,x_i \,\| \,}^2
= {c \,}^2 \cdot {\,\| \oplus _{i \in I} x_i \,\| \,}^2,
\end{align*}
so $\| \oplus _{i \in I} b_i \,\| \leq c$. The converse inequality is obvious.
\pagebreak
\end{proof}

\begin{definition}[direct sum of representations]%
\label{directsumsigma}\index{direct sum!of representations}%
\index{representation!direct sum}%
Consider a \linebreak family $\{ \,\pi_i \,\} _{\,i \in I}$ of representations
of $A$ in Hilbert spaces $H_i$ $(i \in I)$ respectively. If
$\{ \,\| \,\pi _i (a) \,\| \,\} _{\,i \in I}$ is bounded for each $a \in A$, we may
\linebreak define the direct sum representation $\oplus _{i \in I} \,\pi _i$
given by
\begin{alignat*}{2}
\oplus _{i \in I} \,\pi _i : \ & A \ &\to \ & B(\oplus _{i \in I} \,H_i) \\
 \ & a \ &\mapsto \ & \oplus _{i \in I} \pi _i (a).
\end{alignat*}
Hence the direct sum of a family of $\sigma$-contractive
representations of a normed \st-algebra is well-defined and
$\sigma$-contractive.
\end{definition}

\begin{remark}
If the family $\{ \,H_i \,\}_{\,i \in I}$ is a family of pairwise ortho\-gonal
subspaces of a Hilbert space $H$, we imbed the direct sum
$\oplus _{i \in I} H_i$ in $H$.
\end{remark}

It is easily seen that a direct sum of non-degenerate representations
is non-degenerate. In the converse direction we have:

\begin{theorem}\index{direct sum!decomposition}%
\index{representation!direct sum!decomposition}%
\index{decomposition!direct sum}\label{directsumdeco}%
If $\pi$ is non-degenerate and if $H \neq \{ 0 \}$, then
$\pi$ is the direct sum of cyclic subrepresentations.
\end{theorem}

\begin{proof}
Let $Z$ denote the set of all non-empty sets consisting of mutually
orthogonal, cyclic \ref{cyclicsubspace} closed invariant subspaces
of $H$. The set $Z$ is not empty since for $0 \neq c \in H$, we may
consider $M_c := \overline{\pi (A)c}$, so that $\{\, M_c \,\} \in Z$ by
\ref{cyclicsubspace}. Furthermore, $Z$ is inductively ordered by
inclusion. Thus, Zorn's Lemma guarantees the existence of a maximal
set $F$. Consider now $H_1 := \oplus _{M \in F} M$. We claim that
$H_1 = H$. Indeed, if $0 \neq c \in {H_1}^{\perp}$, then
$M_c := \overline{\pi (A)c}$ is a cyclic closed invariant subspace of
$H$ by \ref{cyclicsubspace}, and $F \cup \{ \,M_c \,\}$ would be a set
in $Z$ strictly larger than the maximal set $F$. It follows that
$\pi = \oplus _{M \in F} \,\pi _M$. 
\end{proof}

\medskip
Our next aim is lemma \ref{coeffequal}, which needs some preparation.

\begin{proposition}\label{closureunitary}%
Let $I$ be an isometric linear operator defined on a dense subspace
of a Hilbert space $H_1$, with dense range in a Hilbert space $H_2$.
Then the closure $U$ of $I$ is a unitary operator from $H_1$ onto $H_2$.
\end{proposition}

\begin{proof}
The range of $U$ is a complete, hence closed
subspace of $H_2$, containing the dense range of $I$.
\pagebreak
\end{proof}

\begin{definition}[spatial equivalence]\label{spatequiv}%
\index{spatially equivalent!representations}%
\index{equivalent!spatially!representations}%
Let $\pi_1$, $\pi_2$ be two representations of a \st-algebra $A$
in Hilbert spaces $H_1$, $H_2$  respectively. The representations
$\pi_1$ and $\pi_2$ are called \underline{spatially equivalent}
if there exists a unitary operator $U$ from $H_1$ to $H_2$ such that
\[ U \pi_1 (a) = \pi_2 (a) U \qquad (a \in A). \]
\end{definition}

\begin{definition}[intertwining operators, $C(\pi_1,\pi_2)$]%
\label{intertwinop}\index{intertwining operator}%
\index{operator!intertwining}\index{C45@$C(\pi_1,\pi_2)$}%
Let $\pi_1$, $\pi_2$ be two representations of a \st-algebra $A$
in Hilbert spaces $H_1$, $H_2$ \linebreak respectively. One
says that a bounded linear operator $b$ from $H_1$ to $H_2$
\underline{intertwines} $\pi_1$ and $\pi_2$, if
\[ b \pi_1 (a) = \pi_2 (a) b \qquad  (a \in A). \]
The set of all bounded linear operators intertwining $\pi_1$ and
$\pi_2$ is \linebreak denoted by $C(\pi_1,\pi_2)$. Please note that
$\pi ' = C(\pi, \pi)$, cf.\ \ref{repcommutant}.
\end{definition}

The following result is almost miraculous.

\begin{lemma}\label{coeffequal}%
Let $\pi_1$, $\pi_2$ be cyclic representations of a \st-algebra $A$
in Hilbert spaces $H_1$, $H_2$ with cyclic vectors $c_1$, $c_2$. For
$a \in A$, let
\begin{align*}
\varphi_1 (a) & := \langle \pi_1 (a) c_1, c_1 \rangle, \\
\varphi_2 (a) & := \langle \pi_2 (a) c_2, c_2 \rangle .
\end{align*}
If $\varphi_1 = \varphi_2$, then $\pi_1$ and $\pi_2$ are spatially
equivalent. Furthermore, there then exists a unitary operator in
$C(\pi_1, \pi_2)$ taking $c_1$ to $c_2$. 
\end{lemma}

\begin{proof} One defines an operator $I$ by putting
\[ Iz := \pi_2 (a) c_2 \text{ when } z = \pi_1 (a) c_1. \]
It shall be shown that $I$ is well-defined and
isometric. Indeed, we have
\begin{align*}
\langle Iz, Iz \rangle & = \langle \pi_2 (a) c_2, \pi_2 (a) c_2 \rangle \\
 & = \langle \pi_2 (a^*a) c_2, c_2 \rangle \\
 & = \varphi _2 (a^*a) = \varphi _1 (a^*a) \\
 & = \langle \pi_1 (a^*a) c_1, c_1 \rangle \\
 & = \langle \pi_1 (a) c_1, \pi_1 (a) c_1 \rangle = \langle z, z \rangle.
\end{align*}
It follows that $I$ is well-defined because if for $a, b \in A$ one has
$\pi_1 (a) c_1 = \pi_1 (b) c_1$, then $\pi_1 (a-b) c_1 = 0$, whence
also $\pi_2 (a-b) c_2 = 0$, i.e.\ $\pi_2 (a) c_2 = \pi_2 (b) c_2$.
Furthermore $I$ is densely defined in $H_1$ with range dense in $H_2$.
It follows that the closure $U$ of $I$ is a unitary operator from
$H_1$ onto $H_2$, cf.\ \ref{closureunitary}. It shall be shown that this
operator $U$ intertwines $\pi_1$ and $\pi_2$. Indeed, if
$z = \pi_1 (a) c_1$, we have for all $b \in A$:
\begin{align*}
\pi_2 (b) Iz & = \pi_2 (b)\pi_2 (a) c_2 \\
 & = \pi_2 (ba) c_2 \\
 & = I \pi_1 (ba) c_1 \\
 & = I \pi_1 (b) \pi_1 (a) c_1 = I \pi_1 (b) z.
\end{align*}
By continuity of $U$ and $\pi (b)$ \ref{HellToepl}, it follows that
\[ \pi_2 (b) U = U \pi_1 (b) \qquad (b \in A), \]
so that $U$ intertwines $\pi_1$ and $\pi_2$. It remains to be
shown that $U c_1 = c_2$. For $a \in A$, we have
\[ U \pi_1 (a) c_1 = \pi_2 (a) c_2, \]
whence
\[ \pi_2 (a) U c_1 = \pi_2 (a) c_2. \]
Since $\pi_2$ is cyclic, and thereby non-degenerate, it follows that
$U c_1 = c_2$.
\end{proof}

\clearpage


\chapter{States}

\setcounter{section}{23}


\section{The GNS Construction}

\begin{introduction}\index{GNS construction|(}\index{coefficient}%
\index{representation!coefficient}%
If $\pi$ is a representation of a \st-algebra $A$ in a
Hilbert space $H$, then for $x \in H$, the positive linear functional
\[ a \mapsto \langle \pi (a) x, x \rangle \qquad (a \in A) \]
is called a \underline{coefficient} of $\pi$. (In analogy with matrix
coefficients.)

In this paragraph as well as in the next one, we shall
be occupied with inverting this relationship.

That is, given a normed \st-algebra $A$ and a positive linear functional
$\varphi$ on $A$ satisfying certain conditions, we shall construct a
representation $\pi _{\varphi}$ of $A$ in such a way that $\varphi$ is a
coefficient of $\pi _{\varphi}$. The present paragraph serves as a rough
foundation for the ensuing one which puts the results into sweeter words.
\end{introduction}

\begin{proposition}\label{inducedHf}%
If $\varphi$ is a positive linear functional on a \st-algebra $A$, then
\[ \langle a, b \rangle_\varphi := \varphi (b^*a) \qquad (\,a,b \in A \,) \]
defines a positive Hilbert form $\langle \cdot , \cdot \rangle _\varphi$
on $A$ in the sense of the following definition.
\end{proposition}

\begin{definition}%
[positive Hilbert forms \hbox{\protect\cite[p.\ 349]{Dieu}}]%
\label{Hilbertform}\index{Hilbert form}\index{positive!Hilbert form}%
A \underline{positive} \underline{Hilbert form} on a \st-algebra $A$ is a positive
semidefinite sesquilinear form $\langle \cdot , \cdot \rangle$ on $A$ such that
\[ \langle ab, c \rangle = \langle b, a^*c \rangle \quad \text{for all } a,b,c \in A. \]
Please note that then also
\begin{align*}
\langle a, b \rangle & = \overline{\langle b, a \rangle}, 
\tag*{\textit{(the sesquilinear form is Hermitian)}} \\
{| \,\langle a, b \rangle \,| \,}^2 & \leq \langle a, a \rangle \cdot \langle b, b \rangle.
\pagebreak \tag*{\textit{(Cauchy-Schwarz inequality)}}%
\end{align*}%
\end{definition}

\begin{definition}\label{quotientA}%
\index{isotropic}\index{subspace!isotropic}%
\index{A9@$\underline{A}$}\index{a5@$\underline{a}$}%
Let $\langle \cdot , \cdot \rangle$ be a positive Hilbert form on a
\st-algebra $A$. Denote by $I$ the isotropic subspace of the inner
product space $(A, \langle \cdot, \cdot \rangle)$, that is
$I = \{ \,b \in A : \langle b, c \rangle = 0 \text{ for all } c \in A \,\}$.
The isotropic subspace $I$ is a left ideal in $A$ because
$\langle ab, c \rangle = \langle b, a^* c \rangle$.
Denote by $\underline{A}$ the quotient space $A/I$.
For $a \in A$, let $\underline{a} := a+I \in \underline{A}$.
In $\underline{A}$ one defines an inner product by
\[ \langle \underline{a} , \underline{b} \rangle :=
\langle a, b \rangle \qquad (a, b \in A). \]
In this way $\underline{A}$ becomes a pre-Hilbert space (as
$I = \{ \,b \in A : \langle b, b \rangle = 0 \,\}$ by the Cauchy-Schwarz
inequality).
\end{definition}

\begin{definition}\label{GNS}%
\index{representation!GNS}%
Let $\langle \cdot, \cdot \rangle$ be a positive Hilbert form on a
\st-algebra $A$. For $a \in A$, the translation operator
\begin{align*}
\pi (a) : \underline{A} & \to \underline{A} \\
\underline{b} & \mapsto \pi (a) \underline{b} := \underline{ab}
\end{align*}
is well defined because $I$ is a left ideal. One thus obtains a
representation $\pi$ of $A$ in the pre-Hilbert space $\underline{A}$.
It is called the \underline{representation} \underline{associated with
$\langle \cdot , \cdot \rangle$}. If all of the operators $\pi (a)$ $(a \in A)$ are
bounded, then by continuation we get a representation of $A$
in the completion of $\underline{A}$. This latter representation is called the
\underline{GNS representation}
\underline{associated with $\langle \cdot , \cdot \rangle$}.
\end{definition}

\begin{definition}\label{GNSvarphi}%
\index{H1@$H_{\varphi}$}\index{GNS construction!Hphi@$H _{\varphi}$}%
\index{GNS construction!piphi@$\pi_ {\varphi}$}%
\index{p2@$\pi_ {\varphi}$}\index{representation!GNS}%
Assume that the Hilbert form $\langle \cdot , \cdot \rangle$ derives from
a positive linear functional $\varphi$ on $A$, as in \ref{inducedHf}. One
denotes by $H_{\varphi}$ the completion of $\underline{A}$.
One says that the representation associated with
$\langle \cdot , \cdot \rangle _{\varphi}$ is the representation associated
with $\varphi$. If the operators $\pi (a)$ $(a \in A)$ are bounded, one
denotes by $\underline{\pi _{\varphi}}$ the GNS representation associated
with $\langle \cdot , \cdot \rangle _{\varphi}$. It is then called the
\underline{GNS representation associated with $\varphi$}.
\end{definition}

The acronym GNS stands for Gelfand, Na\u{\i}mark and Segal, the
originators of the theory.

\begin{definition}[weak continuity on $A\sa$]%
\label{weakbdedsadef}%
\index{functional!w@weakly continuous on $A\protect\sa$}%
\index{weakly continuous!As@on $A\protect\sa$}%
A positive linear functional $\varphi$ on a normed \st-algebra
$A$ shall be called \underline{weakly continuous}
\underline{on $A\sa$} if for each $b \in A$, the positive linear
functional $a \mapsto \varphi (b^*ab)$ is continuous on $A\sa$.
\pagebreak
\end{definition}

\begin{theorem}\label{weakbdedsarep}%
Let $\varphi$ be a positive linear functional on a normed \linebreak
\st-algebra $A$. Then $\varphi$ is weakly continuous on $A\sa$ if
and only if the \linebreak representation associated with $\varphi$
is weakly continuous on $A\sa$. In this case, the GNS representation
$\pi _{\varphi}$ exists and is a $\sigma$-contractive representation
of $A$ in $H_{\varphi}$. Cf.\ \ref{weaksa}, \ref{sigmacontr}.
\end{theorem}

\begin{proof}
For $a,b \in A$ we have $\varphi (b^*ab) =
\langle \pi (a) \underline{b} , \underline{b} \rangle _{\varphi}$. 
\end{proof}

\begin{proposition}\label{Banachweakbdedsa}%
Let $A$ be a Banach \st-algebra. Then every positive linear
functional on $A$ is weakly continuous on $A\sa$.
\end{proposition}

\begin{proof}
This follows now from theorem \ref{Banachbounded}.
\end{proof}

\begin{definition}[variation \hbox{\cite[p.\ 307]{Gaal}}]%
\index{variation}\index{v(phi)@$v(\varphi)$}\label{variation}%
\index{GNS construction!v(phi)@$v(\varphi)$}%
A positive linear functional $\varphi$ on a \st-algebra $A$
is said to have \underline{finite variation} if
\[ {| \,\varphi (a) \,| \,}^2 \leq \gamma \,\varphi (a^*a) \qquad (a \in A) \]
for some $\gamma \geq 0$. If this is the case, we shall say that
\[ v ( \varphi ) := \,\inf \,\{ \,\gamma \geq 0 :
{| \,\varphi (a) \,| \,}^2 \leq \gamma \,\varphi (a^*a) \text{ for all } a \in A \,\} \]
is the \underline{variation} of $\varphi$. We then also have
\[ {| \,\varphi (a) \,| \,}^2 \leq v ( \varphi ) \,\varphi (a^*a) \qquad (a \in A). \]
\end{definition}

\begin{proof}
Assume that $\varphi$ has finite variation. With 
\[ S := \{ \,\gamma \geq 0 : {| \,\varphi (a) \,|}^{\,2} \leq \gamma \,\varphi (a^*a)
\ \text{for all} \ a \in A \,\}\]
we have that
\[ v (\varphi) = \inf S = \text{the greatest lower bound of} \ S. \]
Let now $b \in A$ be fixed. We have to prove that
\[ {| \,\varphi (b) \,|}^{\,2} \leq v (\varphi) \,\varphi (b^*b). \]
We know that
\[ {| \,\varphi (b) \,|}^{\,2} \leq \gamma \,\varphi (b^*b) \ \text{for all} \ \gamma \in S. \]
That is:
\[ {| \,\varphi (b) \,|}^{\,2}\ \text{is a lower bound of} \ S \,\varphi (b^*b). \]
Since
\[ v (\varphi) \,\varphi (b^*b) \ \text{is the greatest lower bound of} \ S \,\varphi (b^*b), \]
it follows that
\[ {| \,\varphi (b) \,|}^{\,2} \leq v (\varphi) \,\varphi (b^*b), \]
as was to be shown.
\end{proof}

\begin{proposition}\label{variationinequal}%
Let $\pi$ be a representation of a \st-algebra $A$ in a \linebreak
pre-Hilbert space $H$. For $x \in H$, consider the positive linear
functional $\varphi$ on $A$ defined by
\[ \varphi(a):=\langle \pi(a)x, x \rangle \qquad (a \in A). \]
Then $\varphi$ is Hermitian and has finite variation
$v(\varphi) \leq \langle x, x \rangle$.
\end{proposition}

\begin{proof}
For $a \in A$, we have
\[ {| \,\varphi (a) \,| \,}^2 = {|\,\langle \pi (a) x, x \rangle \,| \,}^2
\leq {\|\,\pi(a)x\,\| \,}^2 \,{ \|\,x\,\| \,}^2 = \varphi(a^*a) \,\langle x, x \rangle, \]
so that $v ( \varphi ) \leq \langle x, x \rangle$. Finally $\varphi$
is Hermitian \ref{Hermfunct}, because for $a \in A$, we have
\[ \varphi(a^*)=\langle\pi(a^*)x,x \rangle=\langle x,\pi(a)x \rangle
=\overline{\langle \pi (a) x,x \rangle} = \overline{\varphi(a)}. \qedhere \]
\end{proof}

\begin{proposition}\label{variationequal}%
Let $\pi$ be a \underline{non-degenerate} representation of a
\linebreak \st-algebra $A$ in a \underline{Hilbert} space $H$.
For $x \in H$, consider the positive \linebreak linear functional
$\varphi$ on $A$ defined by
\[ \varphi (a):= \langle \pi (a) x, x \rangle \qquad (a \in A). \]
The variation of $\varphi$ then is
\[ v(\varphi) = \langle x, x \rangle. \]
\end{proposition}

\begin{proof}
The preceding result shows that $v ( \varphi ) \leq \langle x, x \rangle $.
In order to prove the opposite inequality, consider the closed subspace
$M := \overline{\pi (A) x}$ of $H$. Then $x \in M$ because $\pi$ is
non-degenerate, cf.\ \ref{cyclicsubspace}. There thus exists a sequence
$(a_n)$ in $A$ such that
\[ x = \lim _{n \to \infty} \pi (a_n) x. \]
We then have
\begin{align*}
 & \,\left| \ {| \,\varphi (a_n) \,| \,}^2
 - {\| \,x \,\| \,}^2 \cdot \varphi ({a_n}^*{a_n}) \ \right| \\
 = & \,\left| \ {| \,\langle \pi (a_n) x, x \rangle \,| \,}^2
- {\| \,x \,\| \,}^2 \cdot {\| \,\pi (a_n) x \,\| \,}^2 \ \right| \\
 \leq & \,\left| \ {| \,\langle \pi (a_n) x, x \rangle \,| \,}^2
 - {\| x \| \,}^4 \ \right|
 + {\| \,x \,\| \,}^2 \cdot \left| \ {\| \,\pi (a_n) x \,\| \,}^2
 - {\| \,x \,\| \,}^2 \ \right|
\end{align*}
so that
\[ \lim _{n \to \infty} \left| \ {| \,\varphi (a_n) \,| \,}^2 -
\,{\| \,x \,\| \,}^2 \cdot \varphi ({a_n}^*{a_n}) \ \right| = 0, \]
by
\[ \lim _{n \to \infty} \langle \pi (a_n) x, x \rangle = {\| \,x \,\| \,}^2 \]
and
\[ \lim _{n \to \infty} \| \,\pi (a_n) x \,\| = \| \,x \,\|. \pagebreak \qedhere \]
\end{proof}

\begin{definition}[reproducing vector, $c_{\varphi}$]%
\label{reprovectordef}%
\index{c@$c _{\varphi}$}\index{reproducing vector}%
\index{GNS construction!cphi@$c _{\varphi}$}%
Let $\varphi$ be a positive linear functional on a \st-algebra $A$.
A \underline{reproducing vector for $\varphi$} is a vector $c_{\varphi}$
in $H_{\varphi}$ such that
\[ \varphi (a) = \langle \underline{a}, c_{\varphi} \rangle _{\varphi}
\quad \text{for all } a \in A. \]
It is then uniquely determined by $\varphi$. We shall reserve the
notation $c_{\varphi}$ for the reproducing vector of $\varphi$.
\end{definition}

\begin{proposition}%
Let $\varphi$ be a positive linear functional on a unital \linebreak
\st-algebra $A$. Then $\underline{e}$ is the reproducing vector for
$\varphi$.
\end{proposition}

\begin{proof}
For $a \in A$ we have
\[ \varphi (a) = \varphi (e^*a) =
\langle \underline{a}, \underline{e} \rangle _{\varphi}. \qedhere \]
\end{proof}

\begin{proposition}\label{reprovectorfinvar}%
Let $\varphi$ be a positive linear functional on a \linebreak
\st-algebra $A$. Then $\varphi$ has a reproducing vector if
and only if $\varphi$ has finite variation. In the affirmative
case, we have
\[ v ( \varphi ) = \langle c_{\varphi}, c_{\varphi} \rangle _{\varphi}. \]
\end{proposition}

\begin{proof}
If $\varphi$ has a reproducing vector $c_{\varphi}$, the
Cauchy-Schwarz inequality implies that
\[ {| \,\varphi(a) \,| \,}^2
= {| \,\langle \underline{a}, c_{\varphi} \rangle _{\varphi} \,| \,}^2
\leq {\| \,\underline{a} \,\| \,}^2 \cdot {\| \,c_{\varphi} \,\| \,}^2
= \varphi (a^*a) \cdot {\| \,c_{\varphi}\,\| \,}^2, \]
which shows that
\[ v ( \varphi ) \leq {\| \,c_{\varphi} \,\| \,}^2. \tag*{$(*)$} \]
Conversely assume that $\varphi$ has finite variation. For
$a \in A$ we obtain
\[ | \,\varphi (a) \,| \leq {v ( \varphi ) \,}^{1/2} \cdot {\varphi (a^*a) \,}^{1/2} =
{v ( \varphi ) \,}^{1/2} \cdot \| \,\underline{a} \,\|. \]
Thus $\varphi$ vanishes on the isotropic subspace of
$( A , \langle \cdot , \cdot \rangle _{\varphi} )$, and so $\varphi$
may be considered as a bounded linear functional of norm
$\leq {v(\varphi) \,}^{1/2}$ on the pre-Hilbert space $\underline{A}$.
Hence there exists a vector $c_{\varphi}$ in $H_{\varphi}$ with
\[ \| \,c_{\varphi} \,\| \leq {v ( \varphi ) \,}^{1/2} \tag*{$(**)$} \]
such that
\[ \varphi (a) = \langle \underline{a}, c_{\varphi} \rangle _{\varphi} \]
for all $a \in A$. We have
\[ v ( \varphi ) = \langle c_{\varphi}, c_{\varphi} \rangle _{\varphi} \]
by the inequalities $(*)$ and $(**)$.
\end{proof}

\begin{theorem}\label{finvarcyclic}%
Let $\varphi$ be a positive linear functional of finite vari\-ation on a
normed \st-algebra $A$, which is weakly continuous on $A\sa$.
Let $c_{\varphi}$ be the reproducing vector of $\varphi$. We then have
\[ \pi_{\varphi}(a)c_{\varphi}=\underline{a} \qquad (a \in A). \]
In particular $\pi _{\varphi}$ is cyclic with cyclic vector $c_{\varphi}$
(if $\varphi \neq 0$), and
\[ \varphi(a)=\langle\pi_{\varphi}(a)c_{\varphi},c_{\varphi}\rangle_{\varphi}
\qquad (a \in A).  \]
\end{theorem}

\begin{proof}
For $a, b \in A$, we have
\[ \langle\underline{b},\pi_{\varphi}(a)c_{\varphi}\rangle_{\varphi}
 = \langle\pi_{\varphi}(a^*)\underline{b},c_{\varphi}\rangle_{\varphi}
 = \langle\underline{a^*b},c_{\varphi}\rangle_{\varphi}
 = \varphi(a^*b) = \langle\underline{b},\underline{a}\rangle_{\varphi}, \]
which implies that $\pi_{\varphi}(a)c_{\varphi}=\underline{a}$ for all
$a \in A$. One then inserts this into the defining relation of a
reproducing vector.
\end{proof}

\begin{corollary}\label{Banachfinvarbded}%
Every positive linear functional of finite vari\-ation
on a Banach \st-algebra is continuous.
\end{corollary}

\begin{proof}
Let $\varphi$ be a positive linear functional of finite variation on a Banach
\st-algebra $A$. Then $\varphi$ is weakly continuous on $A\sa$ by
\ref{Banachweakbdedsa}. Hence the GNS representation $\pi _{\varphi}$
exists, cf.\ \ref{weakbdedsarep}, and by \ref{finvarcyclic} we have
$\varphi(a) = \langle \pi _{\varphi} (a) c_{\varphi} , c_{\varphi} \rangle _{\varphi}$.
for all $a \in A$. The continuity of $\pi _{\varphi}$, cf.\ \ref{Banachbounded},
now shows that $\varphi$ is continuous. 
\end{proof}

\begin{theorem}\label{finvarconvex}%
Let $A$ be a normed \st-algebra. The set of positive \linebreak
linear functionals of finite variation on $A$ which are weakly
con\-tinuous on $A\sa$, is a convex cone, and the variation
$v$ is additive and \linebreak $\mathds{R}_+$-homogeneous
on this convex cone.
\end{theorem}

\begin{proof}
It is elementary to prove that $v$ is $\mathds{R}_+$-homogeneous.
Let $\varphi _1$ and $\varphi _2$ be two positive linear
functionals of finite variation on $A$, which are weakly
continuous on $A\sa$. The sum $\varphi := \varphi_1 + \varphi_2$
then also is weakly continuous on $A\sa$. Let $\pi _1$ and $\pi _2$
be the GNS representations associated with $\varphi _1$ and
$\varphi _2$. Consider the direct sum $\pi := \pi _1 \oplus \pi _2$.
It is a non-degenerate representation. Let
$\langle \cdot , \cdot \rangle$ denote the inner product in the
direct sum $H_{\text{\footnotesize{$\varphi_1$}}}
\oplus H_{\text{\footnotesize{$\varphi_2$}}}$.
Let $c _1$ and $c _2$ be the reproducing vectors of
$\varphi _1$ and $\varphi _2$ respectively. For the vector
$c := c _1 \oplus c _2$, one finds
\[ \langle \pi (a) c, c \rangle
= \langle \pi _1 (a) c _1, c _1 \rangle _{\text{\footnotesize{$\varphi_1$}}}
+ \langle \pi _2 (a) c _2, c_2 \rangle _{\text{\footnotesize{$\varphi_2$}}}
= \varphi _1 (a) + \varphi _2 (a) = \varphi (a) \]
for all $a \in A$. Therefore, according to \ref{variationequal}, the
variation of $\varphi$ is
\[ v ( \varphi ) = \langle c, c \rangle
= \langle c _1, c _1 \rangle _{\text{\footnotesize{$\varphi_1$}}}
+ \langle c _2, c _2 \rangle _{\text{\footnotesize{$\varphi_2$}}}
= v ( \varphi _1 ) + v ( \varphi _2 ). \pagebreak \qedhere \]
\end{proof}

\clearpage


\section{States on Normed \texorpdfstring{$*$-}{\052\055}Algebras}

In this paragraph, let $A$ be a normed \st-algebra.

\begin{definition}[states, $S(A)$]\index{state}%
\index{functional!state}\index{S4@$S(A)$}%
A positive linear functional $\psi$ on $A$ shall be called a \underline{state}
if it is weakly continuous on $A\sa$, and has finite variation $v(\psi) = 1$. The
set of states on $A$ is denoted by $\underline{S(A)}$.
\end{definition}

\begin{definition}[quasi-states, $QS(A)$]%
\index{QS(A)@$QS(A)$}\index{quasi-state}\index{functional!quasi-state}%
A \underline{quasi-state} on $A$ shall be a positive linear functional $\varphi$
on $A$ of finite variation $v(\varphi) \leq 1$, which is weakly continuous on
$A\sa$. The set of quasi-states on $A$ is denoted by $\underline{QS(A)}$.
\end{definition}

\begin{theorem}\label{SAconvex}%
The sets $S(A)$ and $QS(A)$ are convex.
\end{theorem}

\begin{proof} \ref{finvarconvex}. \end{proof}

\begin{theorem}\label{staterep}\index{state}%
Let $\pi$ be a non-degenerate $\sigma$-contractive representation
of $A$ in a Hilbert space $H$. If $x$ is a unit vector in $H$, then
\[ \psi (a) := \langle \pi (a) x, x \rangle \qquad (a \in A) \]
defines a state on $A$.
\end{theorem}

\begin{proof} \ref{sigmaAsa}, \ref{variationequal}. \end{proof}

\begin{theorem}\label{repstate}\index{representation!GNS}%
Let $\psi$ be a state on $A$. The GNS representation $\pi _{\psi}$ then is a
$\sigma$-contractive cyclic representation in the Hilbert space $H_{\psi}$.
The reproducing vector $c _{\psi}$ is a cyclic unit vector for $\pi _{\psi}$. We have
\[ \psi (a) = \langle \pi _{\psi} (a) c _{\psi}, c _{\psi} \rangle _{\psi} \qquad (a \in A). \]
\end{theorem}

\begin{proof}
\ref{weakbdedsarep}, \ref{reprovectorfinvar} and \ref{finvarcyclic}.
\end{proof}

\begin{corollary}\label{propquasi}%
A quasi-state $\varphi$ on $A$ is $\sigma$-contractive and Hermitian.
In the auxiliary \st-norm we have $\| \,\varphi \,\| \leq 1$. In particular,
$\varphi$ is contractive on $A\sa$.
\end{corollary}

\begin{proof}
This follows now from \ref{variationinequal} and \ref{sigmaAsa}. \pagebreak
\end{proof}

\begin{definition}\label{topQS}%
We imbed the set $QS(A)$ of quasi-states on $A$ in the unit ball of the dual
space of the real normed space $A\sa$. The space $QS(A)$ is a compact
Hausdorff space in the weak* topology. (By Alaoglu's Theorem, cf.\ the appendix
\ref{Alaoglu}.)
\end{definition}

\begin{theorem}%
Let $\pi$ be a cyclic $\sigma$-contractive representation of $A$ in a Hilbert
space $H$. There then exists a state $\psi$ on $A$ such that $\pi$ is spatially
equivalent to $\pi _{\psi}$. Indeed, each cyclic unit vector $c$ of $\pi$ gives
rise to such a state $\psi$ by putting
\[ \psi (a) := \langle \pi (a) c, c \rangle \qquad (a \in A). \]
There then exists a unitary operator intertwining $\pi$ and $\pi _{\psi}$
 taking $c$ to the reproducing vector $c_{\psi}$.
\end{theorem}

\begin{proof}
This follows from \ref{coeffequal} because
\[ \langle \pi (a) c, c \rangle = \psi (a) =
\langle \pi _{\psi} (a) c_{\psi}, c_{\psi} \rangle _{\psi}. \qedhere \]
\end{proof}

\begin{theorem}\label{stateratio}%
The set of states on $A$ parametrises the class of cyclic $\sigma$-contractive
representations of $A$ in Hilbert spaces up to spatial equivalence. More
precisely, if $\psi$ is a state on $A$, then $\pi _{\psi}$ is a \linebreak
$\sigma$-contractive cyclic representation of $A$. Conversely, if $\pi$ is a
\linebreak $\sigma$-contractive cyclic representation of $A$ in a Hilbert space,
there \linebreak exists a state $\psi$ on $A$ such that $\pi$ is spatially
equivalent to $\pi _{\psi}$.
\end{theorem}

Also of importance is the following uniqueness theorem. (Indeed there is
another way to define the GNS representation associated with a state.)

\begin{theorem}%
The GNS representation $\pi_{\psi}$, associated with a state $\psi$
on $A$, is determined up to spatial equivalence by the condition
\[ \psi (a) = \langle \pi _{\psi} (a) c _{\psi} , c _{\psi} \rangle _{\psi} \qquad (a \in A), \]
in the following sense. If $\pi$ is another cyclic representation of $A$ in a
Hilbert space, with cyclic vector $c$, such that
\[ \psi (a) = \langle \pi (a) c , c \rangle \qquad (a \in A), \]
then $\pi$ is spatially equivalent to $\pi _{\psi}$, and there exists a
unitary operator intertwining $\pi$ and $\pi _{\psi}$, taking $c$ to $c _{\psi}$.
\index{GNS construction|)}
\end{theorem}

\begin{proof} \ref{coeffequal}. \pagebreak \end{proof}

\begin{proposition}\label{plusstate}%
If $\psi$ is a state on $A$, then
\[ \psi (a) \geq 0 \quad \text{ for all } a \in A_+. \]
\end{proposition}

\begin{proof}
For $a \in A_+$, we have $\pi _{\psi} (a) \in B(H) _+$, cf.\ \ref{hompos},
which then has a Hermitian square root $b$ by \ref{C*sqrootrestate}, so
\[ \psi (a) = \langle \pi _{\psi} (a) c _{\psi}, c _{\psi} \rangle _{\psi}
= \langle b c _{\psi}, b c _{\psi} \rangle _{\psi} \geq 0. \qedhere \]
\end{proof}

\medskip
The remaining items \ref{bdedstatesdef} - \ref{kerpsi} give a sort of justification
for our terminology ``state'' for important classes of normed \st-algebras.

\begin{definition}[normed \protect\st-algebras with continuous states]%
\index{algebra!normed s-algebra@normed \protect\st-algebra!with continuous states}%
\index{states!with continuous states}\label{bdedstatesdef}%
\index{continuous states!normed \protect\st-algebra with}%
We shall say that $A$ \underline{has continuous states}, or that $A$ is a normed
\st-algebra \underline{with continuous states}, if every state on $A$ is continuous.
\end{definition}

In practice, $A$ has continuous states as the next proposition
shows. (For a counterexample, see \ref{counterex} below.)

\begin{proposition}\label{bdedstatessuffcond}%
$A$ has continuous states whenever any of the three following conditions is satisfied:
\begin{itemize}
   \item[$(i)$] the involution in $A$ is continuous,
  \item[$(ii)$] $A$ is complete, that is, $A$ is a Banach \st-algebra,
 \item[$(iii)$] $A$ is a commutative \st-subalgebra of a Banach \st-algebra.
\end{itemize}
\end{proposition}

\begin{proof}
\ref{continvol}, \ref{Banachfinvarbded}, and \ref{commsubBsigmacontr}.
\end{proof}

\medskip
We shall prove later that any \st-subalgebra of a Hermitian Banach
\linebreak \st-algebra $A$ has continuous states, cf.\ \ref{Hermcont}.
(The condition that $A$ be Hermitian can actually be dropped,
cf.\ \ref{Palmer}, but we shall not prove this.)

\begin{theorem}\label{bdedstatesthm}%
If $A$ has continuous states, then every $\sigma$-contrac\-tive
representation $\pi$ of $A$ in a Hilbert space $H$ is continuous.
\end{theorem}

\begin{proof} Let $x$ be a vector in the unit ball of $H$.
The positive linear functional
\[ \varphi (a) := \langle \pi (a) x, x \rangle\quad (a \in A) \]
then is a quasi-state (by \ref{sigmaAsa} as well as \ref{variationinequal}),
and hence is continuous. This says that $\pi$ is weakly continuous in the
Hilbert space $H$, and therefore continuous, cf.\ \ref{weakcontHilbert}.
\pagebreak
\end{proof}

\begin{corollary}\label{bdedstatesGNSbded}%
If $A$ has continuous states, then the GNS representations
associated with states on $A$ are continuous.
\end{corollary}

\begin{remark}\label{kerpsi}%
It is well known that a linear functional $\psi$ on a normed space is
continuous if and only if $\ker \psi$ is closed, cf.\ e.g.\ \cite[III.5.3]{Conw},
and thence if and only if $\psi$ has closed graph. 
\end{remark}

\clearpage


\section{Presence of a Unit or Approximate Unit}

\begin{proposition}\label{vare}%
A positive linear functional $\varphi$ on a unital \linebreak \st-algebra
$A$ is Hermitian and has finite variation $v(\varphi) = \varphi(e)$.
\end{proposition}

\begin{proof}
For $a \in A$, we have
\begin{align*}
\varphi (a^*) = \varphi (a^*e) = \langle e, a \rangle _{\varphi}
= \overline{\langle a, e \rangle} _{\varphi} = \overline{\varphi (e^*a)}
= \overline{\varphi (a)},
\end{align*}
as well as
\begin{align*}
{| \,\varphi (a) \,| \,}^2 & = {| \,\varphi (e^*a) \,| \,}^2
= {| \,\langle a, e \rangle _{\varphi} \,| \,}^2 \\
 & \leq \langle a, a \rangle _{\varphi} \,\langle e, e \rangle _{\varphi}
 = \varphi (e^*e) \,\varphi (a^*a) = \varphi (e) \,\varphi (a^*a).\ \qedhere
\end{align*}
\end{proof}

\begin{corollary}%
Let $A$ be a unital normed \st-algebra with continuous states. Then a
state on $A$ can be defined as a continuous positive linear functional
$\psi$ on $A$ with $\psi (e) = 1$.
\end{corollary}

\begin{theorem}\label{eBbded}%
A positive linear functional on a unital Banach \linebreak \st-algebra is
continuous. A state on a unital Banach  \st-algebra can be defined as an
automatically continuous positive linear functional $\psi$ with $\psi (e) = 1$.
\end{theorem}

\begin{proof}
This follows from \ref{vare} and \ref{Banachfinvarbded}.
\end{proof}

\begin{proposition}[extensibility]\label{extendable}
\index{extensible}\index{functional!extensible}%
A positive linear functional $\varphi$ on a \st-algebra
$A$ without unit can be extended to a positive linear functional
$\tld{\varphi}$ on $\tld{A}$ if and only if $\varphi$ is Hermitian
and has finite variation. In this case, $\varphi$ is called
\underline{extensible}. One then can choose for
$\tld{\varphi}(e)$ any real value $\geq v(\varphi)$, and
the choice $\tld{\varphi}(e) = v(\varphi)$ is minimal.
\end{proposition}

\begin{proof}
The ``only if'' part follows from \ref{vare}. Conversely,
let $\varphi$ be Hermitian and have finite variation. Extend
$\varphi$ to a linear functional $\tld{\varphi}$ on $\tld{A}$ with
$\tld{\varphi}(e) = \gamma \geq v (\varphi)$. For $a \in A$ and
$\lambda \in \mathds{C}$, we get (using the Hermitian nature
of $\varphi$):
\begin{align*}
\tld{\varphi}\,\bigl( \,(\lambda e+a)^* (\lambda e+a) \,\bigr)
= & \ \tld{\varphi} \,\bigl( \,{| \,\lambda \,| \,}^2 \,e
+ \overline{\lambda} a + \lambda a^* + a^*a \,\bigr) \\
= & \ \gamma \,{| \,\lambda \,| \,}^2 +
 2 \,\mathrm{Re} \,\bigl( \,\overline{\lambda} \varphi (a) \bigr) + \varphi (a^*a).
\end{align*}
Let now $\lambda =: r \exp (\iu \alpha)$, $a =: \exp (\iu \alpha) \,b$
with $r \in \mathds{R}_+$, $\alpha \in \mathds{R}$, and $b \in A$.
The above expression then becomes
\[ \gamma \,{r \,}^2 + 2 \,r \,\mathrm{Re} \,\bigl( \varphi (b) \bigr) + \varphi (b^*b).
\pagebreak \]
This is non-negative because by construction we have
\[ {\bigl| \,\mathrm{Re} \,\bigl( \varphi (b) \bigr) \,\bigr| \,}^2 \leq {| \,\varphi (b) \,| \,}^2
\leq v (\varphi) \,\varphi (b^*b) \leq \gamma \,\varphi (b^*b). \]
The choice $\tld{\varphi}(e) = v (\varphi)$ is minimal as
a positive linear functional $\psi$ on $\tld{A}$ extending
$\varphi$ satisfies $v(\varphi) \leq v(\psi) = \psi(e)$, cf.\ \ref{vare}.
\end{proof}

\begin{proposition}\label{extquasi}%
Let $\varphi$ be a quasi-state on a normed  \st-algebra $A$
without unit. If $A$ is a pre-C*-algebra and if the completion of $A$
contains a unit, assume that furthermore $\varphi$ is a state.
By extending $\varphi$ to a linear functional $\tld{\varphi}$ on $\tld{A}$
with $\tld{\varphi}(e) := 1$, we obtain a state $\tld{\varphi}$ on $\tld{A}$
which extends $\varphi$. Every state on $\tld{A}$ is of this form.
\end{proposition}

\begin{proof} Recall that a quasi-state is Hermitian by \ref{propquasi},
and thus extensible. We shall show that $\tld{\varphi}$ is continuous on
$\tld{A}\sa$. It is clear that $\tld{\varphi}$ is continuous on $\tld{A}\sa$ with
respect to the norm $| \,\lambda \,| + \| \,a \,\|$, $(\lambda \in \mathds{R}$,
$a \in A\sa)$. Thus, if $A$ is not a pre-C*-algebra, then $\tld{\varphi}$ is
continuous on $\tld{A}\sa$ by definition \ref{normunitis}.
Assume now that $A$ is a pre-C*-algebra. Let $B$ be the completion of $A$.
If $B$ contains no unit, we apply \ref{C*unitis} to the effect that $\tld{\varphi}$
is continuous on $\tld{A}\sa$ with respect to the norm as in \ref{preCstarunitis}.
Also the last statement is fairly obvious in the above two cases.
If $B$ does have a unit, then $A$ is dense in $\tld{A}$ by \ref{complunit}.
We see that the states on $A$ and $\tld{A}$ are related by continuation.
(This uses the continuity of the involution.) This way of extending states
is compatible with the way described in the proposition, by \ref{vare}.
The statement follows immediately.
\end{proof}

\begin{definition}[approximate units]%
\index{approximate unit}\index{unit!approximate}%
Let $A$ be a normed algebra. A \underline{left approximate unit} in $A$
is a net $(e_i) _{i \in I}$ in $A$ such that
\[ a = \lim _{i \in I} e_i a\quad \text{for all } a \in A. \]
Analogously one defines a right approximate unit, a one-sided
approximate unit, as well as a two-sided approximate unit.
A one-sided or two-sided approximate unit $(e_i) _{i \in I}$ in $A$ is
called \underline{bounded} if the net $(e_i) _{i \in I}$ is bounded in $A$.
\end{definition}

Recall that a C*-algebra has a bounded two-sided approximate unit
contained in the open unit ball, cf.\ \ref{C*canapprunit}. \pagebreak

\begin{proposition}\label{bdedfinvar}%
Let $\varphi$ be a continuous positive linear functional on a normed
\st-algebra $A$ with continuous involution and bounded one-sided
approximate unit $(e_i) _{i \in I}$. Then $\varphi$ has finite variation.
\end{proposition}

\begin{proof}
Let $c$ denote a bound of the approximate unit and let $d$ denote
the norm of the involution. The Cauchy-Schwarz inequality implies
\begin{align*}
{| \,\varphi (e_i a) \,| \,}^2
 & = {| \,\langle a, {e_i}^* \rangle _{\varphi} \,| \,}^2
\leq \langle a, a \rangle _{\varphi}
\,\langle {e_i}^*, {e_i}^*  \rangle _{\varphi} \\
 & \leq \varphi (e_i {e_i}^*) \,\varphi (a^*a)
\leq | \,\varphi \,| \,d\,{c\,}^2 \,\varphi (a^*a)
\end{align*}
for all $a \in A$. Hence in the case of a bounded left approximate unit
\[ {| \,\varphi (a) \,| \,}^2 \leq d\,{c\,}^2 \,| \,\varphi \,| \,\varphi (a^*a) \]
for all $a \in A$. If $(e_i)_{i \in I}$ is a bounded right approximate unit,
then $({e_i}^*)_{i \in I}$ is a bounded left approximate unit by continuity
of the involution.
\end{proof}

\begin{theorem}\label{varnorm}\index{v(phi)@$v(\varphi)$}%
Let $A$ be a normed \st-algebra with isometric involution and with
a bounded one-sided approximate unit in the closed unit ball.
If $\varphi$ is a continuous positive linear functional on $A$, then
$\varphi$ has finite variation
\[ v (\varphi) = | \,\varphi \,|. \]
\end{theorem}

\begin{proof} By the preceding proof we have $v (\varphi) \leq | \,\varphi \,|$.
For the converse inequality, please note that $\pi_{\varphi}$ is contractive
by \ref{weakbdedsarep} \& \ref{reprisomcontr}. By \ref{finvarcyclic}:
\[ | \,\varphi (a) \,| =
   | \,\langle \pi_{\varphi} (a) c_{\varphi}, c_{\varphi} \rangle_{\varphi} \,| \]
for all $a \in A$, whence, by \ref{reprovectorfinvar}:
\[ | \,\varphi (a) \,| \leq | \,\pi_{\varphi} (a) \,| \cdot {\| \,c_{\varphi} \,\| \,}^2 \leq
   | \,a \,| \cdot v (\varphi) \qquad (a \in A). \qedhere \]
\end{proof}

\begin{corollary}\label{stateredef}%
Let $A$ be a normed \st-algebra with isometric involution and with
a bounded one-sided approximate unit in the closed unit ball.
Then a state on $A$ can be defined as a continuous positive linear
functional of norm $1$.
\end{corollary}

\begin{proof}
This follows from the fact that $A$ has continuous states,
cf.\ \ref{bdedstatessuffcond} (i).
\end{proof}

\medskip
For the preceding two results, see also \ref{postVarop} below. \pagebreak

\begin{corollary}\label{C*varnorm}
A positive linear functional $\varphi$ on a C*-algebra
is continuous and has finite variation
\[ v (\varphi) = \| \,\varphi \,\|. \]
It follows that a state on a C*-algebra can be defined as an
automatically continuous positive linear functional of norm $1$.
\end{corollary}

\begin{proof}
This follows now from \ref{plfC*bded} and \ref{C*canapprunit}.
\end{proof}

\medskip
Now for the unital case.

\begin{corollary}\label{unitalcase}%
Let $A$ be a normed \st-algebra with isometric \linebreak involution and with
unit $e$ of norm $1$. If $\varphi$ is a continuous positive linear
functional on $A$, then
\[ | \,\varphi \,| = \varphi (e). \]
\end{corollary}

\begin{proof}
This is a trivial consequence of \ref{vare} and \ref{varnorm}.
\end{proof}

\medskip
The next two results are frequently used. 

\begin{theorem}\label{posunitalC*}
A continuous linear functional $\varphi$ on a unital \linebreak
C*-algebra $A$ is positive if and only if
\[ \| \,\varphi \,\| = \varphi (e). \]
\end{theorem}

\begin{proof}
The ``only if'' part follows from the preceding corollary \ref{unitalcase}.
Suppose now that $\| \,\varphi \,\| = \varphi (e)$. We can then assume
that also $\| \,\varphi \,\| = 1$. If $a \in A$ and $\varphi (a^*a)$
is not $\geq 0$, there exists a (large) closed disc
$\{ z \in \mathds{C} : | \,z\0 - z \,| \leq \rho \}$ which contains $\s(a^*a)$,
but does not contain $\varphi (a^*a)$. (Please note that
$\s(a^*a) \subset [0,\infty[$ by the Shirali-Ford Theorem \ref{ShiraliFord}.)
Consequently $\s(z\0 e - a^*a)$ is contained in the disc
$\{ z \in \mathds{C} : | \,z \,| \leq \rho \}$, and so
$\| \,z\0 e -a^*a \,\| = \rlambda (z\0 e - a^*a) \leq \rho$, cf.\ \ref{C*rl}.
Hence the contradiction
\[ | \,z\0 - \varphi (a^*a) \,| = | \,\varphi (z\0 e - a^*a) \,|
\leq \| \,\varphi \,\| \cdot \| \,z\0 e - a^*a \,\| \leq \rho. \qedhere \]
\end{proof}

\begin{corollary}\label{stateunitalC*}
A state on a unital C*-algebra can be defined
as a continuous linear functional $\psi$ such that
\[ \| \,\psi \,\| = \psi (e) = 1. \]
\end{corollary}

\begin{proof}
This follows now from \ref{vare}. \pagebreak
\end{proof}

\begin{theorem}\label{repapprunit}%
Let $A$ be a normed \st-algebra with a bounded left approximate unit
$(e_i) _{i \in I}$. Let $\pi$ be a continuous non-degenerate representation
of $A$ in a Hilbert space $H$. We then have
\[ \lim _{i \in I} \pi (e_i) x = x \quad \text{for all } x \in H. \]
\end{theorem}

\begin{proof}
Suppose first that $x \in H$ is of the form $x = \pi (a)y$ for
some $a \in A$, $y \in H$. By continuity of $\pi$ it follows that
\[ \lim _{i \in I} \| \,\pi (e_i) x - x \,\|
= \lim _{i \in I} \| \,\pi (e_i a-a)y \,\| = 0. \]
Thus $\lim _{i \in I} \pi (e_i)z = z$ for all $z$ in the dense subspace
\[ Z := \mathrm{span} \,\{ \,\pi (a)y \in H : a \in A, y \in H \,\}. \]
Denote by $c$ a bound of the representation, and by $d$ a bound
of the approximate unit. Given $x \in H$, $\varepsilon > 0$, there exists
$z \in Z$ with
\[ \| \,x-z \,\| < \frac{\varepsilon}{3(1+cd)}. \]
Then let $i \in I$ such that for all $j \in I$ with $j \geq i$ one has
\[ \| \,\pi (e_j)z-z \,\| < \frac{\varepsilon}{3}. \]
For all such $j$ we get
\begin{align*}
&\ \| \,\pi (e_j)x-x \,\| \\
\leq &\ \| \,\pi (e_j)x-\pi (e_j)z \,\| + \| \,\pi (e_j)z-z \,\| + \| \,z-x \,\| < \varepsilon.
 \qedhere
\end{align*}
\end{proof}

\medskip
We next give a description of an important class of states obtained
from the Riesz Representation Theorem, cf.\ the appendix \ref{Riesz}.

\begin{lemma}\label{C0Riesz}%
Let $\Omega \neq \varnothing$ be a locally compact Hausdorff space.
For a state $\psi$ on $C\0(\Omega)$ there exists a unique inner regular
Borel probability measure $\mu$ on $\Omega$ such that 
\[ \psi (f) = \int f \,d \mu\quad \text{for all } f \in C\0(\Omega). \]
(For the definition of an ``inner regular Borel probability measure'',
see the appendix \ref{Bairedef} \& \ref{inregBormeas}.)
\end{lemma}

\begin{proof}
This follows from the Riesz Representation Theorem, cf.\ the \linebreak
appendix \ref{Riesz}. Please note \ref{pointwiseorder}, \ref{C*varnorm},
and \ref{Ccdense}. \pagebreak
\end{proof}

\begin{definition}[$M_1(\Omega)$]\index{M1@$M_1(\Omega)$}%
If $\Omega \neq \varnothing$ is a locally compact Hausdorff space,
one denotes by \underline{$M_1(\Omega)$} the set of inner regular
Borel probability measures on $\Omega$, cf.\ the appendix
\ref{Bairedef} \& \ref{inregBormeas}.
\end{definition}

\begin{theorem}\label{predisint}%
For a state $\psi$ on a commutative C*-algebra $A$, there
exists a unique measure $\mu \in M_1\bigl(\Delta(A)\bigr)$ such that
\[ \psi (a) = \int \wht{a} \,d \mu\quad \text{for all } a \in A. \]
\end{theorem}

\begin{proof}
We note that $A \neq \{ 0 \}$ because $\psi$ is a state on $A$. It follows
that the Gelfand transformation $a \mapsto \wht{a}$ establishes an
isomorphism of the C*-algebra $A$ onto the C*-algebra $C\0\bigl(\Delta(A)\bigr)$,
by the Commutative Gelfand-Na\u{\i}mark Theorem \ref{commGN}. By defining
\[ \wht{\psi}(\wht{\vphantom{\psi}a}) := \psi (a)
\quad \text{for all } a \in A, \]
we get a state $\wht{\psi}$ on $C\0\bigl(\Delta(A)\bigr)$. By lemma \ref{C0Riesz}
there exists a measure $\mu \in M_1\bigl(\Delta(A)\bigr)$ such that
\[ \psi (a) = \wht{\psi}(\wht{\vphantom{\psi}a})
= \int \wht{\vphantom{\psi}a} \,d \mu
\quad \text{for all } a \in A. \]
Uniqueness is now clear from lemma \ref{C0Riesz} again.
\end{proof}

\medskip
We shall return to this theme in \ref{Commutativity}.

The impatient reader can go from here to part \ref{part3} of the book, on the
spectral theory of representations.

\clearpage


\section{The Theorem of Varopoulos}%
\label{Varopsect}

Our next objective is the Theorem of Varopoulos \ref{Varopoulos}. As it
is little used, we do not give a complete proof. We state without proof
two closely related Factorisation Theorems of Paul J.\ Cohen and Nicholas
Th.\ Varopoulos. We give them here in a form needed later.

\begin{theorem}[Paul J.\ Cohen]\label{Cohen}%
\index{factorisation!Theorem!Cohen}%
\index{Theorem!Cohen's Factorisation}\index{Cohen's Factorisation Thm.}%
Let $A$ be a Banach algebra with a boun\-ded right approximate unit.
Then every element of $A$ is a product of two elements of $A$.
\end{theorem}

\begin{theorem}[Nicholas Th.\ Varopoulos]\label{Varopoulosfact}%
\index{Theorem!Varopoulos!Factorisation}%
\index{Varopoulos' Theorem!Factorisation}%
\index{factorisation!Theorem!Varopoulos}%
Let $A$ be a  Banach algebra with a bounded right approximate unit. For
a sequence $(a_n)$ in $A$ converging to $0$, there exists a sequence
$(b_n)$ in $A$ converging to $0$ as well as $c \in A$ with $a_n = {b_n}c$
for all $n$.
\end{theorem}

Proofs can be found in \cite[\S\ 11]{BD}, \cite[\S\ 6]{Zel}, or \cite[section 5.2]{Palm}.

\begin{theorem}\label{preVarop1}%
Let $A$ be a Banach \st-algebra with a bounded right approximate unit.
A positive Hilbert form $\langle \cdot , \cdot \rangle$ on $A$, cf.\ \ref{Hilbertform},
is bounded in the sense that there exists $\gamma \geq 0$ such that
\[ | \,\langle a, b \rangle \,| \leq \gamma \cdot | \,a \,| \cdot | \,b \,|
\quad \text{for all } a, b \in A. \]
\end{theorem}

\begin{proof}
It shall first be shown that $\langle \cdot , \cdot \rangle$
is continuous in the first variable. So let $(a_n)$ be a sequence in $A$
converging to $0$. By the Factorisation Theorem of N.\ Th.\ Varopoulos
(theorem \ref{Varopoulosfact} above), there exists a sequence $(b_n)$
converging to $0$ in $A$ as well as $c \in A$, such that $a_n = {b_n}c$
for all $n$. Thus, for arbitrary $d \in A$, we obtain
\[ \langle a_n, d \rangle = \langle {b_n}c, d \rangle
= \langle \pi(b_n) \underline{c}, \underline{d} \rangle \to 0 \]
by continuity of $\pi$ (theorem \ref{Banachbounded}). It follows that
$\langle \cdot , \cdot \rangle$ is separately continuous in both variables
(as $\langle \cdot , \cdot \rangle$ is Hermitian). An application of the
uniform boundedness principle then shows that there exists
$\gamma \geq 0$ such that
\[ | \,\langle a, b \rangle \,| \leq \gamma \cdot | \,a \,| \cdot | \,b \,|
\quad \text{for all } a, b \in A. \pagebreak  \qedhere \]
\end{proof}

\begin{theorem}\label{preVarop2}%
Let $A$ be a normed \st-algebra with a bounded right approximate
unit $(e_i)_{i \in I}$. Let $\langle \cdot , \cdot \rangle$ be a bounded
non-zero positive Hilbert form on $A$. The corresponding GNS
representation $\pi$ then is cyclic. Furthermore, there exists a cyclic
vector $c$ for $\pi$, such that the positive linear functional of finite variation
$a \mapsto \langle \pi (a) c, c \rangle$ $(a \in A)$ induces the positive
Hilbert form $\langle \cdot , \cdot \rangle$. (Cf.\ \ref{inducedHf}.) The
vector $c$ can be chosen to be any adherence point of the net
$(\underline{e_i})_{i \in I}$ in the weak topology. 
\end{theorem}

\begin{proof}
Since the representation associated with $\langle \cdot , \cdot \rangle$
is weakly continuous on $A\sa$, it has an ``extension'' $\pi$ acting on
the completion of $\underline{A}$, cf.\ \ref{weaksa}. Next, by boundedness
of $\langle \cdot , \cdot \rangle$, the net $(\underline{e_i})_{i \in I}$
is bounded and therefore has an adherence point in the weak topology
on the completion of $\underline{A}$. Let $c$ be any such adherence point.
By going over to a subnet, we can assume that $(\underline{e_i})_{i \in I}$
converges weakly to $c$. It shall be shown that $\pi (a)c = \underline{a}$
for all $a \in A$, which implies that $c$ is cyclic for $\pi$. It is enough
to show that for $b \in A$, one has $\langle \pi (a)c, \underline{b} \rangle
= \langle \underline{a}, \underline{b} \rangle$. So one calculates
\begin{align*}
 & \langle \pi (a)c, \underline{b} \rangle
= \langle c, \pi (a)^* \underline{b} \rangle
= \lim _{i \in I} \,\langle \underline{e_i}, \pi (a)^* \underline{b} \rangle \\
= \,&\lim _{i \in I} \,\langle \pi (a)\underline{e_i}, \underline{b} \rangle
= \lim _{i \in I} \,\langle a e_i, b \rangle
= \langle a, b \rangle = \langle \underline{a}, \underline {b} \rangle,
\end{align*}
where we have used the continuity of $\langle \cdot, \cdot \rangle$.
Consider now the positive linear functional $\varphi$ on $A$ defined by
\[ \varphi (a) := \langle \pi (a) c, c \rangle \qquad (a \in A). \]
We have
\[ \varphi (a^*a) = \langle \pi (a^*a) c, c \rangle = {\| \,\pi (a)c \,\| \,}^2
= {\| \,\underline{a} \,\| \,}^2
= \langle \underline{a}, \underline{a} \rangle
= \langle a, a \rangle \]
for all $a \in A$, whence, after polarisation,
\[ \langle a, b \rangle = \varphi(b^*a) \qquad (a, b \in A). \]
That is, $\varphi$ induces the positive Hilbert form $\langle \cdot , \cdot \rangle$,
cf.\ \ref{inducedHf}.
\end{proof}

\medskip
We see from the proofs that so far it is the bounded right
approximate units which are effective. Now for the case
of a bounded two-sided approximate unit. We shall prove
that in the preceding theorem, if we are in presence of a
bounded two-sided approximate unit $(e_i)_{i \in I}$, and
if $A$ is complete, then the cyclic vector $c$ can be chosen
to be the limit in norm of the net $(\underline{e_i})_{i \in I}$.
\pagebreak

\begin{theorem}
Let $A$ be a Banach \st-algebra with a bounded two-sided
approximate unit $(e_i)_{i \in I}$. Let $\langle \cdot , \cdot \rangle$
be a non-zero positive Hilbert form on $A$. The corresponding GNS
representation $\pi$ then is cyclic. Furthermore, there exists a cyclic
vector $c$ for $\pi$, such that the positive linear functional of finite
variation $a \mapsto \langle \pi (a) c, c \rangle$ $(a \in A)$ induces
the positive Hilbert form $\langle \cdot , \cdot \rangle$. The vector
$c$ can be chosen to be the limit in norm of the net
$(\underline{e_i})_{i \in I}$.
\end{theorem}

\begin{proof}
The Hilbert form $\langle \cdot , \cdot \rangle$ is bounded by
\ref{preVarop1}, so \ref{preVarop2} is applic\-able. Only the last
statement needs to be proved. The representation $\pi$ is
continuous by \ref{Banachbounded}, and as in the preceding
proof, we have $\underline{e_i} = \pi (e_i)c$ for all $i \in I$. Since
$(e_i)_{i \in I}$ among other things is a bounded left approximate
unit, \ref{repapprunit} implies that $\underline{e_i} = \pi (e_i)c \to c$
because $\pi$ is cyclic by the preceding theorem, and thus
non-degenerate. 
\end{proof}

\medskip
Next for bounded one-sided approximate units:

\begin{theorem}[Varopoulos, Shirali]\label{Varopoulos}%
\index{Theorem!Varopoulos}\index{Varopoulos' Theorem}%
A positive linear functional on a Banach \st-algebra with a bounded
one-sided approximate unit has finite variation, and thence also is
continuous, cf.\ \ref{Banachfinvarbded}.
\end{theorem}

\begin{proof} In the case of a bounded right approximate unit, the
conclusion follows from the theorems \ref{preVarop1} and
\ref{preVarop2} by an application of Cohen's Factorisation Theorem
\ref{Cohen}. In presence of a bounded left approximate unit, we
change the multiplication to $(a,b) \mapsto ba$, and use the fact
that the functional turns out to be Hermitian, cf.\ \ref{variationinequal}.
\end{proof}

\medskip
This result is due to Varopoulos for the special case of a continuous
involution. It was extended to the general case by Shirali.

\begin{corollary}\index{v(phi)@$v(\varphi)$}\label{postVarop}%
Let $A$ be a Banach \st-algebra with isometric \linebreak involution
and with a bounded one-sided approximate unit in the closed unit ball.
Every positive linear functional $\varphi$ on $A$ is continuous and
has finite variation
\[ v (\varphi) = | \,\varphi \,|. \]
It follows that a state on $A$ can be defined as an automatically
con\-tinuous positive linear functional of norm $1$.
\end{corollary}

\begin{proof}
This follows now from \ref{varnorm} and \ref{stateredef}.
\end{proof}

\medskip
For C*-algebras, we don't need this result, see \ref{C*varnorm}.
\pagebreak

\clearpage


\chapter{The Enveloping C*-Algebra}

\setcounter{section}{27}


\section{C*(A)}

In this paragraph, let $A$ be a normed \st-algebra.

\begin{proposition}%
A state $\psi$ on $A$ satisfies
\[ {\psi (a^*a) \,}^{1/2} \leq \rsigma(a) \qquad \text{for all } a \in A. \]
\end{proposition}

\begin{proof} Theorem \ref{repstate} implies that
\[ {\psi(a^*a) \,}^{1/2}
= {\langle \pi _{\psi}(a^*a) c_{\psi}, c_{\psi} \rangle _{\psi}}^{1/2}
= \| \,\pi _{\psi}(a)c_{\psi} \,\| \leq \rsigma(a). \qedhere \]
\end{proof}

\begin{definition}[$\,\| \cdot \|\0\,$]\index{a6@${"|}{"|}\,a\,{"|}{"|}\0$}%
For $a \in A$, one defines
\[ \| \,a \,\|\0 := \,\sup _{\text{\small{$\psi \in S(A)$}}}
\,{\psi (a^*a) \,}^{1/2} \leq \rsigma(a). \]
If $S(A) = \varnothing$, one puts $\| \,a \,\|\0 := 0$ for all $a \in A$.
\end{definition}

We shall see that $\| \cdot \|\0$ is a seminorm with certain
extremality properties.

\begin{theorem}\label{boundedo}%
If $\pi$ is a representation of $A$ in a pre-Hilbert space $H$,
which is weakly continuous on $A\sa$, then
\[ \| \,\pi (a) \,\| \leq \| \,a \,\|\0 \qquad \text{for all } a \in A. \]
\end{theorem}

\begin{proof}
We may assume that $H$ is a Hilbert space $\neq \{ 0 \}$ and
that $\pi$ is non-degenerate. Let $x \in H$ with $\| \,x \,\| = 1$.
Consider the state $\psi$ on $A$ defined by
\[ \psi (a) := \langle \pi (a) x, x \rangle \qquad (a \in A). \]
We obtain
\[ \| \,\pi (a) x \,\| = {\langle \pi (a) x, \pi (a) x \rangle \,}^{1/2}
= {\psi (a^*a) \,}^{1/2} \leq \| \,a \,\|\0. \pagebreak \qedhere \]
\end{proof}

\begin{definition}[the universal representation, $\pi\univ$]%
\index{representation!universal}\index{universal!representation}%
\index{p3@$\pi_{u}$}\index{H2@$H_{u}$}%
We define
\[ \pi\univ := \oplus \,_{\text{\small{$\psi \in S(A)$}}}
\,\pi _{\text{\small{$\psi$}}}, \qquad
H\univ := \oplus \,_{\text{\small{$\psi \in S(A)$}}}
\,H _{\text{\small{$\psi$}}}. \]
One says that $\pi \univ$ is the \underline{universal representation}
of $A$ and that $H \univ$ is the universal Hilbert space of $A$.
The universal representation is $\sigma$-contractive
by \ref{directsumsigma}. 
\end{definition}

\begin{theorem}\label{GNuniv}%
For $a \in A$, we have
\[ \| \,\pi \univ (a) \,\| = \| \,a \,\|\0. \]
\end{theorem}

\begin{proof}
The representation $\pi \univ$ is $\sigma$-contractive as noted
before, so that, by \ref{boundedo},
\[ \| \,\pi \univ (a) \,\| \leq \| \,a \,\|\0 \]
for all $a \in A$. In order to show the converse inequality, let
$a \in A$ be fixed. Let $\psi \in S(A)$. It suffices then to prove that
$\psi (a^*a) \leq {\| \,\pi \univ (a) \,\| \,}^2$. We have
\begin{align*}
\psi (a^*a)
& = \langle \pi _{\psi} (a^*a) c_{\psi}, c_{\psi} \rangle _{\psi}
= \langle \pi _{\psi} (a) c_{\psi}, \pi _{\psi} (a) c_{\psi} \rangle _{\psi} \\
& = {\| \,\pi _{\psi} (a) c_{\psi} \,\| \,}^2
\leq {\| \,\pi _{\psi} (a) \,\| \,}^2 \leq {\| \,\pi \univ (a) \,\| \,}^2. \qedhere
\end{align*}
\end{proof}

\begin{definition}[the Gelfand-Na\u{\i}mark seminorm]%
\index{Gelfand-Naimark@Gelfand-Na\u{\i}mark!seminorm}%
\label{GNSseminorm}%
We see from the preceding theorem that $\| \cdot \| \0$ is a
seminorm on $A$. The seminorm $\| \cdot \| \0$ is called the
\underline{Gelfand-Na\u{\i}mark seminorm}. It satisfies
\[ \| \,a \,\|\0 \leq \rsigma(a) \leq \| \,a \,\| \qquad (a \in A), \]
where $\| \cdot \|$ denotes the auxiliary \st-norm on $A$.
\end{definition}

\begin{proposition}\label{bstGNc}%
If $A$ has continuous states, then the Gelfand-Na\u{\i}mark seminorm
is continuous on $A$. If $A$ furthermore is commutative, then the
Gelfand-Na\u{\i}mark seminorm is dominated by the norm in $A$.
\end{proposition}

\begin{proof}
This follows from \ref{bdedstatesthm} and \ref{commbdedcontr}.
\end{proof}

\medskip
A basic object in representation theory is the following enveloping
\linebreak C*-algebra. It is useful for pulling down results. \pagebreak

\begin{definition}[the enveloping C*-algebra, $C^*(A)$]%
\index{s035@\protect\st-rad(A)}\index{s04@\protect\st-radical}%
\index{C5@$C^*(A)$}\index{j@$j$}\index{C6@C*-algebra!enveloping}%
\index{enveloping C*-algebra}\index{algebra!C*-algebra!enveloping}%
The \underline{\st-radical} of $A$ is defined to be the \st-stable
\ref{selfadjointsubset} two-sided ideal on which $\| \cdot \|\0$ vanishes.
It shall be denoted by \st-$\mathrm{rad}(A)$.
We introduce a norm $\| \cdot \|$ on $A/$\st-$\mathrm{rad}(A)$
by defining
\[ \| \,a+\text{\st-}\mathrm{rad}(A) \,\| := \| \,a \,\|\0 \text{ for all } a \in A. \]
Then $(A/$\st-$\mathrm{rad}(A), \| \cdot \,\|)$ is isomorphic as a normed
\st-algebra to the range of the universal representation, and thus is a
pre-C*-algebra. Let $(C^*(A),\| \cdot \|)$ denote the completion. It is
called the \underline{enveloping} \underline{C*-algebra} of $A$. We put
\[ j : A \to C^*(A),\quad a \mapsto a+\text{\st-}\mathrm{rad}(A). \]
\end{definition}

\begin{proposition}\label{jdense}%
The \st-algebra homomorphism $j : A \to C^*(A)$ is \linebreak
$\sigma$-contractive and its range is dense in $C^*(A)$.
\end{proposition}

\begin{proposition}%
If $A$ has continuous states, then the mapping \linebreak
$j : A \to C^*(A)$ is continuous. If furthermore $A$ is commutative,
then $j$ is contractive.
\end{proposition}

\begin{proof}
This follows from \ref{bstGNc}.
\end{proof}

\medskip
We now come to a group of statements which collectively might be
called the ``\underline{universal property}'' of the enveloping C*-algebra.
\index{universal!property!of enveloping C*-algebra}

\begin{proposition}\index{p4@$\pi \0$}%
Let $H$ be a Hilbert space. The assignment
\[ \pi := \pi \0 \circ j \]
establishes a bijection from the set of all representations $\pi \0$ of
$C^*(A)$ in $H$ onto the set of all $\sigma$-contractive representations
$\pi$ of $A$ in $H$.
\end{proposition}

\begin{proof} Let $\pi \0$ be a representation of $C^*(A)$
in $H$. It then is contractive, and it follows that $\pi := \pi \0 \circ j$ is
a $\sigma$-contractive representation of $A$ in $H$. Conversely, let $\pi$
be a $\sigma$-contractive representation of $A$ in $H$. We then have
$\| \,\pi (a) \,\| \leq \| \,a \,\|\0$ by \ref{boundedo}. It follows that $\pi$
vanishes on \st-$\mathrm{rad}(A)$ and thus induces a representation
$\pi \0$ of $C^*(A)$ such that $\pi = \pi \0 \circ j$. To show injectivity, let
$\rho \0$, $\sigma \0$ be two representations of $C^*(A)$ in $H$ with
$\rho \0 \circ j = \sigma \0 \circ j$. Then $\rho \0$ and $\sigma \0$
agree on the dense set $j(A)$, hence everywhere. \pagebreak
\end{proof}

\begin{proposition}\label{piopi}%
Let $\pi$ be a representation of $A$ in a Hilbert space $H$, which is
weakly continuous on $A\sa$. Let $\pi \0$ be the corresponding
representation of $C^*(A)$. Then $R(\pi \0)$ is the closure of $R(\pi)$.
It follows that
\begin{itemize}
   \item[$(i)$] $\pi$ is non-degenerate if and only if $\pi \0$ is so,
  \item[$(ii)$] $\pi$ is cyclic if and only if $\pi \0$ is so,
 \item[$(iii)$] ${\pi \0}' = {\pi}'.$
\end{itemize}
\end{proposition}

\begin{proof}
The range of $\pi \0$ is a C*-algebra by \ref{rangeC*}.
Furthermore $R(\pi)$ is dense in $R(\pi \0)$ by \ref{jdense},
from which it then follows that $R(\pi \0)$ is the closure of $R(\pi)$.
\end{proof}

\begin{proposition}%
If $\psi$ is a state on $A$, then
\[ | \,\psi(a) \,| \leq \| \,a \,\|\0 \qquad \text{for all } a \in A. \]
\end{proposition}

\begin{proof} We have
\[ | \,\psi (a) \,| \leq {\psi (a^*a) \,}^{1/2} \leq
\sup _{\varphi \in S(A)} {\varphi (a^*a) \,}^{1/2} = \| \,a \,\|\0. \qedhere \]
\end{proof}

\begin{corollary}\label{linboundedo}%
If $\varphi$ is a positive linear functional of  finite variation
on $A$, which is weakly continuous on $A\sa$, then
\[ | \,\varphi (a) \,| \leq v( \varphi ) \,\| \,a \,\|\0 \qquad \text{for all } a \in A. \]
\end{corollary}

\begin{proposition}\index{f1@$\varphi \0$}%
By putting
\[ \varphi := \varphi \0 \circ j \]
we establish a bijection from the set of positive linear functionals
$\varphi \0$ on $C^*(A)$ onto the set of positive linear functionals
$\varphi$ of finite variation on $A$, which are weakly continuous on
$A\sa$. Then $v ( \varphi ) = v ( \varphi \0 )$.
\end{proposition}

\begin{proof} Let $\varphi \0$ be a positive linear functional on $C^*(A)$.
Then $\varphi \0$ is continuous by \ref{plfC*bded} and consequently has
finite variation by \ref{bdedfinvar} as $C^*(A)$ has a bounded two-sided
approximate unit \ref{C*canapprunit}. It follows that $\varphi := \varphi \0 \circ j$
is a positive linear functional on $A$, which is weakly continuous on $A\sa$
because $j$ is contractive on $A\sa$. For $a \in A$,\,we have
\begin{gather*}
{| \,\varphi (a) \,| \,}^2 = {| \,( \varphi \0 \circ j ) (a) \,| \,}^2, \\
\varphi (a^*a) = ( \varphi \0 \circ j ) (a^*a),
\end{gather*}
so that $v ( \varphi ) = v ( \varphi \0 |_{j(A)})$. By continuity of $\varphi \0$
on $C^*(A)$ and by continuity of the involution in $C^*(A)$, it follows that
$v ( \varphi ) = v ( \varphi \0 )$. \pagebreak

Conversely, let $\varphi$ be a positive linear functional of finite variation
on $A$, which is weakly continuous on $A\sa$. By \ref{linboundedo},
$\varphi$ also is bounded with respect to $\| \cdot \|\0$. It therefore
induces a positive \linebreak linear functional $\varphi \0$ on $C^*(A)$
such that $\varphi = \varphi \0 \circ j$. To show \linebreak injectivity, let
$\varphi \0$, ${\varphi \0}'$ be positive linear functionals on $C^*(A)$ with
$\varphi \0 \circ j = {\varphi \0}' \circ j$. Then $\varphi \0$ and ${\varphi \0}'$
coincide on the dense set $j(A)$, hence \linebreak everywhere, by automatic
continuity of positive linear functionals on \linebreak C*-algebras \ref{plfC*bded}. 
\end{proof}

\begin{reminder}%
We imbed the set $QS(A)$ of quasi-states on $A$ in the unit ball of the dual
space of the real normed space $A\sa$. The space $QS(A)$ is a compact
Hausdorff space in the weak* topology. (By Alaoglu's Theorem, cf.\ the appendix
\ref{Alaoglu}.)
\end{reminder}

\begin{theorem}\label{affinehomeom}%
The assignment
\[ \varphi := \varphi \0 \circ j \]
is an \underline{affine homeomorphism} from the set $QS\bigl(C^*(A)\bigr)$ of
quasi-states $\varphi \0$ on $C^*(A)$ onto the set $QS(A)$ of quasi-states
$\varphi$ on $A$, and hence from the set $S\bigl(C^*(A)\bigr)$ onto the set $S(A)$.
\end{theorem}

\begin{proof}
This assignment is continuous by the universal property of the weak* topology,
cf.\ the appendix \ref{weak*top}. The assignment, being a continuous bijection
from the compact space $QS\bigl(C^*(A)\bigr)$ to the Hausdorff space $QS(A)$,
is a homeomorphism, cf.\ the appendix \ref{homeomorph}. 
\end{proof}

\medskip
The next two results put the ``extension processes'' for states
and their GNS representations in a useful relationship.

\begin{proposition}\label{piopsio}%
Let $\psi$ be a state on $A$, and let $\psi\0$ be the corresponding state
on $C^*(A)$. Let $(\pi _{\psi}) \0$ denote the representation of $C^*(A)$
corresponding to the GNS representation $\pi _{\psi}$. We then have
\[ \psi \0 (b) = \langle (\pi _{\psi}) \0 (b) c _{\psi},c _{\psi} \rangle _{\psi}
\quad \text{ for all } b \in C^*(A). \]
\end{proposition}

\begin{proof}
Consider the state $\psi _1$ on $C^*(A)$ defined by
\[ \psi _1 (b) := \langle (\pi _{\psi}) \0 (b) c _{\psi}, c _{\psi} \rangle _{\psi}
\qquad \bigl(b \in C^*(A)\bigr). \]
We get
\[ \psi _1 \bigl(j(a)\bigr)
= \langle (\pi _{\psi}) \0 \bigl(j(a)\bigr) c _{\psi}, c_{\psi} \rangle _{\psi}
= \langle \pi _{\psi} (a)c _{\psi}, c_{\psi} \rangle _{\psi} = \psi (a), \]
for all $a \in A$, so that $\psi _1 = \psi \0$. \pagebreak
\end{proof}

\begin{corollary}\label{oequiv}%
For a state $\psi$ on $A$, the representation $\pi _{(\psi \0)}$ is spatially
equivalent to $(\pi _{\psi}) \0$, and there exists a unitary operator intertwining
$\pi _{(\psi \0)}$ and $(\pi _{\psi}) \0$ taking $c _{(\psi \0)}$ to $c _{\psi}$.
\end{corollary}

\begin{proof} \ref{coeffequal}. \end{proof}

\begin{proposition}\label{reform2}%
If $A$ has continuous states, then \st-$\mathrm{rad}(A)$ is closed in $A$.
\end{proposition}

\begin{proof}
This is a consequence of \ref{bstGNc}. For Banach \st-algebras,
this also follows from \ref{kerclosed}.
\end{proof}

\medskip
We close this paragraph with a discussion
of \st-semisimple normed \linebreak \st-algebras.

\begin{definition}[\st-semisimple normed \st-algebras]%
\index{s05@\protect\st-semisimple}%
We  will say that $A$ is \underline{\st-semisimple} if
\st-$\mathrm{rad}(A) = \{0\}$.
\end{definition}

Please note that $A$ is \st-semisimple if and only
if its universal representation is faithful.

\begin{proposition}\label{faithsemi}%
If $A$ has a faithful representation in a pre-Hilbert space, which is
weakly continuous on $A\sa$, then $A$ is \st-semi\-simple.
\end{proposition}

\begin{proof}
A representation which is weakly continuous on
$A\sa$ vanishes on the \st-radical by \ref{boundedo}.
\end{proof}

\begin{proposition}\label{stsemipreC*}%
If $A$ is \st-semisimple, then $(A, \| \cdot \|\0)$
is a pre-C*-algebra and $C^*(A)$ is its completion. 
\end{proposition}

\begin{proposition}\label{reform1}%
If $A$ is \st-semisimple, then the involution in $A$ has closed
graph. In particular, it then follows that $A\sa$ is closed in $A$.
\end{proposition}

\begin{proof}
This follows from \ref{injclosgraph}.
\end{proof}

\begin{corollary}\label{stsemiscont}%
The involution in a \st-semisimple Banach \st-al\-gebra is continuous.
\pagebreak
\end{corollary}

\clearpage


\section{The Theorems of Ra\texorpdfstring{\u{\i}}{i}kov and of %
Gelfand \texorpdfstring{\&}{\046} Na\texorpdfstring{\u{\i}}{i}mark}

Our next objective is a result of Ra\u{\i}kov \ref{Raikov}, of which the
Gelfand-Na\u{\i}mark Theorem \ref{GelfandNaimark} is a special case.

\begin{theorem}\label{Krein}%
Let $C$ be a convex cone in a real vector space $B$. Let $f\0$ be a linear
functional on a subspace $M\0$ of $B$ such that $f\0(x) \geq 0$ for all
$x \in M\0 \cap C$. Assume that $M\0+C = B$. The functional $f\0$
then has a linear extension $f$ to $B$ such that $f(x) \geq 0$ for all $x \in C$.
\end{theorem}

\begin{proof}
Assume that $f\0$ has been extended to a linear functional $f_1$ on a
subspace $M_1$ of $B$ such that $f_1(x) \geq 0$ for all $x \in M_1 \cap C$.
Assume that $M_1 \neq B$. Let $y \in B$, $y \notin M_1$ and put
$M_2 := \mathrm{span} \,( M_1 \cup \{ \,y \,\} )$. Define
\[ E' := M_1 \cap (y-C), \quad E'' := M_1 \cap (y+C). \]
The sets $E'$ and $E''$ are non-void because $y, -y \in M\0+C$
(by the assumption that $M\0+C = B$).
If $x' \in E'$, and $x'' \in E''$, then $y-x' \in C$ and $x''-y \in C$,
whence $x''-x' \in C$, so that $f_1(x') \leq f_1(x'')$. It follows that with
\[ a := \sup _{x' \in E'} f_1(x'),\quad b := \inf _{x'' \in E''} f_1(x'') \]
one has $-\infty < a \leq b < \infty$. With $a \leq c \leq b$ we then have
\[ f_1(x') \leq c \leq f_1(x'') \qquad (x' \in E', x'' \in E''). \]
Define a linear functional $f_2$ on $M_2 = \mathrm{span} \,( M_1 \cup \{ \,y \,\} )$ by
\[ f_2(x+\alpha y) := f_1(x)+\alpha c\quad (x \in M_1, \alpha \in \mathds{R}). \]
It shall be shown that $f_2(x) \geq 0$ for all $x \in M_2 \cap C$. So let
$x \in M_1$, $\alpha \in \mathds{R} \setminus \{0\}$. There then exist $\beta > 0$
and $z \in M_1$ such that either $x+\alpha y = \beta (z+y)$ or $x+\alpha y =
\beta (z-y)$. Thus, if $x+\alpha y \in C$, then either $z+y \in C$ or $z-y
\in C$. If $z+y \in C$, then $-z \in E'$, so that $f_1(-z) \leq c$, or
$f_1(z) \geq -c$, whence $f_2(z+y) \geq 0$. If $z-y \in C$, then $z \in E''$,
so that $f_1(z) \geq c$, whence $f_2(z-y) \geq 0$. Thus $f_1$ has been
extended to a linear functional $f_2$ on the subspace $M_2$ properly
containing $M_1$, such that $f_2(x) \geq 0$ for all $x \in M_2 \cap C$.

One defines
\begin{align*}
Z := \{ \,(M,f) : \ & M \text{ is a subspace of } B \text{ containing } M\0, \\
 & f \text{ is a linear functional on } M \text{ extending } f\0 \\
 & \text{such that } f(x) \geq 0 \text{ for all } x \in M \cap C \,\}. 
\end{align*}
Applying Zorn's Lemma yields a maximal element $(M,f)$ of $Z$ (with the
obvious ordering). It follows from the above that $M = B$. \pagebreak
\end{proof}

\begin{theorem}\label{extstate}%
If $A$ is a Hermitian Banach \st-algebra, then a state on a closed
\st-subalgebra of $A$ can be extended to a state on all of $A$.
\end{theorem}

\begin{proof} Let $D$ be a closed \st-subalgebra of the Hermitian Banach
\linebreak
\st-algebra $A$. Put $E := \mathds{C}e+D$, where $e$ is the unit in $\tld{A}$.
Then $E$ is a Hermitian Banach \st-algebra\vphantom{$\tld{A}$}, cf.\ \ref{Herminher}
as well as the proof of \ref{Banachunitis}. Let $\psi$ be a state on $D$. Putting
$\tld{\psi}(e) := 1$ and extending linearly, we get a state $\tld{\psi}$ on $E$ which
extends $\psi$, cf.\ the proof of \ref{extendable}, as well as \ref{eBbded}. We want
to apply the preceding theorem with $B := \tld{A}\sa$, $C := \tld{A}_+$,
$M\0 := E\sa$, $f\0 := \tld{\psi}|_{\text{\small{$E\sa$}}}$. The assumptions are
satisfied because $\tld{\psi}$ is Hermitian by \ref{propquasi}, $\tld{A}_+$ is convex
by \ref{plusconvexcone}, $\tld{\psi}(b) \geq 0$ for all $b \in E_+$ by \ref{plusstate},
$E_+ = \tld{A}_+ \cap E\sa$ by \ref{plussubal}, and $E\sa+\tld{A}_+ = \tld{A}\sa$
because for $a \in \tld{A}\sa$ one has $\rlambda(a)e+a \geq 0$. By the preceding
extension theorem, it follows that $\tld{\psi}|_{\text{\small{$E\sa$}}}$ has an
extension to a linear functional $f$ on $\tld{A}\sa$ such that $f(a) \geq 0$ for all
$a \in \tld{A}_+$. Let $\tld{\varphi}$ be the linear extension of $f$ to $\tld{A}$. We
then have $\tld{\varphi} \,(a^*a) \geq 0$ for all $a \in \tld{A}$ because
$a^*a \in \tld{A}_+$ by the Shirali-Ford Theorem \ref{ShiraliFord}. Moreover,
$\tld{\varphi}$ is a state on $\tld{A}$ because $\tld{\varphi}(e) = \tld{\psi}(e) = 1$,
cf.\ \ref{eBbded}. Let $\varphi$ be the restriction of $\tld{\varphi}$ to $A$. Then
$\vphantom{\tld{A}}\varphi$ is a quasi-state extending $\psi$, and thus a
state\vphantom{$\tld{A}$}.
\end{proof}

\begin{theorem}\label{posequiv}%
Let $A$ be a Hermitian Banach \st-algebra, and let $a$ be a
Hermitian element of $A$. The following properties are equivalent.
\begin{itemize}
   \item[$(i)$] $a \geq 0$,
  \item[$(ii)$] $\pi \univ (a) \geq 0$,
 \item[$(iii)$] $\psi(a) \geq 0$\quad for all $\psi \in S(A)$.
\end{itemize}
\end{theorem}

\begin{proof}
The implications (i) $\Rightarrow$ (ii) $\Leftrightarrow$ (iii) hold for
normed \st-algebras in general. Indeed (i) $\Rightarrow$ (ii) follows
from \ref{hompos} while (ii) $\Leftrightarrow$ (iii) follows from
$\pi \univ (a) \geq 0 \Leftrightarrow \langle \pi\univ (a)x,x\rangle \geq 0$
for all $x \in H\univ$, cf.\ \ref{posop}. The closed subalgebra $B$ of $A$
generated by $a$ is a commutative Hermitian Banach \st-algebra
by \ref{Herminher}. One notes that every multiplicative linear functional
$\tau$ on $B$ is a state. Indeed, as $\tau$ is Hermitian \ref{mlfHerm},
we have $| \,\tau (b) \,|\,^2 = \tau (b^*b) $ for all $b \in B$,
so that $\tau$ is positive and has variation $1$. Furthermore, $\tau$ is
continuous by \ref{mlfbounded}. The proof of (iii) $\Rightarrow$ (i)
goes as follows. If $\psi(a) \geq 0$ for all $\psi$ in $S(A)$ then by the
preceding theorem also $\psi(a) \geq 0$ for all $\psi$ in $S(B)$. In
particular $\tau(a) \geq 0$ for all $\tau \in \Delta (B)$. It follows from
\ref{rangeGT} that $a \in B_+$, whence $a \in A_+$ by \ref{specsubalg}.
\pagebreak
\end{proof}

\begin{theorem}[Ra\u{\i}kov's Criterion]\label{Raikov}%
\index{Theorem!Ra\u{\i}kov}\index{Ra\u{\i}kov's Criterion}%
For a Banach \st-algebra $A$, the following statements are equivalent.
\begin{itemize}
   \item[$(i)$] $A$ is Hermitian,
  \item[$(ii)$] $\| \,\pi \univ (a) \,\| = \rsigma(a)$\quad for all $a \in A$,
 \item[$(iii)$] $\| \,a \,\|\0 = \rsigma(a)$\quad for all $a \in A$.
\end{itemize}
\end{theorem}

\begin{proof}
(ii) $\Leftrightarrow$ (iii): We know that $\| \,a \,\|\0 = \| \,\pi \univ (a) \,\|$,
cf.\ \ref{GNuniv}.
(ii) $\Rightarrow$ (i) follows from \ref{rsHerm} and \ref{C*rs}.
(i) $\Rightarrow$ (ii): It is easily seen that the Gelfand-Na\u{\i}mark
seminorm on $\tld{A}$ coincides on $A$ with $\| \cdot \|\0$,
cf.\ \ref{extquasi}. One can thus assume that $A$ is unital. Clearly
$\pi \univ$ then is unital \ref{unitalhom}. If $A$ is Hermitian,
then the preceding theorem and lemma \ref{posposrs} imply
that for $a \in A$, one has $\rsigma \bigl(\pi \univ(a)\bigr) = \rsigma(a)$,
or $\| \,\pi \univ (a) \,\| = \rsigma(a)$, cf.\ \ref{C*rs}.%
\end{proof}

\begin{theorem}[the Gelfand-Na\u{\i}mark Theorem]%
\index{Theorem!Gelfand-Na\u{\i}mark}\label{GelfandNaimark}%
\index{Gelfand-Naimark@Gelfand-Na\u{\i}mark!Theorem}%
If $(A, \| \cdot \|)$ is a C*-algebra, then
\[  \| \,\pi \univ (a) \,\| = \| \,a \,\|\0 = \| \,a \,\| \qquad \text{for all } a \in A, \]
so that the universal representation of $A$ establishes an isomorphism
of \linebreak C*-algebras from $A$ onto a C*-algebra of bounded linear
operators on a Hilbert space. Furthermore $C^*(A)$ can be identified with $A$.
\end{theorem}

\begin{proof}
This follows now from $\rsigma(a) = \| \,a \,\|$, cf.\ \ref{C*rs}.
\end{proof}

\begin{proposition}%
A Banach \st-algebra $A$ is Hermitian if and only if
the Pt\'{a}k function $\rsigma$ is a seminorm.
\end{proposition}

\begin{proof}
The ``only if'' part is clear by $\rsigma(a) = \| \,a \,\|\0$ if $A$ is Hermitian.
Conversely assume that $\rsigma$ is a seminorm. In order to show that
$A$ is Hermitian, we can go over to the closed subalgebra generated by
a Hermitian element, cf.\ \ref{specsubalg}. Hence we may assume that
$A$ is commutative. By assumption we have
$\rsigma(a+b) \leq \rsigma(a)+\rsigma(b)$, whence also
$| \,\rsigma(a) - \rsigma(b) \,| \leq \rsigma(a-b)$ for all $a,b \in A$. We also have
$\rsigma(a^*a) = {\rsigma(a) \,}^2$ and $\rsigma(a^*) = \rsigma(a)$ (by \ref{rsC*}),
as well as $\rsigma(ab) \leq \rsigma(a) \,\rsigma(b)$ (by \ref{commrs})
for all $a,b \in A$. It follows that the set $I$ on which $\rsigma$
vanishes is a \st-stable two-sided ideal in $A$. Furthermore, it follows that
$ \| \,a+I \,\| := \rsigma(a)$ $(a \in A)$ is well-defined and makes $A/I$ into a
pre-C*-algebra. Upon completing $A/I$, we get a \st-algebra homomorphism
$\pi$ from $A$ to a C*-algebra such that $\rsigma(a) = \| \,\pi(a) \,\|$ for all
$a \in A$. Now \ref{rsHerm} implies that $A$ is Hermitian. \pagebreak
\end{proof}

\begin{proposition}\label{weakPalmer}%
On a Hermitian Banach \st-algebra, the Pt\'{a}k function $\rsigma$ is
uniformly continuous. 
\end{proposition}

\begin{proof}
\ref{Raikov}, \ref{bdedstatessuffcond} (ii), \ref{bstGNc}.
\end{proof}

\medskip
As an immediate consequence, we obtain:

\begin{corollary}\label{Hermcont}%
If $A$ is a \st-subalgebra of a Hermitian Banach \linebreak \st-algebra,
then every $\sigma$-contractive linear map from $A$ to a normed space
is continuous. In particular, $A$ has continuous states.
\end{corollary}

\begin{proof}
\ref{continvol} (iv) $\Rightarrow$ (vi).
\end{proof}

\medskip
(See also \ref{Palmer} as well as the remark following
\ref{bdedstatessuffcond}.)

Now for the commutative case.

\begin{theorem}\label{maincommHerm}%
If $A$ is a commutative Hermitian Banach \st-algebra with
$\Delta (A) \neq \varnothing$, then for $a \in A$, we have
\[ \| \,a \,\|\0 = | \,\wht{a} \,|_{\infty}. \]
\end{theorem}

\begin{proof}
$\rsigma(a) = \rlambda(a) = | \,\wht{a} \,|_{\infty}$
by \ref{fundHerm} (iv) \& \ref{normGTrl}.
\end{proof}

\begin{corollary}%
If $A$ is a commutative Hermitian Banach \st-al\-gebra
with $\Delta (A) \neq \varnothing$, then $A$ is \st-semisimple if
and only if the Gelfand transformation $A \to C\0\bigl(\Delta(A)\bigr)$
is injective.
\end{corollary}

\begin{corollary}%
If $A$ is a commutative Hermitian Banach \st-al\-gebra
with $\Delta (A) \neq \varnothing$, then $C^*(A)$ is isomorphic
as a C*-algebra to $C\0\bigl(\Delta(A)\bigr)$.
\end{corollary}

\begin{proof}
By theorem \ref{maincommHerm},
$\bigl(A/$\st-$\mathrm{rad}(A),\| \cdot \|\bigr)$
is isomorphic as a pre-C*-algebra to
$\{ \,\wht{a} \in C\0\bigl(\Delta(A)\bigr) : a \in A \,\}$,
and the latter set is dense in $C\0\bigl(\Delta(A)\bigr)$,
cf.\ \ref{Hermdense}.
\end{proof}

\clearpage


\section{Commutativity}%
\label{Commutativity}

In this paragraph, let $A$ be a \underline{commutative} normed
\st-algebra.

\begin{definition}[$\Delta^*\protect\bsa(A)$]\label{mlfHbsa}%
\index{D(A)bsa@$\Delta^*\protect\bsa(A)$}%
We denote by $\Delta^*\bsa(A)$ the set of \linebreak Hermitian
multiplicative linear functionals on $A$, which are bounded on $A\sa$.
Please note that
\[ \Delta^*\bsa(A) = S(A) \cap \Delta (A). \]
\end{definition}

As in \ref{Gtopo} we obtain:

\begin{definition}%
We equip the set $\Delta^*\bsa(A)$ with the topology \linebreak
inherited from QS(A). It becomes in this way a locally compact
\linebreak Hausdorff space.
\end{definition}

\begin{proposition}%
If $A$ is a Hermitian commutative Banach \linebreak \st-algebra, then
$\Delta^*\bsa(A) = \Delta(A)$, also in the sense of topological spaces.
\end{proposition}

\begin{proposition}%
\index{a07@$\protect\wht{a}$}\label{commdense}%
If $\Delta^*\bsa(A) \neq \varnothing$, the Gelfand transformation
\begin{align*}
A & \to C\0\bigl(\Delta^*\bsa(A)\bigr) \\
a & \mapsto \wht{a}
\end{align*}
then is a $\sigma$-contractive \st-algebra homomorphism,
whose range is dense in $C\0\bigl(\Delta^*\bsa(A)\bigr)$.
(We shall abbreviate
$\wht{a}|_{\Delta^*_{\text{\small{bsa}}}(A)}$ to $\wht{a}$).
\end{proposition}

\begin{proof}
One applies the Stone-Weierstrass Theorem, cf.\ the appendix
\linebreak \ref{StW}.
\end{proof}

\begin{proposition}%
The homeomorphism $\varphi\0 \mapsto \varphi\0 \circ j$ \ref{affinehomeom}
from \linebreak $S\bigl(C^*(A)\bigr)$ onto $S(A)$ restricts to a
homeomorphism from $\Delta\bigl(C^*(A)\bigr)$ onto $\Delta^*\bsa(A)$.
\end{proposition}

\begin{corollary}%
If $\Delta^*\bsa(A) \neq \varnothing$, then $C^*(A)$ is isomorphic as a \linebreak
C*-algebra to $C\0\bigl(\Delta^*\bsa(A)\bigr)$. Furthermore, for $a \in A$, we have
\[ \| \,a \,\|\0 = \bigl| \,\wht{a}|_{\Delta^*_{\text{\small{bsa}}}(A)} \,\bigr|_{\infty} . \]
It follows that $A$ is \st-semisimple if and only if the Gelfand trans\-formation
$A \to C\0\bigl(\Delta^*\bsa(A)\bigr)$ is injective. \pagebreak
\end{corollary}

\begin{theorem}[the abstract Bochner Theorem]%
\label{Bochner}\index{Theorem!abstract!Bochner}%
\index{abstract!Bochner Theorem}\index{Bochner Thm., abstr.}%
The formula
\[ \psi (a) = \int \wht{a} \,d \mu \qquad (a \in A) \]
establishes an affine bijection from $M_1\bigl(\Delta^*\bsa(A)\bigr)$\,onto $S(A)$.
\end{theorem}

\begin{proof} Let $\psi$ be a state on $A$. Let $\psi \0$ be the
corresponding state on $C^*(A)$. By \ref{predisint}, there exists a
unique measure $\mu \0 \in M_1\Bigl(\Delta\bigl(C^*(A)\bigr)\Bigr)$ such that
\[ \psi \0(b) = \int \wht{b} \,d \mu \0
\quad \text{for all } b \in C^*(A). \]
Let $\mu$ be the image measure of $\mu \0$ under the homeomorphism
$\tau \0 \mapsto \tau \0 \circ j$ from $\Delta\bigl(C^*(A)\bigr)$ onto
$\Delta^*\bsa(A)$. (For the concept of image measure,
see the appendix \ref{imagebegin} - \ref{imageend}.) For $a \in A$, we get
\begin{align*}
\psi (a) & = \psi \0\bigl(j(a)\bigr) \\
& = \int \wht{j(a)} \,d \mu \0 \\
& = \int \wht{j(a)}(\tau \0) \,d \mu \0 (\tau \0) \\
& = \int \tau \0 \bigl(j(a)\bigr) \,d \mu \0 (\tau \0) \\
& = \int \wht{a} (\tau \0 \circ j) \,d \,\mu \0 (\tau \0) \\
& = \int \wht{a} (\tau ) \,d \mu (\tau ) = \int \wht{a} \,d \mu.
\end{align*}
Uniqueness of $\mu$ follows from the fact that the functions
$\wht{a}$ $(a \in A)$ are dense in $C\0\bigl(\Delta^*\bsa(A)\bigr)$,
cf.\ \ref{commdense}.

Conversely, let $\mu$ be a measure in $M_1\bigl(\Delta^*\bsa(A)\bigr)$
and consider the linear functional $\varphi$ on $A$ defined by
\[ \varphi (a) := \int \wht{a} \,d \mu \qquad (a \in A). \]
Then $\varphi$ is a positive linear functional on $A$ because for
$a \in A$, we have
\[ \varphi (a^*a) = \int {| \,\wht{a} \,| \,}^2 \,d \mu \geq 0. \pagebreak \]
The positive linear functional $\varphi$ is contractive on $A\sa$
because all the $\tau \in \Delta^*\bsa(A)$ are contractive on $A\sa$
as they are states. Indeed
\[ | \,\varphi (b) \,| \leq \int | \,\tau(b) \,| \,d\mu (\tau) \leq | \,b \,| 
\quad \text{for all } b \in A\sa. \]
Furthermore, we have $v (\varphi ) \leq 1$ by the Cauchy-Schwarz inequality:
\[ {| \,\varphi (a) \,| \,}^2 = {\biggl| \,\int \wht{a} \,d\mu \,\biggr| \,}^2
\leq \int {| \,\wht{a} \,| \,}^2 \,d\mu \cdot \int {1 \,}^2 \,d\mu = \varphi (a^*a)
\quad \text{for all } a \in A . \]
Thus, it follows that $\varphi = \lambda \psi$ with $\psi$ a state on $A$
and $0 \leq \lambda \leq 1$. By the preceding, there exists a unique measure
$\nu \in M_1\bigl(\Delta^*\bsa(A)\bigr)$ such that
\[ \psi (a) = \int \wht{a} \,d \nu\quad (a \in A). \]
We obtain
\[ \varphi (a) = \int \wht{a} \,d(\lambda \nu ). \]
as well as
\[ \varphi (a) = \int \wht{a} \,d \mu. \]
By density of the functions $\wht{a}$ in
$C\0\bigl(\Delta^*\bsa(A)\bigr)$, it follows
$\mu = \lambda \nu$, whence $\lambda = 1$,
so that $\varphi = \psi$ is a state on $A$.
\end{proof}

\medskip
See also \ref{Bochnerbis} \& \ref{Bremark} below.

\begin{counterexample}\label{counterex}%
There exists a commutative unital nor\-med \st-algebra $A$ such that
each Hermitian multiplicative linear \linebreak functional on $A$ is
contractive on $A\sa$, yet discontinuous, and thereby a discontinuous state.
\end{counterexample}

\begin{proof} We take $A := \mathds{C}[\mathds{Z}]$ as
\st-algebra, cf.\ \ref{CG}. (We shall define a norm later on.)
Let $\tau$ be a Hermitian multiplicative linear functional on $A$.
If we put $u := \tau(\delta _1)$, then
\[ \tau(a) = \sum _{n \in \mathds{Z}} \,a(n) \,u^{\,n}
\quad \text{for all } a \in A. \]
With $e := \delta_0$ the unit in $A$, we have
\begin{align*}
{| \,u \,| \,}^2 & = u \,\overline{u} = \tau(\delta_1) \,\overline{\tau(\delta _1)}
= \tau(\delta_1) \,\tau({\delta _1}^*) \\
 & = \tau(\delta_1) \,\tau(\delta _{-1}) = \tau(\delta_1 \delta_{-1}) = \tau (e) = 1,
\pagebreak
\end{align*}
which makes that $u = \tau(\delta _1)$ belongs to the unit circle.
Conversely, if $u$ is in the unit circle then
\[ \tau(a) := \sum _{n \in \mathds{Z}} \,a(n) \,u^{\,n} \quad (a \in A)  \]
defines a Hermitian multiplicative linear functional on $A$.
In other words, the Hermitian multiplicative linear functionals $\tau$
on $A$ cor\-respond to the points $u := \tau (\delta _1)$ in the unit circle.

We shall now introduce a norm on $A$. Let $\gamma > 1$ be fixed. Then
\[ | \,a \,| := \sum _{n \in \mathds{Z}} \,| \,a(n) \,| \,{\gamma \,}^n \]
defines an algebra norm on $A$. (We leave the verification to the reader.)
If $a \in A\sa$, then $a(-n) = \overline{a(n)}$
for all $n \in \mathds{Z}$, so for $a \in A\sa$ we have
\[ | \,a \,| = | \,a(0) \,| + \sum _{n \geq 1} \,| \,a(n) \,|
\,\Bigl( {\gamma \,}^n + \frac{1}{{\gamma \,}^n} \Bigr). \]
Let $\tau$ be a Hermitian multiplicative linear functional on $A$.
Then $\tau(\delta _n) = \bigl(\tau(\delta _1)\bigr)^n$ is in the unit circle
for all $n \in \mathds{Z}$, so for $a \in A\sa$ we find
\[ | \,\tau(a) \,| \leq | \,a(0) \,| + 2 \,\sum _{n \geq 1} \,| \,a(n) \,| . \]
We see that if $a \in A\sa$, then $| \,\tau(a) \,| \leq | \,a \,|$, so $\tau$
is contractive on $A\sa$. However
\begin{alignat*}{2}
 & | \,\delta _n \,|  = {\gamma \,}^n \to 0 & &
\quad \text{for } n \to -\infty, \\
 & | \,\tau(\delta _n) \,|  = 1 & &
\quad \text{for all } n \in \mathds{Z},
\end{alignat*}
which shows that $\tau$ is discontinuous.
\end{proof}

\begin{remark}\label{counterexbis}
This gives an example of a commutative  normed \linebreak \st-algebra which
cannot be imbedded in a Banach  \st-algebra, see \ref{bdedstatessuffcond}
(iii). See also \ref{Brlsymm}, \ref{Bnormalrsleqrl}, and \ref{Palmer}. The above
normed \st-algebra also serves as a counterexample for the latter
statements with normed \st-algebras instead of Banach \st-algebras. Please
note that the above commutative unital normed \st-algebra does
not carry any continuous Hermitian multiplicative linear functional.
\end{remark}

\clearpage


\chapter{Irreducible Representations}

\setcounter{section}{30}


\section{General \texorpdfstring{$*$-}{\052\055}Algebras: Indecomposable Functionals}

Let throughout $A$ be a \st-algebra and $H \neq \{0\}$ a Hilbert space.

\begin{definition}[irreducible representations]%
\index{representation!irreducible}\index{irreducible}%
A representation of $A$ in $H$ is called \underline{irreducible}
if it is non-zero and if it has no closed reducing subspaces other
than $\{0\}$ and $H$.
\end{definition}

\begin{proposition}%
A non-zero representation $\pi$ of $A$ in $H$ is irreducible
if and only if every non-zero vector in $H$ is cyclic for $\pi$.
\end{proposition}

\begin{proof} Assume that $\pi$ is irreducible. Then $\pi$ must be
non-degenerate, cf.\ \ref{nondegenerate}. Let $x$ be a non-zero vector
in $H$ and consider $M = \overline{R(\pi)x}$. Then $x \in M$, because
$\pi$ is non-degenerate, cf.\ \ref{cyclicsubspace}. In particular, $M$ is a
closed reducing subspace $\neq \{0\}$, so that $M = H$. The proof of the
converse implication is left to the reader.
\end{proof}

\begin{lemma}[Schur's Lemma]\index{Theorem!Schur's Lemma}%
\index{Schur's Lemma}\index{Lemma!Schur's}\label{Schur}%
A non-zero representation $\pi$ of $A$ in $H$
is irreducible if and only if $\pi' =\mathds{C}\mathds{1}$.
\end{lemma}

\begin{proof} Assume that $\pi' = \mathds{C}\mathds{1}$. Let $M$
be a closed reducing subspace of $\pi$ and let $p$ be the projection on
$M$. Then $p$ belongs to $\pi'$ by \ref{invarcommutant}, whence either
$p = 0$ or $p = \mathds{1}$. Assume conversely that $\pi$ is irreducible.
Then $\pi'$ is a unital C*-subalgebra of $B(H)$. Consider a Hermitian
element $a = a^*$ of $\pi'$ and let $C$ be the C*-subalgebra of $\pi'$
generated by $a$ and $\mathds{1}$. In order to prove the lemma, it suffices
to show that $\Delta (C)$ contains only one point. So assume that $\Delta (C)$
contains more than one point. There then exist non-zero $c_1, c_2 \in C$ with
$c_1 c_2 = 0$ (by Urysohn's Lemma). Choosing $x \in H$ with $c_2 x \neq 0$,
we obtain the contradiction $c_1 H = c_1 \overline{\pi (A)c_2x} \subset
\overline{\pi (A)c_1 c_2 x} = \{0\}$. \pagebreak
\end{proof}

\begin{corollary}\label{commirroned}%
If $A$ is commutative, then a non-zero representation $\pi$ of
$A$ in $H$ is irreducible if and only if $H$ is one-dimensional.
\end{corollary}

\begin{proof} We have $R(\pi) \subset \pi'$. \end{proof}

\medskip
An irreducible representation of a non-commutative \st-algebra
need not be one-dimensional, cf.\ e.g.\ \ref{BHirred} below.

\begin{proposition}\label{spherical}%
Let $\pi$ be an irreducible representation of $A$ in $H$.
Assume that $x, y$ are vectors in $H$ with
\[ \langle \pi (a) x, x \rangle = \langle \pi (a) y, y \rangle
\quad \text{for all } a \in A. \]
There then exists a complex number $u$ of modulus $1$ with $x = uy$.
\end{proposition}

\begin{proof} Assume that both $x$ and $y$ are
non-zero. Then $x,y$ are cyclic vectors for $\pi$ and there exists a
unitary operator $U \in \pi' = C(\pi,\pi)$ taking $x$ to $y$, cf.\ \ref{coeffequal}.
Since $\pi$ is irreducible, it follows from Schur's Lemma \ref{Schur} that
$U = u\mathds{1}$ with $u \in \mathds{C}, | \,u \,| = 1$. Assume next that
$x = 0$, $y \neq 0$. The equation
\[ 0=\langle \pi (a)x,x \rangle =\langle\pi (a)y,y \rangle \qquad (a \in A) \]
would then imply that the vector $y \neq 0$ was orthogonal to all the vectors
$\pi (a)y$ $(a \in A)$, so that $y$ could not be a cyclic vector.
\end{proof}

\begin{lemma}\label{posoprep}%
Let $\varphi$ be a positive semidefinite sesquilinear form on
$H$ with $\sup _{\,\| \,x \,\| \leq 1} \varphi (x, x) < \infty$. There
then exists a unique linear operator $a \in B(H)_+$ such that
\[ \varphi (x, y) = \langle ax, y \rangle\quad (x,y \in H). \]
\end{lemma}

\begin{proof} We have $\sup _{\,\| \,x \,\|, \| \,y \,\| \leq 1} | \,\varphi (x,y) \,| < \infty$
by the Cauchy-Schwarz \linebreak inequality. It follows from a well-known
theorem that there exists a unique bounded linear operator $a$ on $H$
with $\langle ax,y \rangle = \varphi (x,y)$ for all $x,y \in H$. The operator
$a$ is positive because $\varphi$ is positive semidefinite, cf.\ \ref{posop}.
\end{proof}

\begin{definition}[subordination]%
\index{subordinate}\index{functional!subordinate}%
Let $\varphi_1, \varphi_2$ be positive linear functionals on $A$.
One says that $\varphi_1$ is \underline{subordinate} to $\varphi_2$ if
there is $\lambda \geq 0$ such that $\lambda \varphi_2 - \varphi_1$ is
a positive linear functional. \pagebreak
\end{definition}

\begin{theorem}\label{subordrep}%
Let $\pi$ be a cyclic representation of $A$ in $H$ and let $c$ be a cyclic
vector for $\pi$. Consider the positive linear functional $\varphi$ on $A$
defined by
\[ \varphi (a) := \langle \pi (a) c, c \rangle\quad (a \in A). \]
The equation
\[ \varphi_1 (a) = \langle b \pi(a) c, c \rangle\quad (a \in A) \]
establishes a bijection between operators $b \in {\pi'}_+$
and Hermitian positive linear functionals $\varphi_1$ on $A$
of finite variation which are subordinate to $\varphi$. 
\end{theorem}

\begin{proof}
Let $b \in {\pi'}_+$. It is immediate that
\[ \varphi_1 (a) := \langle b \pi (a)c,c \rangle
= \langle \pi(a) b^{1/2} c, b^{1/2} c \rangle \quad (a \in A) \]
defines a Hermitian positive linear functional of finite variation on $A$,
cf.\ \ref{C*sqrootrestate} and \ref{variationinequal}. To show that $\varphi_1$
is subordinate to $\varphi$, we note that $b \leq \| \,b \,\| \,\mathds{1}$, whence
\[ \| \,b \,\| \,\langle x,x \rangle - \langle bx,x \rangle \geq 0 \]
for all $x \in H$. With $x := \pi (a) c$ it follows that
\[ \| \,b \,\| \,\langle \pi (a^*a)c,c \rangle
- \langle b \pi (a^*a)c,c \rangle \geq 0, \]
which says that $\| \,b \,\| \,\varphi - \varphi_1$ is a positive linear
functional on $A$.

Conversely, let $\varphi_1$ be a Hermitian positive linear functional
of finite variation subordinate to $\varphi$. There then exists
$\lambda \geq 0$ such that
\[ 0 \leq \varphi_1 (a^*a) \leq \lambda \,\varphi (a^*a)  \tag*{(\st)} \]
for all $a \in A$.

For $x = \pi (f)c, y = \pi (g)c$, $(f,g \in A)$, we put
\[ \alpha (x,y) := \varphi_1 (g^*f). \]
This expression depends only on $x,y$. Indeed, let for example also
$x = \pi (h)c$, i.e.\ $\pi (f-h)c = 0$. We obtain
\[ \varphi \,\bigl((f-h)^*(f-h)\bigr) = \langle \pi (f-h)c, \pi (f-h)c \rangle = 0, \]
so that by (\st),
\[ \varphi_1 \bigl((f-h)^*(f-h)\bigr) = 0. \]
From the Cauchy-Schwarz inequality, it follows that
\[ \varphi_1 \bigl(g^*(f-h)\bigr) = 0,  \]
so that
\[ \varphi_1 (g^*f) = \varphi_1 (g^*h), \pagebreak \]
as claimed. Denote by $H\0$ the set of all vectors $x \in H$ of the form
$x = \pi (f)c$ $(f \in A)$. The vector space $H\0$ is dense in $H$ because
$c$ is cyclic for $\pi$. The positive semidefinite sesquilinear form $\alpha$
is bounded on $H \0$ since (\st) may be rewritten as
\[ \alpha (x, x) \leq \lambda \,\langle x, x \rangle \tag*{(\st\st)} \]
for all $x \in H\0$. Therefore $\alpha$ has a unique continuation to a bounded
positive semidefinite sesquilinear form on $H$. The inequality (\st\st) then
is valid for all $x \in H$. By \ref{posoprep} there is a unique operator
$b \in B(H)_+$ such that
\[ \alpha (x,y) = \langle bx,y \rangle\quad (x,y \in H). \]
It shall be shown that $b \in \pi'$. For $x = \pi (f)c, y = \pi (g)c$, and $a \in A$,
we have
\[ \langle b \pi (a)x,y \rangle = \langle b \pi (af)c, \pi (g)c \rangle
= \alpha \bigl(\pi(af)c, \pi(g)c\bigr) = \varphi_1 (g^*af), \]
as well as
\begin{align*}
\langle \pi (a)bx,y \rangle & = \langle bx,\pi (a^*)y \rangle
= \langle b \pi (f)c,\pi (a^*g)c \rangle = \alpha \bigl(\pi(f)c, \pi(a^*g)c\bigr) \\
 & = \varphi_1 \bigl((a^*g)^*f\bigr) = \varphi_1 (g^*af),
\end{align*}
so that
\[ \langle b \pi (a)x,y \rangle = \langle \pi (a)bx,y \rangle \]
for all $x,y \in H\0$ and hence for all $x,y \in H$, so that $b \in \pi'$.
In order to prove that
\[ \varphi_1 (a) = \langle b \pi (a)c,c \rangle\quad (a \in A), \]
one first notices that
\[ {| \,\varphi_1 (a) \,| \,}^2 \leq v(\varphi_1) \,\varphi_1(a^*a)
\leq v(\varphi_1) \,\lambda \,\varphi (a^*a)
= \lambda \,v(\varphi_1) \,{\| \,\pi (a)c \,\| \,}^2, \]
so that $\pi (a)c \mapsto \varphi_1 (a)$ is a bounded linear functional on
$H\0$, hence extends to a bounded linear functional on $H$. (The expression
$\varphi_1 (a)$ depends only on $\pi (a)c$ as is seen as follows: if
$\pi (f)c = \pi (g)c$, we have
\[ {| \,\varphi_1 (f) - \varphi_1 (g) \,| \,}^2
= {| \,\varphi_1 (f-g) \,| \,}^2 \leq \lambda \,v(\varphi_1 ) \,{\| \,\pi (f-g)c \,\| \,}^2
= 0.) \]
There thus exists a vector $d \in H$ with
\[ \varphi_1 (a) = \langle \pi (a)c,d \rangle\quad (a \in A).  \]
We get
\begin{align*}
\langle \pi (f)c,b \pi (g)c \rangle & = \langle b \pi (f)c, \pi (g)c \rangle
= \alpha (\pi (f)c,\pi (g)c) \\
 & = \varphi_1 (g^*f) = \langle \pi (g^*f)c,d \rangle
= \langle \pi (f)c, \pi (g) d \rangle,
\end{align*}
so that
\[ b \pi (g)c = \pi (g)d\quad (g \in A). \pagebreak \]
It follows that
\[ \langle b \pi (a)c,c \rangle = \langle \pi (a)d,c \rangle
= \overline{\langle \pi (a^*)c,d \rangle}
= \overline{\varphi_1(a^*)} = \varphi_1(a). \]
The uniqueness statement follows from the identity
\[ \varphi_1 (g^*f) = \langle b \pi (f)c,\pi (g)c \rangle, \]
and the fact that $H\0$ is dense in $H$. 
\end{proof}

\begin{definition}[indecomposable positive linear functionals]%
\index{indecomposable}\index{functional!indecomposable}%
Let $\varphi$ be a Hermitian positive linear functional
of finite variation on $A$. The functional $\varphi$ is called
\underline{indecomposable} if every other Hermitian positive
linear functional of finite variation on $A$, which is
subordinate to $\varphi$, is a multiple of $\varphi$.
\end{definition}

\begin{theorem}\label{indecequiv}%
Let $\pi$ be a cyclic representation of $A$ in $H$ and let $c$ be a cyclic
vector for $\pi$. Consider the positive linear functional $\varphi$ on $A$
defined by
\[ \varphi (a) := \langle \pi (a)c,c \rangle\quad (a \in A). \]
Then $\pi$ is irreducible if and only if $\varphi$ is indecomposable.
\end{theorem}

\begin{proof} Assume that $\pi$ is irreducible and let $\varphi_1$ be
a Hermitian positive linear functional of finite variation subordinate to
$\varphi$. There exists $b \in {\pi'}_+$ with
$\varphi_1 (a) = \langle b \pi (a) c,c \rangle$ for all $a \in A$. Since
$\pi$ is irreducible, it follows that $b = \lambda \mathds{1}$ for some
$\lambda \geq 0$, so that $\varphi_1 = \lambda \varphi$. Conversely,
let $\varphi$ be indecomposable and let $p$ denote the projection
on a closed invariant subspace of $\pi$. We then have $p \in {\pi'}_+$,
cf.\ \ref{invarcommutant}, so that the Hermitian positive linear functional
$\varphi_1$ of finite variation defined by
$\varphi_1(a) := \langle p \pi (a)c,c \rangle$ $(a \in A)$, is subordinate
to $\varphi$, whence $\varphi_1 = \lambda \varphi$ for some
$\lambda \in \mathds{R}$. Since $p$ is uniquely determined by
$\varphi_1$, it follows that $p = \lambda \mathds{1}$, but then
$\lambda = 0$ or $1$ as $p$ is a projection, and so $\pi$ is irreducible.
\end{proof}

\clearpage


\section{Normed \texorpdfstring{$*$-}{\052\055}Algebras: Pure States}

Throughout this paragraph, let $A$ be a normed \st-algebra,
and \linebreak assume that $\psi$ is a state on $A$.

\begin{definition}[pure states, $PS(A)$]%
\index{pure state}\index{state!pure}%
\index{P9@$PS(A)$}\index{functional!state!pure}%
One says that $\psi$ is a \linebreak \underline{pure state} on $A$
if it is an extreme point of the convex set $S(A)$ of states on $A$,
cf.\ \ref{SAconvex}. The set of pure states on $A$ is denoted by
$\underline{PS(A)}$.
\end{definition}

\begin{lemma}\label{subordo}%
Let $\psi \0$ be the state on $C^*(A)$ corresponding to $\psi$.
Then $\psi$ is indecomposable precisely when $\psi \0$ is so.
\end{lemma}

\begin{proof}
The representation $\pi _{(\psi \0)}$ is spatially equivalent to
$(\pi _{\psi}) \0$, cf.\ \ref{oequiv}. Therefore $\pi _{(\psi \0)}$ is
irreducible precisely when $(\pi _{\psi}) \0$ is. Now $(\pi _{\psi}) \0$
is irreducible if and only if $\pi _{\psi}$ is irreducible, as is shown
by \ref{piopi} (iii).
\end{proof}

\begin{theorem}\label{indecopure}%
The state $\psi$ is indecomposable if and only if it is a pure state.
\end{theorem}

\begin{proof} We shall first prove this if $A$ has continuous
involution and a bounded left approximate unit. Let the state
$\psi$ be indecomposable. Assume that
$\psi = \lambda _1 \psi _1 + \lambda _2 \psi _2$ where
$\lambda _1, \lambda _2 > 0, \lambda _1 + \lambda _2 = 1$ and
$\psi _1, \psi _2 \in S(A)$. Since $\lambda _k \psi _k$ is subordinate to
$\psi$, we have $\lambda _k \psi _k = \mu _k \psi$ for some $\mu _k \geq 0$.
We obtain $\lambda _k = v(\lambda _k \psi _k) = v(\mu _k \psi ) = \mu _k$
so that $\psi _k = \psi$ for $k = 1,2$. Therefore $\psi$ is an extreme point
of $S(A)$. Let next the state $\psi$ be pure. Consider a Hermitian
positive linear functional $\psi _1$ of finite variation on $A$ subordinate
to $\psi$. Then $\psi _1$ is bounded on $A\sa$ by applying
\ref{subordrep} to $\pi_{\psi}$. But then $\psi _1$ and $\psi$ are bounded
as the involution in $A$ is continuous by assumption.
It must be shown that $\psi _1$ is a multiple of $\psi$. Assume without loss
of generality that $\psi _1 \in S(A)$. There then exists $\lambda > 0$ such
that $\psi - \lambda \psi _1 =: \psi'$ is a positive linear functional. This
functional is continuous and thus has finite variation because by assumption
$A$ has continuous involution and a bounded left approximate unit,
cf.\ \ref{bdedfinvar}. It must be shown that $\psi _1$ is a multiple of $\psi$.
Assume that $\psi' \neq 0$. Put $\psi _2 := \psi'/\mu$ where $\mu := v(\psi') > 0$
as $\psi' \neq 0$. Then $\psi _2$ is a state on $A$ with
$\psi = \lambda \psi _1 + \mu \psi _2$. It follows that $1 = \lambda + \mu$
by additivity of the variation, cf.\ \ref{finvarconvex}. Since $\lambda, \mu > 0$,
it follows that $\psi _1 = \psi _2 = \psi$ because $\psi$ is an extreme point
of $S(A)$, and in particular $\psi _1$ is a multiple of  $\psi$.
This proves the theorem in presence of a continuous involution and a
bounded left approximate unit. We shall pull down the general case from
$C^*(A)$. Since the bijection $\psi \mapsto \psi \0$ is affine, $\psi$ is a pure
state precisely when $\psi \0$ is a pure state. Moreover, by \ref{subordo},
$\psi$ is indecomposable if and only if $\psi \0$ is indecomposable.
\end{proof}

\begin{theorem}\label{irredindecpure}%
The following conditions are equivalent.
\begin{itemize}
   \item[$(i)$] $\pi _{\psi}$ is irreducible,
  \item[$(ii)$] $\psi$ is indecomposable,
 \item[$(iii)$] $\psi$ is a pure state.
\end{itemize}
\end{theorem}

\begin{proof}
\ref{indecequiv} and \ref{indecopure}.
\end{proof}

\begin{theorem}\label{cPSD}%
If $A$ is \underline{commutative}, then $PS(A) = \Delta^*\bsa(A)$.
\end{theorem}

\begin{proof} \ref{commirroned} and \ref{mlfHbsa}. \end{proof}

\begin{proposition}\label{homPS}%
The affine homeomorphism $\psi \mapsto \psi \0$ \ref{affinehomeom}
from $S(A)$ onto $S\bigl(C^*(A)\bigr)$ restricts to a homeomorphism
from $PS(A)$ onto $PS\bigl(C^*(A)\bigr)$.
\end{proposition}

\begin{proof}
\ref{subordo} and \ref{indecopure}. (Recall \ref{affinehomeom}.)
\end{proof}

\begin{theorem}%
The set of pure states on $A$ parametrises the class of irreducible
$\sigma$-contractive representations of $A$ in Hilbert spaces up to
spatial equivalence. More precisely, if $\psi$ is a pure state on $A$,
then $\pi _{\psi}$ is an irreducible $\sigma$-contractive representation
of $A$. Conversely, if $\pi$ is an irreducible $\sigma$-contractive
representation of $A$ in a Hilbert space, there exists a pure state $\psi$
on $A$ such that $\pi$ is spatially equivalent to $\pi _{\psi}$.
(Cf.\ \ref{stateratio}.)
\end{theorem}

Remark however that \ref{spherical} implies that if $\pi$ is a
multi-dimensional irreducible $\sigma$-contractive representation of a
(necessarily non-commut\-ative) normed \st-algebra, there exists a
continuum of pure states $\varphi$ such that $\pi$ is spatially
equivalent to $\pi _{\varphi}$.

\begin{proposition}%
The set of extreme points of the non-empty compact convex set $QS(A)$
is $PS(A) \cup \{0\}$. Cf.\ \ref{SAconvex} \& \ref{topQS}. \pagebreak
\end{proposition}

\begin{proof}
It shall first be shown that $0$ is an extreme point of $QS(A)$. So
suppose that $0 = \lambda \,\varphi _1 +(1-\lambda) \,\varphi _2$
with $0 < \lambda < 1$ and $\varphi _1, \varphi _2 \in QS(A)$. For
$a \in A$, we then have $\varphi _k (a^*a) = 0$ by $\varphi _k (a^*a) \geq 0$,
and so $\varphi _k (a) = 0$ by
$| \,\varphi _k (a) \,|^2 \leq v(\varphi _k) \,\varphi _k (a^*a)$.
It shall next be shown that each $\psi \in PS(A)$ is an
extreme point of $QS(A)$. So let $\psi \in PS(A)$,
$\psi = \lambda \,\varphi _1 +(1-\lambda ) \,\varphi _2$, where
$0 < \lambda < 1$ and $\varphi _1$, $\varphi _2 \in QS(A)$. By
additivity of the variation \ref{finvarconvex}, we have
$1 = \lambda \,v(\varphi _1)+( 1 - \lambda ) \,v(\varphi _2)$,
whence $v(\varphi _1) = v(\varphi _2) = 1$, so that $\varphi _1$
and $\varphi _2$ are states, from which it follows that
$\psi = \varphi _1 = \varphi _2$. It remains to be shown that
$\varphi \in QS(A)$ with $0 < v(\varphi ) < 1$ is not an extreme point of
$QS(A)$. This is so because then
$\varphi = v(\varphi) \cdot \bigl(v(\varphi)^{-1}
\varphi \bigr)+\bigl(1-v(\varphi)\bigr) \cdot 0$. 
\end{proof}

\begin{corollary}%
The non-empty compact convex set $QS(A)$ is the  closed convex hull of
$PS(A) \cup \{0\}$.
\end{corollary}

\begin{proof}
The Kre\u{\i}n-Milman Theorem says that a non-empty compact convex set
in a separated locally convex topological vector space is the closed convex
hull of its extreme points.
\end{proof}

\begin{corollary}\label{existPS}%
If a normed \st-algebra has a state, then it has a pure state.
\end{corollary}

\begin{theorem}%
The state $\psi$ is in the closed convex hull of $PS(A)$.
\end{theorem}

\begin{proof} The state $\psi$ is in the closed convex hull of $PS(A) \cup \{0\}$,
so that $\psi$ is the limit of a net $(\varphi _i)_{i \in I}$ where each
$\varphi _i$ $(i \in I)$ is of the form
\[ \lambda \cdot 0 + \sum _{j=1}^{n} \,\lambda _j \cdot \psi _j \]
with
\[ \lambda + \sum _{j=1}^{n} \,\lambda _j =1,
\ \lambda \geq 0, \ \lambda _j \geq 0,
\ \psi _j \in PS(A)\ (j = 1, \ldots ,n). \]
We then have
\[ v(\varphi _i) = \sum _{j=1}^{n} \,\lambda _j \leq 1, \]
so $\limsup v(\varphi _i) \leq 1$. We shall show that
$\liminf v(\varphi _i) \geq 1$, from which $\lim v(\varphi _i) = 1$.
Assume that $\alpha := \liminf v(\varphi _i) < 1$. 
By going over to a subnet, we can then assume that
$\alpha = \lim v(\varphi _i)$. For $a \in A$, one would have
\[ {| \,\psi (a) \,| \,}^2 = \lim {| \,\varphi _i (a) \,| \,}^2
\leq \lim v(\varphi _i) \,\varphi _i(a^*a) = \alpha \,\psi (a^*a), \]
which would imply $v(\psi ) \leq \alpha < 1$, a contradiction. We have
thus shown that $\lim v(\varphi _i) = 1$. In particular, by going over to
a subnet, we may assume that $v(\varphi _i) \neq 0$ for all $i \in I$.
The net $(v(\varphi _i)^{-1} \varphi _i) _{i \in I}$ then converges to
$\psi$, and consists of convex combinations of pure states.
\end{proof}

\begin{definition}[the spectrum $\wht{A}$]%
\index{spectrum!of an algebra}\index{AA@$\wht{A}$}%
We shall say that two states on $A$ are spatially equivalent if their
corresponding GNS representations are spatially equivalent. This
defines an equivalence relation on $S(A)$. We define
\underline{the spectrum $\wht{A}$} of $A$ as the set of spatial
equivalence classes of $PS(A)$.
\end{definition}

Please note that if $A$ is commutative, then $\wht{A}$ may be
identified with $\Delta^*\bsa(A)$, cf.\ \ref{cPSD}. 

\begin{definition}%
We shall say that
$(\psi _{\text{\small{$\lambda$}}}) _{\text{\small{$\lambda \in \wht{A}$}}}$
is a choice function if $\psi _{\text{\small{$\lambda$}}} \in \lambda$ for each
$\lambda \in \wht{A}$.
There exists$\vphantom{ _{\text{\small{$\lambda\in \wht{A}$}}}}$
a choice function
$(\psi _{\text{\small{$\lambda$}}}) _{\text{\small{$\lambda \in \wht{A}$}}}$
by the axiom of choice$\vphantom{\tld{A}}$ (cf.\ \ref{existPS}).
\end{definition}

\begin{definition}[the reduced atomic representation, $\pi \ra$]%
\index{p5@$\pi_{ra}$}\index{H3@$H_{ra}$}%
\index{representation!reduced atomic}%
If \linebreak
$(\psi _{\text{\small{$\lambda$}}}) _{\text{\small{$\lambda \in \wht{A}$}}}$
is a choice function, we define
\[ \pi \ra := \oplus \,_{\text{\small{$\lambda \in \wht{A}$}}}
\ \pi _{\text{\small{$\psi _{\text{\small{$\lambda$}}}$}}}. \]
One says that $\pi \ra$ is a reduced atomic representation of $A$. Any two
reduced atomic representations of $A$ are spatially equivalent, so that
one speaks of ``the'' \underline{reduced atomic representation} of $A$.
\end{definition}

Our next aim is theorem \ref{ranorm}.

\begin{theorem}%
If $A$ is a Hermitian Banach \st-algebra and if $B$ be a closed \st-subalgebra
of $A$, then every pure state on $B$ can be extended to a pure state on $A$.
\pagebreak
\end{theorem}

\begin{proof}
Let $\varphi$ be a pure state on $B$. Let $K$ be the set of states on $A$ which
extend $\varphi$, cf.\ \ref{extstate}. Then $K$ is a non-empty compact convex
subset of $S(A)$. (Indeed every adherence point of $K$ is a quasi-state extending
$\varphi$, hence a state too.) It follows that $K$ has an extreme point, $\psi$,
say. It shall be shown that  $\psi$ is a pure state. So let
$\psi = \lambda \psi _1+(1-\lambda )\psi _2$ with
$0 < \lambda < 1$ and $\psi _1, \psi _2 \in S(A)$. We then also have
$\psi |_B = \lambda \psi _1|_B + (1-\lambda ) \psi _2|_B$, so that
$\psi _1|_B = \psi _2|_B = \psi |_B = \varphi$. It follows that
$\psi _1, \psi _2$ belong to $K$. As $\psi$ is an extreme point of $K$, one has
$\psi _1 = \psi _2$, and so $\psi$ is a pure state.
\end{proof}

\begin{theorem}%
If $A$ is a Hermitian Banach \st-algebra and if $a$ is a Hermitian element
of $A$, then the following properties are equivalent.
\begin{itemize}
   \item[$(i)$] $a \geq 0$,
  \item[$(ii)$] $\pi \ra (a) \geq 0$,
 \item[$(iii)$] $\psi (a) \geq 0$\quad for all $\psi \in PS(A)$.
\end{itemize}
\end{theorem}

\begin{proof}
This follows in the same way as \ref{posequiv}.
\end{proof}

\begin{theorem}\label{ranorm}%
For all $a \in A$ we have
\[ \| \,\pi \ra (a) \,\| = \| \,a \,\|\0. \]
\end{theorem}

\begin{proof} Assume first that $A$ is a Hermitian Banach \st-algebra.
As in \ref{Raikov}, it follows that $\| \,\pi \ra (a) \,\| = \| \,a \,\|\0 = \rsigma(a)$.
If $A$ is merely a normed \st-algebra, we pull down the result
from the enveloping \linebreak C*-algebra. For $a \in A$, we have
\[ \| \,a \,\|\0 = \| \,j(a) \,\| = \| \,j(a) \,\|\0 = \| \,\pi \ra \bigl(j(a)\bigr) \,\| = \| \,\pi \ra (a) \,\|, \]
where we have used the Gelfand-Na\u{\i}mark Theorem \ref{GelfandNaimark}
as well as \ref{subordo}.
\end{proof}

\begin{theorem}%
The following statements are equivalent.
\begin{itemize}
   \item[$(i)$] $A$ is \st-semisimple,
  \item[$(ii)$] $\pi\ra$ is faithful,
 \item[$(iii)$] the irreducible $\sigma$-contractive representations of $A$ in Hilbert
                       spaces separate the points of $A$.
\end{itemize}
\end{theorem}

\begin{proof} The proof is left to the reader. \pagebreak \end{proof}

\clearpage


\section{The Left Spectrum in a Hermitian Banach \texorpdfstring{$*$-}{\052\055}Algebra}%
\label{leftspectrum}

In this paragraph, let $A$ be a \underline{Hermitian} Banach \st-algebra.

\begin{proposition}\label{leftidstate}%
If $A$ is unital, then for a proper left ideal $I$ in $A$, there exists a state
$\varphi$ on $A$ such that
\[ \varphi(a^*a) = 0 \qquad \text{for all } a \in I. \]
\end{proposition}

\begin{proof}
Let $M\0$ denote the real subspace of $A\sa$ spanned by $e$ and the
Hermitian elements of $I$. Define a linear functional $f\0$ on $M\0$ by
\[ f\0 (\lambda e + y) := \lambda \qquad (\lambda \in \mathds{R},\ y^* = y \in I). \]
Please note here that $e \notin I$, cf.\ \ref{notinvideal} (i). Let $C$ denote the
convex cone $A_+$ in $A\sa$, cf.\ \ref{plusconvexcone}. We have $M\0 + C = A\sa$
because for $a \in A\sa$ one has $\rlambda(a)e+a \geq 0$. Furthermore $f\0(x) \geq 0$
for all $x \in M\0 \cap C$. Indeed, if $x = \lambda e + y \in M\0 \cap C$ with $y^* = y \in I$,
then $\lambda e - x = -y$ is a Hermitian element of $I$, which is a proper left ideal,
so that $\lambda e -x$ is not left invertible by \ref{notinvideal} (ii), which implies
that $f\0 (x) = \lambda \in \s(x) \subset [0, \infty[$. It follows from \ref{Krein} that $f\0$
has a linear extension $f$ to $A\sa$ such that $f(x) \geq 0$ for all $x \in C = A_+$.
Then $f(a^*a) \geq 0$ for all $a \in A$ by the Shirali-Ford Theorem \ref{ShiraliFord}.
The linear extension $\varphi$ of $f$ from $A\sa$ to $A$ satisfies the requirement,
cf.\ \ref{eBbded}.
\end{proof}

\begin{proposition}\label{ispropleft}%
If $\varphi$ is a state on a unital normed  \st-algebra $B$, then
\[ M := \{ \,b \in B : \varphi(b^*b) = 0 \,\} \]
is a proper left ideal of $B$.
\end{proposition}

\begin{proof}
Let $\langle \cdot , \cdot \rangle _{\varphi}$ denote the positive Hilbert
form induced by $\varphi$, cf.\ \ref{inducedHf}. An element $b$ of $B$
is in $M$ if and only if $\langle b , c \rangle _{\varphi} = 0$ for all $c \in B$
by the Cauchy-Schwarz inequality. So $M$ is the isotropic subspace
of $\langle \cdot , \cdot \rangle _{\varphi}$, and thus a left ideal in $B$,
cf.\ \ref{quotientA}. The left ideal $M$ is proper as
$\varphi(e^*e) = \varphi(e) = 1$, cf.\ \ref{vare}.
\end{proof}

\begin{proposition}\label{maxleftstate}%
If $A$ is unital, then for a maximal left ideal $M$ in $A$ there exists a state
$\varphi$ on $A$ such that
\[ M = \{ \,a \in A : \varphi(a^*a) = 0 \,\}. \]
\end{proposition}

\begin{proof}
There exists a state $\varphi$ on $A$ such that $\varphi(a^*a) = 0$ for all
$a \in M$, cf.\ \ref{leftidstate}. The set $\{ \,a \in A : \varphi(a^*a) = 0 \,\}$
is a proper left ideal in $A$ containing $M$, cf.\ \ref{ispropleft} and thus
equal to $M$ by maximality of $M$.\pagebreak%
\end{proof}

\medskip
See also \ref{maxleftpurestate} below. 

\begin{theorem}\label{invstate}%
If $A$ is unital, then an element $a$ of $A$ is left inver\-tible if and only if
\[ \varphi(a^*a) > 0 \qquad \text{for all states } \varphi \text{ on } A. \]
\end{theorem}

\begin{proof}
Let $a \in A$. If $a$ is not left invertible, then $a$ lies in a proper
left ideal by \ref{notinvideal} (ii), which is contained in a maximal left
ideal by \ref{inmaxideal}, so that $\varphi(a^*a) = 0$ for some state
$\varphi$ on $A$ by proposition \ref{maxleftstate}. Conversely, if
$\varphi(a^*a) = 0$ for some state $\varphi$ on $A$, then $a$ lies
in the proper left ideal $\{ \,b \in A : \varphi(b^*b) = 0 \,\}$, cf.\ proposition
\ref{ispropleft}, so that $a$ cannot be left invertible by \ref{notinvideal} (ii).
\end{proof}

\begin{theorem}\label{maxleftpurestate}%
If $A$ is unital, then for a maximal left ideal $M$ in $A$, there exists a pure
state $\psi$ on $A$ such that
\[ M = \{ \,a \in A : \psi(a^*a) = 0 \,\}. \]
\end{theorem}

\begin{proof}
The set $K$ of those states $\varphi$ on $A$ with
\[ M = \{ \,a \in A : \varphi(a^*a) = 0 \,\} \]
is non-empty by \ref{maxleftstate}. Furthermore, $K$ is a compact
convex subset of the state space of $A$. (Indeed, if $\ell$ is in the
closed convex hull of $K$, then $\ell$ is a state (as $A$ is unital) with
\[ M \subset \{ \,a \in A : \ell \,(a^*a) = 0 \,\}. \]
But the set on the right hand side is a proper left ideal, and therefore
equal to $M$, by maximality of $M$.) Thus $K$ has an extreme point,
$\psi$ say. It remains to be shown that the state $\psi$ is pure. So let
\[ \psi = \lambda \varphi_1 + \mu \varphi_2 \]
with $\varphi_1$, $\varphi_2$ states on $A$ and $\lambda$,
$\mu \in \ ]0, 1[$, $\lambda + \mu = 1$. For $i \in \{ \,1, 2 \,\}$, we then have
\[ M \subset \{ \,a \in A : \varphi_i(a^*a) = 0 \,\} \]
by positivity of $\varphi_i$. As above, it follows that both sets are equal,
i.e.\ $\varphi_i \in K$, which implies $\psi = \varphi_1 = \varphi_2$.
\end{proof}

\begin{theorem}\label{invpstate}
If $A$ is unital, then an element $a$ of $A$ is left inver\-tible if and only if
\[ \psi(a^*a) > 0 \qquad \text{for all pure states } \psi \text{ on } A. \]
\end{theorem}

\begin{proof}
In the proof of \ref{invstate}, replace ``state'' by ``pure state''. \pagebreak
\end{proof}

A similar characterisation for the right invertible elements holds
as well, and thus we obtain: 

\begin{corollary}
A normal element of a unital Hermitian Banach \linebreak \st-algebra is
left invertible if and only if it it is right invertible.
\end{corollary}

\begin{definition}[$E(b)$]\label{Ea}\index{E(a)@$E(a)$}%
Let $B$ be a normed \st-algebra, and let $b \in B$. We shall consider
the set $E(b)$ defined as the set of those pure states $\psi$ on $B$,
for which $c_{\psi}$ is an eigenvector of $\pi_{\psi}(b)$ to the
eigenvalue $\psi(b)$. That is,
\[ E(b) := \{ \,\psi \in PS(B) : \pi_{\psi}(b)c_{\psi} = \psi(b)c_{\psi} \,\}. \]
\end{definition}

We have the following characterisation of $E(b)$:

\begin{proposition}\label{Eachar}%
For an element $b$ of a normed \st-algebra $B$, we have
\[ E(b) = \{ \,\psi \in PS(B) : {| \,\psi(b) \,| \,}^2 = \psi(b^*b) \,\}. \]
\end{proposition}

\begin{proof}
For a state $\psi$ on $B$, we calculate
\begin{align*}
&\ {\| \,\bigl( \pi_{\psi}(b) - \psi(b) \bigr) c_{\psi} \,\| \,}^2 \\
= &\ \langle \bigl( \pi_{\psi}(b) - \psi(b) \bigr) c_{\psi} ,
\bigl( \pi_{\psi}(b) - \psi(b) \bigr) c_{\psi} \rangle_{\psi} \\
= &\ {\| \,\pi_{\psi}(b)c_{\psi} \,\| \,}^2 + {\| \,\psi(b)c_{\psi} \,\| \,}^2
- 2 \,\mathrm{Re} \,\bigl( \overline{\psi(b)} \langle \pi_{\psi(b)}c_{\psi},
c_{\psi} \rangle_{\psi} \bigr) \\
= &\ \psi(b^*b) - {| \,\psi(b) \,| \,}^2. \qedhere
\end{align*}
\end{proof}

\begin{definition}[the left spectrum, $\mathrm{lsp}_B(b)$]%
\index{L15@$\mathrm{lsp}(a)$}\index{spectrum!of an element!left}%
If $b$ is an element of an algebra $B$, we define the
\underline{left spectrum} of $b$ as the set $\mathrm{lsp}_B(b)$
of those complex numbers $\lambda$ for which $\lambda e - b$
is not left invertible in $\tld{B}$. We shall also abbreviate
$\mathrm{lsp}(b) := \mathrm{lsp}_B(b)$.
\end{definition}

\begin{theorem}\label{leftirred}%
If $A$ is unital, then for an element $a$ of $A$, we have
\[ \mathrm{lsp}_A(a) = \{ \,\psi(a) : \psi \in E(a) \,\}. \]
In particular, each point in the left spectrum is an eigenvalue
in some irreducible representation. \pagebreak
\end{theorem}

\begin{proof}
Let $a \in A$, $\lambda \in \mathds{C}$ and put $b := \lambda e - a$.
For a pure state $\psi$ on $A$, we compute
\begin{align*}
\psi(b^*b) & = {| \,\lambda \,| \,}^2 + \psi(a^*a)
- 2 \,\mathrm{Re} \,\bigl( \overline{\lambda}\psi(a) \bigr) \\
 & = {| \,\lambda - \psi(a) \,| \,}^2 + \psi(a^*a) - {| \,\psi(a) \,| \,}^2.
\end{align*}
Now
\[ \psi(a^*a) - {| \,\psi(a) \,| \,}^2 \geq 0 \]
as $v(\psi) = 1$. Therefore $\psi(b^*b)$ vanishes if and only if both
$\lambda = \psi(a)$ and ${| \,\psi(a) \,| \,}^2 = \psi(a^*a)$. The statement
follows now from the results \ref{invpstate} and \ref{Eachar}.
\end{proof}

\begin{theorem}\label{lspC}%
If $A$ is unital, then for an element $a$ of $A$, we have
\[ \mathrm{lsp}_A(a) = \mathrm{lsp}_{C^*(A)} \bigl( j(a) \bigr). \]
Here, the left spectrum can be replaced with the right spectrum
and with the spectrum.
\end{theorem}

\begin{proof}
The bijection $\psi\0 \mapsto \psi\0 \circ j$ from $PS\bigl(C^*(A)\bigr)$
to $PS(A)$, cf.\ \ref{homPS}, restricts to a bijection from $E\bigl( j(a) \bigr)$
to $E(a)$. The statement concerning the right spectrum follows by
changing the multiplication in $A$ to $(a,b) \mapsto ba$.
The statement concerning the spectrum follows from \ref{leftrightinv}.
\end{proof}

\medskip
Now for non-unital algebras.

\begin{theorem}\label{leftirrednu}%
For an element $a$ of $A$, we have
\[ \mathrm{lsp}_A (a) \setminus \{0\}
= \{ \,\psi (a) : \psi \in E(a) \,\} \setminus \{0\}. \]
\end{theorem}

\begin{proof}
The inclusion ``$\supset$'' follows from \ref{extquasi}, \ref{leftirred},
and \ref{specunit}. Conversely, let $0 \neq \lambda \in \mathrm{lsp}_A (a)$.
We then also have $0 \neq \lambda \in \mathrm{lsp}_{\tld{A}} (a)$ by
\ref{specunit}. Hence there exists a pure state $\tld{\psi}$ on $\tld{A}$
with $0 \neq \lambda = \tld{\psi}(a)$ and ${| \,\tld{\psi} (a) \,| \,}^2 = \tld{\psi} (a^*a)$,
cf.\ \ref{leftirred} \& \ref{Eachar}. The restriction $\psi$ of $\tld{\psi}$ to $A$
then is a quasi-state with $0 \neq \lambda = \psi (a)$, and
${| \,\psi (a) \,| \,}^2 = \psi (a^*a)$. The fact that
$0 \neq {| \,\psi (a) \,| \,}^2 = \psi (a^*a)$ implies that $\psi$ actually is a state.
The state $\psi$ is pure because every state $\varphi$ on $\tld{A}$ satisfies
$\varphi (e) = 1$, cf. \ref{vare}. This makes that $\psi \in E(a)$, cf.\ \ref{Eachar}.
\pagebreak
\end{proof}

\begin{theorem}\label{lspCw}%
For an element $a$ of $A$, we have
\[ \mathrm{lsp}_A(a) \setminus \{0\}
= \mathrm{lsp}_{C^*(A)} \bigl( j(a) \bigr) \setminus \{0\}. \]
Here, the left spectrum can be replaced with the right spectrum
and with the spectrum.
\end{theorem}

\begin{proof}
The proof is the same as that of \ref{lspC}.
\end{proof}

\medskip
It is not very far from here to Ra\u{\i}kov's Criterion \ref{Raikov}.

\begin{theorem}
For an element $a$ of $A$, we have
\[ \mathrm{lsp}_A(a) \setminus \{0\}
= \mathrm{lsp} \bigl( \pi\univ (a) \bigr) \setminus \{0\}
= \mathrm{lsp} \bigl( \pi\ra (a) \bigr) \setminus \{0\}, \]
and the latter two sets consist entirely of eigenvalues.
\end{theorem}

\begin{proof}
This follows now from \ref{leftirrednu}, \ref{Ea}, and from the fact that
\[ \mathrm{lsp} \bigl( \pi (a) \bigr) \setminus \{0\}
\subset \mathrm{lsp} ( a ) \setminus \{0\} \]
for an algebra homomorphism $\pi$, cf.\ the proof of \ref{spechom}.
\end{proof}

\clearpage


\part{Spectral Theory of Representations}\label{part3}

\chapter{Representations by Normal Bounded Operators}

\setcounter{section}{33}


\section{Spectral Measures}

\begin{definition}[spectral measures]%
\index{spectral!measure}\index{measure!spectral}%
Let $H \neq \{0\}$ be a Hilbert space, and let $\mathcal{E}$ be
a $\sigma$-algebra on a set $\Omega \neq \varnothing$. Then a
\underline{spectral measure} defined on $\mathcal{E}$ and acting
on $H$ is a map
\[ P : \mathcal{E} \to B(H) \]
such that
\begin{itemize}
   \item[$(i)$] $P(\omega)$ is a projection in $H$ for each $\omega \in \mathcal{E}$,
  \item[$(ii)$] $P(\Omega) = \mathds{1}$,
 \item[$(iii)$] for all $x \in H$, the function
\end{itemize}
\begin{align*}
     \quad \langle P x, x \rangle : \mathcal{E} & \to \mathds{R}_+ \\
      \omega & \mapsto \langle P(\omega) x, x \rangle  = {\| \,P(\omega) x \,\| \,}^2
\end{align*}
\qquad\quad\ \ \ is a measure.

One notes that if $x$ is a unit vector then $\langle P x, x \rangle$ is a
probability measure. For $x$, $y \in H$ arbitrary, one defines a bounded
complex measure $\langle P x, y \rangle$ on $\Omega$ by polarisation:
\[ \langle Px,y \rangle \,:=
\,\frac14 \,\sum _{k=1}^{4} \,\iu ^k \,\langle P(x + \iu ^ky),(x + \iu ^ky) \rangle. \]
We then have $\langle P x, y \rangle (\omega) = \langle P(\omega) x, y \rangle$
for all $\omega \in \mathcal{E}$, and all $x,y \in H$.

By abuse of notation, one also puts
\[ P := \{ \,P (\omega) \in B(H) : \omega \in \mathcal{E} \,\}. \]
\end{definition}

For the remainder of this paragraph, let $H \neq \{0\}$ denote a Hilbert space,
let $\mathcal{E}$ denote a $\sigma$-algebra on a set $\Omega \neq \varnothing$,
and let $P$ denote a spectral measure defined on $\mathcal{E}$ and acting on $H$.
\pagebreak

\begin{proposition}%
The following statements hold.
\begin{itemize}
   \item[$(i)$] $P$ is finitely additive, i.e.\ if $D$ is a finite collection
of mutually disjoint sets in $\mathcal{E}$, then
\[ P\Bigl(\bigcup _{\text{\footnotesize{$\omega \in D$}}} \omega\Bigr)
= \sum _{\text{\footnotesize{$\omega \in D$}}} P(\omega), \]
  \item[$(ii)$] $P$ is orthogonal, i.e.\ if $\omega_1, \omega_2 \in \mathcal{E}$
and $\omega_1 \cap \omega_2 = \varnothing$, then
\[ P(\omega_1) P(\omega_2) = 0, \]
 \item[$(iii)$] the events in $\mathcal{E}$ are independent, i.e.\ for
$\omega_1, \omega_2 \in \mathcal{E}$, we have
\[ P(\omega_1 \cap \omega_2) = P(\omega_1)P(\omega_2). \]
In particular $P(\omega_1)$ and $P(\omega_2)$ commute.
\end{itemize}
\end{proposition}

\begin{proof}
(i): For all $x, y \in H$, we have
\[ \langle P\Bigl(\bigcup _{\text{\footnotesize{$\omega \in D$}}} \omega \Bigr) \,x, y \rangle
= \sum _{\text{\footnotesize{$\omega \in D$}}} \langle P(\omega)x,y \rangle
= \langle \sum _{\text{\footnotesize{$\omega \in D$}}} P(\omega)x,y \rangle. \]
(ii): Let $\omega_1,\omega_2 \in \mathcal{E}$
with $\omega_1 \cap \omega_2 = \varnothing$.
By (i), we have
\[ P(\omega_1) + P(\omega_2) = P(\omega_1 \cup \omega_2), \]
which says that $P(\omega_1) + P(\omega_2)$ is a projection. It follows
\begin{align*}
     &\ P(\omega_1)+P(\omega_2) = {\bigl(P(\omega_1)+P(\omega_2)\bigr)}^2 \\
 = &\ P(\omega_1)+P(\omega_2)+P(\omega_1)P(\omega_2)+P(\omega_2)P(\omega_1).
\end{align*}
This implies
\[ P(\omega_1)P(\omega_2) + P(\omega_2)P(\omega_1) = 0. \]
Multiplication from the left and from the right with $P(\omega_2)$ yields
\[ P(\omega_2)P(\omega_1)P(\omega_2)
+ P(\omega_2)P(\omega_1)P(\omega_2) = 0, \]
and so
\[ P(\omega_2)P(\omega_1)P(\omega_2) = 0. \]
The C*-property \ref{preC*alg} implies
\[ {\| \,P(\omega_1)P(\omega_2) \,\| \,}^2
= \| \,P(\omega_2)P(\omega_1)P(\omega_2) \,\| = 0, \]
that is
\[ P(\omega_1)P(\omega_2) = 0. \pagebreak \]
(iii): This follows from (i) and (ii) in the following way.
\begin{align*}
& \ P(\omega_1)P(\omega_2) \\
= & \ P\bigl((\omega_1 \cap \omega_2) \cup (\omega_1 \setminus \omega_2)\bigr)
\cdot P\bigl((\omega_2 \cap \omega_1) \cup (\omega_2 \setminus \omega_1)\bigr) \\
= & \ \bigl(P(\omega_1 \cap \omega_2)+P(\omega_1 \setminus \omega_2)\bigr)
\cdot \bigl(P(\omega_2 \cap \omega_1)+P(\omega_2 \setminus \omega_1)\bigr) \\
= & \ P(\omega_1 \cap \omega_2). \qedhere%
\end{align*}%
\end{proof}%

\begin{definition}[$I_P$]\index{I1@$I_P$}%
Let $\underline{I_P}$ be the linear map on the complex vector space
of $\mathcal{E}$-step functions with
$I_P(1_{\text{\small{$\omega$}}}) = P(\omega)$
for all $\omega \in \mathcal{E}$. That is, if
\[ h = \sum _{\text{\footnotesize{$\omega \in D$}}}
\alpha_{\text{\footnotesize{$\omega$}}} \,1_{\text{\footnotesize{$\omega$}}} \]
is an $\mathcal{E}$-step function with $D \subset \mathcal{E}$ a finite
partition of $\Omega$, then
\[ I_P (h) := \sum _{\text{\footnotesize{$\omega \in D$}}}
\alpha_{\text{\footnotesize{$\omega$}}} \,P(\omega). \]
This map $I_P$ is a representation of the \st-algebra of $\mathcal{E}$-step
functions in $H$ by the preceding proposition. Here $I_P$ stands for
``integral''.
\end{definition}

\begin{proposition}\label{stepint}%
For an $\mathcal{E}$-step function $h$, we have
\[ \langle I_P (h) x, x \rangle = \int h \,d \langle P x, x \rangle. \]
\end{proposition}

\begin{proof}
Indeed, if
\[ h = \sum _{\text{\footnotesize{$\omega \in D$}}}
\alpha_{\text{\footnotesize{$\omega$}}} \,1_{\text{\footnotesize{$\omega$}}} \]
is a $\mathcal{E}$-step function with $D \subset \mathcal{E}$ a finite
partition of $\Omega$, then
\begin{align*}
   \langle I_P (h) x, x \rangle
   & = \langle \sum _{\text{\footnotesize{$\omega \in D$}}}
          \alpha_{\text{\footnotesize{$\omega$}}} \,P(\omega) \,x, x \,\rangle \\
   & = \sum _{\text{\footnotesize{$\omega \in D$}}}
          \alpha_{\text{\footnotesize{$\omega$}}} \,\langle P(\omega) \,x, x \rangle \\
   & = \sum _{\text{\footnotesize{$\omega \in D$}}}
          \alpha_{\text{\footnotesize{$\omega$}}}
          \int \,1_{\text{\footnotesize{$\omega$}}} \,d \,\langle P x, x \rangle
       = \int h \,d \langle P x, x \rangle. \qedhere
\end{align*}
\end{proof}

\begin{definition}[$\mathcal{M}_b(\mathcal{E})$]%
\index{M2@$\mathcal{M}_b(\mathcal{E})$}%
We shall denote by $\underline{\mathcal{M}_b(\mathcal{E})}$ the
C*-al\-gebra of bounded complex-valued $\mathcal{E}$-measurable
functions equipped with the supremum norm $| \cdot | _{\infty}$,
cf.\ the appendix \ref{imagebegin}. \pagebreak
\end{definition}

\begin{definition}[spectral integrals]%
Let $f \in \mathcal{M}_b(\mathcal{E})$ and $a \in B(H)$. One writes
\[ a = \int f \,dP \qquad \text{(weakly)} \]
if
\[ \langle ax,x \rangle = \int f \,d \langle Px,x \rangle
\quad \text{for all } x \in H. \]
We then also have
\[ \langle ax,y \rangle = \int f \,d \langle Px,y \rangle
\quad \text{for all } x,y \in H, \]
so that then $a$ is uniquely determined by $f$ and $P$.

We shall write
\[ a = \int f \,dP \qquad \text{(in norm)} \]
if for every $\mathcal{E}$-step function $h$ one has
\[ \| \,a - I_P(h) \,\| \leq | \,f-h \,|_{\infty}. \]
\end{definition}

\begin{proposition}%
Let $f \in \mathcal{M}_b(\mathcal{E})$ and $a \in B(H)$. If
\[ a = \int f \,dP  \qquad \text{(in norm)} \]
then also
\[ a = \int f \,dP \qquad \text{(weakly)}. \]
In particular, $a$ then is uniquely determined by $f$ and $P$.
\end{proposition}

\begin{proof}
Since the $\mathcal{E}$-measurable function $f$ is bounded, there exists
a sequence $(h_n)$ of $\mathcal{E}$-step functions converging uniformly
to $f$. From
\[ a = \int f \,dP \qquad \text{(in norm)} \]
we have
\[ \| \,a - I_P ( h_n ) \,\| \leq | \,f - h_n \,|_{\infty} \to 0, \]
that is,
\[ a = \lim _{n \to \infty} I_P(h_n).  \]
With \ref{stepint} we get for all $x \in H$:
\[ \langle ax,x \rangle
   = \lim _{n \to \infty} \langle I_P(h_n)x,x \rangle
   = \lim _{n \to \infty} \int h_n \,d \langle Px,x \rangle
   = \int f \,d \langle Px,x \rangle, \pagebreak \]
which is the same as
\[ a = \int f \,dP \qquad \text{(weakly)}. \qedhere \]
\end{proof}

\begin{theorem}[$\pi_{P}(f)$]\index{p6@$\pi_P$}%
For $f \in \mathcal{M}_b(\mathcal{E})$ there exists a (necessarily
\linebreak unique) operator $\underline{\pi_{P}(f)} \in B(H)$ such that
\[ \pi _P(f) = \int f \,dP \qquad \text{(in norm)}. \]
The map
\begin{align*}
\pi _P : \mathcal{M}_b(\mathcal{E}) & \to B(H) \\
             f & \mapsto \pi_{P}(f)
\end{align*}
then defines a representation $\pi_P$ of $\mathcal{M}_b(\mathcal{E})$ in $H$.
\end{theorem}

\begin{proof}
Let $S(\mathcal{E})$ denote the complex vector space of $\mathcal{E}$-step
functions. For $h \in S(\mathcal{E})$, we may write
\[ h = \sum _{\text{\footnotesize{$\omega \in D$}}}
\alpha_{\text{\footnotesize{$\omega$}}} \,1_{\text{\footnotesize{$\omega$}}} \]
where $D \subset \mathcal{E}$ is a finite partition of $\Omega$. We obtain
\[ I_P(h) = \sum _{\text{\footnotesize{$\omega \in D$}}}
\alpha_{\text{\footnotesize{$\omega$}}} \,P(\omega). \]
In order to show that the representation $I_P$ is contractive:
\[ \| \,I_P(h) \,\| \leq | \,h \,|_{\infty}, \]
we note that for $x \in H$, we have
\begin{align*}
{\| \,I_P(h)x \,\| \,}^2 & = \sum _{\text{\footnotesize{$\omega \in D$}}}
{| \,\alpha_{\text{\footnotesize{$\omega$}}} \,| \,}^2 \ {\| \,P(\omega)x \,\| \,}^2 \\
 & \leq {| \,h \,|_{\infty} \,}^2 \sum _{\text{\footnotesize{$\omega \in D$}}}
{\| \,P(\omega)x \,\| \,}^2 = {| \,h \,|_{\infty} \,}^2 \ {\| \,x \,\| \,}^2.
\end{align*}
Let $\overline{S(\mathcal{E})}$ denote the closure of $S(\mathcal{E})$ in
the C*-algebra ${\ell\,}^{\infty}(\Omega)$. We have $\overline{S(\mathcal{E})}
\subset \mathcal{M}_b(\mathcal{E})$ as pointwise limits of sequences of
$\mathcal{E}$-measurable functions are $\mathcal{E}$-measurable. Also,
$\mathcal{M}_b(\mathcal{E}) \subset \overline{S(\mathcal{E})}$ because
every bounded $\mathcal{E}$-measurable function is the uniform limit of
a sequence of $\mathcal{E}$-step functions. We obtain the equality
$\mathcal{M}_b(\mathcal{E}) = \overline{S(\mathcal{E})}$.
It then follows that the mapping
\[ I_P : h \mapsto I_P(h) \qquad \bigl( h \in S(\mathcal{E}) \bigr) \pagebreak \]
has a unique extension to a continuous mapping
\[ \pi _P : \mathcal{M}_b(\mathcal{E}) \to B(H) \]
This mapping then is a contractive representation
(since $I_P$ is so), from which it follows that
\[ \pi _P(f) = \int f \,dP  \qquad \text{(in norm)} \]
for each $f \in \mathcal{M}_b(\mathcal{E})$.
\end{proof}

\begin{proposition}\label{weaklyinnorm}%
Let $a$ be a bounded linear operator in $H$ such that
\[ a = \int f \,dP \qquad \text{(weakly)} \]
for some $f \in \mathcal{M}_b(\mathcal{E})$. We then also have
\[ a = \int f \,dP  \qquad \text{(in norm)}. \]
\end{proposition}

\begin{proof}
We have $a = \pi _P(f)$ by $\pi _P(f) = \int f \,dP$ (weakly) as well.
\end{proof}

\begin{proposition}\label{normpointwise}%
For $f \in \mathcal{M}_b(\mathcal{E})$ and $x \in H$, we have
\[ \| \, \pi_P (f) \,x \,\| 
= \| \,f \,\|_{\text{\small{$\langle Px,x \rangle,2$}}}
= {\biggl(\int {| \,f \,| \,}^2 \ d \langle Px,x \rangle \biggr)}^{1/2}. \]
\end{proposition}

\begin{proof} One calculates
\[ {\| \,\pi _P(f) \,x \,\| \,}^2 = \langle \pi_P\bigl({ \,| \,f \,| \,}^2 \,\bigr) \,x, x \rangle
= \int {| \,f \,| \,}^2 \,d \langle Px,x \rangle. \qedhere \]
\end{proof}

\begin{theorem}\label{spanP}%
We have
\[ R(\pi _P) = \overline{\mathrm{span}}(P). \]
It follows that $R(\pi _P)$ is the C*-subalgebra of $B(H)$ generated by $P$.
\end{theorem}

\begin{proof}
By \ref{rangeC*} it follows that $R(\pi _P)$ is a
C*-algebra and so a closed subspace of $B(H)$,
containing $\mathrm{span}(P)$ as a dense subset.
\end{proof}

\begin{corollary}\label{spprojcomm}%
We have ${\pi _P}' = P'$, cf.\ \ref{repcommutant}.
\end{corollary}

\begin{definition}[$P$-a.e.]\index{P95@$P$-a.e.}%
A property applicable to the  points in $\Omega$ is said
to hold \underline{$P$-almost everywhere} ($P$-a.e.)
if it holds $\langle Px, x \rangle$-almost everywhere for
all $x \in H$. \pagebreak
\end{definition}

\begin{definition}[$N(P)$]\index{N(P)@$N(P)$}%
We denote by $N(P)$ the closed \st-stable ideal in
$\mathcal{M}_b(\mathcal{E})$ consisting of the functions
in $\mathcal{M}_b(\mathcal{E})$ which vanish \linebreak
$P$-almost everywhere.
\end{definition}

\begin{proposition}\label{kerspectralint}%
We have $N(P) = \ker \pi _P$.
\end{proposition}

\begin{proof}
This follows from \ref{normpointwise}.
\end{proof}

\begin{theorem}\label{spprojform}%
Every projection in $R(\pi _P)$ is of the form
$P(\omega)$ for some $\omega \in \mathcal{E}$.
\end{theorem}

\begin{proof} Let $f \in \mathcal{M}_b(\mathcal{E})$ be such
that $\pi _P(f) = {\pi _P(f) \,}^2$. We obtain $f = {f \,}^2$ $P$-a.e.
It follows that $f$ is $P$-a.e.\ equal to either $0$ or $1$.
Then, with $\omega = f^{-1}(1) \in \mathcal{E}$, we have
$f =1_{\textstyle\omega}$ $P$-a.e.\ and thus
$\pi _P(f) = \pi _P(1_{\textstyle\omega}) = P(\omega)$.
\end{proof}

\begin{theorem}\label{piPpos}%
For $f \in \mathcal{M}_b(\mathcal{E})$, we have
\[ \pi _P(f) \geq 0 \Leftrightarrow f \geq 0 \ P\text{-a.e.} \]
\end{theorem}

\begin{proof} We have
\begin{align*}
\pi _P(f) \geq 0 & \Leftrightarrow \pi _P(f) = \bigl| \,\pi _P(f) \,\bigr| \\
& \Leftrightarrow \pi _P(f) = \pi _P\bigl( \,| \,f \,| \,\bigr)
\Leftrightarrow f = | \,f \,| \ P \text{-a.e.}
\end{align*}
We used that $\pi _P\bigl( \,| \,f \,| \,\bigr) = \bigl| \,\pi _P(f) \,\bigr|$,
cf.\ \ref{homabs} \& the proof of \ref{Coabs}.
\end{proof}

\begin{theorem}[the Dominated Convergence Theorem]%
\label{domconvergence}%
If $(f_n)$ is a norm-bounded sequence in $\mathcal{M}_b(\mathcal{E})$
which converges pointwise $P$-a.e.\ to a function
$f \in \mathcal{M}_b(\mathcal{E})$, then
\[ \lim _{n \to \infty} \pi _P(f_n) \,x = \pi _P(f) \,x
\quad \text{for all } x \in H.  \]
\end{theorem}

\begin{proof}
This follows from Lebesgue's Dominated Convergence Theorem
because by \ref{normpointwise} we have
\[ \| \,\pi _P(f) \,x - \pi _P(f_n) \,x \,\|
= {\biggl( \int {| \,f-f_n \,| \,}^2 \,d \langle Px,x \rangle \biggr)}^{1/2} \to 0.
\pagebreak \qedhere \]
\end{proof}

\begin{corollary}[strong $\sigma$-additivity]%
If $(\omega _n)$ is a sequence of pairwise disjoint sets in $\mathcal{E}$, then
\[ P\Bigl(\,\bigcup _n \omega _n \Bigr) \,x = \sum _n P(\omega _n) \,x
\quad \text{for all } x \in H. \]
\end{corollary}

\begin{corollary}[strong $\sigma$-continuity]\label{strsgcont}%
If $(\omega _n)$ is an increasing sequence of sets in $\mathcal{E}$, then
\[ P\Bigl(\,\bigcup _n \omega _n \Bigr) \,x = \lim _{n \to \infty} P(\omega _n) \,x
\quad \text{for all } x \in H. \]
\end{corollary}

\begin{theorem}[the Monotone Convergence Theorem]%
\label{monconvthm}%
If $(f_n)$ is a \linebreak norm-bounded sequence in
$\mathcal{M}_b(\mathcal{E})\sa$ such that
$f_n \leq f_{n+1}$ $P$-a.e.\ for all $n$, then
\[ \sup _n \int f_n \,dP = \int \sup _n f_n \,dP. \]
\end{theorem}

\begin{proof} This follows from \ref{domconvergence}
and \ref{ordercompl}. \end{proof}

\begin{definition}[$L^{\infty}(P)$]\index{L2@$L^{\infty}(P)$}%
We denote
\[ L^{\infty}(P) := \mathcal{M}_b(\mathcal{E})/N(P), \]
which is a commutative C*-algebra, cf.\ \ref{C*quotient}.
\end{definition}

\begin{proposition}\label{LPmain}%
The representation $\pi _P$ factors to an isomorphism
of C*-algebras from $L^{\infty}(P)$ onto $R(\pi _P)$.
\end{proposition}

\begin{proof}
This follows from \ref{kerspectralint} and \ref{rangeC*}.
\end{proof}

\begin{definition}[image of a spectral measure]\label{imagespdef}%
\index{spectral!measure!image}\index{image!spectral measure}%
Let $\mathcal{E}'$ be a \linebreak $\sigma$-algebra on a set
$\Omega' \neq \varnothing$, and let $f : \Omega \to \Omega'$ be an
$\mathcal{E}$-$\mathcal{E}'$ measurable function. Then
\begin{alignat*}{2}
 f(P) : \mathcal{E}' & \to           &\ & B(H)_+ \\
                 \omega' & \mapsto &\ & P\bigl(f^{-1}(\omega')\bigr)
\end{alignat*}
defines a spectral measure $f(P) = P \circ f^{-1}$ defined on $\mathcal{E}'$,
called the \underline{image} of $P$ under $f$. For
$g \in \mathcal{M}_b(\mathcal{E'})$, one has
\[ \int g\ d\,f (P) = \int (g \circ f) \,dP. \]
See the appendix \ref{imagebegin} - \ref{imageend}. \pagebreak
\end{definition}

\clearpage


\section{Spectral Theorems}

\begin{definition}[resolution of identity]%
\index{resolution of identity}\index{spectral!measure!resolution of id.}%
Let $\Omega \neq \varnothing$ be a locally compact Hausdorff space.
Then a \underline{resolution of identity} on $\Omega$ is a \linebreak
spectral measure $P$, defined on the Borel $\sigma$-algebra of $\Omega$
and acting on a Hilbert space $H \neq \{0\}$, such that for all unit vectors $x$
in $H$, the Borel probability measure $\langle Px, x \rangle$ is inner regular,
cf.\ the appendix \ref{Bairedef} \& \ref{inregBormeas}.
\end{definition}

\begin{definition}[the support]\index{S5@$\mathrm{supp}(P)$}%
\index{support}\index{resolution of identity!support}%
Let $P$ be a resolution of identity on a locally compact Hausdorff space
$\Omega \neq \varnothing$. If there exists a largest open subset $\omega$
of $\Omega$ with $P(\omega) = 0$, then $\Omega \setminus \omega$ is
called the \underline{support} of $P$. The support of $P$ exists, as the
union $O$ of all open subsets $\omega$ of $\Omega$ with $P(\omega) = 0$
satisfies $P(O) = 0$, by a compactness argument and by inner regularity.
The support of $P$ is denoted by $\underline{\mathrm{supp}(P)}$.
\end{definition}

Please note that if the support of $P$ is all of $\Omega$, then
$C_b(\Omega)$ is imbedded in $L^{\infty}(P)$ as two continuous functions
differ on an open subset. This imbedding is isometric by \ref{C*injisometric}.

\begin{theorem}[the Spectral Theorem, archetypal form]\label{spthmC*}%
\index{Theorem!Spectral!commutative C*-algebra}%
\index{spectral!theorem!commutative C*-algebra}%
\index{spectral!resolution!commutative C*-algebra}%
For a commutative C*-algebra $B$ of bounded linear operators
acting non-degenerately on a Hilbert space $H \neq \{0\}$, there exists a
unique \linebreak resolution of identity $P$ on $\Delta(B)$, acting on $H$,
such that
\[ b = \int \wht{b} \,dP\quad \text{(in norm)}
\qquad \text{for all } b \in B. \]
One says that $P$ is the \underline{spectral resolution} of $B$.
We have
\[ B' = P', \]
and the support of $P$ is all of $\Delta(B)$.
\end{theorem}

\begin{proof}
If $x$ is a unit vector in $H$, then
\[ b \mapsto \langle bx,x \rangle \]
is a state on $B$, cf.\ \ref{staterep}. It follows that there is a unique inner regular
Borel probability measure $\langle Px,x \rangle$ on $\Delta(B)$ such that
\[ \langle bx,x \rangle = \int \wht{b} \,d \langle Px,x \rangle
\qquad \text{for all } b \in B, \pagebreak \]
cf.\ \ref{predisint}. Polarisation yields complex measures $\langle Px,y \rangle$
on $\Delta(B)$ such that
\[ \langle bx,y \rangle = \int \wht{b} \,d \langle Px,y \rangle \]
for all $x,y \in H$ and all $b \in B$. For a Borel set $\omega$ in
$\Delta(B)$, consider the positive semidefinite sesquilinear form
$\varphi$ defined by
\[ \varphi (x,y) := \int 1_{\textstyle\omega} \,d \langle Px,y \rangle
\qquad (x,y \in H). \]
Since $\varphi (x,x) \leq 1$ for all $x$ in the unit ball of $H$, there
exists a unique operator $P(\omega) \in B(H)_+$ such that
\[ \varphi (x,y) = \langle P(\omega)x,y \rangle \quad \text{for all } x, y \in H, \]
cf.\ \ref{posoprep}. In other words
\[ \langle P(\omega)x,y \rangle = \int 1_{\textstyle\omega} \,d \langle Px,y \rangle
\qquad (x,y \in H). \]
In order to show that $P : \omega \mapsto P(\omega)$ is a spectral measure
and hence a resolution of identity, it suffices to prove that each
$P(\omega)$ is an idempot\-ent. We shall prove more.

Let $\omega, \omega'$ be Borel subsets of $\Delta(B)$. We shall show that
$P(\omega) P(\omega') = P(\omega \cap \omega')$. Let $x$, $y \in H$.
For $a$, $b \in B$, we have
\[ \ \int \wht{\vphantom{b}a} \,\wht{b} \,d \langle Px,y \rangle
= \ \langle abx,y \rangle
= \ \int \wht{\vphantom{b}a} \,d \langle Pbx,y \rangle. \tag*{$(\st)$} \]
Now $\{ \,\wht{a} : a \in B \,\} = C\0\bigl(\Delta(B)\bigr)$, see the Commutative
Gelfand-Na\u{\i}mark Theorem \ref{commGN}. It follows that
\[ \wht{b} \cdot \langle Px,y \rangle = \langle Pbx,y \rangle \]
in the sense of equality of complex measures.
The integrals in $(\st)$ therefore remain equal if $\wht{a}$ is replaced by
$1_{\textstyle\omega}$. Hence
\begin{align*}
\int 1_{\textstyle\omega} \,\wht{b} \,d \langle Px,y \rangle
= & \ \int 1_{\textstyle\omega} \,d \langle Pbx,y \rangle \\
= & \ \langle P(\omega) bx,y \rangle \\
= & \ \langle bx,P(\omega)y \rangle
=     \int \wht{b} \,d \langle Px,P(\omega)y \rangle. \pagebreak \tag*{$(\st\st)$}
\end{align*}
The same reasoning as above shows that the integrals in $(\st\st)$
remain equal if $\wht{b}$ is replaced by $1_{\textstyle\omega'}$. Consequently
we have
\begin{align*}
\langle P(\omega \cap \omega') x,y \rangle
= & \ \int 1_{\textstyle{\omega \cap \omega'}} \,d \langle Px,y \rangle \\
= & \ \int 1_{\textstyle{\omega\vphantom{'}}} 1_{\textstyle{\omega'}} \,d \langle Px,y \rangle \\
= & \ \int 1_{\textstyle{\omega'}} \,d \langle Px,P(\omega)y \rangle \\
= & \ \langle P(\omega')x,P(\omega)y \rangle
=     \langle P(\omega)P(\omega')x,y \rangle,
\end{align*}
so that $P(\omega \cap \omega') = P(\omega)P(\omega')$ indeed, which
finishes the proof that $P$ is a resolution of identity. 

From
\[ \langle bx,y \rangle = \int \wht{b} \,d \langle Px,y \rangle \]
it follows that
\[ b = \int \wht{b} \,dP \qquad \text{(weakly)}, \]
whence also
\[ b = \int \wht{b} \,dP \qquad \text{(in norm)}, \]
cf.\ \ref{weaklyinnorm}.

Uniqueness follows from $\{ \,\wht{b}: b \in B \,\} = C\0\bigl(\Delta(B)\bigr)$,
cf.\ the Commutative Gelfand-Na\u{\i}mark Theorem \ref{commGN}.

It shall now be shown that the support of $P$ is all of $\Delta(B)$,
i.e.\ $\varnothing$ is the largest open subset $\omega$ of $\Delta(B)$
with $P(\omega) = 0$. Let $\omega$ be a non-empty open subset of
$\Delta(B)$ and let $x \in \omega$. There then exists $f \in C_c(\Delta(B)$
such that
\[ 0 \leq f \leq 1,\ f(x) = 1,\ \mathrm{supp} (f) \subset \omega. \]
Let $b \in B$ with $\wht{b} = f$. We then have $b \neq 0$, as well as $b \geq 0$
by \ref{pointwiseorder}. Since $f \leq 1_{\textstyle\omega}$, it follows by \ref{piPpos}:
\[ 0 \lneq b = \pi _P \bigl( \,\wht{b} \,\bigr) = \pi _P(f) \leq \pi _P(1_{\textstyle\omega}) = P(\omega), \]
whence $P(\omega) \neq 0$. \pagebreak

It shall next be shown that $B' = P'$. For $a \in B(H)$, the following
statements are equivalent.
\begin{align*}
& a \in P', \\
& P(\omega)a = aP(\omega) \qquad \text{ for all Borel sets } \omega, \\
& \langle P(\omega) ax,y \rangle = \langle aP(\omega)x,y \rangle
\qquad \text{ for all Borel sets } \omega, \text{ for all } x, y \in H, \\
& \langle P(\omega)ax,y \rangle = \langle P(\omega)x,a^*y \rangle
\qquad \text{ for all Borel sets } \omega, \text{ for all } x, y \in H, \\
& \langle Pax,y \rangle = \langle Px,a^*y \rangle
\qquad \text{ for all } x, y \in H, \\
& \int \wht{b} \,d \langle Pax,y \rangle
= \int \wht{b} \,d \langle Px,a^*y \rangle
\qquad \text{ for all } b \in B, \text{ for all } x, y \in H, \\
& \langle bax,y \rangle = \langle bx,a^*y \rangle
\qquad \text{ for all } b \in B, \text{ for all } x, y \in H, \\
& \langle bax,y \rangle = \langle abx,y \rangle
\qquad \text{ for all } b \in B, \text{ for all } x, y \in H, \\
& ba = ab \qquad \text{ for all } b \in B, \\
& a \in B'.
\end{align*}
The proof is complete. 
\end{proof}

\medskip
Our next aim is the Spectral Theorem for a normal
bounded linear operator \ref{spthmnormalbded}.
We shall need the following result.

\begin{theorem}[Fuglede - Putnam - Rosenblum]%
\index{Theorem!Fuglede-Putnam-Rosenblum}%
\index{Fuglede-Putn.-Ros.\ Thm.}
Let $A$ be a \linebreak pre-C*-algebra. Let $n_1, n_2$ be normal
elements of $A$. Let $a \in A$ and assume that $an_1 = n_2a$.
We then also have $a{n_1}^* = {n_2}^*a$.

In particular, if $n$ is a normal bounded linear operator on a
Hilbert space, then $\{n\}' = \{n^*\}'$.
\end{theorem}

\begin{proof}
We may assume that $A$ is a C*-algebra.
Please note that it follows from the hypothesis that $a{n_1}^k = {n_2}^ka$
for all integers $k \geq 0$, so if $p \in \mathds{C}[z]$, then
$a p(n_1) = p(n_2)a$. It follows that
\[ a \,\exp(\iu \overline{z}n_1) = \exp(\iu \overline{z}n_2) \,a \]
for all $z \in \mathds{C}$, or equivalently
\[ a = \exp(\iu \overline{z}n_2) \,a \,\exp(-\iu \overline{z}n_1). \]
Since $\exp(x+y) = \exp(x)\exp(y)$ when $x$ and $y$ commute, \pagebreak
the normality of $n_1$ and $n_2$ implies
\begin{align*}
 & \exp\bigl(\iu z{n_2}^*\bigr) \,a \,\exp\bigl(-\iu z{n_1}^*\bigr) \\
= & \ \exp\bigl(\iu z{n_2}^*\bigr)\exp\bigl(\iu \overline{z}n_2\bigr)
\,a \,\exp\bigl(-\iu \overline{z}n_1\bigr)\exp\bigl(-\iu z{n_1}^*\bigr) \\
= & \ \exp\bigl(\iu (z{n_2}^*+\overline{z}n_2)\bigr)
\,a \,\exp\bigl(-\iu (\overline{z}n_1+z{n_1}^*)\bigr)
\end{align*}
The expressions $z{n_2}^*+\overline{z}n_2$ and
$\overline{z}n_1+z{n_1}^*$ are Hermitian, and hence \linebreak
$\exp\bigl(\iu (z{n_2}^*+\overline{z}n_2)\bigr)$ and
$\exp\bigl(-\iu (\overline{z}n_1+z{n_1}^*)\bigr)$ are unitary. This implies
\[ \| \,\exp\bigl(\iu z{n_2}^*) \,a \,\exp(-\iu z{n_1}^*\bigr) \,\| \leq \| \,a \,\| \]
for all $z \in \mathds{C}$. However
\[ z \mapsto \exp(\iu z{n_2}^*) \,a \,\exp(-\iu z{n_1}^*) \]
is an entire function, hence constant by Liouville's Theorem, so that
\[ \exp\bigl(\iu z{n_2}^*\bigr) \,a = a \,\exp\bigl(\iu z{n_1}^*\bigr) \]
for all $z \in \mathds{C}$. By equating the coefficients for $z$
we obtain ${n_2}^*a = a{n_1}^*$.
\end{proof}

\medskip
We remark that the above result immediately carries over to
\st-semi\-simple normed \st-algebras, cf.\ \ref{stsemipreC*}.

\begin{theorem}%
[the Spectral Theorem for normal bounded linear operators]%
\label{spthmnormalbded}%
\index{Theorem!Spectral!normal bounded operator}%
\index{spectral!theorem!normal bounded operator}%
\index{spectral!resolution!normal bounded operator}%
Let $b$ be a normal bounded linear operator on a Hilbert space
$H \neq \{0\}$. There then is a unique resolution of identity $P$
on $\s(b)$, acting on $H$, such that
\[ b = \int _{\text{\small{$\s(b)$}}} \id _{\text{\small{$\s(b)$}}} \,dP
\quad \text{(in norm)}. \]
This is usually written as
\[ b = \int z \,dP(z). \]
One says that $P$ is the \underline{spectral resolution} of $b$. We have
\[ \{ \,b \,\}' = P', \]
and the support of $P$ is all of $\s(b)$. \pagebreak
\end{theorem}

\begin{proof} Let $B$ be the C*-subalgebra of $B(H)$ generated by
$b$, $b^*$, and $\mathds{1}$. It is commutative as $b$ is normal. By
\ref{spechomeom}, we may identify $\Delta(B)$ with $\s(b)$ in such a way
that $\wht{b}$ becomes $\id _{\text{\small{$\s(b)$}}}$, and existence follows.
On the other hand, if $Q$ is another resolution of identity on $\s(b)$ with
\[ b = \int z \,dQ(z), \]
then
\[ p(b,b^*) = \int p(z,\overline{z}) \,dQ(z) \]
for all $p \in \mathds{C}[z,\overline{z}]$. Uniqueness follows then from the
Stone-Weierstrass \linebreak Theorem, cf.\ the appendix \ref{StW}. The
statement concerning the \linebreak commutants follows from the preceding
theorem.
\end{proof}

\begin{addendum}%
Let $b$ be a normal bounded linear operator on a Hilbert space
$H \neq \{0\}$. Let $P$ be the spectral resolution of $b$. The map
\begin{align*}
L^{\infty}(P) & \to B(H) \\
f+N(P) & \mapsto \pi _P(f)
\end{align*}
extends the operational calculus for $b$ \ref{opcalc}:
\begin{align*}
C\bigl(\s(b)\bigr) & \to           B(H) \\
                            f & \mapsto f(b).
\end{align*}
\end{addendum}

\begin{remark}\label{rangespres}%
The range of the operational calculus for $b$ is the \linebreak
C*-subalgebra of $B(H)$ generated by $b$, $b^*$, and $\mathds{1}$.
We shall give a similar characterisation of the range of the extended
map, see \ref{PBorel} below.
\end{remark}

\begin{proposition}\label{eival}%
Let $b$ be a normal bounded linear operator on a Hilbert space $H \neq \{0\}$.
Let $P$ be its spectral resolution. A complex number $\lambda \in \s(b)$
is an eigenvalue of $b$ if and only if $P\bigl(\{ \,\lambda \,\}\bigr) \neq 0$, i.e.\ if
$P$ has an atom at $\lambda$, in which case the range of
$P\bigl(\{ \,\lambda \,\}\bigr)$ is the eigenspace corresponding to the
eigenvalue $\lambda$.

Hence an isolated point of the spectrum of $b$ is an eigenvalue of $b$.%
\pagebreak%
\end{proposition}

\begin{proof}
Assume first that $P\bigl(\{ \,\lambda \,\}\bigr) \neq 0$, and let
$0 \neq x \in H$ with $P\bigl(\{ \,\lambda \,\}\bigr) \,x = x$. For such $x$, we have
\begin{align*}
b\,x & = \int z \,dP(z) \,x
= \int z \,dP(z) \,P(\{\,\lambda \,\}) \,x \\
 & = \left[ \int _{\text{\small{$\s(b) \setminus \{ \,\lambda \,\}$}}} z \,dP(z)
+ \int _{\text{\small{$\{ \,\lambda \,\}$}}} z \,dP(z) \right]
\cdot P(\{ \,\lambda \,\}) \,x \\
 & = \lambda \,P\bigl(\{ \,\lambda \,\}\bigr)^2\,x
= \,\lambda \,P\bigl(\{ \,\lambda \,\}\bigr)\,x = \lambda \,x,
\end{align*}
so that $x$ is an eigenvector of $b$ to the eigenvalue $\lambda$.
Conversely, let $x \neq 0$ be an eigenvector of $b$ to the eigenvalue $\lambda$.
For $n \geq 1$, consider the function $f_n$ given by
$f_n(z) = (\lambda - z)^{-1}$ if $| \,\lambda -z \,| > \frac1n$,
and $f_n(z) = 0$ if $| \,\lambda - z \,| \leq \frac1n$. We get
\[ \pi _P(f_n)(\lambda \mathds{1} - b) = \pi _P \bigl(f_n(\lambda - \id)\bigr)
= P\bigl(\bigl\{ \,z : | \,z - \lambda \,| > \textstyle{\frac1n} \,\bigr\}\bigr). \]
Since $(\lambda \mathds{1} - b)x = 0$, we have
\[ P\bigl(\bigl\{ \,z : | \,z - \lambda \,| > \textstyle{\frac1n} \,\bigr\}\bigr) \,x = 0. \]
Letting $n \to \infty$, we obtain $P\bigl(\{ \,z : z \neq \lambda \,\}\bigr) \,x = 0$,
cf.\ \ref{strsgcont}, whence $P(\{ \,\lambda \,\}) \,x = x$.
The last statement follows from the fact that the support of $P$ is all of $\s(b)$.
\end{proof}

\begin{definition}[$\Delta^*(A)$]\index{D(A)@$\Delta^*(A)$}%
If $A$ is a \st-algebra, we denote the set of Hermitian multiplicative
linear functionals on $A$ by \underline{$\Delta^*(A)$} and equip it
with the weak* topology, cf.\ the appendix \ref{weak*top}.
\end{definition}

\begin{definition}[$\s(\pi), \pi^*$]%
\label{specpi}\index{s3@$\protect\s(\pi)$}%
Let $A$ be a commutative  \st-algebra. Let $\pi$ be a non-zero
representation of $A$ in a Hilbert space $H \neq \{0\}$.
Let $B$ denote the closure of $R(\pi)$ in $B(H)$. We define
\[ \s(\pi) := \{ \,\tau \circ \pi \in \mathds{C}^{\,\text{\Small{$A$}}} : \tau \in \Delta(B) \,\}. \]
Then $\s(\pi)$ is a subset of $\Delta^*(A)$ with $\s(\pi) \cup \{0\}$
weak* compact. Hence $\s(\pi)$ is a closed and locally compact
subset of $\Delta^*(A)$. The functions $\wht{a}| _{\text{\small{$\s(\pi)$}}}$
$(a \in A)$ are dense in $C\0\bigl(\s(\pi)\bigr)$ by the Stone-Weierstrass
\linebreak Theorem, cf.\ the appendix \ref{StW}. The ``adjoint'' map
\begin{alignat*}{2}\index{p7@$\pi^*$}
\pi^* : \Delta(B) & \to &\ & \s(\pi) \\
             \tau & \mapsto &\ & \tau \circ \pi
\end{alignat*}
is a homeomorphism. \pagebreak
\end{definition}

\begin{proof}
For $\tau \in \Delta(B)$, the functional $\tau \circ \pi$ is non-zero
by density of $R(\pi)$ in $B$. We have $\s(\pi) \subset \Delta^*(A)$
by \ref{mlfHerm}. Extend the ``adjoint'' map $\pi^*$ to a map
$\Delta(B) \cup \{0\} \to \s(\pi) \cup \{0\}$, taking $0$ to $0$. The
extended map is bijective and continuous by the universal property
of the weak* topology, cf.\ the appendix \ref{weak*top}. Since the
domain $\Delta(B) \cup \{0\}$ is compact \ref{mlf0compact}, and
the range $\s(\pi) \cup \{0\}$ is Hausdorff, it follows that the extended
map is a homeomorphism, cf.\ the appendix \ref{homeomorph}.
In particular, $\s(\pi) \cup \{0\}$ is compact, whence the statement.
\end{proof}

\begin{remark}\label{sppisa}%
Assume that furthermore $A$ is a normed \st-algebra. Since a
multiplicative linear functional on a Banach algebra is contractive
\ref{mlfbounded}, we have the following facts. If $\pi$ is weakly
continuous on $A\sa$, then each $\tau \in \s(\pi)$ is contractive on
$A\sa$, cf.\ \ref{weaksa}, \ref{sigmaAsa}. Also, if $\pi$ is weakly
continuous, then each $\tau \in \s(\pi)$ is contractive, cf.\ \ref{weakcontHilbert}.
\end{remark}

\begin{proposition}%
Let $A$ be a commutative \st-algebra. Let $\pi$ be a non-zero
representation of $A$ in a Hilbert space. Then
\[ \s\bigl(\pi(a)\bigr) \setminus \{0\} = \wht{a}\,\bigl(\s(\pi)\bigr) \setminus \{0\}
\qquad \text{for all } a \in A, \]
whence, by \ref{C*rl},
\[ \| \,\pi(a) \,\| = \bigl| \,\wht{a}| _{\text{\small{$\s(\pi)$}}} \,\bigr|_{\infty}
\qquad \text{for all } a \in A. \]
\end{proposition}

\begin{proof}
If $B$ denotes the closure of $R(\pi)$, and $a \in A$, then
\begin{align*}
\s\bigl(\pi(a)\bigr) \setminus \{0\}
& = \wht{\pi (a)}\,\bigl(\Delta(B)\bigr) \setminus \{0\} \qquad \text{by \ref{rangeGT}} \\
& = \wht{a}\,\bigl(\s(\pi)\bigr) \setminus \{0\}. \qedhere
\end{align*}
\end{proof}

\begin{theorem}[the Spectral Theorem for representations]%
\label{spthmrep}%
\index{Theorem!Spectral!representation}%
\index{spectral!theorem!representation}%
\index{spectral!resolution!representation}%
Let $\pi$ be a non-degenerate representation of a commutative \st-algebra
$A$ in a \linebreak Hilbert space $H \neq \{0\}$. There then exists a unique
resolution of identity $P$ on $\s(\pi)$, acting on $H$, such that
\[ \pi(a) = \int _{\text{\small{$\s(\pi)$}}} \wht{a}|_{\text{\small{$\s(\pi)$}}} \,dP
\quad \text{(in norm)} \qquad \text{for all } a \in A. \]
One says that $P$ is the \underline{spectral resolution} of $\pi$. We have
\[ \pi' = P', \]
and the support of $P$ is all of $\s(\pi)$. \pagebreak
\end{theorem}

\begin{proof} Uniqueness follows from the fact that the functions
$\wht{a}| _{\text{\small{$\s(\pi)$}}}$ $(a \in A)$ are dense in
$C\0\bigl(\s(\pi)\bigr)$, cf.\ \ref{specpi}. For existence, consider the closure
$B$ of $R(\pi)$. It is a commutative C*-algebra acting non-degenerately
on $H$. The spectral resolution $Q$ of $B$ is a resolution of identity on
$\Delta(B)$ such that
\begin{gather*}
b = \int \wht{b} \,dQ\quad \text{for all } b \in B, \\
B' = Q',
\end{gather*}
and such that the support of $Q$ is all of $\Delta(B)$. Let $P$ be the
image of $Q$ under the ``adjoint'' map
\begin{align*}
\pi^* : \Delta(B) & \to \s(\pi) \\
                \tau & \mapsto \tau \circ \pi,
\end{align*}
which is a homeomorphism, cf.\ \ref{specpi}. (For the notion of an
image of a spectral measure, recall \ref{imagespdef}, and see the
appendix \ref{imagebegin} - \ref{imageend}.) For $a \in A$, we have
\[ \int_{\text{\small{$\s(\pi)$}}} \wht{a}|_{\text{\small{$\s(\pi)$}}} \,dP
= \int \wht{\pi(a)} \,dQ = \pi(a), \]
cf.\ e.g.\ the appendix \ref{imageend}. Since $R(\pi)$ is dense in $B$,
we have
\[ \pi' = B' = Q' = P'. \]
Since the ``adjoint'' map $\pi^*$ is a homeomorphism, it is clear that
the support of $P$ is all of $\s(\pi)$.
\end{proof}

\medskip
The uniqueness statement above presupposes
that the resolution of identity lives on $\s(\pi)$.
A stronger uniqueness statement goes like this:

\begin{addendum}%
Let $\pi$ be a non-degenerate representation of a commutative \st-algebra
$A$ in a Hilbert space $H \neq \{0\}$. Let $P$ be the spectral resolution
of $\pi$. Assume that $S$ is a subset of $\Delta^*(A)$ with $S \cup \{0\}$ weak*
compact, and that $Q$ is a resolution of identity on $S$, acting on $H$, with
\[ \pi(a) = \int _{\text{\small{$S$}}} \wht{a}|_{\text{\small{$S$}}} \,dQ
\quad \text{(in norm)} \qquad \text{for all } a \in A. \]
Assume also that the support of $Q$ is all of $S$. Then $S = \s(\pi)$ and $Q = P$.
\pagebreak
\end{addendum}

\begin{proof}
Let $a \in A$. Then $\| \,\pi(a) \,\| = \| \,\pi_Q ( \,\wht{a}|_{\text{\small{$S$}}} ) \,\|$
equals the quotient norm of $\wht{a}|_{\text{\small{$S$}}}$ in the quotient
C*-algebra $L^{\infty}(Q)$, cf.\ \ref{LPmain}. Hence, for any $\alpha > \| \,\pi(a) \,\|$,
we have $Q$-almost everywhere that
$\bigl| \ \wht{a}|_{\text{\small{$S$}}} \,\bigr| \leq \alpha$,
by definition of the quotient norm, see the appendix \ref{quotspace}.
(It is used that for $X \subset \mathds{R}$ and $\alpha \in \mathds{R}$ one has
$\inf X < \alpha$ if and only if there exists $x \in X$ with $x < \alpha$.
The ``if'' part is trivial. Conversely, if $\alpha \leq x$ for all $x \in X$, then
$\alpha$ is a lower bound of $X$, and therefore $\alpha \leq \inf X$
as $\inf X$ is the greatest lower bound of $X$.)
Since $\wht{a}|_{\text{\small{$S$}}}$ is continuous, the set
\[ T := \bigl\{ \,\sigma \in S : \bigl| \ \wht{a}|_{\text{\small{$S$}}} (\sigma) \,\bigr|
> \| \,\pi(a) \,\| \,\bigr\} \]
is an open subset of $S$ with $Q(T) = 0$, by \ref{strsgcont}. Since by assumption
the support of $Q$ is all of $S$, it follows that actually $T = \varnothing$. In other
words
\[ \bigl| \,\sigma (a) \,\bigr| = \bigl| \ \wht{a}|_{\text{\small{$S$}}} (\sigma) \,\bigr|
\leq \| \,\pi(a) \,\| \tag*{$(\st)$} \]
holds for all $\sigma \in S$ and all $a \in A$.
Consider any $\sigma \in S$. It is a Hermitian multiplicative linear functional
on $A$, which vanishes on $\ker \pi$, as is seen from $(\st)$. Thus $\sigma$
factors to a Hermitian multiplicative linear functional $\varphi$ on $R(\pi)$ such
that $\varphi\bigl(\pi(a)\bigr) = \sigma (a)$ for all $a \in A$. (The functional
$\varphi$ is not identically zero because $\sigma$ isn't.) The inequality $(\st)$
implies that $\varphi$ is contractive, and so $\varphi$ has a unique continuation
to a Hermitian multiplicative linear functional $\tau$ on the closure of $R(\pi)$ in
$B(H)$. It is immediate that $\tau \circ \pi = \sigma$, and consequently
$\sigma \in \s(\pi)$. We have shown that $S \subset \s(\pi)$. We also have that
$S$ is a closed subset of $\s(\pi)$. We thus may consider $Q$ as a resolution
of identity on $\s(\pi)$, whose support is $S$. Then $Q = P$ by the fact that the
functions $\wht{a}| _{\text{\small{$\s(\pi)$}}}$ $(a \in A)$ are dense in
$C\0\bigl(\s(\pi)\bigr)$, cf.\ \ref{specpi}. Hence also
$S = \mathrm{supp}(Q) = \mathrm{supp}(P) = \s(\pi)$.
\end{proof}

\medskip
We remark that we cannot take all of $\Delta ^* (A)$ as the space of
decomposition, as it can fail to be locally compact.

The hurried reader can go from here to chapter \ref{unbdself}
on unbounded self-adjoint operators.

\clearpage


\section{Von Neumann Algebras}

Let throughout $H$ denote a Hilbert space.

We now come to commutants and second commutants in $B(H)$,
which we encountered in the preceding paragraphs, 
especially \ref{questimbed} and \ref{CyclicNonDeg}.

\begin{definition}\index{commutant}%
Let $S \subset B(H)$. The set of all operators in $B(H)$ which commute
with every operator in $S$ is called the \underline{commutant} of $S$
and is denoted by $S'$. Please note that $S \subset S''$.
\end{definition}

\begin{proposition}%
For two subsets $S$, $T$ of $B(H)$, we have
\[ T \subset S \Rightarrow S' \subset T'. \]
\end{proposition}

\begin{theorem}%
For a subset $S$ of $B(H)$, we have
\[ S''' = S'. \]
\end{theorem}

\begin{proof}
From $S \subset S''$ it follows that $(S'')' \subset S'$.
Also $S' \subset (S')''$.
\end{proof}

\medskip
Recall that we called a subset of a \st-algebra \underline{\st-stable}
(or self-adjoint), if with each element $a$, it also contains the adjoint
$a^*$, cf.\ \ref{selfadjointsubset}.

\begin{proposition}%
Let $S$ be a subset of $B(H)$. If $S$
is \st-stable, then $S'$ is a C*-algebra.
\end{proposition}

\begin{definition}%
\index{algebra!von Neumann}\index{von Neumann!algebra}%
\index{algebra!W*-algebra}\index{W*-algebra}%
A \underline{von Neumann algebra} or \underline{W*-algebra}
on $H$ is a subset of $B(H)$ which is the commutant of a
\st-stable subset of $B(H)$.
\end{definition}

\begin{proposition}%
If $H \neq \{0\}$, then each von Neumann algebra
on $H$ is a unital C*-subalgebra of $B(H)$.
\end{proposition}

\begin{example}
$B(H)$ is a von Neumann algebra because
$B(H) = (\mathds{C}\mathds{1})'$.
If $\pi$ is a representation of a \st-algebra in a Hilbert space,
then $\pi'$ is a von Neumann algebra.
\end{example}

\begin{remark}
The term ``W*-algebra'' is also used in a wider sense for a
C*-algebra which is isomorphic (as a C*-algebra) to a von
\linebreak Neumann algebra, see
\cite[vol.\ I, Def.\ III.3.1, p.\ 130]{Tak}. \pagebreak
\end{remark}

Von Neumann algebras come in pairs:

\begin{theorem}%
If $\mathscr{M}$ is a von Neumann algebra, then
\begin{itemize}
  \item[$(i)$] $\mathscr{M}'$ is a von Neumann algebra,
 \item[$(ii)$] $\mathscr{M}'' = \mathscr{M}$.
\end{itemize}
\end{theorem}

The property in the following theorem is often taken as the
definition of von Neumann algebras. 

\begin{theorem}\index{commutant!second}%
A C*-subalgebra $A$ of $B(H)$ is a von
Neumann \linebreak algebra if and only if $A'' = A$.
\end{theorem}

For two vectors $x$, $y \in H$, we shall denote by $x \odot y$
the bounded operator of rank one on $H$ given by
$(x \odot y) \,z := \langle z, y \rangle \,x$ $(z \in H)$.

\begin{proposition}%
Let $\{ \,x_i \,\}_{i \in I}$ be an orthonormal basis of $H$. If
$a \in B(H)$ commutes with all the bounded operators
$x_i \odot x_j$ $(i, j \in I)$, then $a$ is a scalar multiple of
the unit operator on $H$.
\end{proposition}

\begin{proof}
If $a$ commutes with each $x_i \odot x_j$ $(i, j \in I)$, we get
\[ \langle a z, x_j \rangle \,x_i = \langle z, x_j \rangle \,a x_i \quad
\text{for all } z \in H, \ i, j \in I. \]
Especially for $z = x_k$ $(k \in I)$, we obtain
\[ \langle a x_k, x_j \rangle \,x_i
= \langle x_k ,x_j \rangle \,a x_i = \delta_{kj} a x_i. \]
In particular we have
\[ \langle a x_k, x_j \rangle = 0 \quad
\text{for all } k, j \in I, \ k \neq j \]
and
\[ \langle a x_j, x_j \rangle \,x_i = a x_i \quad
\text{for all } i, j \in I. \]
It follows that there exists some $\lambda \in \mathds{C}$ such that
\[ \langle a x_j, x_j \rangle = \lambda \quad \text{for all } j \in I. \]
Thus
\[ \langle a x_i, x_j \rangle = \lambda \,\delta_{ij} \quad
\text{for all } i, j \in I. \qedhere \]
\end{proof}

\begin{corollary}\label{BHirred}%
We have $B(H)' = \mathds{C}\mathds{1}$. In particular,
$\mathds{C}\mathds{1}$ is a von Neumann algebra. \pagebreak
\end{corollary}

\begin{proposition}%
If $S$ is a \st-stable subset of $B(H)$, then the set $S''$
is the smallest von Neumann algebra containing $S$.
\end{proposition}

\begin{proof}
By the above, $S''$ is a von Neumann algebra containing $S$.
If $\mathscr{M}$ is a von Neumann algebra containing $S$,
then $S \subset \mathscr{M}$ whence $\mathscr{M}' \subset S'$
and thus $S'' \subset \mathscr{M}'' = \mathscr{M}$.
\end{proof}

\begin{definition}\index{W*(pi)@$W^*(\pi)$}%
If $S$ is an arbitrary subset of $B(H)$, we denote by $\underline{W^*(S)}$
the smallest von Neumann algebra containing $S$. Obviously
$W^*(S) = (S \cup S^*)''$. One says that $W^*(S)$ is the
\underline{von Neumann algebra} \underline{generated by $S$}.
If $\pi$ is a representation of a \st-algebra in a Hilbert space,
we denote $\underline{W^*(\pi)} := W^*\bigl(R(\pi)\bigr)$, so that
$W^*(\pi) = R(\pi)'' = (\pi')'$.%
\end{definition}

\begin{proposition}\label{Wstarcomm}%
The von Neumann algebra generated by a \linebreak normal
subset \ref{selfadjointsubset} of $B(H)$ is commutative.
\end{proposition}

\begin{proof}
See \ref{scndcomm}.
\end{proof}

\begin{proposition}\label{spprojsubset}%
Let $\pi$ be a non-degenerate representation of a commutative \st-algebra
in a Hilbert space $\neq \{0\}$. Let $P$ be the spectral resolution of $\pi$.
We then have $R(\pi _P) \subset W^*(\pi)$. In particular, for the spectral
projections, we have $P \subset W^*(\pi)$.
\end{proposition}

\begin{proof}
We have
\[ {\pi _P}' = P' = {\pi}', \]
see \ref{spprojcomm} \& \ref{spthmrep}, whence
\[ R(\pi _P) \subset R(\pi _P)'' = ({\pi _P}')' = (\pi ')' = W^*(\pi). \qedhere \]
\end{proof}

\begin{proposition}\label{gensame}%
Let $\pi$ be a non-degenerate representation of a commutative \st-algebra
in a Hilbert space $\neq \{0\}$. Let $P$ be the spectral \linebreak resolution
of $\pi$. Then $\pi _P$ and $\pi$ generate the same von Neumann
\linebreak algebra. Furthermore we have $W^*(\pi) = W^*(P)$.
\end{proposition}

\begin{proof} This follows again from ${\pi _P}' = P' = \pi'$. \end{proof}

\medskip
Please note that $R(\pi _P)$ is merely the C*-algebra generated
by $P$, cf.\ \ref{spanP}. On a non-separable Hilbert space, it can
happen that $W^*(\pi)$ is strictly larger than $R(\pi _P)$,
cf.\ \ref{nonsep} below. \pagebreak

\begin{theorem}\label{Abelian}%
Let $\mathscr{M}$ be a commutative von Neumann  algebra
on a Hilbert space $\neq \{0\}$. Let $P$ be the spectral 
resolution of $\mathscr{M}$. We then have $\mathscr{M} = R(\pi _P)$.
\end{theorem}

\begin{proof} By the Spectral Theorem we have
$\mathscr{M} \subset R(\pi _P)$, and by \ref{spprojsubset} we also have
$R(\pi _P) \subset W^*(\mathscr{M}) = \mathscr{M}$.
\end{proof}

\begin{definition}[$\mathscr{P}(\mathscr{M})$]%
For a von Neumann algebra $\mathscr{M}$, we denote by
$\mathscr{P}(\mathscr{M})$ the set of projections in $\mathscr{M}$.
\end{definition}

\begin{proposition}%
Let $\mathscr{M}$ be a commutative von Neumann algebra on a Hilbert
space $\neq \{0\}$. Let $P$ be the spectral resolution of $\mathscr{M}$.
We then have $\mathscr{P}(\mathscr{M}) = P$ and
$\mathscr{M} = \overline{\mathrm{span}}\bigl(\mathscr{P}(\mathscr{M})\bigr)$.
\end{proposition}

\begin{proof}
This follows from $\mathscr{M} = R(\pi _P) = \overline{\mathrm{span}}(P)$,
cf.\ \ref{Abelian} \& \ref{spanP}, and from the fact that every projection
in $R(\pi _P)$ is of the form $P(\omega)$, cf.\ \ref{spprojform}. 
\end{proof}

\begin{theorem}\label{W*main}%
A von Neumann algebra $\mathscr{M}$ satisfies
\[ \mathscr{M} = \overline{\mathrm{span}}\bigl(\mathscr{P}(\mathscr{M})\bigr). \]
For emphasis: a von Neumann algebra is the closed linear span of its
\linebreak projections.
\end{theorem}

\begin{proof}
Let $c \in \mathscr{M}$. Then $c$ can be decomposed as $c = a + \iu b$ with
$a,b \in \mathscr{M}\sa$. It now suffices to consider the commutative von
Neumann algebras $W^*(a), W^*(b) \subset \mathscr{M}$, cf.\ \ref{Wstarcomm}.
\end{proof}

\medskip
The above theorem is of fundamental importance for representation
\linebreak theory. With its help, questions concerning a von
Neumann algebra $\mathscr{M}$ can be transferred to questions
concerning the set $\mathscr{P}(\mathscr{M})$ of projections in
$\mathscr{M}$. E.g.\ Schur's Lemma \ref{Schur} appears trivial in
this light.

\begin{corollary}%
A von Neumann algebra $\mathscr{M}$ satisfies
\[ \mathscr{M}' = \mathscr{P}(\mathscr{M})', \]
whence also $\mathscr{M} = W^*\bigl(\mathscr{P}(\mathscr{M})\bigr)$.
\end{corollary}

A nice observation is the following one: \pagebreak

\begin{proposition}\label{disconn}%
Let $\mathscr{M}$ be a commutative von Neumann algebra on a
Hilbert space $\neq \{0\}$. Let $P$ be the spectral resolution of $\mathscr{M}$.
The imbedding of $C\bigl(\Delta(\mathscr{M})\bigr)$ into
$L^{\infty}(P)$ then is an isomorphism of C*-algebras. In particular, for a
bounded Borel function $f$ on $\Delta(\mathscr{M})$, there exists a unique
continuous function which is $P$-a.e.\ equal to $f$.
\end{proposition}

\begin{proof}
The mapping in question is an imbedding indeed as the support of $P$
is all of $\Delta(\mathscr{M})$, cf.\ \ref{spthmC*}. It is surjective by
$\mathscr{M} = R(\pi _P)$, cf.\ \ref{Abelian}. It thus is a \st-algebra
isomorphism, and therefore isometric, cf.\ \ref{Cstarisom}.
\end{proof}

\smallskip
We shall need the following technical lemma.

\begin{lemma}\label{lemmabic}%
Let $\pi$ be a representation of a \st-algebra in $H$. Let $b \in W^*(\pi)$.
We then also have $\oplus _n \,b \in W^*(\oplus _n \,\pi)$. Furthermore, if
$x \in \oplus _n \,H$, and $M := \overline{R \,(\oplus _n \,\pi) \,x}$, then the
closed subspace $M$ is invariant under $\oplus _n \,b$, and
$(\oplus _n \,b)|_M \in W^*\bigl((\oplus _n \,\pi)|_M\bigr)$.
(Here, $\oplus _n$ is shorthand for $\oplus _{n=1} ^{\infty}$.)
\end{lemma}

\begin{proof}
The first statement follows from the fact that each
element $c$ of $(\oplus _n \,\pi)'$ has a ``matrix representation''
$(c_{mn})_{m,n}$ with each $c_{mn} \in \pi'$. Clearly $M$ is a closed subspace
of $\oplus _n \,H$, invariant under $\oplus _n \,\pi$. With $p$ denoting the
projection on $M$, it follows that $p \in (\oplus _n \,\pi)'$, cf.\ \ref{invarcommutant}.
From the first statement, we have $\oplus _n \,b \in W^*(\oplus _n \,\pi)$, so
that $\oplus _n \,b$ commutes with $p$, which in turn implies that $M$ is
invariant under $\oplus _n \,b$, by \ref{invarcommutant} again. But then also
$(\oplus _n \,b)|_M \in W^*\bigl((\oplus _n \,\pi)|_M\bigr)$.
\end{proof}

\smallskip
The remainder of this paragraph won't be used in the sequel,
though it is of considerable importance. We shall need the
following somewhat lengthy definition.

\begin{definition}%
\index{topology!weak operator}\index{topology!strong operator}%
The \underline{strong operator topology} on $B(H)$ is the locally convex
topology induced by the seminorms
\[ a \mapsto \| \,ax \,\| \qquad (x \in H). \]

The \underline{weak operator topology} on $B(H)$
is the locally convex topology induced by the linear functionals
\[ a \mapsto \langle ax,y \rangle \qquad (x, y \in H). \]

The weak operator topology on $B(H)$ is weaker than the strong operator
topology on $B(H)$ which in turn is weaker than the norm topology.
\pagebreak

The \underline{$\sigma$-strong topology} on $B(H)$ is the locally convex
topology \linebreak induced by the seminorms
\[ a \,\mapsto \,{\biggl( \,\sum _n {\,\| \,ax_n \,\|\,}^2 \biggr)}^{1/2} \]
with $(x_n)$ any sequence in $H$ such that
$\bigl( \,\| \,x_n \,\| \,\bigr) \in {\ell\,}^2$.
\index{topology!sigma-weak@$\sigma$-weak}%
\index{topology!sigma-strong@$\sigma$-strong}%

The \underline{$\sigma$-weak topology} on $B(H)$ is the locally
convex topology \linebreak induced by the linear functionals
\[ a \,\mapsto \,\sum _n \,\langle ax_n,y_n \rangle \]
with $(x_n)$, $(y_n)$ any sequences in $H$
such that $\bigl( \,\| \,x_n \,\| \,\bigr), \bigl( \,\| \,y_n \,\| \,\bigr) \in {\ell\,}^2$.

The $\sigma$-weak topology is weaker than the $\sigma$-strong topology
which in turn is weaker than the norm topology.

Also, the weak operator topology is weaker than the $\sigma$-weak
topology, and similarly the strong operator topology is weaker that the
$\sigma$-strong topology.
\end{definition}

\begin{theorem}[von Neumann's Bicommutant Theorem]%
\index{Theorem!von Neumann!Bicommutant}\label{bicomm}%
\index{Bicommutant Theorem}\index{von Neumann!Theorem!Bicommutant}%
Let $\pi$ be a non-degenerate representation of a \st-algebra in $H$.
Then $W^*(\pi)$ is the \linebreak closure of $R(\pi)$ in any of the four
topologies: strong operator topology, weak operator topology, $\sigma$-strong
topology, and $\sigma$-weak topology. 
\end{theorem}

\begin{proof}
It is easily seen that a von Neumann algebra is closed in the weak
operator topology, and hence in all of the four topologies, because
the weak operator topology is the weakest among them. Since the \linebreak
$\sigma$-strong topology is the strongest of these topologies, it
suffices to prove that $R(\pi)$ is $\sigma$-strongly dense in $W^*(\pi)$.
Let $b \in W^*(\pi)$. Let $\varepsilon >  0$ and let $(x_n)$ be a sequence
in $H$ with $\bigl( \,\| \,x_n \,\| \,\bigr) \in {\ell\,}^2$. It suffices to show that
there exists an element $a$ of the \st-algebra such that
\[ {\biggl( \,\sum _n \,{\bigl\| \,\bigl( \,b - \pi(a) \,\bigr) \,x_n \,\bigr\|\,}^2 \,\biggl)}^{\,1/2}
< \varepsilon. \]
Consider $x := \oplus _n \,x_n \in \oplus _n \,H$ and put
$M := \overline{R \,( \oplus _n \,\pi ) \,x}$. Now $M$ is a closed subspace
invariant and cyclic under $\oplus _n \,\pi$, with cyclic vector $x$,
cf.\ \ref{cyclicsubspace}. The closed subspace $M$ is invariant under
$\oplus _n \,b$ by lemma \ref{lemmabic}. Since then
$(\oplus _n \,b) \,x \in M$, there exists by definition of
$M$ an element $a$ of the \st-algebra with
\[ \bigl\| \,( \,\oplus _n \,b \,) \,x - \bigl( \,\oplus _n \,\pi \,(a) \,\bigr) \,x \,\bigr\|
< \varepsilon. \pagebreak \qedhere \]
\end{proof}

\begin{corollary}%
Let $A$ be a \st-subalgebra of $B(H)$ acting non-degen\-erately on $H$.
Then $W^*(A)$ is the closure of $A$ in any of the four topo\-logies: strong
operator topology, weak operator topology, $\sigma$-strong topology,
and $\sigma$-weak topology.
\end{corollary}

\begin{corollary}%
The von Neumann algebras on $H$ are precisely the C*-subalgebras of
$B(H)$ acting non-degenerately on $H$ which are closed in any, hence all,
of the four topologies: strong operator topology, weak operator topology,
$\sigma$-strong topology, and $\sigma$-weak topology.
\end{corollary}

\begin{remark}%
Some authors define von Neumann algebras as \linebreak C*-algebras
of bounded operators on a Hilbert space, which are closed in the strong
operator topology, say. This is a more inclusive definition than ours, but
one can always put oneself in our situation by considering a smaller
Hilbert space, cf.\ \ref{nondegenerate}.
\end{remark}

We mention here that among the four topologies considered above,
the $\sigma$-weak topology is of greatest interest, because
it actually is a weak* topology resulting from a predual of $B(H)$, namely
the Banach \linebreak \st-algebra of trace class operators on $H$,
equipped with the trace-norm.

\clearpage


\chapter{Multiplication Operators and their Applications}

\setcounter{section}{36}


\section{Multiplication Operators}

\begin{definition}[$M_1(\Omega)$]\index{M1@$M_1(\Omega)$}%
If $\Omega \neq \varnothing$ is a locally compact Hausdorff space,
one denotes by \underline{$M_1(\Omega)$} the set of inner regular
Borel probability measures on $\Omega$, cf.\ the appendix
\ref{Bairedef} \& \ref{inregBormeas}.
\end{definition}

\begin{theorem}[multiplication operators]%
\index{M15@$M_{\mu}(g)$}\index{D(g)@$\mathcal{D}(g)$}%
Let $\Omega \neq \varnothing$ be a locally  compact Hausdorff
space, and let $\mu \in M_1(\Omega)$. If $g$ is a fixed
$\mu$-measur\-able function on $\Omega$, one defines
\[ \mathcal{D}(g) := \{ \,f \in L^2(\mu) : gf \in L^2(\mu) \,\} \]
as well as
\begin{align*}
M\smu(g) : \mathcal{D}(g) & \to L^2(\mu) \\
f & \mapsto gf.
\end{align*}
One says that $M\smu(g)$ is the \underline{multiplication operator}
with $g$. Then $M\smu(g)$ is bounded if and only if
$g \in {\mathscr{L} \,}^{\infty}(\mu)$, in which case
\[ \| \,M\smu(g) \,\| = \| \,g \,\|_{\,\text{\small{$\mu$}},\infty}. \]
For $M\smu(g)$ to be bounded, it suffices to be bounded on
$C_c(\Omega)$.
\end{theorem}

\begin{proof}
If $g \in {\mathscr{L} \,}^{\infty}(\mu)$, one sees that
$\| \,M\smu(g) \,\| \leq \| \,g \,\|_{\,\text{\small{$\mu$}},\infty}$. The
converse inequality is established as follows. Assume that
\[ \| \,gh \,\|_{\,\text{\small{$\mu$}},2}
\leq c \,\| \,h \,\|_{\,\text{\small{$\mu$}},2}
\quad \text{for all } h \in C_c(\Omega). \]
It shall be proved that $\| \,g \,\|_{\,\text{\small{$\mu$}},\infty} \leq c$.
Let $A := \{ \,t \in \Omega : | \,g(t) \,| > c \,\}$, which is a
$\mu$-measurable set, and thus $\mu$-integrable. We have to show
that $A$ is a $\mu$-null set. For this it is enough to show that each
compact subset of $A$ is a $\mu$-null set. (By inner regularity of the
completion of $\mu$.) So let $K$ be a compact subset of $A$. For
$h \in C_c(\Omega)$ with $1_K \leq h \leq 1$, we get
\begin{align*}
\int {| \,g \,| \,}^2 \,1_K \,d \mu
& \,\leq \int {| \,g \,{h \,}^{1/2} \,| \,}^2 \,d \mu
= {\| \,g \,{h \,}^{1/2} \,\|_{\,\text{\small{$\mu$}},2} \,}^2 \\
& \,\leq {c \,}^2 \,{\| \,{h \,}^{1/2} \,\|_{\,\text{\small{$\mu$}},2} \,}^2
= {c \,}^2 \int h \,d \mu.
\end{align*}
Thus, by taking the infimum over all such functions $h$, we obtain
\[ \int {| \,g \,| \,}^2 \,1_K \,d \mu \,\leq {c \,}^2 \,\mu(K). \]
Since $| \,g \,| > c$ on $K$, the set $K$ must be $\mu$-null, whence
$A$ is a $\mu$-null set. We have shown that
$\| \,g \,\|_{\,\text{\small{$\mu$}},\infty} \leq c$.
\end{proof}

\begin{lemma}%
Let $\mu \in M_1(\Omega)$, where $\Omega \neq \varnothing$ is a
locally compact Hausdorff space. If a bounded linear operator $b$
on $L^2(\mu)$ commutes with every $M\smu(f)$ with $f \in C_c(\Omega)$,
then $b = M\smu(g)$ for some $g \in {\mathscr{L} \,}^{\infty}(\mu)$.
\end{lemma}

\begin{proof}
Please note first that $1_{\Omega} \in {\mathscr{L} \,}^2(\mu)$.
For every $f \in C_c(\Omega)$ we have
\[ bf = b \bigl(M\smu(f)1_{\Omega} \bigr)
= M\smu(f)(b1_{\Omega}) = M\smu(f)g \]
with $g := b1_{\Omega} \in {\mathscr{L} \,}^2(\mu)$. We continue to compute
\[ M\smu(f)g = fg = gf = M\smu(g)f. \]
Thus the operators $b$ and $M\smu(g)$ coincide on $C_c(\Omega)$.
Then $M\smu(g)$ is bounded by the last statement of the preceding
theorem, and so $g \in {\mathscr{L} \,}^{\infty}(\mu)$. By density
of $C_c(\Omega)$ in $L^2(\mu)$, cf.\ \ref{CcdenseinLp}, it follows that
$b$ and $M\smu(g)$ coincide everywhere.
\end{proof}

\begin{definition}%
[maximal commutative \protect\st-subalgebras]%
\index{algebra!von Neumann!maximal commutative}%
\index{von Neumann!algebra!maximal commutative}%
\index{maximal!commutative}%
If $H$ is a Hilbert space, then
a \underline{maximal commutative} \st-subalgebra of $B(H)$
is a commutative \st-subalgebra of $B(H)$ which is not properly
contained in any other commutative \st-subalgebra of $B(H)$.
\end{definition}

\begin{proposition}\label{maxcommprime}%
If $H$ is a Hilbert space, then a \st-subalgebra $A$ of $B(H)$
is maximal commutative if and only if $A' = A$.
\end{proposition}

\begin{proof}
For $a = a^* \in A'$, the \st-subalgebra generated by
$A$ and $a$ is a commutative \st-subalgebra of $B(H)$.
\end{proof}

\medskip
In particular maximal commutative \st-subalgebras of $B(H)$ are
von Neumann algebras. \pagebreak

\begin{corollary}\index{L3@$L^{\infty}(\mu)$}\label{maxcommL}%
Let $\mu \in M_1(\Omega)$, with $\Omega \neq \varnothing$
a locally compact Hausdorff space. The mapping
\begin{align*}
L^{\infty}(\mu) & \to B\bigl(L^2(\mu)\bigr) \\
g & \mapsto M\smu(g)
\end{align*}
establishes an isomorphism of C*-algebras onto a maximal
commutative \linebreak \st-subalgebra of $B\bigl(L^2(\mu)\bigr)$.
In short, one speaks of the maximal \linebreak commutative von
Neumann algebra $L^{\infty}(\mu)$. 
\end{corollary}

\clearpage


\section{The Spectral Representation}

In this paragraph, we shall consider a fixed \underline{cyclic} representation
$\pi$ of a \underline{commutative} \st-algebra $A$ in a Hilbert
space $H \neq \{0\}$. We shall denote by $P$ the spectral resolution of $\pi$.

\begin{theorem}\label{Bochnerbis}\index{representation!spectral}%
\index{p8@$\pi_{\mu}$}\index{spectral!representation}%
If $c$ is a unit cyclic vector for $\pi$, there exists a unique
measure $\mu \in M_1\bigl(\s(\pi)\bigr)$ such that
\[ \langle \pi(a)c,c \rangle = \int \wht{a} \,d \mu
\qquad \text{for all } a \in A, \]
namely $\mu = \langle Pc,c \rangle$. The representation $\pi$ then is spatially
equivalent to the so-called \underline{spectral representation} $\pi \smu$, given by
\begin{align*}
\pi \smu : A & \to B\bigl(L^2(\mu)\bigr) \\
a & \mapsto M\smu(\wht{a}).
\end{align*}
The philosophy is that the spectral representation $\pi \smu$ is
particularly simple as it acts through multiplication operators.
\end{theorem}

\begin{proof}
Uniqueness follows from the fact that
$\{ \,\wht{a}| _{\text{\small{$\s(\pi)$}}} : a \in A \,\}$
is dense in $C\0\bigl(\s(\pi)\bigr)$, cf.\ \ref{specpi}.
Let $\psi$ be the positive linear functional on $A$ defined by
\[ \psi (a) := \langle \pi (a)c,c \rangle\quad (a \in A). \]
Putting $\mu := \langle Pc,c \rangle$, it is clear that
\[ \psi (a) =
\int _{\text{\small{$\s(\pi)$}}} \wht{a}|_{\text{\small{$\s(\pi)$}}} \,d \mu
\quad \text{for all } a \in A. \]
Now consider the dense subspace
\[ H\0 := \{ \,\pi (a)c \in H : a \in A \,\} \]
of $H$. An isometry $V : H\0 \to L^2(\mu)$ is constructed
in the following way: For $x = \pi (a)c \in H\0$, define
$Vx := \wht{a} \in L^2(\mu)$. It has to be verified that this is
well-defined. So let $\pi (b)c = \pi (a)c$. For $d \in A$ we obtain
\[ \psi (db) = \langle \pi (d) \pi (b)c,c \rangle
= \langle \pi (d) \pi (a)c,c \rangle = \psi (da), \]
whence
\[ \int \wht{d} \ \wht{b} \,d \mu
= \int \wht{d} \ \wht{\vphantom{b}a} \,d \mu
\qquad \text{for all } d \in A,  \]
which implies
\[ \int f \,\wht{b} \,d \mu = \int f \,\wht{\vphantom{b}a} \,d \mu
\qquad \text{for all } f \in C\0\bigl(\s(\pi)\bigr). \pagebreak \]
This means that the measures
$\wht{b} \cdot \mu$ and $\wht{\vphantom{b}a} \cdot \mu$
are equal, whence
\[ \wht{b} = \wht{\vphantom{b}a} \ \mu \text{-a.e.} \]
It shall next be shown that the mapping $V$
defined so far is isometric. For $a \in A$ we have
\[ {\| \,\pi (a)c \,\| \,}^2 = \psi (a^*a) =
\int {| \,\wht{a} \,| \,}^2 \,d \mu = {\| \,\wht{a} \,\|_{\,\text{\small{$\mu$}},2} \,}^2. \]
It shall now be proved that
\[ V \pi (a) = \pi \smu(a) V\quad \text{on } H\0 \]
for all $a \in A$. So let $a \in A$ be fixed. For $x = \pi (b)c$ we have
\[ \pi \smu(a) Vx = M\smu\bigl(\,\wht{\vphantom{b}a}\,\bigr) Vx
= M\smu\bigl(\,\wht{\vphantom{b}a}\,\bigr) \,\wht{b} = \wht{\vphantom{b}a} \,\wht{b} \]
as well as
\[ V \pi (a)x = V \pi (ab)c = \wht{(ab)} = \wht{\vphantom{b}a} \,\wht{b}. \]
The closure of $V$ then is a unitary operator which intertwines
$\pi$ and $\pi \smu$, cf.\ \ref{closureunitary}.
\end{proof}

\begin{proposition}\index{Theorem!abstract!Bochner}\label{Bremark}%
\index{abstract!Bochner Theorem}\index{Bochner Thm., abstr.}%
If $A$ furthermore is a commutative \underline{normed} \linebreak \st-algebra,
if $\psi$ is a state on $A$, and if in the above we choose $\pi := \pi _{\psi}$,
$c := c _{\psi}$, then the measure obtained is  essentially the same
as the one provided by the abstract Bochner Theorem \ref{Bochner},
cf.\ \ref{sppisa}.
\end{proposition}

\begin{theorem}\index{P95@$P_{\mu}$}\label{exspectres}%
Let $c$ be a unit cyclic vector for $\pi$, and let
$\mu = \langle Pc, c \rangle \in M_1\bigl(\s(\pi)\bigr)$.
For a Borel set $\omega$ in $\s(\pi)$ one defines
\[ P\smu(\omega) := M\smu(1_{\textstyle\omega}) \in B\bigl(L^2(\mu)\bigr). \]
The map
\[ P \smu : \omega \mapsto P \smu(\omega) = M \smu(1_{\textstyle\omega}) \]
then is a resolution of identity on $\s(\pi)$. Indeed we have
\[ \langle P \smu x,y \rangle = x\, \overline{\vphantom{b}y}
\cdot \mu \qquad \text{for all } x,y \in L^2(\mu). \]
The resolution of identity $P \smu$ is the spectral resolution of $\pi \smu$.
For every bounded Borel function $f$ on $\s(\pi)$, we have
\[ \pi _{\text{\footnotesize{$P$}} \smu}(f) = M \smu(f).  \]
It follows that the C*-algebra $R(\pi _{\text{\footnotesize{$P$}} \smu})$ is
equal to the maximal commutative von Neumann algebra $L^{\infty}(\mu)$
\ref{maxcommL}. Indeed a function $f \in {\mathscr{L} \,}^{\infty}(\mu)$ is
$\mu$-a.e.\ equal to a bounded Borel (even Baire) function, cf.\ the appendix
\ref{Baire}. \pagebreak
\end{theorem}

\begin{proof}
For a Borel set $\omega$ in $\s(\pi)$ and $x,y$ in $L^2(\mu)$, we have
\[ \langle P \smu (\omega) \,x, y \rangle
= \langle M\smu(1_{\textstyle\omega}) \,x, y \rangle
= \int 1_{\textstyle\omega} \,x\, \overline{\vphantom{b}y} \,d \mu, \]
whence
\[ \langle P \smu \,x, y \rangle = x\, \overline{\vphantom{b}y} \cdot \mu. \]
It follows that for each $x \in H$, the measure $\langle P \smu x, x \rangle$
is an inner regular Borel measure on $\s(\pi)$, so that $P \smu$ is a resolution
of identity on $\s(\pi)$. Next, for $a \in A$, we have
\[ \langle \pi \smu (a) \,x, y \rangle
= \langle M\smu\bigl(\,\wht{\vphantom{b}a}\,\bigr) \,x, y \rangle
= \int \wht{\vphantom{b}a} \,x\, \overline{\vphantom{b}y} \,d \mu
= \int \wht{\vphantom{b}a} \,d \langle P \smu x, y \rangle \]
for all $x,y \in L^2(\mu)$, which says that
\[ \pi \smu (a) = \int _{\text{\small{$\s(\pi)$}}} \wht{a}|_{\text{\small{$\s(\pi)$}}} \,d P \smu
= \int _{\text{\small{$\s(\pi\smu)$}}} \wht{a}|_{\text{\small{$\s(\pi\smu)$}}} \,d P \smu
\quad \text{(weakly).} \]
(Please note here that $\s(\pi \smu) = \s(\pi)$ as both representations are spatially
equivalent.) This implies that $P \smu$ is the spectral resolution of $\pi \smu$.
For a bounded Borel function $f$ on $\s(\pi)$, we have
\[ \pi _{\text{\footnotesize{$P$}} \smu}(f) = \int f \,d P \smu\quad \text{(weakly),} \]
so that for $x,y \in L^2(\mu)$ we get
\[ \langle \pi _{\text{\footnotesize{$P$}} \smu} (f) \,x, y \rangle
= \int f \,d \langle P \smu x, y \rangle
= \int f \,x\, \overline{\vphantom{b}y} \,d \mu = \langle M \smu(f) \,x, y \rangle. \]
This shows that
\[ \pi _{\text{\footnotesize{$P$}} \smu} (f) = M\smu (f). \qedhere \]
\end{proof}

\begin{theorem}\label{spectralrepmain}%
Let $Q$ be the restriction of $P$ to the Baire sets. Then
\[ R(\pi _Q) = R(\pi _P) = W^*(\pi) = \pi' \]
is a maximal commutative von Neumann algebra. 
\end{theorem}

\begin{proof}
Let $c$ be a cyclic unit vector for $\pi$ and let $\mu$ be the measure
$\langle Pc,c \rangle$. Then $\pi$ is spatially equivalent to $\pi \smu$,
cf.\ \ref{Bochnerbis}. Let $S$ denote the restriction of $P \smu$ to the
Baire sets. The equalities
\[ R(\pi _{\text{\footnotesize{$S$}}}) = R(\pi _{\text{\footnotesize{$P$}} \smu})
= W^*(\pi \smu) = {\pi \smu}' \]
follow from the fact that $R(\pi _{\text{\footnotesize{$S$}}})$ is a maximal
commutative von Neumann algebra, cf.\ the last two statements of the 
preceding theorem \ref{exspectres}. \pagebreak
\end{proof}

\medskip
Please note that if $\pi$ is not necessarily cyclic, we merely have the inclusions
\[ R(\pi _Q) \subset R(\pi _P) \subset W^*(\pi) \subset \pi'. \]
The last inclusion holds by commutativity of $W^*(\pi)$, cf.\ \ref{Wstarcomm}.

\bigskip\begin{corollary}[reducing subspaces]%
\index{reducing}\index{subspace!reducing}%
\index{invariant}\index{subspace!invariant}%
A closed subspace $M$ of $H$ \linebreak reduces $\pi$ if and only if $M$
is the range of some $P(\omega)$, where $\omega$ is a Borel (or even Baire)
subset of $\s(\pi)$.
\end{corollary}

\begin{proof}
Let $Q$ be the restriction of $P$ to the Baire sets. Let $M$ be a closed
subspace of $H$ and let $p$ denote the projection on $M$. The subspace
$M$ reduces $\pi$ if and only if $p \in \pi' = R(\pi _{Q})$,
cf.\ \ref{invarcommutant} and the preceding theorem. However a projection
in $R(\pi _{Q})$ is of the form $Q(\omega)$ for some Baire subset
$\omega$ of $\s(\pi)$, cf.\ \ref{spprojform}.
\end{proof}

\medskip
Please note that the above statement does not explicitly involve any specific
cyclic vector.

\begin{remark}
Let $\pi$ be a non-degenerate representation of a commutative
\st-algebra in a Hilbert space $\neq \{0\}$. Then $\pi$ is the direct
sum of cyclic subrepresentations as has been shown in
\ref{directsumdeco}. Hence $\pi$ is spatially equivalent to a direct
sum of spectral representations. Using some care, one can even
give the operators in $R(\pi)$ a representation as multiplication
operators, but then these multiplication operators will in general
live on a space larger than $\s(\pi)$. What is more, the measure
won't be bounded in general. We refer the reader to Mosak
\cite[Thm.\ (7.11) p.\ 104]{MosS}. We shall bypass this construction,
and show instead that in the case of a separable Hilbert space
(the case of interest), we can give the C*-algebra $W^*(\pi)$ a
faithful representation by multiplication operators living on $\s(\pi)$,
see \ref{decospace}. However this new representation won't be
implemented by a unitary operator like in a spatial equivalence;
it will only be a C*-algebra isomorphism. This may not be acceptable
from the operator theoretic point of view, but it may still be so from the
algebraic point of view.
\end{remark}

\clearpage


\section{Separability: the Range}

(This paragraph and the next one can be skipped as a whole.)

In this paragraph, let $H$ be a Hilbert space $\neq \{0\}$.

\begin{introduction}\label{decospace}%
Let $\pi$ be a non-degenerate representation of a commutative
\st-algebra in $H$, and let $P$ be the spectral resolution of $\pi$.
The representation $\pi _P$ factors to an isomorphism of C*-algebras
from $L^{\infty} (P)$ onto $R(\pi _P)$, cf.\ \ref{LPmain}. It shall be
shown in this paragraph that if $H$ is separable, then (on the side of
the \textit{range} of the factored map), the \linebreak C*-algebra $R(\pi _P)$ is
all of the von Neumann algebra $W^*(\pi)$, cf.\ \ref{spprojsubset}.
It shall be shown in the next paragraph that if $H$ is separable, then
(on the side of the \textit{domain} of the factored map), the C*-algebra
$L^{\infty} (P)$ can be identified with the maximal commutative von
Neumann algebra $L^{\infty} ( \langle Px, x \rangle )$ for some unit
vector $x \in H$. Thus, in the case of a separable Hilbert space, the
representation $\pi _P$ factors to a C*-algebra isomorphism of von
Neumann algebras. This C*-algebra isomorphism
$L^{\infty}( \langle P x, x \rangle ) \to W^*(\pi)$ is called the
\underline{$L^{\infty}$ functional calculus}. The inverse of this map
faithfully represents $W^*(\pi)$ by multiplication operators living on
$\s(\pi)$. (So much for the titles of these two paragraphs.)
\end{introduction}

\begin{theorem}%
Let $\pi$ be a non-degenerate representation of a commutative
\st-algebra in $H$. Let $P$ be the spectral resolution of $\pi$,
and let $b \in W^*(\pi)$. Then for any sequence $(x_n)$ in $H$,
there exists a bounded Baire function $f$ on $\s(\pi)$ with
\[ b \,x_n = \pi _P (f) \,x_n \qquad \text{for all } n. \]
\end{theorem}

\begin{proof}
We can assume that $\bigl( \,\| \,x_n \,\| \,\bigr) \in {\ell\,}^2$. Let
$x := \oplus _n x_n \in \oplus _n H$, and put
$M := \overline{R(\oplus _n \pi) \,x}$. Then $M$ is a closed subspace
invariant and cyclic under $\oplus _n \pi$, with cyclic vector $x$,
cf.\ \ref{cyclicsubspace}. Let $Q$ denote the spectral resolution of the cyclic
representation $(\oplus _n \pi )|_M$. Let $S$ denote the restriction of $Q$
to the Baire sets. We have $W^*\bigl((\oplus _n \pi)|_M\bigr) = R(\pi _S)$ by
\ref{spectralrepmain}. By Lemma \ref{lemmabic}, the closed subspace $M$
is invariant under $\oplus _n b$, and
$(\oplus _n b)|_M \in W^*\bigl((\oplus _n \pi)|_M\bigr)$. It follows that
$(\oplus _n b)|_M \in R(\pi _S)$, and so there exists a bounded Baire
function $f$ on $\s(\pi)$ with $(\oplus _n b)|_M = \pi _S (f)$. Since
$\pi _Q = (\oplus _n \pi _P)|_M$ by the uniqueness property of the spectral
resolution, it follows that $b \,x_n = \pi _P (f) \,x_n$ for all $n$.
\end{proof}

\medskip
An immediate consequence is: \pagebreak

\begin{theorem}\label{sepmain}%
Let $\pi$ be a non-degenerate representation of a commutative
\st-algebra in $H$. Let $P$ be the spectral resolution of $\pi$,
and let $Q$ be the restriction of $P$ to the Baire sets. If $H$ is
\underline{separable}, then
\[ W^*(\pi) = R(\pi_P) = R(\pi_Q). \]
\end{theorem}

\begin{definition}[function of a normal bounded operator]%
\index{function!of an operator}\label{Borelfunct}%
Let $b$ be a normal bounded linear operator on $H$. Let $P$ be its
spectral resolution. If $f$ is a bounded Borel function on $\s(b)$, one
writes
\[ f(b) := \pi _P (f) = \int f \,d P \quad \text{(in norm)}, \]
and one says that $f(b)$ is a \underline{bounded Borel function of $b$}.

Please note that the spectral projections $P(\omega)$ with $\omega$
a Borel subset of $\s(b)$, are bounded Borel functions of $b$, as
$P(\omega) = \pi _P ( 1_{\omega} ) = 1_{\omega}(b)$.
\end{definition}

\begin{theorem}[von Neumann]\label{functions}%
\index{Theorem!von Neumann}\index{von Neumann!Theorem}%
Let $b$ be a normal bounded linear operator on $H$. If $H$ is
\underline{separable}, then a bounded linear operator on $H$ is a
bounded Borel function of $b$ if and only if it commutes with
every bounded linear operator in $H$ which commutes with $b$.
\end{theorem}

With a little effort, we can improve on \ref{sepmain}, see \ref{sepspm}
below.

\begin{definition}[Borel \protect\st-algebras \hbox{\cite[4.5.5]{PedC}}]
\index{Borel \protect\st-algebra}\index{algebra!Borel \protect\st-algebra}%
A \underline{Borel \st-algebra} on $H$ is a C*-subalgebra $\mathscr{B}$ of
$B(H)$ such that $\mathscr{B}\sa$ contains the supremum in $B(H)\sa$ of
each increasing sequence in $\mathscr{B}\sa$ which is upper bounded in
$B(H)\sa$. (Cf.\ \ref{ordercompl}.)
\end{definition}

\begin{example}%
When $P$ is a spectral measure, then $R(\pi _P)$ is a commutative Borel
\st-algebra containing $\mathds{1}$. (Use \ref{ordernorm} \& \ref{monconvthm}.)
\end{example}

Since the intersection of Borel \st-algebras on $H$ is a Borel \st-algebra,
questions can be raised on the Borel \st-algebra generated by a set of
operators.

\begin{theorem}\label{QBorel}%
Let $A$ be a commutative C*-subalgebra of $B(H)$ containing $\mathds{1}$.
Let $P$ be the spectral resolution of $A$, and let $Q$ be the restriction of $P$
to the Baire sets. Then $R(\pi _Q)$ is the Borel \st-algebra generated by $A$.
\pagebreak
\end{theorem}

\begin{proof}
By the preceding example, $R(\pi _Q)$ is a Borel
\st-algebra containing $A$. It shall be shown that $R(\pi _Q)$ is the
smallest Borel \st-algebra containing $A$. For this purpose, let
$\mathscr{C}$ be any other Borel \st-algebra containing $R(\pi)$. It then
suffices to show that $\pi _P(f) \in \mathscr{C}$ whenever $f$ is a bounded
Baire function. It is enough to show this for $f$ of the form
$f = 1_{\textstyle{\omega}}$ where $\omega$ is a Baire set. (This is so because
every bounded Baire function is the uniform limit of Baire step functions.) But
then it suffices to show that $\pi _P(f) \in \mathscr{C}$ whenever $f$ is
of the form $f = 1_{\text{\small{$\{ \,g = 0 \,\}$}}}$ for some
$g \in C\bigl(\Delta(A)\bigr)$. This is so because the Baire $\sigma$-algebra
is generated by the sets of the form $\{ \,g = 0 \,\}$ where $g$ runs through
$C\bigl(\Delta(A)\bigr)$, cf.\ the appendix \ref{Bairedef}. One then concludes
with
\[ 1_{\text{\small{$\{ \,g = 0 \,\}$}}} =
1 - \vee _{n \geq 1} \,\bigl( \,1 \wedge n \,| \,g \,| \,\bigr).  \qedhere \]
\end{proof}

\begin{corollary}%
A commutative Borel \st-algebra $\mathscr{B}$ containing $\mathds{1}$
and acting on a separable Hilbert space is a von Neumann algebra.
\end{corollary}

\begin{proof}
Let $P$ be the spectral resolution of $\mathscr{B}$ and let
$Q$ be the restriction of $P$ to the Baire sets. We then have
$W^*(\mathscr{B}) = R(\pi _Q) = \mathscr{B}$,
cf.\ \ref{sepmain} and \ref{QBorel}.
\end{proof}

\begin{theorem}\label{sepspm}%
If $P$ is a spectral measure acting on a separable Hilbert space,
then $R(\pi _P)$ is a von Neumann algebra.
\end{theorem}

\begin{remark}\label{nonsep}
The theorem of von Neumann \ref{functions}, and with it \ref{sepmain}
and \ref{sepspm}, fail to hold in non-separable Hilbert spaces, cf.\ Folland
\cite[page 29]{Foll}.
\end{remark}

Next three miscellaneous results on Borel \st-algebras.

\begin{theorem}\label{PBorel}%
Let $b$ be a normal bounded linear operator on $H$. Let $P$ be its
spectral resolution. Then $R(\pi _P)$ is the Borel  \st-algebra
generated by $b$, $b^*$ and $\mathds{1}$. See also \ref{rangespres}.
\end{theorem}

\begin{proof}
As $\s(b)$ is a metric space, the Borel and Baire $\sigma$-algebras
on $\s(b)$ coincide, cf.\ the appendix \ref{metric}.
\end{proof}

\medskip
The next corollary explains the name ``Borel \st-algebra''. \pagebreak

\begin{corollary}%
Let $\mathscr{B}$ be a Borel \st-algebra containing $\mathds{1}$.
If $b$ is a normal operator in $\mathscr{B}$, then every
bounded Borel function of $b$ lies in $\mathscr{B}$. In particular
the spectral resolution of $b$ takes values in $\mathscr{B}$.
\end{corollary}

\begin{proof}
By the preceding theorem, a bounded Borel function of $b$
lies in the Borel \st-algebra generated by $b,b^*$ and $\mathds{1}$,
which is contained in $\mathscr{B}$.
\end{proof}

\medskip
One says that Borel \st-algebras containing
$\mathds{1}$ have a Borel functional calculus.

\begin{corollary}%
A Borel \st-algebra containing $\mathds{1}$ is
the closed \linebreak linear span of its projections.
\end{corollary}

\begin{proof}
This is proved like \ref{W*main}. 
\end{proof}

\clearpage


\section{Separability: the Domain}%
\label{separating}

(This paragraph can be skipped at first reading.)

Let throughout $H$ denote a Hilbert space $\neq \{0\}$,
and let $A$ be a \linebreak \st-subalgebra of $B(H)$.

\begin{definition}[separating vectors and cyclic vectors]%
\index{vector!separating}\index{vector!cyclic}\index{separating!vector}%
A vector $x \in H$ is called \underline{separating} for $A$,
whenever $ax = 0$ for some $a \in A$ implies $a = 0$. A vector $x \in H$
is called \underline{cyclic} under $A$ if $Ax$ is dense in $H$, and in this
case the \st-algebra $A$ is called cyclic.
\end{definition}

\begin{proposition}%
A vector $x \in H$ is cyclic under $A$ if and only if it is separating for $A'$.
\end{proposition}

\begin{proof}
Assume that $x$ is cyclic under $A$, and let $b \in A'$ with $bx = 0$.
Then for $a \in A$, we have $bax = abx = 0$, whence $b = 0$ as
$Ax$ is dense in $H$. Next assume that $x$ is separating for $A'$.
With $p$ denoting the projection on the closed subspace $\overline{Ax}$,
we have $p \in A'$ by \ref{invarcommutant}, as well as $px = x$ by
\ref{cyclicsubspace}. Thus $\mathds{1} - p \in A'$ and
$( \mathds{1} - p ) x = 0$, which implies $\mathds{1} - p = 0$ as $x$ is
separating for $A'$. Hence $p = \mathds{1}$, so $x$ is cyclic under $A$.
\end{proof}

\begin{corollary}
If $A$ is a commutative cyclic \st-subalgebra of $B(H)$ with
cyclic vector $c$, then $c$ also is a separating vector for $A$.
\end{corollary}

\begin{proof}
Since $A$ is commutative, we have $A \subset A'$.
\end{proof}

\begin{theorem}\label{commsepsep}%
If $A$ is commutative and if $H$ is separable, then $A$ has a separating
vector.
\end{theorem}

\begin{proof}
Let $A$ be commutative and $H$ separable. Then $H$ is the
\linebreak countable direct sum of closed subspaces $H_n$
$(n \geq 1)$ invariant \linebreak under $A$, and such that for
each $n \geq 1$, the subspace $H _n$ contains a unit cyclic
vector $c_n$ under $A|_{H_n}$ cf.\ \ref{directsumdeco}. It shall
be shown that the vector $x := \sum _{n \geq 1} {2}^{-n/2} c_n \in H$
is a separating vector for $A$. So let $a \in A$ with $ax = 0$.
With $p_n$ denoting the projection on the closed invariant subspace
$H _n$, we have $p_n \in A'$, cf.\ \ref{invarcommutant}, and so
\[ a A c_n = a A p_n x = A a p_n x = A p_n a x= \{0\}, \]
which implies that $a$ vanishes on $H _n$, hence everywhere.
We have shown that $x$ is separating for $A$. \pagebreak
\end{proof}

\begin{theorem}\label{sepmaxcomm}%
Let $\pi$ be a non-degenerate representation of a \linebreak \st-algebra
in $H$. Assume that the \st-algebra is commutative, and that $H$ is
separable. Then $\pi$ is cyclic if and only if $\pi'$ is commutative.

This gives a characterisation of cyclicity in terms of only the \linebreak
algebraic structure of the commutant.
\end{theorem}

\begin{proof}
The ``only if'' part holds without the assumption of separability on $H$,
cf.\ \ref{spectralrepmain}. The ``if'' part holds without the assumption
of commutativity on the \st-algebra, as we shall show next. If $\pi'$ is
\linebreak
commutative, then since $H$ is separable, $\pi'$ has a separating
vector, which will be cyclic under $\pi$.
\end{proof}

\begin{definition}[spatial equivalence]%
\index{equivalent!spatially!von Neumann algebra}%
\index{spatially equivalent!von Neumann algebra}%
Two von Neumann algebras $\mathscr{M}_1$ and $\mathscr{M}_2$
acting on Hilbert spaces $H_1$ and $H_2$ respectively, are called
\underline{spatially equivalent}, if there exists a unitary operator
$U : H_1 \to H_2$ such that $U \mathscr{M}_1 U^{-1} = \mathscr{M}_2$.
\end{definition}

\begin{theorem}\label{maxcommchar}%
For a commutative von Neumann algebra $\mathscr{M}$ on a
\underline{separable} Hilbert space $\neq \{0\}$, the following
statements are equivalent.
\begin{itemize}
   \item[$(i)$] $\mathscr{M}$ is maximal commutative,
  \item[$(ii)$] $\mathscr{M}$ is cyclic,
 \item[$(iii)$] $\mathscr{M}$ is spatially equivalent to the maximal
                       commutative von Neumann algebra $L^{\infty}(\mu)$
                       for some inner regular Borel probability measure $\mu$
                       on a locally compact Hausdorff space $\neq \varnothing$.
\end{itemize}
\end{theorem}

\begin{proof}
(i) $\Rightarrow$ (ii): \ref{maxcommprime} and \ref{sepmaxcomm}.
(ii) $\Rightarrow$ (iii): \ref{Abelian}, \ref{Bochnerbis} and \ref{exspectres}.
(iii) $\Rightarrow$ (i): \ref{maxcommL}.
\end{proof}

\medskip
We see that on separable Hilbert spaces $\neq \{0\}$, the
prototypes of the maximal commutative von Neumann algebras
are the maximal commutative von Neumann algebras $L^{\infty}(\mu)$
with $\mu$ an inner regular Borel probability measure on a locally
compact Hausdorff space.

\begin{definition}[scalar spectral measure]%
\index{measure!scalar spectral}\index{spectral!measure!scalar}%
\index{scalar spectral measure}%
Let $P$ be a spectral measure, defined on a $\sigma$-algebra
$\mathcal{E}$. Then a probability measure $\mu$, defined on
$\mathcal{E}$, is called a \underline{scalar spectral measure}
for $P$, if for $\omega \in \mathcal{E}$, one has $P(\omega) = 0$
precisely when $\mu(\omega) = 0$. That is, $P$ and $\mu$ should
be mutually absolutely continuous. \pagebreak
\end{definition}

We restrict the definition of scalar spectral measures to probability
measures for the reason that the applications - foremostly the crucial
item \ref{identify} below - require bounded measures. 

The following two results deal with the existence of scalar spectral measures.

\begin{proposition}[separating vectors for spectral measures]%
\index{separating!vector}\index{vector!separating}%
Let $P$ be a spectral measure, defined on a $\sigma$-algebra
$\mathcal{E}$, and acting on $H$. Then for a unit vector $x \in H$,
the probability measure $\langle P x, x \rangle$ is a scalar spectral
measure for $P$ if and only if $x$ is \underline{separating} for $P$.
(In the obvious extended sense that for $\omega \in \mathcal{E}$,
one has $P(\omega)x = 0$ only when $P(\omega) = 0$.)
\end{proposition}

\begin{proof}
For $x \in H$ and $\omega \in \mathcal{E}$, we have
\[ {\| \,P(\omega)x \,\| \,}^2 = \langle P x , x \rangle (\omega), \]
so $P(\omega)x = 0$ if and only if $\langle P x , x \rangle (\omega) = 0$.
Therefore $x$ is separating for $P$ if and only if $\omega \in \mathcal{E}$
and $\langle P x , x \rangle (\omega) = 0$ imply $P(\omega) = 0$,
i.e.\ if and only if $\langle P x , x \rangle$ is a scalar spectral measure. 
\end{proof}

\begin{theorem}\label{exscalspm}%
Let $P$ be a spectral measure acting on $H$. If $H$ is separable,
then there exists a unit vector $x \in H$ such that the probability
measure $\langle P x, x \rangle$ is a scalar spectral measure for $P$.
The vector $x$ can be chosen to be any unit separating vector for $P$,
for example any unit separating vector for $R(\pi_P)$.
\end{theorem}

\begin{proof}
The range $R(\pi_P)$ has a separating vector by \ref{commsepsep}.
\end{proof}

\medskip
The next two results exhibit the purpose of scalar spectral measures.

\begin{theorem}\label{identify}%
Let $P$ be a spectral measure, defined on a $\sigma$-algebra
$\mathcal{E}$. If $P$ admits a probability measure $\mu$, defined
on $\mathcal{E}$, as a scalar spectral measure, then we may identify
the C*-algebras $L^{\infty}(P)$ and $L^{\infty}(\mu)$.
\end{theorem}

\begin{proof}
Assume that some probability measure $\mu$, defined on
$\mathcal{E}$, is a scalar spectral measure for $P$. A function $f$
in ${\mathscr{L} \,}^{\infty}(\mu)$ then is \linebreak $\mu$-a.e.\ equal
to a function $g \in \mathcal{M}_b(\mathcal{E})$, cf.\ the appendix
\ref{general}, and the function $f$ is $\mu$-a.e.\ zero if and only if
the function $g$ is $P$-a.e.\ zero.\pagebreak%
\end{proof}

\begin{theorem}\label{scalpurpose}%
Let $P$ be a spectral measure, defined on a $\sigma$-algebra
$\mathcal{E}$. If $P$ admits a probability measure $\mu$,
defined on $\mathcal{E}$, as a scalar spectral measure, then
$\pi_P$ factors to a C*-algebra isomorphism from $L^{\infty}(\mu)$
onto $R(\pi_P)$.
\end{theorem}

\begin{proof}
This follows now from \ref{LPmain}.
\end{proof}

\begin{theorem}%
Let $P$ be a spectral measure acting on $H$. If $H$ is separable,
there exists a unit vector $x \in H$ such that $\pi_P$ factors to a
C*-algebra isomorphism from $L^{\infty} ( \langle P x, x \rangle )$
onto $W^*(P)$. The vector $x$ can be chosen to be any unit vector
in $H$, such that the probability measure $\langle P x, x \rangle$ is
a scalar spectral measure for $P$, for example any unit separating
vector for $P$, or even $W^*(P)$.
\end{theorem}

\begin{proof}
Since $H$ is separable, we know from \ref{sepspm} that $R(\pi_P)$
is a von Neumann algebra, namely $W^*(\pi_P)$, which equals
$W^*(P)$ by \ref{spprojcomm}. The statement follows now from
\ref{exscalspm} and \ref{scalpurpose}.
\end{proof}

\begin{theorem}[the $L^{\infty}$ functional calculus]%
Let $\pi$ be a non-degen\-erate representation of a
commutative \st-algebra in $H$, and let $P$ be the spectral
resolution of $\pi$. If $H$ is separable, there exists a unit
vector $x \in H$ such that the probability measure
$\langle P x, x \rangle$ is a scalar spectral measure for $P$,
and in this case $\pi_P$ factors to a C*-algebra isomorphism
from $L^{\infty}(\langle P x, x \rangle)$ onto $W^*(\pi)$. The
vector $x$ can be chosen to be any unit separating vector for
$P$, for example any unit separating vector for $W^*(\pi)$.
\end{theorem}

\begin{proof}
This follows from \ref{sepmain}, \ref{exscalspm} and \ref{scalpurpose}.
The result also follows from the preceding one together with \ref{gensame}.
\end{proof}

\medskip
We also obtain the following representation theorem
for commut\-ative von Neumann algebras on separable
Hilbert spaces $\neq \{0\}$.

\begin{theorem}\label{vNCstar}%
A commutative von Neumann algebra $\mathscr{M}$ on a
\linebreak separable Hilbert space $\neq \{0\}$ is isomorphic
as a C*-algebra to $L^{\infty}(\mu)$ for some inner regular
Borel probability measure $\mu$ on the compact Hausdorff
space $\Delta(\mathscr{M})$.
The measure $\mu$ can be chosen to be any inner regular
Borel probability measure $\mu$ on $\Delta(\mathscr{M})$
which is a scalar spectral measure for the spectral resolution
of $\mathscr{M}$. In this case, the support of $\mu$ is all of
$\Delta(\mathscr{M})$, and every function in
${\mathscr{L} \,}^{\infty}(\mu)$ is $\mu$-a.e.\ equal to a unique
continuous function on $\Delta(\mathscr{M})$. \pagebreak
\end{theorem}

\begin{proof}
See \ref{disconn} and the appendix \ref{Baire}.
\end{proof}

\medskip
Please note that here we only have a C*-algebra isomorphism,
while in \ref{maxcommchar}, we have a spatial equivalence
of von Neumann algebras, which is a stronger condition.

\clearpage


\chapter{Single Unbounded Self-Adjoint Operators}%
\label{unbdself}

\setcounter{section}{40}


\section{Spectral Integrals of Unbounded Functions}

In this paragraph, let $\mathcal{E}$ be a $\sigma$-algebra on a set
$\Omega \neq \varnothing$, let $H \neq \{0\}$ be a Hilbert space, and let
$P$ be a spectral measure defined on $\mathcal{E}$ and acting on $H$.

We first have to extend the notion of ``linear operator''.

\begin{definition}[linear operator in $H$]\index{operator}%
By a \underline{linear operator in $H$} we mean a linear operator
defined on a subspace of $H$, taking values in $H$. We stress that
we then say: linear operator \underline{in} $H$; not \underline{on} $H$.
The domain (of definition) of such a linear operator $a$ is denoted
by $D(a)$.\index{D2@$D(a)$}
\end{definition}

\begin{definition}[$\mathcal{M(E)}$]%
\index{M3@$\mathcal{M(E)}$}%
We shall denote by $\underline{\mathcal{M(E)}}$ the set of
complex-valued $\mathcal{E}$-measurable functions on $\Omega$,
cf.\ the appendix, \ref{imagebegin}.
\end{definition}

We shall associate with each $f \in \mathcal{M(E)}$ a linear operator
in $H$, defined on a dense subspace $\mathcal{D}_f$ of $H$.
Next the definition of $\mathcal{D}_f$.

\begin{proposition}[$\mathcal{D}_f$]\index{D3@$\mathcal{D}_f$}%
Let $f \in \mathcal{M(E)}$. One defines
\[ \underline{\mathcal{D}_f} :=
\{ \,x \in H : f \in {\mathscr{L} \,}^{2}(\langle Px,x \rangle) \,\}. \]
The set $\mathcal{D}_f$ is a dense subspace of $H$.
\end{proposition}

\begin{proof}
For $x,y \in H$ and $\omega \in \mathcal{E}$, we have
(the parallelogram identity):
\[ {\| \,P(\omega)(x+y) \,\| \,}^2 + { \,\| \,P(\omega)(x-y) \,\| \,}^2
= 2 \,\bigl( \,{\| \,P(\omega)x \,\| \,}^2 + {\| \,P(\omega)y \,\| \,}^2 \,\bigr), \]
whence
\[ \langle P(\omega)(x+y),(x+y) \rangle
\leq 2 \,\bigl( \,\langle P(\omega)x,x \rangle + \langle P(\omega)y,y \rangle \,\bigr),
\pagebreak \]
so that $\mathcal{D}_f$ is closed under addition. It shall next be
proved that $\mathcal{D}_f$ is dense in $H$. For $n \geq 1$ let
$\omega_n := \{ \,| f \,| < n \,\}$. We have $\Omega = \cup _n \omega _n$
as $f$ is complex-valued. Let $\omega \in \mathcal{E}$.
For $z$ in the range of $P(\omega _n)$, we have
\[ P(\omega) z = P(\omega) P(\omega _n) z = P(\omega \cap \omega_n) z, \]
so that
\[ \langle Pz,z \rangle (\omega) = \langle Pz,z \rangle (\omega \cap \omega _n), \]
and therefore
\[ \int {| \,f \,| \,}^2 \,d \langle Pz,z \rangle
= \int {| \,f \,| \,}^2 \,1_{\textstyle\omega _n} \,d \langle Pz,z \rangle
\leq {n \,}^2 \,{\| \,z \,\| \,}^2 < \infty. \]
This shows that the range of $P(\omega _n)$ is contained in
$\mathcal{D}_f$. We have
\[ y = \lim _{n \to \infty} P(\omega _n) y \]
for all $y \in H$ by $\Omega = \cup _n \omega _n$ and by \ref{strsgcont}.
This says that $\mathcal{D}_f$ is dense in $H$. 
\end{proof}

\begin{definition}[spectral integrals]%
Consider a linear  operator $a : H \supset D(a) \to H$, defined
on a subspace $D(a)$ of $H$. Let $f \in \mathcal{M(E)}$. We write
\[ a = \int f \,d P \qquad \text{(pointwise)} \]
if $D(a) = \mathcal{D}_f$ and if for all $x \in \mathcal{D}_f$ and all
$h \in \mathcal{M}_b(\mathcal{E})$ one has
\[ \| \,\bigl(a - \pi _P(h)\bigr) x \,\| 
= \| \,f-h \,\|_{\text{\small{$\langle Px,x \rangle$,2}}}
= {\biggl(\int {| \,f-h \,| \,}^2 \ d \langle Px,x \rangle \biggr)}^{1/2}. \]
One writes
\[ a = \int f \,d P \qquad \text{(weakly)} \]
if $D(a) = \mathcal{D}_f$ and
\[ \langle ax,x \rangle = \int f \,d \langle Px,x \rangle \]
for all $x \in \mathcal{D}_f$. Then $a$ is
uniquely determined by $f$ and $P$. \pagebreak
\end{definition}

\begin{proposition}%
Let $f \in \mathcal{M(E)}$ and let $a$ be a linear operator in $H$. Assume that
\[ a = \int f \,d P \qquad \text{(pointwise)}. \]
We then also have
\[ a = \int f \,d P \qquad \text{(weakly)}. \]
In particular, $a$ is uniquely determined by $f$ and $P$.
\end{proposition}

\begin{proof}
For $n \geq 1$, put
$f_n := f \,1_{\text{\small{$\{ \,| \,f \,| < n \,\}$}}}
\in \mathcal{M}_b(\mathcal{E})$.
For $x \in \mathcal{D}_f$, we get
\[ ax = \lim _{n \to \infty} \pi _P(f_n)x \]
by the Lebesgue Dominated Convergence Theorem,
whence, by the same theorem,
\[ \langle ax,x \rangle
= \lim _{n \to \infty} \langle \pi _P(f_n)x,x \rangle
= \lim _{n \to \infty} \int f_n \,d \langle Px,x \rangle
= \int f \,d \langle Px,x \rangle. \qedhere \]
\end{proof}

\begin{theorem}[$\Psi_P(f)$]\index{P99@$\Psi _P(f)$}%
For $f \in \mathcal{M(E)}$ there exists a (necessarily unique)
linear operator $\underline{\Psi _P(f)}$ in $H$ such that
\[ \Psi _P(f) = \int f \,d P \qquad \text{(pointwise)}. \]
\end{theorem}

\begin{proof}
Let $x \in \mathcal{D}_f$ be fixed. Look at
$\mathcal{M}_b(\mathcal{E})$ as a subspace of \linebreak
${\mathscr{L} \,}^{2}(\langle Px,x \rangle)$.
Consider the map
\[ I_x : {\mathscr{L} \,}^{2}(\langle Px,x \rangle) \supset
\mathcal{M}_b(\mathcal{E}) \to H,\quad h \mapsto \pi _P(h) x. \]
Since $\mathcal{M}_b(\mathcal{E})$ is dense in
${\mathscr{L} \,}^{2}(\langle Px,x \rangle)$, and since $I_x$ is an
isometry by \ref{normpointwise}, it follows
that $I_x$ has a unique continuation to an iso\-metry
$L^2(\langle Px,x \rangle) \to H$, which shall also be denoted
by $I_x$. Putting \linebreak $\Psi _P(f)x := I_x(f)$, we obtain that
$\Psi _P(f)$ is an operator in $H$ with domain $\mathcal{D}_f$.
It shall be proved that
\[ \Psi _P(f) = \int f \,d P \qquad \text{(pointwise)}. \]
Let $x \in \mathcal{D}_f$. For $h \in \mathcal{M}_b(\mathcal{E})$,
we then have
\[ \| \,\bigl( \Psi _P(f) - \pi _P(h)\bigr)x \,\|
= \| \,I_x(f-h) \,\|
= \| \,f-h \,\|_{\text{\small{$\langle Px,x \rangle,2$}}} \]
because $I_x$ is an isometry. Also $\Psi _P(f)$ is linear because it is
pointwise approximated by linear operators $\pi _P(h)$. \pagebreak
\end{proof}

\begin{corollary}%
Let $f \in \mathcal{M(E)}$. Let $a$ be a linear operator in $H$ such that
\[ a = \int f \,d P \qquad \text{(weakly)}. \]
We then also have
\[ a = \int f \,d P \qquad \text{(pointwise)}. \]
\end{corollary}

\begin{proof}
We have $a = \Psi _P(f)$ by $\Psi _P(f) = \int f \,dP$ (weakly) as well.
\end{proof}

\begin{proposition}\label{heqnull}%
Let $f \in \mathcal{M(E)}$. For $x \in \mathcal{D}_f$, we then have
\[ \| \,\Psi _P (f) \,x \,\| 
= \| \,f \,\|_{\text{\small{$\langle Px,x \rangle,2$}}}
= {\biggl(\int {| \,f \,| \,}^2 \ d \langle Px,x \rangle \biggr)}^{1/2}. \]
\end{proposition}

\begin{proof}
Put $h := 0$ in $\Psi _P(f) = \int f \,d P$ (pointwise).
\end{proof}

\begin{lemma}\label{Dlemma}%
Let $g \in \mathcal{M(E)}$ and $x \in \mathcal{D}_g$.
For $y := \Psi _P(g)x$, we then have
\[ \langle Py,y \rangle = {| \,g \,| \,}^2 \cdot \langle Px,x \rangle. \]
Moreover, for $h \in \mathcal{M}_b(\mathcal{E})$, we have
\[ \pi _P(h) y = \Psi _P(hg)x. \]
\end{lemma}

\begin{proof}
With $g_n := g \,1_{\text{\small$\{ \,| \,g \,| < n \,\}$}}$
for all $n \geq 1$, we have
\[ y = \Psi _P(g)x = \lim _{n \to \infty} \pi _P(g_n)x. \]
For $h \in \mathcal{M}_b(\mathcal{E})$, we obtain
\begin{align*}
\pi _P(h)y = \pi _P(h) \Psi _P(g)x
& = \lim _{n \to \infty} \pi _P(h) \pi _P(g_n)x \\
& = \lim _{n \to \infty} \pi _P(h g_n)x = \Psi _P(hg)x.
\end{align*}
The last equality follows from $\mathcal{D}_{hg} \supset \mathcal{D}_g$
by boundedness of $h$. In particular, from \ref{heqnull} we get
\[ \int {| \,h \,| \,}^2 \,d \langle Py,y \rangle = {\| \,\pi_P(h)y \,\| \,}^2
= {\| \,\Psi_P(hg)x \,\| \,}^2 = \int {| \,hg \,| \,}^2 \,d \langle Px,x \rangle \]
for all $h \in \mathcal{M}_b(\mathcal{E})$. This means that
\[ \langle Py,y \rangle
= {| \,g \,| \,}^2 \cdot \langle Px,x \rangle. \pagebreak \qedhere \]
\end{proof}

We need some notation:

\begin{definition}[extension]%
If $a$ and $b$ are linear operators in $H$ with respective
domains $D(a)$ and $D(b)$, one writes $a \subset b$,
if $b$ is an \underline{extension} of $a$, i.e.\ if $D(a) \subset D(b)$
and $ax = bx$ for all $x \in D(a)$.
\end{definition}

\begin{theorem}[Addition and Multiplication Theorems]\label{AddMult}%
\index{Theorem!Addition \& Multiplication}\index{Addit.\ \& Multiplicat.\ Thm.}%
For $f,g$ in $\mathcal{M(E)}$ we have
\begin{gather*}
\Psi _P(f) + \Psi _P(g) \subset \Psi _P(f+g), \\
\Psi _P(f) \,\Psi _P(g) \subset \Psi _P(fg).
\end{gather*}
The precise domains are given by
\begin{gather*}
D\bigl(\Psi _P(f) + \Psi _P(g)\bigr) = \mathcal{D}_{| \,f \,| + | \,g \,|}, \\
D\bigl(\Psi _P(f) \,\Psi _P(g)\bigr) = \mathcal{D}_{fg} \cap \mathcal{D}_g.
\end{gather*}
\end{theorem}

\begin{proof}
We shall first prove the part of the statement relating to addition.
Let $f,g \in \mathcal{M(E)}$. If $x \in D\bigl(\Psi _P(f) + \Psi _P(g)\bigr)
:= D\bigl(\Psi _P(f)\bigr) \cap D\bigl(\Psi _P(g)\bigr)$, then
$f, g \in {\mathscr{L} \,}^{2}(\langle Px,x \rangle)$, whence
$f + g \in {\mathscr{L} \,}^{2}(\langle Px,x \rangle)$, so that
$x \in \mathcal{D}_{f+g}$. It follows that
\[ \Psi _P(f+g)x = \Psi _P(f)x+\Psi _P(g)x \]
because $\Psi _P$ is pointwise approximated by the linear map $\pi _P$.
We have
\[ D\bigl(\Psi _P(f)+\Psi _P(g)\bigr) = \mathcal{D}_{| \,f \,|+| \,g \,|}. \]
Indeed, if $x \in H$, then for the complex-valued $\mathcal{E}$-measurable
functions $f$ and $g$, we have that both
$f, g \in {\mathscr{L} \,}^{2}(\langle Px,x \rangle)$
if and only if $| \,f \,|+| \,g \,| \in {\mathscr{L} \,}^{2}(\langle Px,x \rangle)$.

We shall now prove the part of the statement relating to composition.
Let $g \in \mathcal{M(E)}$, $x \in \mathcal{D}_g$, and put $y := \Psi _P(g)x$.
Let $f \in \mathcal{M(E)}$. Lemma \ref{Dlemma} shows that
$y \in \mathcal{D}_f$ if and only if $x \in \mathcal{D}_{fg}$, and so
\[ D\bigl(\Psi _P(f) \,\Psi _P(g)\bigr) = \mathcal{D}_{fg} \cap \mathcal{D}_g
\qquad (f,g \in \mathcal{M(E)}). \]
For $n \geq 1$ define $f_n := f \,1_{\text{\small{$\{ \,| \,f \,| < n \,\}$}}}$.
If $x \in \mathcal{D}_{fg} \cap \mathcal{D}_g$, then
$f_n g \to fg$ in ${\mathscr{L} \,}^{2}(\langle Px,x \rangle)$ and so
$f_n \to f$ in ${\mathscr{L} \,}^{2}(\langle Py,y \rangle)$ by \ref{Dlemma}.
Applying \ref{Dlemma} with $h := f_n$, we obtain
\begin{align*}
\Psi _P(f) \Psi _P(g) x = \Psi _P(f) y & = \lim _{n \to \infty} \pi _P(f_n) y \\
& = \lim _{n \to \infty} \Psi _P(f_n g) x = \Psi _P(fg) x. \pagebreak \qedhere
\end{align*}
\end{proof}

\begin{definition}[self-adjoint linear operators]%
\index{self-adjoint!operator}\index{operator!self-adjoint}%
\index{adjoint!operator}\index{operator!adjoint}%
\index{symmetric operator}\index{operator!symmetric}%
Let $a$ be a linear operator in $H$ defined on a dense
subspace $D(a)$. One defines
\begin{align*}
D(a^*) := \{ \,y \in H : \ & \text{there exists } y^* \in H \text{ with} \\
                           & \langle ax,y \rangle = \langle x,y^* \rangle
                             \text{ for all } x \in D(a) \,\},
\end{align*}
which is a subspace of $H$. The associated $y^*$ is determined
uniquely because $D(a)$ is dense in $H$. Thus
\[ a^*y := y^* \qquad \bigl(y \in D(a^*)\bigr) \]
determines a linear operator $a^*$ with domain $D(a^*)$ in $H$.
The operator $a^*$ is called the \underline{adjoint} of $a$. The
operator $a$ is called \underline{self-adjoint} if $a = a^*$, and
\underline{symmetric} if $a \subset a^*$.
\end{definition}

\begin{definition}[closed operators]%
\index{closed operator}\index{operator!closed}%
A linear operator $a$ in $H$ with domain $D(a) \subset H$,
is called \underline{closed}, if its graph is closed in $H \times H$.
\end{definition}

\begin{proposition}\label{adjclosed}%
The adjoint of a densely defined linear operator in $H$ is closed.
\end{proposition}

\begin{proof}
Let $a$ be a densely defined linear operator in $H$.
Let $(y_n)$ be a sequence in $D(a^*)$ such that
$y_n \to y$, and $a^*y_n \to z$ in $H$. For $x \in D(a)$ we get
\[ \langle ax,y \rangle = \lim _{n \to \infty} \langle ax,y_n \rangle
   = \lim _{n \to \infty} \langle x,a^*y_n \rangle = \langle x,z \rangle, \]
whence $y \in D(a^*)$ and $a^*y = z$.
\end{proof}

\begin{corollary}\label{selfadjclosed}%
A self-adjoint linear operator in $H$ is closed.
\end{corollary}

\begin{proposition}\label{submult}%
If $a,b$ are densely defined linear operators in $H$,
and if $ab$ is densely defined as well, then
\[ b^* a^* \subset (ab)^*. \]
\end{proposition}

\begin{proof}
Let $x \in D(ab)$ and $y \in D(b^*a^*)$. We get
\[ \langle abx,y \rangle = \langle bx,a^*y \rangle \]
because $bx \in D(a)$ and $y \in D(a^*)$. Similarly we have
\[ \langle bx,a^*y \rangle = \langle x,b^*a^*y \rangle \]
because $x \in D(b)$ and $a^*y \in D(b^*)$. It follows that
\[ \langle abx,y \rangle = \langle x,b^*a^*y \rangle, \]
whence $y \in D\bigl((ab)^*\bigr)$ and $(ab)^*y = b^*a^*y$.
\pagebreak
\end{proof}

\begin{theorem}%
For $f \in \mathcal{M(E)}$, the following statements hold:
\begin{itemize}
   \item[$(i)$] ${\Psi _P(f) \,}^* = \Psi _P(\overline{f})$,
  \item[$(ii)$] $\Psi _P(f)$ is closed,
 \item[$(iii)$] ${\Psi _P(f) \,}^* \,\Psi _P(f) = \Psi _P \bigl( \,{| \,f \,| \,}^2 \,\bigr)
= \Psi _P(f) \,{\Psi _P(f) \,}^*$.
\end{itemize}
\end{theorem}

\begin{proof}
(i): For $x \in \mathcal{D}_f$, and
$y \in \mathcal{D}_f = \mathcal{D}_{\overline{f}}$ we have
\begin{align*}
\langle \Psi _P(f)x,y \rangle = \int f \,d \langle Px,y \rangle
& = \overline{\int \overline{f} \,d \langle Py,x \rangle} \\
& = \overline{\langle \Psi _P(\overline{f}) y,x \rangle}
= \langle x,\Psi _P(\overline{f})y \rangle,
\end{align*}
so that $y \in D\bigl({\Psi _P(f) \,}^*\bigr)$ and
$\Psi _P(\overline{f}) \subset {\Psi _P(f) \,}^*$.
In order to prove the reverse inclusion, consider the truncations
$h_n := 1_{\text{\small{$\{ \,| \,f \,| < n \,\}$}}}$ for $n \geq 1$.
The Multiplication Theorem \ref{AddMult} gives
\[ \Psi _P(f) \,\pi _P(h_n) = \pi _P (f h_n). \]
One concludes with the help of \ref{submult} that
\[ \pi _P(h_n) \,{\Psi _P(f) \,}^*
\subset {[ \,\Psi _P(f) \,\pi _P(h_n) \,] \,}^*
= {\pi _P(f h_n) \,}^* = \pi _P(\overline{f} h_n). \]
For $z \in D \,\bigl( \,{\Psi _P(f) \,}^* \,\bigr)$ and $v = {\Psi _P(f) \,}^* \,z$ it follows
\[ \pi _P(h_n)v = \pi _P(\overline{f} h_n)z. \]
Hence
\[ \int {| \,\overline{f} h_n \,| \,}^2 \,d \langle Pz,z \rangle
= \int h_n \,d \langle Pv,v \rangle \leq \langle v,v \rangle \]
for all $n \geq 1$, so that $z \in \mathcal{D}_{\overline{f}}$.

(ii) follows from (i) together with \ref{adjclosed}.

(iii) follows now from the Multiplication Theorem  \ref{AddMult}
because $\mathcal{D}_{\overline{f}f} \subset \mathcal{D}_f$
by the Cauchy-Schwarz inequality:
\[ {\left| \,\int {| \,f \,| \,}^2 \,d \langle Px,x \rangle \,\right| \,}^2
\leq \int {| \,f \,| \,}^4 \,d \langle Px,x \rangle
\cdot \int {1 \,}^2 \,d \langle Px,x \rangle. \qedhere \]
\end{proof}

\bigskip\begin{corollary}%
If $f$ is a real-valued function in $\mathcal{M(E)}$,
then $\Psi _P(f)$ is a self-adjoint linear operator in $H$.
\end{corollary}

\clearpage


\section{The Spectral Theorem for Self-Adjoint Operators}

In this paragraph, let $a : H \supset D(a) \to H$ be a
self-adjoint linear operator in a Hilbert space $H \neq \{0\}$,
defined on a dense subspace $D(a)$ of $H$.

\begin{definition}[the spectrum]\index{spectrum!of an operator}%
One says that $\lambda \in \mathds{C}$ is a regular value of the
self-adjoint linear operator $a$, if $\lambda \mathds{1}- a$ is injective
as well as surjective onto $H$, and if its left inverse is bounded.
The \underline{spectrum} $\s(a)$ of $a$ is defined as the set of
those complex numbers which are not regular values of $a$.
\end{definition}

\begin{theorem}\label{specisreal}%
The spectrum of the self-adjoint operator $a$ is real.
\end{theorem}

\begin{proof}
Let $\lambda \in \mathds{C} \setminus \mathds{R}$
and let $\lambda =: \alpha + \iu \beta$ with $\alpha$ and $\beta$ real.
For $x \in D(a)$, we get
\[ {\| \,\bigl( ( \alpha \mathds{1} - a ) + \iu \beta \mathds{1} \bigr) x \,\| \,}^2
= {\| \,( \alpha \mathds{1} - a ) x \,\| \,}^2 + {\beta \,}^2 \,{\| \,x \,\| \,}^2
\geq {\beta \,}^2 \,{\| \,x \,\| \,}^2, \]
which shows that $\lambda \mathds{1} - a$ has a bounded left inverse. It remains
to be shown that this left inverse is everywhere defined. It shall first be shown that
the range $R(\lambda \mathds{1} -a)$ is closed. Let $(y_n)$ be a sequence in
$R(\lambda \mathds{1} - a)$ which converges to an element $y \in H$. There then
exists a sequence $(x_n)$ in $D(a)$ such that $y_n = (\lambda \mathds{1} - a)x_n$
for all $n$. Next we have
\[ \| \,x_n-x_m \,\| \,\leq \,{| \,\beta \,| \,}^{-1}
\cdot \,\| \,( \lambda \mathds{1} - a ) (x_n-x_m) \,\|
\,= \,{| \,\beta \,| \,}^{-1} \cdot \,\| \,y_n-y_m \,\|, \]
so that $(x_n)$ is Cauchy and thus converges to some $x \in H$. Since $a$ is closed
\ref{selfadjclosed}, this implies that $x \in D(a)$ and $y = ( \lambda \mathds{1} - a ) x$,
so that $y \in R(\lambda \mathds{1}-a)$, which shows that $R(\lambda \mathds{1}-a)$
is closed. It shall be shown that $R(\lambda \mathds{1} - a) = H$. Let $y \in H$ be
orthogonal to $R(\lambda \mathds{1} - a)$. We have to show that $y = 0$. We get
\[ \langle (\lambda \mathds{1}-a)x,y \rangle = 0 = \langle x,0 \rangle
\quad \text{for all } x \in D(a), \]
which means that $(\lambda \mathds{1} - a)^*y = (\overline{\lambda}\mathds{1}-a)y$
exists and equals zero. It follows that $y = 0$ because
$\overline{\lambda}\mathds{1}-a$ is injective by the above.
\end{proof}

\begin{definition}[the Cayley transform]\index{Cayley transform}%
The operator
\[ u := (a - \iu \mathds{1}) (a + \iu \mathds{1})^{-1} \]
is called the \underline{Cayley transform} of $a$. \pagebreak
\end{definition}

\begin{proposition}\label{onenoteival}%
The Cayley transform $u$ of $a$ is unitary and $1$ is not an
eigenvalue of $u$. 
\end{proposition}

\begin{proof} For $x \in D(a)$ we have
\[ {\| \,(a - \iu \mathds{1})x \,\| \,}^2 = {\| \,ax \,\| \,}^2 + {\| \,x \,\| \,}^2
= {\| \,(a + \iu \mathds{1})x \,\| \,}^2, \]
so that
\[ \| \,(a - \iu \mathds{1})(a + \iu \mathds{1})^{-1}y \,\| = \| \,y \,\|
\quad \text{for all } y \in D(u), \]
which shows that $u$ is isometric. It follows from the preceding theorem
that $u$ is bijective, hence unitary. Next, let $x \in D(a)$ and
$y := (a + \iu \mathds{1})x$. We then have $uy = (a - \iu \mathds{1})x$.
By subtraction we obtain $(\mathds{1}-u)y = 2 \iu x$, so that if
$(\mathds{1}-u)y = 0$, then $x = 0$, whence $y = 0$, which shows
injectivity of $\mathds{1}-u$.
\end{proof}

\begin{proposition}\label{spau}%
Let
\[ u := (a - \iu \mathds{1}) (a + \iu \mathds{1})^{-1} \]
be the Cayley transform of $a$.
For $\lambda \in \mathds{R}$, consider
\[ \nu := (\lambda - \iu) (\lambda + \iu)^{-1}. \]
Then $\lambda \in \s(a)$ if and only if $\nu \in \s(u)$.
\end{proposition}

\begin{proof}
We have
\begin{align*}
 &\ u-\nu \mathds{1}\\
= &\ \left[ (a - \iu \mathds{1})-\frac{\lambda -\iu}{\lambda + \iu} \,(a + \iu \mathds{1}) \right]
\cdot (a + \iu \mathds{1})^{-1} \\
= &\ \frac{1}{\lambda + \iu} \,\Bigl[ \,(\lambda + \iu) \,(a - \iu \mathds{1})
    -(\lambda - \iu) \,(a + \iu \mathds{1}) \,\Bigr] \cdot (a + \iu \mathds{1})^{-1} \\
= &\ \frac{2 \iu}{\lambda + \iu} \ (a - \lambda \mathds{1}) \cdot (a + \iu \mathds{1})^{-1}.
\end{align*}
This shows that $u-\nu \mathds{1}$ is injective precisely when
$a-\lambda \mathds{1}$ is so. Also, for the ranges, we have
$R(u - \nu \mathds{1}) = R(a - \lambda \mathds{1})$,
so that the statement follows from the Closed Graph Theorem.
(A left invertible operator is closed if and only if its left inverse is closed.)
\end{proof}

\begin{corollary}%
The spectrum of $a$ is a non-empty closed subset of the real line.
In particular, it is locally compact. \pagebreak
\end{corollary}

\begin{proof}
The Moebius transformation $\lambda \mapsto \nu$ of the preceding
proposition \ref{spau} maps $\mathds{R} \cup \{\infty\}$ homeomorphically
onto the unit circle, the point $\infty$ being mapped to $1$. Now the spectrum
of a unitary operator in a Hilbert space is a non-empty compact subset of
the unit circle, cf.\ \ref{C*unitary} and \ref{specradform}. If $a$ had empty
spectrum, then the spectrum of its Cayley transform would consist of $1$
alone, which would therefore be an eigenvalue, cf.\ \ref{eival}, in
contradiction with \ref{onenoteival}. 
\end{proof}

\begin{proposition}\label{swapadj}%
Let $b, c$ be two densely defined linear operators in $H$ such that
$b \subset c$. Then $c^* \subset b^*$.
\end{proposition}

\begin{proof}
For $x \in D(b)$, and $y \in D(c^*)$, we have
\[ \langle bx,y \rangle = \langle x,c^*y \rangle \]
because $b \subset c$. It follows that $y \in D(b^*)$
and $b^*y = c^*y$, so $c^* \subset b^*$.
\end{proof}

\begin{definition}[maximal symmetric operators]%
\index{operator!symmetric!maximal}%
\index{symmetric operator!maximal}%
\index{maximal!symmetric}%
A symmetric \linebreak operator in a Hilbert space is called
\underline{maximal symmetric} if it has no proper symmetric extension.
\end{definition}

\begin{proposition}\label{ismaxsymm}%
The self-adjoint linear operator $a$ is maximal symmetric.
\end{proposition}

\begin{proof}
Suppose that $b$ is a symmetric extension of $a$. By proposition
\ref{swapadj} above, we then have
\[ b \subset b^* \subset a^* = a \subset b, \]
so that $a = b$.
\end{proof}

\begin{definition}\label{commute}%
\index{operator!commuting}\index{commuting operator}%
A \underline{bounded} linear operator $c$ on $H$ is said to
\underline{commute} with $a$ if $D(a)$ is invariant under
$c$ and $cax = acx$ for all $x \in D(a)$.
\end{definition}

\begin{proposition}%
A bounded linear operator $c$ on $H$ commutes with $a$
if and only if $ca \subset ac$.
\end{proposition}

\begin{definition}\index{a65@$\{a\}'$}%
The set of bounded linear operators on $H$ \linebreak
commuting with $a$ is denoted by $\{a\}'$. \pagebreak
\end{definition}

\begin{theorem}%
[the Spectral Theorem for unbounded  self-adjoint operators]%
\index{Theorem!Spectral!self-adjoint unb.\ operator}%
\index{spectral!theorem!self-adjoint unb.\ operator}%
\index{spectral!resolution!self-adjoint unb.\ operator}%
\label{specthmunbded}
For the self-adjoint linear operator $a$ in $H$, there exists a unique
resolution of identity $P$ on $\s(a)$, acting on $H$, such that
\[ a = \int _{\text{\small{$\s(a)$}}} \id_{\text{\small{$\s(a)$}}} \,dP
\qquad \text{(pointwise)}. \]
This is usually written as
\[ a = \int \lambda \,dP(\lambda). \]
One says that $P$ is the \underline{spectral resolution} of $a$. We have
\[ \{a\}' = P', \]
and the support of $P$ is all of $\s(a)$. 
\end{theorem}

We split the proof into two lemmata.

\begin{lemma}\label{firstlemmaspecthm}%
Let $u$ denote the Cayley transform of $a$.
Let $\Omega := \s(u) \setminus \{1\}$. The function
\[ \kappa : \lambda \mapsto (\lambda - \iu)(\lambda + \iu)^{-1} \]
maps $\s(a)$ homeomorphically onto $\Omega$, cf.\ \ref{spau}.
Let $Q$ be the spectral resolution of $u$. Then $Q(\{1\}) = 0$, so that also
\[ u = \int _{\Omega} \id_{\Omega} \,dQ = \int \id \,d(Q|_{\Omega})
       \qquad \text{(pointwise)}. \]
Let $P := \kappa^{-1}(Q|_{\Omega})$ be the image of $Q|_{\Omega}$
under $\kappa^{-1}$. Then $P$ is a spectral resolution on $\s(a)$ such that
\[ a = \int _{\text{\small{$\s(a)$}}} \id_{\text{\small{$\s(a)$}}} \,dP
      \qquad \text{(pointwise)}. \]
Let conversely $P$ be a resolution of identity on $\s(a)$ such that
\[ a = \int _{\text{\small{$\s(a)$}}} \id_{\text{\small{$\s(a)$}}} \,dP
      \qquad \text{(pointwise)}. \]
The image $\kappa(P)$ of $P$ under $\kappa$ then is $Q|_{\Omega}$.
\pagebreak
\end{lemma}

\begin{proof}
We have $Q(\{1\}) = 0$ as $1$ is not an eigenvalue of $u$, see
\ref{onenoteival} and \ref{eival}. Hence also
\[ u = \int \id \,dQ = \int _{\Omega} \id_{\Omega} \,dQ = \int \id \,d(Q|_{\Omega}). \]
As in \ref{imagespdef} and the appendix \ref{imagedef} below, let
\[ P := \kappa^{-1}(Q|_{\Omega}) \]
be the image of $Q|_{\Omega}$ under $\kappa^{-1}$ and consider
\[ b := \Psi _P(\id) = \int \id \,dP = \int \id \,d\bigl(\kappa^{-1}(Q|_{\Omega})\bigr). \]
By the appendix \ref{imageend} below, we have
\[ b = \int _{\Omega}\kappa^{-1} \,dQ = \int _{\Omega} \iu \,\frac{1+\id}{1-\id} \,dQ. \]
It shall be shown that $b = a$. Since
\[ \mathds{1}-u = \int _{\Omega} (1-\id) \,dQ, \]
the Multiplication Theorem \ref{AddMult} gives
\[ b \,(\mathds{1}-u) = \int \iu \,(1+\id) \,dQ = \iu \,(\mathds{1}+u), \tag*{$(i)$} \]
and in particular
\[ R (\mathds{1}-u) \subset D(b). \tag*{$(ii)$} \]
On the other hand, for $x \in D(a)$, put
\[ (a + \iu \mathds{1}) \,x =: z. \]
We then have
\[ (a - \iu \mathds{1}) \,x = uz, \]
so
\begin{align*}
(\mathds{1}-u) \,z & = 2 \iu x \tag*{$(iii)$} \\
(\mathds{1}+u) \,z & = 2ax,
\end{align*}
whence
\[ a \,(\mathds{1}-u) \,z = a2 \iu x = \iu 2ax = \iu \,(\mathds{1}+u) \,z \tag*{$(iv)$} \]
for all $z \in R(a + \iu \mathds{1})$, which is equal to $H$ as $\s(a)$ is real,
cf.\ \ref{specisreal}. From $(ii)$ and $(iii)$ we have
\[ D(a) = R(\mathds{1}-u) \subset D(b). \tag*{$(v)$}  \]
From $(i)$ and $(iv)$ it follows now that $b$ is a self-adjoint extension of $a$.
Since $a$ is maximal symmetric \ref{ismaxsymm}, this implies that
$a = b = \Psi _P(\id)$.\pagebreak

Let conversely $P$ be a resolution of identity on $\s(a)$ such that
\[ a = \int _{\text{\small{$\s(a)$}}} \id_{\text{\small{$\s(a)$}}} \,dP. \]
Let $R : = \kappa(P)$ be the image of $P$ under $\kappa$ and put
\[ v := \Psi _{R} (\id) = \int \id \,dR = \int \id \,d\bigl(\kappa(P)\bigr). \]
By the appendix \ref{imageend} again, we find
\[ v= \int \kappa \,dP = \int \frac{\id - \iu}{\id + \iu} \,dP, \]
so that the Multiplication Theorem \ref{AddMult} gives
\[ v \,(a + \iu \mathds{1}) = a - \iu \mathds{1}, \]
or
\[ v = (a - \iu \mathds{1})(a + \iu \mathds{1})^{-1} = u. \]
By the uniqueness of the spectral resolution, we must have
\[ \kappa(P) = R = Q|_{\Omega}. \qedhere \]
\end{proof}

\begin{lemma}%
If $u$ is the Cayley transform of $a$, then
\[ \{a\}' = \{u\}'. \]
\end{lemma}

\begin{proof}
If $c \in \{a\}'$ then
\[ c \,(a \pm \iu \mathds{1}) \subset (a \pm \iu \mathds{1}) \,c, \]
whence
\[ (a \pm \iu \mathds{1})^{-1} \,c = c \,(a \pm \iu \mathds{1})^{-1}. \]
It follows that $c\,u = u\,c$.

Conversely, let $c \in \{u\}'$. From the formulae $(iv)$ and $(v)$
of the preceding lemma, $\mathds{1}-u$ has range $D(a)$, and
\[ a \,(\mathds{1}-u) = \iu \,(\mathds{1}+u). \]
Therefore we have
\begin{align*}
c \,a \,(\mathds{1}-u) & = \iu \,c \,(\mathds{1}+u) \\
 & = \iu \,(\mathds{1}+u) \,c \\
 & = a \,(\mathds{1}-u) \,c = a \,c \,(\mathds{1}-u),
\end{align*}
whence $c \in \{a\}'$. \pagebreak
\end{proof}

\begin{definition}[function of a self-adjoint operator]%
\index{function!of an operator}\label{unbdedf}%
Let $P$ be the spectral resolution of the self-adjoint linear
operator $a$. If $f$ is a complex-valued Borel function on $\s(a)$,
one writes
\[ f(a) := \Psi_P (f) = \int f \,d P \quad \text{(pointwise)}, \]
and one says that $f(a)$ is a \underline{function} of $a$.
\end{definition}

We have actually shown the following:

\begin{addendum}
The Cayley transform
\[ u := (a - \iu \mathds{1}) (a + \iu \mathds{1})^{-1} \]
is the function $\kappa (a)$ with
\begin{alignat*}{2}
\kappa : \s(a)       & &\ \to            &\ \mathds{C} \\
               \lambda & &\ \mapsto  &\ (\lambda - \iu)(\lambda + \iu)^{-1}.
\end{alignat*}
\end{addendum}

\begin{proof}
See the proof of the converse part of lemma \ref{firstlemmaspecthm}.
\end{proof}

\clearpage


\section{Application: an Initial Value Problem}%
\label{scopug}

In this paragraph, let $H$ be a Hilbert space $\neq \{0\}$.

\begin{definition}%
A function
\begin{alignat*}{2}
\mathds{R} & \to           &\ & B(H) \\
                   t & \mapsto &\ & U_t
\end{alignat*}
is called a \underline{strongly continuous one-parameter unitary group}, if
\begin{itemize}
   \item[$(i)$] $U_t$ is unitary for all $t \in \mathds{R}$,
  \item[$(ii)$] $U_{t+s} = U_t \,U_s$ for all $s$, $t \in \mathds{R}$,
                      \quad (the group property)
 \item[$(iii)$] for all $x \in H$, the function $t \mapsto U_t \,x$
                       is continuous on $\mathds{R}$.
\end{itemize}
\end{definition}

\begin{definition}[infinitesimal generator]%
\index{infinitesimal generator}%
If $t \mapsto U_t$ is a strongly continuous one-parameter unitary
group in $H$, one says that the \linebreak
\underline{infinitesimal generator} of $t \mapsto U_t$ is the linear
operator $a$ in $H$, with domain $D(a) \subset H$, characterised by
\[ x \in D(a),\ y = a x \quad \Leftrightarrow \quad
\ - \frac{1}{\iu} \frac{d}{dt} \Bigl|_{\text{\small{$t=0$}}} \,U_t \,x
\text{\ exists and equals\ } y. \]
This is also written as
\[ - \frac{1}{\iu} \frac{d}{dt} \Bigl|_{\text{\small{$t=0$}}} \,U_t
= a\quad \text{(pointwise)}. \]
\end{definition}

\begin{theorem}
Let $a$ be a self-adjoint linear operator in $H$, and let $P$ be its
spectral resolution. For $t \in \mathds{R}$, we denote by $U_t$ the
bounded linear operator on $H$ given by
\[ U_t := \exp(- \iu t a) = \int _{\text{\small{$\s(a)$}}} \exp(- \iu t \lambda) \,dP (\lambda)
\quad \text{(pointwise)}. \]
(Cf.\ \ref{unbdedf}.) Then $t \mapsto U_t$ is a strongly continuous one-parameter
unitary group with infinitesimal generator $a$. One actually has
\[ - \frac{1}{\iu} \frac{d}{dt} \Bigl|_{\text{\small{$t=s$}}} \,U_t \,x
= a \,U_s \,x \]
for all $s \in \mathds{R}$ and all $x \in D(a)$, i.e.\ $t \mapsto U_t \,x$ is a
solution to the initial value problem
\[ - \frac{1}{\iu} \frac{d}{dt} \,u(t) = a \,u(t),\quad u(0) = x.  \]
Please note that $t \mapsto U_t$ is the Fourier transform of $P$. \pagebreak
\end{theorem}

\begin{proof}
The spectral resolution $P$ of $a$ lives on the spectrum of $a$,
which is a subset of the real line, cf.\ \ref{specisreal}. The
fact that $\pi_P$ is a representation in $H$ implies therefore that
$U_t$ is unitary for all $t \in \mathds{R}$, and that $U_{t+s} = U_t \,U_s$
for all $s$, $t \in \mathds{R}$. Furthermore, the function $t \mapsto U_t \,x$
is continuous on $\mathds{R}$ for all $x \in H$ by the Dominated
Convergence Theorem \ref{domconvergence}. Thus $t \mapsto U_t$
$(t \in \mathds{R})$ is a strongly continuous one-parameter unitary group.
We shall show that $a$ is the infinitesimal generator of $t \mapsto U_t$.
For $x \in D(a)$ and $t \in \mathds{R}$, we have
\[ {\Bigl\| \ {- \frac{1}{\iu t} \,(U_t \,x - x) - ax} \ \Bigr\| \,}^2
= \int {| \,f_t \,| \,}^2 \,d \langle Px, x \rangle, \]
with
\[ f_t (\lambda) = - \frac{1}{\iu t} \,\bigl(\exp(- \iu t \lambda) - 1 + \iu t \lambda \bigr)
\qquad (\lambda \in \mathds{R}). \]
By power series expansion of $\exp(- \iu t \lambda)$, one sees that $f_t \to 0$
pointwise as $t \to 0$. We have $| \,\exp( \iu s) -1 \,| \leq | \,s \,|$ for all real
$s$ (as chord length is smaller than arc length). It follows that
$| \,f_t (\lambda) \,| \leq 2 \,| \,\lambda \,|$ for all $t \in \mathds{R}$. Since
$\lambda \mapsto 2 \,| \,\lambda \,|$ is a function in
${\mathscr{L} \,}^2( \langle Px, x \rangle )$, the Lebesgue Dominated
Convergence Theorem implies that
\[ \lim _{t \to 0} - \frac{1}{\iu t} \,(U_t \,x - x) = a x\quad \text{for all } x \in D(a). \]
Thus, if $b$ denotes the infinitesimal generator of $t \mapsto U_t$, we
have shown that $b$ is an extension of $a$. Now $b$ is symmetric
because for $x$, $y \in D(b)$ one has
\[ \langle \,- \frac{1}{\iu t} \,(U_t \,x - x), y \,\rangle
= \langle \,- \frac{1}{\iu t} \,(U_t - \mathds{1}) \,x, y \,\rangle
= \langle \,x, \frac{1}{\iu t} \,(U_{-t} \,y - y) \,\rangle. \]
It follows that $b = a$ because $a$ is maximal symmetric, cf.\ \ref{ismaxsymm}.
The second statement follows from the group property of $t \mapsto U_t$
and from the fact that $U_t$ commutes with $a$ by \ref{specthmunbded}
and \ref{spprojcomm}. Please note also that $D(a)$ is invariant under $U_t$,
cf.\ \ref{commute}
\end{proof}

\smallskip
This result is of prime importance for quantum mechanics,
where the initial value problem is the Schr\"odinger equation
for the evolution in time of the quantum mechanical state
described by the vector $U_t \,x$.

\clearpage


\phantomsection

\addcontentsline{toc}{part}{Appendix}


\begin{center} \textbf{\huge{Appendix}} \end{center}


\setcounter{section}{44}

\phantomsection

\addcontentsline{toc}{chapter}{\texorpdfstring{\ \ \ \S\ 44.\quad Quotient Spaces}%
{\textsection\ 44. Quotient Spaces}}

\pagestyle{myheadings}
\markboth{\textsc{APPENDIX}}{\S\ 44. \textsc{QUOTIENT SPACES}}

\smallskip
\begin{center} \textbf{\S\ 44.\:\:Quotient Spaces} \end{center}

\setcounter{theorem}{0}

\begin{reminder}\index{quotient!space}\label{quotspace}%
Let $(A,| \cdot |)$ be a normed space and $B$
a \underline{closed} subspace of $A$. The quotient space $A/B$
is defined as
\[ A/B := \{ a+B \subset A : a \in A \}. \]
The canonical projection (here denoted by an underscore)
\begin{alignat*}{2}
\_ \,: A\ & & \to     &\ A/B \\
       a\ & & \mapsto &\ a+B =: \underline{a}
\end{alignat*}
is a linear map from $A$ onto $A/B$. The definition
\[ | \,\underline{a}\, | := \inf \,\{ \,|\,c\,| : c \in \underline{a} \,\}
\tag*{$( \,\underline{a} \in A/B \,) $} \]
then defines a norm on A/B. It is called the
\underline{quotient norm} on $A/B$.
\end{reminder}

\begin{theorem}\label{quotspaces}%
Let $(A,| \cdot |)$ be a normed space and let $B$
be a \underline{closed} subspace of $A$. If $A$ is a Banach space,
then $A/B$ is a Banach space. If both $A/B$ and $B$ are Banach
spaces, then $A$ is a Banach space.
\end{theorem}

\begin{proof}
Assume first that $A$ is complete. Let
$(\underline{c_k})_{k \geq 0}$ be a Cauchy sequence in $A/B$. We can then
find a subsequence $(\underline{d_n})_{n \geq 0}$ of
$(\underline{c_k})_{k \geq 0}$ such that
\[ | \,\underline{d_{n+1}} - \underline{d_n} \,| < 2^{-(n+1)} \]
for all $n \geq 0$. We construct a sequence $(a_n)_{n \geq 0}$ in $A$ such
that
\[ \underline{a_n} = \underline{d_n} \]
and
\[ | \,a_{n+1} -a_n \,| < 2^{-(n+1)} \]
for all $n \geq 0$. This is done as follows. Let $a\0 \in \underline{d\0}$.
Suppose that $a_n$ has been constructed and let
$d_{n+1} \in \underline{d_{n+1}}$. We obtain
\[ \inf \,\{ \,| \,d_{n+1} -a_n +b \,| : b \in B \,\}
= | \,\underline{d_{n+1}} - \underline{a_n} \,|
= | \,\underline{d_{n+1}} - \underline{d_n} \,| < 2^{-(n+1)} \]
so that there exists $b_{n+1} \in B$ such that
\[ | \,d_{n+1} -a_n +b_{n+1} \,| < 2^{-(n+1)}.  \]
We then put $a_{n+1} := d_{n+1} + b_{n+1} \in \underline{d_{n+1}}$.
This achieves the construction of $(a_n)_{n \geq 0}$ as required.
But then $(a_n)_{n \geq 0}$ is a Cauchy sequence in $A$ and converges
in $A$ by assumption. By continuity of the projection,
$(\underline{d_n})_{n \geq 0}$ converges in $A/B$, so that also
$(\underline{c_k})_{k \geq 0}$ converges in $A/B$. \pagebreak

Assume now that both $A/B$ and $B$ are complete. Let $(a_n)_{n \geq 0}$
be a Cauchy sequence in $A$. Then $(\underline{a_n})_{n \geq 0}$
is a Cauchy sequence in $A/B$ because the projection is contractive. By
assumption, there exists $a \in A$ such that
\[ \underline{a} = \lim _{n \to \infty} \underline{a_n}. \]
For $n \geq 0$, we can find $b_n \in B$ such that
\[ |\,(a-a_n)+b_n\,|
< | \,\underline{a}-\underline{a_n} \,| + \frac{1}{n}. \]
But then
\begin{align*}
| \,b_m-b_n \,| & = \bigl|\,[ \,(a-a_m) + b_m \,]
                  - [ \,(a-a_n) + b_n \,] + (a_m-a_n) \,\bigr| \\
              & \leq |\,\underline{a}-\underline{a_m}\,| + \frac{1}{m} +
                     |\,\underline{a}-\underline{a_n}\,| + \frac{1}{n} +
                     |\,a_m-a_n\,|,
\end{align*}
from which we see that $(b_n)_{n \geq \0}$ is a Cauchy sequence in
$B$ and consequently convergent to some element $b$ in $B$. Finally we have
\[ |\,(a+b)-a_n\,| \leq |\,(a-a_n)+b_n\,| + |\,b-b_n\,| \]
which implies that
\[ \lim _{n \to \infty} a_n = a+b. \qedhere \]
\end{proof}

\clearpage


\setcounter{section}{45}

\phantomsection

\addtocontents{toc}{\protect\vspace{-0.65em}}

\addcontentsline{toc}{chapter}{\texorpdfstring{\ \ \ \S\ 45.\quad Reminder on Topology}%
{\textsection\ 45. Reminder on Topology}}

\pagestyle{myheadings}
\markboth{\textsc{APPENDIX}}{\S\ 45. \textsc{REMINDER ON TOPOLOGY}}

\bigskip
\begin{center} \textbf{\S\ 45.\:\:Reminder on Topology} \end{center}
\medskip

\setcounter{theorem}{0}

We gather here the topological ingredients of the Commutative
Gelfand-Na\u{\i}mark Theorem.
\hypertarget{remtop}{}%

\begin{theorem}[the weak topology]%
Let $X$ be a set. Let $\{ Y_i \}_{i \in I}$ be a family of topological spaces,
and let $\{ f_i \}_{i \in I}$ be a family of functions $f_i : X \to Y_i$ $(i \in I)$.
Among the topologies on $X$ for which all the functions $f_i$ $(i \in I)$
are continuous, there exists a coarsest topology. It is the topology
generated by the base consisting of all the sets of the form
$\bigcap_{t \in T} {f_t}^{-1} (O_t)$, where $T$ is a finite subset of $I$ and
$O_t$ is an open subset of $Y_t$ for all $t \in T$. This topology is called
the \underline{weak topology induced by $\{ f_i \}_{i \in I}$}. It is also called
the initial topology induced by $\{ f_i \}_{i \in I}$.
\end{theorem}

\begin{proof}
We refer the reader to \cite[Prop.\ 1.4.8 p.\ 31]{Eng}.
\end{proof}

\medskip
The philosophy is, that the weaker a topology is, the more quasi-compact
subsets it has.

\begin{theorem}[the universal property]%
Let $X$ be a topological space carrying the weak topology induced by
a family $\{ f_i \}_{i \in I}$. The \linebreak following property is called the
\underline{universal property} of the weak topology induced by
$\{ f_i \}_{i \in I}$. A function $g$ from a topological space to $X$ is
\linebreak continuous if and only if the compositions $f_i \circ g$ are
continuous for all $i \in I$.
\end{theorem}

\begin{proof}
We refer the reader to \cite[Prop.\ 1.4.9 p.\ 31]{Eng}.
\end{proof}

\begin{definition}[the weak* topology]%
\label{weak*top}\index{topology!weak*}%
\index{universal!property!w@of weak* topology}%
Let $A$ be a real or complex vector space, and let $A^*$ be its dual space.
The \underline{weak* topology} on $A^*$ is the weak topology induced by
the evaluations $\wht{a}$ at elements $a \in A$:
\begin{alignat*}{2}
\wht{a} : A^* &           \to & & \ \mathds{R} \text{ or } \mathds{C} \\
              \ell \,& \mapsto & & \ \wht{a}(\ell) := \ell (a).
\end{alignat*}
In particular, the evaluations $\wht{a}$ $(a \in A)$ are continuous in the weak*
topology, and the weak* topology is the coarsest topology on $A^*$ such that
these evaluations are continuous. The \underline{universal property} of the
weak* topology says that a function $g$ from a topological space to $A^*$ is
continuous with respect to the weak* topology on $A^*$, if and only if the
compositions $\wht{a} \circ g$ $(a \in A)$ are continuous. \pagebreak
\end{definition}

\begin{definition}[separating the points]\index{separating!the points}%
Let $\mathcal{F}$ be a set of functions defined on a set $\Omega$.
The set $\mathcal{F}$ is said to \underline{separate the points} of
$\Omega$, if for $x$, $y \in \Omega$ with $x \neq y$, there exists a
function $f \in \mathcal{F}$ such that $f(x) \neq f(y)$.
\end{definition}

\begin{proposition}
Let $X$ be a topological space carrying the weak topology induced
by a family $\{ f_i \}_{i \in I}$ of functions $f_i: X \to Y_i$ $(i \in I)$.
If the set $\{ \,f_i : i \in I \,\}$ separates the points of $X$, and if all the
spaces $Y_i$ $(i \in I)$ are Hausdorff, then $X$ is Hausdorff. 
\end{proposition}

\begin{proof}
The proof is very easy and left as an exercise. Otherwise, see
e.g.\ \cite[Prop.\ 1.5.3 pp.\ 23 f.]{PedT}.
\end{proof}

\medskip
In particular, the weak* topology is Hausdorff.

\medskip
The importance of the weak* topology stems for a large part from the
following compactness result.

\begin{theorem}[Alaoglu's Theorem]\index{Alaoglu's Theorem}%
\label{Alaoglu}\index{Theorem!Alaoglu}%
Let $A$ be a (possibly real) \linebreak normed space, and let $A'$ be
its dual normed space. The closed unit ball in $A'$ then is a compact
Hausdorff space in the weak* topology.
\end{theorem}

\begin{proof}
See e.g.\ \cite[Thm.\ 2.5.2 p.\ 69]{PedT}. 
\end{proof}

\begin{theorem}\label{homeomorph}%
Let $f$ be a continuous function from a quasi-{\linebreak}compact
space $\Omega$ to a Hausdorff space. Then $f$ is a closed function
(in the sense that $f$ maps closed sets to closed sets). If $f$ is injective,
then $f$ is an imbedding. If $f$ is bijective, then $f$ is a homeomorphism.
\end{theorem}

\begin{proof}
A closed subset $\omega$ of the quasi-compact space $\Omega$ is
quasi-compact, and so is its continuous image $f(\omega)$. The set
$f(\omega)$ then is closed, as a quasi-compact subset of a Hausdorff
space is closed. If $f$ is injective, then $f : \Omega \to f(\Omega)$ is an
open function (in the sense that $f$ maps open sets to open sets).
Indeed, if $\omega$ is an open subset of $\Omega$, then
$f(\Omega) \setminus f(\omega) = f(\Omega \setminus \omega)$
is closed in $f(\Omega)$. This implies that $f : \Omega \to f(\Omega)$
is a homeomorphism. That is, $f$ is an imbedding.
\end{proof}

\begin{corollary}
The topology of a compact Hausdorff space cannot be strengthened
without losing the quasi-compactness property, and it cannot be
weakened without losing the Hausdorff property. \pagebreak
\end{corollary}

\begin{definition}[vanishing nowhere]\index{vanishing nowhere}%
Let $\mathcal{F}$ be a set of functions defined on a set $\Omega$.
The set $\mathcal{F}$ is said to \underline{vanish nowhere} on $\Omega$,
if the functions in $\mathcal{F}$ have no common zero on $\Omega$.
\end{definition}

\begin{theorem}[the Stone-Weierstrass Theorem]%
\label{StW}\index{Theorem!Stone-Weierstrass}%
\index{Stone-Weierstrass Theorem}%
Let $\Omega$ be a \linebreak locally compact Hausdorff space.
A \st-subalgebra of $C\0(\Omega)$, which \linebreak separates
the points of $\Omega$, and which vanishes nowhere
on $\Omega$, is dense in $C\0(\Omega)$.
\end{theorem}

\begin{proof}
We refer the reader to \cite[Cor.\ 4.3.5 p.\ 146]{PedT}.
\end{proof}

\medskip
Please note that we need not only a subalgebra, but a
\st-subalgebra.

\clearpage


\setcounter{section}{46}

\phantomsection

\addtocontents{toc}{\protect\vspace{-0.65em}}

\addcontentsline{toc}{chapter}%
{\texorpdfstring{\ \ \ \S\ 46.\quad Complements to Integration Theory}%
{\textsection\ 46. Complements to Integration Theory}}

\pagestyle{myheadings}
\markboth{\textsc{APPENDIX}}%
{\S\ 46. \textsc{COMPLEMENTS TO INTEGRATION THEORY}}

\medskip
\begin{center} \textbf{\S\ 46.\:\:Complements to Integration Theory} \end{center}
\setcounter{theorem}{0}

\begin{definition}[the Borel and Baire $\sigma$-algebras]\label{Bairedef}%
Let $\Omega$ be a \linebreak Hausdorff space. The $\sigma$-algebra of
\underline{Borel sets} is defined as the $\sigma$-algebra generated by the
open subsets of $\Omega$. The $\sigma$-algebra of \underline{Baire sets}
is \linebreak defined as the $\sigma$-algebra generated by the sets
$\{\,f=0\,\}$, with $f \in C(\Omega)$.
\end{definition}

\begin{proposition}\label{metric}%
Every Baire set is a Borel set. On a metric space, the Baire and Borel sets
coincide. (Distance function from a point to a fixed closed set.)
\end{proposition}

\begin{definition}[inner regular Borel probability measure]%
\index{measure!inner regular}\index{measure!Borel}\label{inregBormeas}%
Let $\Omega$ be a Hausdorff space $\neq \varnothing$. A
\underline{Borel probability measure} on $\Omega$ is a probability
measure defined on the Borel sets of $\Omega$. A Borel probability
measure $\mu$ on $\Omega$ is called \underline{inner regular}, if
for every Borel set $\omega \subset \Omega$, one has
\[ \mu(\omega) = \sup \,\{ \,\mu(K) : K \text{ is a compact subset of } \omega \,\}. \]
\end{definition}

We shall use the Riesz Representation Theorem in the following form.
The proof follows easily from other versions of the Theorem.

\begin{theorem}[the Riesz Representation Theorem]%
\label{Riesz}\index{Theorem!Riesz Representation}%
\index{Riesz Representation Thm.}%
Let $\Omega \neq \varnothing$ be a locally compact Hausdorff space. Let
$\ell$ be a linear functional on $C_c(\Omega)$, which is positive in the
sense that $\ell(f) \geq 0$ for every $f \in C_c(\Omega)$ with $f \geq 0$
pointwise. Assume that $\ell$ is bounded with norm $1$ on the
pre-C*-algebra $C_c(\Omega)$. There then exists a unique inner regular
Borel probability measure $\mu$ on $\Omega$ such that
\[ \ell (f) = \int f \,d \mu\quad \text{for every } f \in C_c(\Omega). \]
\end{theorem}

\begin{definition}[measurability]%
\index{measurable}\label{imagebegin}%
Let $\mathcal{E}$, $\mathcal{E}'$ be two $\sigma$-algebras on two
non-empty sets $\Omega,$ $\Omega'$ respectively. A function
$f : \Omega \to \Omega'$ is called \linebreak
\underline{$\mathcal{E}$-\,$\mathcal{E}'$ measurable} if
$f^{-1}(\omega') \in \mathcal{E}$ for all $\omega' \in \mathcal{E}'$.
A complex-valued function on $\Omega$ is called
\underline{$\mathcal{E}$-measurable}, if it is $\mathcal{E}$-Borel
measurable. A complex-valued function on $\Omega$ is
$\mathcal{E}$-measurable if and only if its real and imaginary parts are,
cf.\ \cite[Cor.\ 2.5 p.\ 45]{FollInt}. A real-valued function $f$ on $\Omega$ is
$\mathcal{E}$-measurable if and only if $\{ \,f < \alpha \,\} \in \mathcal{E}$
for all $\alpha \in \mathds{R}$, cf.\ \cite[Prop.\ 2.3 p.\ 44]{FollInt}.
A complex-valued function on a Hausdorff space is called a
\underline{Borel function}, if it is Borel-measurable, and it is called a
\underline{Baire function}, if it is Baire-measurable. \pagebreak
\end{definition}

\begin{proposition}\label{continuous}%
Every continuous complex-valued function is a Baire function.
Every Baire function is a Borel function. 
\end{proposition}

\begin{proof}
The first statement follows from
\[ \,\{ \,\mathrm{Re}\,(f) < \alpha \,\} = \{ \,(\mathrm{Re}\,(f) - \alpha)_+ = 0 \,\}
\setminus \{ \,\mathrm{Re}\,(f) - \alpha = 0 \,\}. \qedhere \]
\end{proof}

\begin{definition}[image measure]%
\index{image!measure}\index{measure!image}\label{imagedef}%
Let $(\Omega, \mathcal{E}, \mu)$ be a probability space. Let $\mathcal{E}'$
be a $\sigma$-algebra on a set $\Omega'$. Let $f : \Omega \to \Omega'$
be an $\mathcal{E}$-\,$\mathcal{E}'$ measurable function. The
\underline{image measure} $f(\mu)$ is defined as the measure
\begin{alignat*}{2}
 f(\mu) : \mathcal{E}' & \to           &\ & [0, 1] \\
                     \omega' & \mapsto &\ & \mu\bigl( f^{-1}(\omega') \bigr).
\end{alignat*}
That is, $f(\mu) = \mu \circ f^{-1}$.
If $f$ has the meaning of a random variable on $(\Omega, \mathcal{E}, \mu)$,
then $f(\mu)$ is the distribution of the random variable $f$.
\end{definition}

\begin{theorem}\label{imageend}%
Let $(\Omega, \mathcal{E}, \mu)$ be a probability space. Let
$\mathcal{E}'$ be a \linebreak $\sigma$-algebra on a set
$\Omega'$. Let $f : \Omega \to \Omega'$ be an
$\mathcal{E}$-\,$\mathcal{E}'$ measurable function. For an
$\mathcal{E}'$-measurable function $g$, we have
\begin{alignat*}{2}
g \in {\mathscr{L} \,}^{1} \bigl(f(\mu)\bigr)
 & \ \Leftrightarrow & & \ g \circ f \in {\mathscr{L} \,}^{1}(\mu) \\
 & \ \Rightarrow & & \ \int g\ d\,f (\mu)
= \int (g \circ f) \,d \mu.
\end{alignat*}
\end{theorem}

\begin{proof}
Approximate positive $\mathcal{E}'$-measurable functions
by increasing sequences of $\mathcal{E}'$-step functions
and apply the Beppo Levi Principle.
\end{proof}

\medskip
The next result is well known.

\begin{theorem}\label{CcdenseinLp}%
Let $\mu$ be an inner regular Borel probability measure on a locally
compact Hausdorff space $\Omega$. For a function
$f \in {\mathscr{L} \,}^{p} \,(\mu)$ with $1 \leq p < \infty$,
there exists a sequence $\{\,f_n\,\}$ in $C_c(\Omega)$
which converges to $f$ in ${\mathscr{L} \,}^{p} \,(\mu)$
and $\mu$-almost everywhere.
\end{theorem}

For a proof, see e.g.\ \cite[Thm.\ 6.15 p.\ 204]{FW}.

\begin{corollary}\label{Baire}%
Let $\mu$ be an inner regular Borel probability \linebreak measure
on a locally compact Hausdorff space $\Omega$. Then a function in
${\mathscr{L} \,}^{\infty}(\mu)$ is $\mu$-a.e.\ equal to a bounded Baire
function on $\Omega$. \pagebreak
\end{corollary}

\begin{proof}
Let $f \in {\mathscr{L} \,}^{\infty}(\mu)$. As $\mu$ is a probability measure,
we have ${\mathscr{L} \,}^{\infty}(\mu) \subset {\mathscr{L} \,}^{1}(\mu)$, so
we apply the preceding theorem with $p = 1$. \linebreak Let $\{\,f_n\,\}$ be
a sequence as described. The trimmed function
\[ \limsup_{n \to \infty} \ \Bigl\{ \,\bigl[ \,\mathrm{Re}\,(f_n)
\wedge \| \,\mathrm{Re}\,(f) \,\|_{\,\text{\small{$\mu$}},\infty} \,1_{\,\Omega} \,\bigr]
\,\vee \,\bigl[ \,- \| \,\mathrm{Re}\,(f) \,\|_{\,\text{\small{$\mu$}},\infty} \,1_{\,\Omega} \,\bigr]
\,\Bigr\} \]
is a bounded Baire function which is $\mu$-a.e.\ equal to $\mathrm{Re}\,(f)$.
\end{proof}

\medskip
We also have the similar result \ref{general}. For the proof, we need
the following theorem, which is well known.

\begin{theorem}\label{pless}%
Let $(\Omega, \mathcal{E}, \mu)$ be a probability space.
Then a function in ${\mathscr{L} \,}^{p} \,(\mu)$ with $1 \leq p < \infty$
is $\mu$-a.e.\ equal to a complex-valued $\mathcal{E}$-measurable
function on $\Omega$.
\end{theorem}

For a proof, see e.g. \cite[Thm.\ 3.2.13 p.\ 480]{CFW}.

\begin{corollary}\label{general}%
Let $(\Omega, \mathcal{E}, \mu)$ be a probability space.
Then a function in ${\mathscr{L} \,}^{\infty} \,(\mu)$ is
$\mu$-a.e.\ equal to a bounded $\mathcal{E}$-measurable
function on $\Omega$.
\end{corollary}

\begin{proof}
Let $f \in {\mathscr{L} \,}^{\infty}(\mu)$. As $\mu$ is a probability measure,
we have ${\mathscr{L} \,}^{\infty}(\mu) \subset {\mathscr{L} \,}^{1}(\mu)$, so
we apply the preceding theorem with $p = 1$.  There exists a real-valued
$\mathcal{E}$-measurable function $g$ on $\Omega$, which is \linebreak
$\mu$-a.e.\ equal to $\mathrm{Re} \,(f)$. The trimmed function
\[ \bigl[ \,g \wedge \| \,\mathrm{Re}\,(f) \,\|_{\,\text{\small{$\mu$}},\infty} \,1_{\,\Omega} \,\bigr]
\,\vee \,\bigl[ \,- \| \,\mathrm{Re}\,(f) \,\|_{\,\text{\small{$\mu$}},\infty} \,1_{\,\Omega} \,\bigr] \]
then is a bounded $\mathcal{E}$-measurable function on $\Omega$,
which is $\mu$-a.e.\ equal to $\mathrm{Re} \,(f)$.
\end{proof}

\clearpage


\pagestyle{headings}


\backmatter


\phantomsection

\addcontentsline{toc}{part}{Back Matter}

\phantomsection

\addtocontents{toc}{\protect\vspace{-1.4ex}}

\clearpage


\phantomsection

\addtocontents{toc}{\protect\vspace{-1.4ex}}

\printindex

\clearpage





\begin{thebibliography}{99}


%


\vspace{8.5ex}

\centerline{\textbf{Algebra}}

\vspace{1.4ex}


\bibitem{Kost} A.\ I.\ Kostrikin. \emph{Introduction to Algebra}.
(Universitext), Springer-Verlag, New York, 1982


\vspace{1ex}

\begin{center} \hspace{-2.3em} \textbf{Analysis} \end{center}

\vspace{1ex}


\bibitem{RuPr} W.\ Rudin. \emph{Principles of Mathematical Analysis}.
Third Edition, (international student edition), McGraw-Hill International
Book Company, 1976


\vspace{1ex}

\begin{center} \hspace{-2.3em} \textbf{Complex Analysis} \end{center}

\vspace{1ex}


\bibitem{Hen} P.\ Henrici. \emph{Applied and Computational Complex
Analysis}. vol.\ I, Wiley-Interscience, 1974


\bibitem{RCAC} W.\ Rudin. \emph{Real and Complex Analysis}. Third Edition,
(international edition), Mathematics Series, McGraw-Hill Book Company, 1987


\vspace{1ex}

\begin{center} \hspace{-2.3em} \textbf{Topology} \end{center}

\vspace{1ex}


\bibitem{EngOut} R.\ Engelking. \emph{Outline of General Topology}.
American Elsevier Publishing Company, New York, 1968


\bibitem{Eng} R.\ Engelking. \emph{General Topology}. Revised and
completed edition, Sigma Series in Pure Mathematics, vol.\ 6, Heldermann
Verlag, Berlin, 1989


\bibitem{TopFoll} G.\ B.\ Folland. \emph{Real Analysis: Modern Techniques
and Their Applications}. 2nd edition, Wiley Interscience, 1999


\bibitem{PedT} G.\ K.\ Pedersen. \emph{Analysis Now}. GTM 118,
Springer-Verlag, New York, 1988


\bibitem{Top} B.\ von Querenburg. \emph{Mengentheoretische Topologie}.
Dritte neu bearbeitete und erweiterte Auflage, Springer-Verlag, Berlin
Heidelberg New York, 2001 


\vspace{1ex}

\begin{center} \hspace{-2.3em} \textbf{Integration Theory} \end{center}

\vspace{1ex}


\bibitem{Bau} H.\ Bauer. \emph{Measure and Integration Theory}.
Studies in Mathematics 26, Walter de Gruyter, Berlin, New York, 2001


\bibitem{CFW} C.\ Constantinescu, W.\ Filter, K.\ Weber in collaboration
with A.\ Sontag. \emph{Advanced Integration Theory}. Mathematics and
its Applications, Kluwer Academic Publishers, 1998


\bibitem{FW} W.\ Filter, K.\ Weber. \emph{Integration Theory}. Mathematics
Series, Chapman \& Hall, 1997 \pagebreak


\bibitem{Flo} K.\ Floret. \emph{Ma\ss- und Integrationstheorie}. Teubner
Studienb\"ucher Mathematik, B.\ G.\ Teubner, Stuttgart, 1981


\bibitem{FollInt} G.\ B.\ Folland. \emph{Real Analysis: Modern Techniques
and Their Applications}. 2nd edition, Wiley Interscience, 1999


\bibitem{Koe} H. K\"{o}nig. \emph{Measure and Integration, An Advanced
Course in Basic Procedures and Applications}. Springer-Verlag, Berlin
Heidelberg New York, 1997


\vspace{1ex}

\begin{center} \hspace{-2.3em} \textbf{Functional Analysis} \end{center}

\vspace{1ex}


\bibitem{Conw} J.\ B.\ Conway. \emph{A Course in Functional Analysis, second
edition}. GTM 96, Springer-Verlag, New York Berlin Heidelberg Tokyo, 1990


\bibitem{Krey} E.\ Kreyszig. \emph{Introductory Functional Analysis with
Applications}. Wiley Classics Library Edition, John Wiley \& Sons, 1989


\bibitem{PedB} G.\ K.\ Pedersen. \emph{Analysis Now}. GTM 118,
Springer-Verlag, New York, 1988


\bibitem{ReSi1} M.\ Reed, B.\ Simon. \emph{Functional Analysis}. Methods of
Modern Mathematical Physics, vol.\ I.\ Academic Press, New York, 1972


\bibitem{RFAF} W.\ Rudin. \emph{Functional Analysis}. International Series in
Pure and Applied Mathematics, McGraw-Hill, TMH Edition 1982, 2nd Edition 1991


\bibitem{Sche} M.\ Schechter. \emph{Principles of Functional Analysis}.
2nd Edition. Graduate Studies in Mathematics vol.\ 36, American Mathematical
Society, 2002 


\vspace{1ex}

\begin{center} \hspace{-2.3em} \textbf{Banach Algebras} \end{center}

\vspace{1ex}


\bibitem{Aup} B.\ Aupetit. \emph{A Primer on Spectral Theory}.
Universitext, Springer-Verlag, New York, 1991


\bibitem{BD} F.\ F.\ Bonsall, J.\ Duncan. \emph{Complete Normed Algebras}.
Ergebnisse der Mathematik und Ihrer Grenzgebiete, Band 80, Springer-Verlag,
Berlin Heidelberg New York, 1973


\bibitem{BB} N.\ Bourbaki. \emph{Th\'{e}ories spectrales, Chapitres 1 et 2}.
\'{E}l\'{e}ments de Math\'{e}matique, R\'{e}impression inchang\'{e}e de
l'\'{e}dition originale de 1967, Springer-Verlag, 2007


\bibitem{GaalB} S.\ A.\ Gaal. \emph{Linear Analysis and Representation Theory}.
Die Grundlehren der math.\ Wissenschaften in Einzeldarstellungen,
Band 198, Springer-Verlag, New York, 1973


\bibitem{GRSh} I.\ Gelfand, D.\ Ra\u{\i}kov, G.\ Shilov. \emph{Commutative
Normed Rings}. Chelsea, New York, 1964


\bibitem{HelBA} A.\ Ya.\ Helemskii. \emph{Banach and Locally Convex Algebras}.
Oxford University Press, 1993


\bibitem{LarBA} R.\ Larsen. \emph{Banach Algebras, an introduction}.
Marcel Dekker, New York, 1973


\bibitem{LoBA} L.\ H.\ Loomis. \emph{An Introduction to Abstract Harmonic Analysis}.
van Nostrand, 1953


\bibitem{Mos} R.\ D.\ Mosak. \emph{Banach Algebras}. Chicago-London,
University of Chicago Press, 1975


\bibitem{NaiBA} M.\ A.\ Na\u{\i}mark. \emph{Normed Algebras}. Third American
edition, translated from the second Russian edition, Wolters-Noordhoff Publishing,
Groningen, 1972


\bibitem{Neu} M.\ A.\ Neumark. \emph{Normierte Algebren}. Translation of the
second Russian edition, Verlag Harri Deutsch, Thun \& Frankfurt am Main, 1990


\bibitem{Palm} T.\ W.\ Palmer. Encyclopedia of Mathematics and its
Applications, vol.\ 49: \emph{Algebras and Banach Algebras}, \& vol.\ 79:
\emph{Banach Algebras and the General Theory of \st-Algebras}. Cambridge
University Press, 1994, 2001 \pagebreak


\bibitem{Ri} C.\ E.\ Rickart. \emph{General Theory of Banach Algebras}. Van
Nostrand, Princeton, 1960


\bibitem{Zel} W.\ \.{Z}elazko. \emph{Banach Algebras}. Elsevier, New York, 1973


\bibitem{ZhuB} K.\ Zhu. \emph{An Introduction to Operator Algebras}.
CRC Press, 1993


\vspace{1ex}

\begin{center}
\hspace{-2.3em} \textbf{Representations and Positive Linear Functionals}
\end{center}

\vspace{1ex}


\bibitem{Dieu} J.\ Dieudonn\'{e}. \emph{Treatise on Analysis}. vol.\ II,
Enlarged and Corrected Printing, Academic Press, 1976 


\bibitem{DixP} J.\ Dixmier. \emph{Les C*-alg\`{e}bres et leurs
repr\'{e}sentations}. 2\`eme \'ed.\ Gauthier-Villars, 1969. R\'eimpression
1996 \'editions Jacques Gabay


\bibitem{Gaal} S.\ A.\ Gaal. \emph{Linear Analysis and Representation Theory}.
Die Grundlehren der math.\ Wissenschaften in Einzeldarstellungen,
Band 198, Springer-Verlag, New York, 1973


\bibitem{HelRP} A.\ Ya.\ Helemskii. \emph{Banach and Locally Convex Algebras}.
Oxford University Press, 1993


\bibitem{MosR} R.\ D.\ Mosak. \emph{Banach Algebras}. Chicago-London,
University of Chicago Press, 1975


\bibitem{NaiRP} M.\ A.\ Na\u{\i}mark. \emph{Normed Algebras}. Third American
edition, translated from the second Russian edition, Wolters-Noordhoff Publishing,
Groningen, 1972


\bibitem{NeuR} M.\ A.\ Neumark. \emph{Normierte Algebren}. Translation of the
second Russian edition, Verlag Harri Deutsch, Thun \& Frankfurt am Main, 1990


\bibitem{RiR} C.\ E.\ Rickart. \emph{General Theory of Banach Algebras}. Van
Nostrand, Princeton, 1960


\vspace{1ex}

\begin{center} \hspace{-2.3em} \textbf{Spectral Theory of Representations} \end{center}

\vspace{1ex}


\bibitem{AcGl} N.\ I.\ Achieser, I.\ M.\ Glasmann. \emph{Lineare Operatoren
im Hilbert-Raum}. Verlag Harri Deutsch, 8-te erweiterte Auflage, Thun, 1981


\bibitem{AkGl} N.\ I.\ Akhiezer, I.\ M.\ Glazman. \emph{Theory of Linear
Operators in Hilbert space}. Frederick Ungar, vol.\ I, 1961, vol.\ 2, 1963.
Reprinted 1993,\ two volumes bound as one, Dover


\bibitem{Conw2} J.\ B.\ Conway. \emph{A Course in Functional Analysis, second
edition}. GTM 96, Springer-Verlag, New York Berlin Heidelberg Tokyo, 1990


\bibitem{HelST} A.\ Ya.\ Helemskii. \emph{Banach and Locally Convex Algebras}.
Oxford University Press, 1993


\bibitem{MosS} R.\ D.\ Mosak. \emph{Banach Algebras}. Chicago-London,
University of Chicago Press, 1975


\bibitem{NaiST} M.\ A.\ Na\u{\i}mark. \emph{Normed Algebras}. Third American
edition, translated from the second Russian edition, Wolters-Noordhoff Publishing,
Groningen, 1972


\bibitem{NeuS} M.\ A.\ Neumark. \emph{Normierte Algebren}. Translation of the
second Russian edition, Verlag Harri Deutsch, Thun \& Frankfurt am Main, 1990


\bibitem{PedST} G.\ K.\ Pedersen. \emph{Analysis Now}. GTM 118,
Springer-Verlag, New York, 1988


\bibitem{RFA} W.\ Rudin. \emph{Functional Analysis}. International Series in
Pure and Applied Mathematics, McGraw-Hill, TMH Edition 1982, 2nd Edition 1991


\bibitem{SeKu} I.\ E.\ Segal, R.\ A.\ Kunze. \emph{Integrals and Operators}.
TMH Edition, Tata McGraw-Hill Publishing Company, 1968


\bibitem{Weid} J.\ Weidmann. \emph{Lineare Operatoren in Hilbertr\"{a}umen}.
Math.\ Leitf\"aden, B.\ G.\ Teubner, Stuttgart, vol.\ I 2002, vol.\ II 2003


\bibitem{ZhuST} K.\ Zhu. \emph{An Introduction to Operator Algebras}.
CRC Press, 1993 \pagebreak



\vspace{1ex}

\begin{center} \hspace{-2.3em} \textbf{Further Material:} \end{center}

\vspace{1ex}

\begin{center}
\hspace{-2.3em} \textbf{C*-Algebras and von Neumann Algebras (advanced)}
\end{center}

\vspace{2ex}



\bibitem{BrRo} O.\ Bratteli, D.\ Robinson. \emph{Operator Algebras and
Quantum Statistical Mechanics}. TMP, Springer-Verlag, vol.\ I (2nd ed.) 1989,
vol.\ II 1981


\bibitem{DixvN} J.\ Dixmier. \emph{Les alg\`{e}bres d'op\'{e}rateurs dans
l'espace hilbertien (alg\`{e}bres de von Neumann)}. 2\`eme \'ed.\
Gauthier-Villars, 1969. R\'eimpression 1996 \'editions Jacques Gabay


\bibitem{DixC} J.\ Dixmier. \emph{Les C*-alg\`{e}bres et leurs
repr\'{e}sentations}. 2\`eme \'ed.\ Gauthier-Villars, 1969. R\'eimpression
1996 \'editions Jacques Gabay


\bibitem{GaalC} S.\ A.\ Gaal. \emph{Linear Analysis and Representation Theory}.
Die Grundlehren der math.\ Wissenschaften in Einzeldarstellungen,
Band 198, Springer-Verlag, New York, 1973


\bibitem{KaRi} R.\ V.\ Kadison, J.\ R.\ Ringrose. \emph{Fundamentals of the
Theory of Operator Algebras}, vol.\ I, II: Academic Press, 1983, 1986,
vol.\ III, IV: Birkh\"{a}user, 1991, 1992


\bibitem{PedC} G.\ K.\ Pedersen. \emph{C*-Algebras and their
Automorphism Groups}. London Mathematical Society Monographs
No.\ 14, Academic Press, London, 1979


\bibitem{Sak} S.\ Sakai. \emph{C*-Algebras and W*-Algebras}. Ergebnisse der
Mathematik und Ihrer Grenzgebiete, vol.\ 60, Springer-Verlag, 1971. Reprinted
1998, Classics in Mathematics


\bibitem{Tak} M.\ Takesaki. \emph{Theory of Operator Algebras}.
Encyclopaedia of Mathematical Sciences, 3 volumes, 2nd printing
of the 1979 First Edition, 2001


\vspace{1.5ex}

\begin{center} \hspace{-2.3em} \textbf{Basic Abstract Harmonic Analysis} \end{center}

\vspace{1.5ex}


\bibitem{BH} N.\ Bourbaki. \emph{Th\'{e}ories spectrales, Chapitres 1 et 2}.
\'{E}l\'{e}ments de Math\'{e}matique, R\'{e}impression inchang\'{e}e de
l'\'{e}dition originale de 1967, Springer-Verlag, 2007


\bibitem{Foll} G.\ B.\ Folland. \emph{A Course in Abstract Harmonic
Analysis}. CRC Press, 1995


\bibitem{LarAHA} R.\ Larsen. \emph{Banach Algebras, an introduction}.
Marcel Dekker, New York, 1973


\bibitem{LoAHA} L.\ H.\ Loomis. \emph{An Introduction to Abstract Harmonic Analysis}.
van Nostrand, 1953


\bibitem{NaiHA} M.\ A.\ Na\u{\i}mark. \emph{Normed Algebras}. Third American
edition, translated from the second Russian edition, Wolters-Noordhoff Publishing,
Groningen, 1972


\bibitem{NeuH} M.\ A.\ Neumark. \emph{Normierte Algebren}. Translation of the
second Russian edition, Verlag Harri Deutsch, Thun \& Frankfurt am Main, 1990


\vspace{1.5ex}

\begin{center} \hspace{-2.3em} \textbf{Unbounded Operator Algebras} \end{center}

\vspace{1.5ex}


\bibitem{KoS} K.\ Schm\"udgen. \emph{Unbounded Operator Algebras
and Representation Theory}. Operator Theory, Vol.\ 37, Birkh\"auser
Verlag, 1990


\vspace{1.5ex}

\begin{center} \hspace{-2.3em} \textbf{Positive Definite Functions on Semigroups} \end{center}

\vspace{1.5ex}


\bibitem{BCR} C.\ Berg, J.\ P.\ R.\ Christensen, P.\ Ressel.
\emph{Harmonic Analysis on Semigroups, Theory of Positive
Definite and Related Functions}. GTM 100, Springer, 1984


\end{thebibliography}
\end{document}